\documentclass[a4paper,12pt,reqno,noindent]{amsart}
\usepackage[body={15cm,21.5cm},centering]{geometry} 
\usepackage{amsmath,amsfonts,amssymb,amsthm}
\usepackage{graphics}
\usepackage{epsfig}
\usepackage{esint}
\usepackage{color}
\newcommand{\ignore}[1]{}
\newcommand{\ho}{\mathrm{hom}}
\newcommand{\e}{\varepsilon}

\newcommand{\R}{{\mathbb R}}

\newcommand{\dps}{\displaystyle}
\newcommand{\norm}{|\!|\!|}
\newcommand{\N}{\mathcal{N}}

\newtheorem{proposition}{Proposition}
\newtheorem{theorem}{Theorem}
\newtheorem{remark}{Remark}
\newtheorem{lemma}{Lemma}
\newtheorem{corollary}{Corollary}

\newcommand{\step}[1]{\noindent \textit{Step} #1.}
\newcommand{\substep}[1]{\noindent \textit{Substep} #1.}
\newcommand{\calR}{\mathcal{R}}

\title[The corrector in stochastic homogenization]{The corrector in stochastic homogenization:
optimal rates, stochastic integrability, and fluctuations}

\author[A. Gloria \& F. Otto]{Antoine Gloria \& Felix Otto}
\date{\today}
\address[Antoine Gloria]{Universit\'e Libre de Bruxelles (ULB) \\ Brussels, Belgium \\ and Inria  \\ Villeneuve d'Ascq, France}
\email{agloria@ulb.ac.be}
\address[Felix Otto]{Max-Planck-Institut f\"ur Mathematik in den Naturwissenschaften \\ Leipzig, Germany}
\email{otto@mis.mpg.de}
\begin{document}
\maketitle

\begin{center}
\begin{minipage}{13cm}
\small{
\noindent {\bf Abstract.}
We consider uniformly elliptic coefficient fields that are randomly
distributed according to a stationary ensemble of a finite range of dependence.
We show that the gradient and flux $(\nabla\phi,a(\nabla \phi+e))$ of the corrector $\phi$, when spatially
averaged over a scale $R\gg 1$ decay like the CLT scaling $R^{-\frac{d}{2}}$.
We establish this optimal rate on the level of sub-Gaussian bounds in terms of
the stochastic integrability, and also establish a suboptimal rate on the level of optimal 
Gaussian bounds in terms of the stochastic integrability. 
The proof unravels and exploits the self-averaging property of the associated
semi-group, which provides a natural and convenient disintegration of scales, and culminates
in a propagator estimate with strong stochastic integrability.
As an application, we characterize the fluctuations of the homogenization commutator, 
and prove sharp bounds on the spatial growth of the corrector, a quantitative two-scale expansion,
and several other estimates of interest in homogenization.
\vspace{10pt}


\vspace{6pt}
}

\end{minipage}
\end{center}

\bigskip

\tableofcontents

\section{Introduction and context}

We are interested in uniformly elliptic coefficient fields $a$ that are randomly
distributed according to a stationary ensemble $\langle \cdot \rangle$ of a finite range of dependence. 
By the qualitative stochastic homogenization theory of Kozlov \cite{Kozlov} and Papanicolaou \& Varadhan \cite{PapanicolaouVaradhan}, the behavior of the inverse operator $(-\nabla \cdot a \nabla)^{-1}$ at large scales is described by that of 
the constant-coefficient operator $(-\nabla \cdot a_\ho \nabla)^{-1}$, where $a_\ho$ are the homogenized coefficients, characterized in direction $e$
by the corrector $\phi$,  unique (up to additive constant) sublinear solution at infinity of
$$
-\nabla \cdot a (\nabla \phi+e)=0 \quad \text{in }\mathbb{R}^d,
$$
via the formula
$$
a_\ho e \,=\,\langle a(\nabla \phi+e)\rangle.
$$
The corrector $\phi$ is the key to quantitative homogenization properties, and the main goal of this contribution is to understand how 
quantitative ergodic properties of the coefficient field $a$ are transmitted to $\phi$.
We shall show that the gradient $\nabla\phi$ and the flux $a(\nabla \phi+e)$ of the corrector $\phi$ display, when spatially averaged, 
stochastic cancellations at almost the rate as if $\nabla\phi$ was a local function of the coefficient field.
More precisely, we prove  
\begin{itemize}
\item[(i)] that the (relative) error between the homogenized semi-group and the heterogenenous semi-group is small, which we quantify 
via a propagator estimate with strong stochastic integrability, cf Theorem~\ref{Tau};
\item[(ii)] that the rescaling $R^\frac{d}{2}\Xi(R\cdot)$ of the homogenization commutator $\Xi$, defined in direction $e$ by $\Xi e:=a(\nabla \phi+e)-a_\ho (\nabla \phi+e)$, 
is close to a Gaussian white noise as $R\uparrow \infty$ as a random Schwartz distribution, cf Theorem~\ref{Tau2};
\item[(iii)] several optimal estimates of interest in stochastic homogenization, cf Corollaries~\ref{cor:corrector}--\ref{cor:r*},
and Theorem~\ref{Tsys}.

\end{itemize}

Let us start with the by-now well-developed literature on the subject.
Optimal rates in stochastic homogenization like the CLT-scaling implicit in (ii) have been captured, at first for small ellipticity
contrast in the pioneering work by Naddaf and Spencer \cite{NaddafSpencer}, and later on by Conlon and Naddaf \cite{ConlonNaddaf} and Conlon and Spencer \cite{ConlonSpencer}, 
and more recently in the general case by the authors in \cite{GloriaOtto1,GloriaOtto2,GloriaOttoCorrectorEquation} and by Neukamm and the authors in 
\cite{GloriaNeukammOttoInventiones} for (iii).  Cancellations in $\nabla\phi$ that lead to (ii) have
been established in \cite{GloriaNeukammOtto}. However, these rates have been only captured with a {\it suboptimal}
integrability, the best being the exponential moments in \cite[Theorem~2]{GloriaNeukammOtto}.
On the other hand, suboptimal rates but with the Gaussian bounds have been established by Armstrong and Smart in \cite[Theorem~3.1]{ArmstrongSmart},
and recently by Armstrong, Kuusi, and Mourrat in \cite[Theorem~1.1]{ArmstrongKuusiMourrat} in a much more quantified way.
The merit of the present work is to capture optimal rates with nearly optimal
stochastic integrability on the one hand and nearly-optimal rates with optimal stochastic integrability on the other hand.
The present contribution additionally contains the series of estimates (iii) that are straightforward consequences
of the control of the semi-group.

\medskip

The Gaussianity of the fluctuations of the energy of the corrector was first identified by Nolen \cite{N} using the Stein method, see also \cite{GloriaNolen} by Nolen and the first author for a CLT.
For the corrector, the covariance structure was identified by Mourrat and the second author in \cite{MO}, see also \cite{MN} for the normality of fluctuations.
Similar results were then obtained by Gu and Mourrat for the fluctuations of the solution of the heterogeneous PDE, see \cite{GuM}.
The complete structure of fluctuations based on the notion of homogenization commutator (that characterizes the fluctuations of the solution, and of the flux and of the gradient of the corrector) was recently unravelled by Duerinckx and the authors in \cite{DGO}, based on functional inequalities. Announced in \cite{DGO}, Theorem~\ref{Tau2} identifies
the fluctuations of the homogenization commutator under a finite range of dependence assumption, cf (ii).

\medskip

The present work borrows both in philosophy and tools quite a bit from earlier works.
On the one hand, we work under the assumption of {\it finite range} (like Armstrong and Smart in
\cite{ArmstrongSmart}, which was extended to mixing conditions in \cite{ArmstrongMourrat} by Armstrong and Mourrat in the spirit of the early work \cite{Yurinskii} by Yurinskii) 
and just use stochastic cancellations which come from summing independent random variables
(rather: spatially averaging stationary random fields that are approximately local).
In particular we do not appeal to a more general concentration of measure property,
which captures CLT cancellations in random variables of a more complicated structure. 
Typically these arguments build on an underlying product or Gaussian structure,
and pass via assumptions on the ensemble like the Spectral Gap (SG) introduced into stochastic homogenization
by Naddaf and Spencer in \cite{NaddafSpencer} and extensively used by the authors in \cite{GloriaOtto1,GloriaOtto2,GloriaOttoCorrectorEquation}
and Neukamm and the authors in \cite{GloriaNeukammOttoInventiones}
or the Logarithmic Sobolev Inequality (LSI) introduced by Marahrens and the second author in \cite{MarahrensOtto} and refined by Neukamm and the authors  in \cite{GloriaNeukammOtto}, Fischer and the second author in \cite{FischerOtto2}, and Duerinckx and the first author in \cite{DuerinckxGloria}.

\medskip

On the other hand, we completely bypass the variational and subadditive arguments used in \cite{ArmstrongSmart}
(and adapted to a non-symmetric situation in \cite{ArmstrongMourrat}), the large-scale regularity theory
for elliptic equations with random coefficients (in the form of \cite{ArmstrongSmart} or \cite{GloriaNeukammOtto}), and thus the related part of the strategy laid
out in \cite{ArmstrongKuusiMourrat}.
Instead we use the parabolic approach as in \cite{GloriaNeukammOttoInventiones}, which yields a
convenient disintegration of scales and more flexible results, based on the homogenization error at the level of the semi-group. However, like
the pioneering works by Avellaneda \& Lin \cite{AvellanedaLin87} in the periodic and \cite{ArmstrongSmart} in the random case, we compare the actual solution
to the homogenized solution on all scales. When it comes to estimating the homogenization error,
we use the tools introduced in \cite{GloriaNeukammOtto}, namely the vector potential $\sigma$ for the harmonic coordinates.
Here, in order to be able to buckle, we use the ``modified'' (our language for introducing a massive
term as an infra-red regularization) version $(\phi_T,\sigma_T)$ of this augmented corrector 
$(\phi,\sigma)$, which necessitates to introduce a third field $g_T$. Another important technical element
in order to get a small {\it relative} homogenization error is a parabolic version of a novel (deterministic)
inner regularity estimate introduced by Bella and the second author in \cite[Proof of Theorem~2, Step~6]{BellaOtto} and refined by Bella, Giunti, and the second author in \cite[Lemma~4]{BellaGiuntiOtto}.
Note that since our arguments only rely on PDE analysis, the statement of the results and their proofs are oblivious to whether $a$ is symmetric or not.

\medskip

We refer the reader to the introduction of \cite{GloriaNeukammOtto} for a more thorough discussion of improved regularity theory in the large for random
elliptic operators in divergence form, subject only touched upon in the present article.
Indeed, the present contribution significantly differs from \cite{GloriaNeukammOtto} for it does not rely on the large-scale regularity theory
which is at the core of  \cite{GloriaNeukammOtto}.
In a nutshell, the approach of \cite{GloriaNeukammOtto} consists in combining a  $C^{1,1-}$-regularity theory in the large
with a sensitivity calculus (based on functional inequalities) in the spirit of our earlier works to prove quantitative estimates ;
it allowed us to cover both for the first time
non-symmetric systems and non-integrable correlations of the coefficient field.
Whereas the main effort in \cite{GloriaNeukammOtto} is on the regularity theory (homogenization errors like the quantitative two-scale expansion and other quantitative estimates
are obtained as by-products using the functional inequality), the homogenization error at the level of the semi-group is the driving quantity in the present approach, cf Theorem~\ref{Tau}.
In particular, our analysis only exploits cancellations that occur by summing approximately local random variables.
As announced in  \cite{DGO}, this allows us to characterize the fluctuations  of the homogenization commutator, see Theorem~\ref{Tau2}, using softer arguments than in \cite{DGO}.

\medskip

To conclude this introduction, let us mention that before completing the revision of this manuscript, Armstrong, Kuusi, and Mourrat independently obtained similar results as ours in \cite{ArmstrongKuusiMourrat2,ArmstrongKuusiMourrat3}.
Whereas our method is based on a semi-group approach, their method is based on variational techniques (completing the program they initiated in~\cite{ArmstrongKuusiMourrat}).


\section{Notation, objects, and statement of results}

\subsection{Notation and assumptions}

We say that a coefficient field $a=a(x)$ on $d$-dimensional Euclidean space $\mathbb{R}^d$
is uniformly $\lambda$-elliptic provided 
\begin{equation}\label{1.3}
\xi\cdot a(x)\xi\ge\lambda|\xi|^2\quad\mbox{and}\quad
\xi\cdot a(x)\xi\ge|a(x)\xi|^2
\end{equation}
for all points $x\in\mathbb{R}^d$ and tangent vectors $\xi$, where the
ellipticity ratio $\lambda>0$ is fixed once for all. In the case of symmetric $a$, (\ref{1.3}) is equivalent to
$\lambda|\xi|^2\le\xi\cdot a(x)\xi\le|\xi|^2$; in general, the second condition in (\ref{1.3}),
which is equivalent to $\xi\cdot a(x)^{-1}\xi\ge|\xi|^2$ and thus also invariant under transposition,
yields the more standard upper bound $|a(x)\xi|\le|\xi|$; but it is in the form of (\ref{1.3})
that the constant in the upper bound is preserved under homogenization, see (\ref{Laux4.3b}) in Lemma \ref{Laux4}.
While we use scalar notation and language like above, $a$ is allowed to be the coefficient field of
an elliptic system. We stress that $a$ needs not be symmetric, in fact, no iota in the proof
would change for asymmetric $a$.

\medskip

We consider an ensemble of (that is, a probability measure on the space of) 
$\lambda$-uniformly elliptic coefficient fields $a$ and use the physicists notation $\langle\cdot\rangle$
to address both the probability measure and the expectation. We assume that $\langle\cdot\rangle$
is {\it stationary} in the sense that for any integrable random variable $F=F(a)$ and every
shift vector $z\in\mathbb{R}^d$, which acts on coefficient fields via $a(\cdot+z)(x)=a(x+z)$
and thus on random variables via $F^{z}(a)=F(a(\cdot+z))$, we have $\langle F\rangle=\langle F^z\rangle$.
Moreover, we assume that $\langle\cdot\rangle$ has finite range, in fact, has a {\it unity range of dependence}
by which we mean that for two square integrable random variables $F_i=F_i(a)$, $i=1,2$, that have
the property that $F_i$ depends on $a$ only through $a_{|D_i}$ for two open sets $D_i$ of distance
larger than unity, we have $\langle F_1F_2\rangle=\langle F_1\rangle\langle F_2\rangle$.
We will be cavalier about
measurability issues; however, in case of the qualitative
anchoring of our quantitative estimate stated in Lemma~\ref{Laux17}, it is advantageous to be
precise: Taking inspiration from the work of Dal Maso \& Modica \cite{DalMasoModica}, where
the space $\Omega$ of $\lambda$-uniformly coefficient fields was endowed with the topology
coming from Spagnolo's G-convergence (which in their nonlinear context is seen
as $\Gamma$-convergence), we take the one coming from Murat \& Tartar's H-convergence \cite{MuratTartar} instead, which extends the ideas from \cite{DalMasoModica}
to our non-symmetric context and also makes $\Omega$ compact.
We consider probability measures $\langle\cdot\rangle$ which respect this topology.
It turns out that the qualitative homogenization result, even uniformly in stationary
ensembles $\langle\cdot\rangle$
of range unity, follows handily thanks to this natural choice of topology, cf Lemma~\ref{Laux17}.

\medskip

It is a classical result in qualitative stochastic homogenization that under these
conditions (in fact: finite range might be replaced by mere qualitative ergodicity, cf \cite{Kozlov,PapanicolaouVaradhan}), 
for a direction $e$ (a unit vector in $\mathbb{R}^d$),
there exists a stationary gradient field $\nabla\phi=\nabla\phi(a,x)$
(where stationary means shift-covariant in the sense of $\nabla\phi(a,x+z)=\nabla\phi(a(\cdot+z),x)$)
that is $\langle\cdot\rangle$-square integrable, of vanishing expectation,
and such that for $\langle\cdot\rangle$-a.\ e.\ realization $a$ we have
\begin{equation}\label{o56}
-\nabla\cdot a(\nabla\phi+e)=0.
\end{equation}
Its associated flux will play a crucial role
\begin{equation}\nonumber
q:=a(\nabla\phi+e).
\end{equation}
In view of (\ref{o56}), $\phi$ itself, which can be constructed as a non-stationary random field
that is non-unique up to a random additive constant, 
might be seen as the correction to the scalar potential of the curl-free harmonic
vector field $\nabla\phi+e$, a closed and thus exact $1$-form. In view of (\ref{o56}), also $q$ is a closed and thus
exact $(d-1)$-form and hence admits a ``vector potential'', that is, a $(d-2)$-form $\sigma$, which can be represented
by a skew-symmetric tensor field $\sigma=\{\sigma_{jk}\}_{j,k=1,\cdots,d}$:
\begin{equation}\label{o5800}
q-\langle q\rangle=\nabla\cdot\sigma, \quad \sigma \text{ skew}.
\end{equation}
Clearly, this $(d-2)$ form is
non-unique up to a random $(d-3)$-form. The natural choice of gauge is given by
\begin{equation}\label{o58}
-\triangle\sigma_{jk}=\partial_jq_k-\partial_kq_j.
\end{equation}
It was recently established in \cite[Lemma~1]{GloriaNeukammOtto} that under the conditions of stationarity and ergodicity
there indeed exists a stationary gradient field $\nabla\sigma$ that is square integrable, of vanishing expectation,
and such
that (\ref{o58}) holds almost-surely wrt $\langle\cdot\rangle$. 
Since within the framework of this work, we may recover $(\nabla\phi,q,\nabla\sigma)$
by an approximation via a massive term in (\ref{o56}) and (\ref{o58}), see (\ref{c62}) and (\ref{o17}) (and \cite[Lemma~2.7]{GloriaOttoCorrectorEquation}),
which is in fact how we establish the bounds,
we ask the reader not to worry about how to construct $(\nabla\phi,\nabla\sigma,q)$ as random objects
and how they are uniquely characterized.
For simplicity, we don't indicate the --- linear --- dependence of $(\nabla\phi,\nabla\sigma,q)$ 
on $e$ by our notation.

\medskip

It is convenient
to define spatial averages $f_R$ on scale $R$ of a field $f=f(x)$ by convolution with the Gaussian 
$\frac{1}{R^d}\sqrt{2\pi}^{-d}\exp(-\frac{1}{2}|\frac{x}{R}|^2)$. This is convenient mostly
because of its connection to the constant-coefficient heat kernel and in particular
its semi-group property $(f_{R})_r=f_{\sqrt{R^2+r^2}}$.
If not specified otherwise, we don't distinguish between $f_R(0)$ and $f_R$ in our notation, 
which is an acceptable abuse of language
since the convolution will only be applied to random fields $f(a,x)$ that are stationary
in the above sense, which implies that the law of $f_R(y)$ under the stationary 
$\langle\cdot\rangle$ does not depend on $y$.

\medskip

We shall approximate the corrector $\phi$ and the flux $q$ by the semi-group as follows. 
Let $u$ be the solution of
\begin{eqnarray}
\partial_\tau u-\nabla \cdot a \nabla u\,=\,0,&& \tau>0,\label{e.1}\\
u\,=\,\nabla \cdot (ae),&& \tau=0,\label{e.2}
\end{eqnarray}
we define $\phi(t):=\int_0^t  u d\tau$ and $q(t)=a(e+\nabla \phi(t))$.
We formally expect that $\lim_{t\uparrow \infty} \nabla \phi(t)=\nabla \phi$ and $\lim_{t\uparrow \infty} q(t)=q$.

\subsection{Propagator estimate}

The main result compares the semi-group $S$ based on the spatially variable coefficients $a$ and the 
semi-group $S^{\ho}$ based on the constant coefficients $a_{\ho}$. By both we understand the
semi-groups acting on (spatial) vector fields $q_0$ (fluxes) and which are defined as follows: 
For two times $t\le T$, the propagator $S_{t\rightarrow T}$ applied to some $q_0$ is defined via
\begin{align}\label{au37}
S_{t\rightarrow T}q_0=q_0+\int_t^Ta\nabla v d\tau,
\end{align}
where the scalar time-space field $v$ is the solution of the initial value problem
\begin{align}\label{au38}
\partial_\tau v-\nabla\cdot a\nabla v=0\;\;\mbox{for}\;\tau>t,\quad
v=\nabla\cdot q_0\;\;\mbox{for}\;\tau=t.
\end{align}
The propagator $S^{\ho}_{t\rightarrow T}q_0$ is defined analogously while replacing $a$ by $a_{\ho}$ in
(\ref{au37}) and (\ref{au38}). Since neither $a$ nor $a_{\ho}$ depend on time,
we clearly have $S_{t\rightarrow T}=S_{0\rightarrow T-t}$ and
$S^{\ho}_{t\rightarrow T}=S^{\ho}_{0\rightarrow T-t}$. We stick to the more general notation since
we will need it in case of an intermediate propagator $S^h_{t\rightarrow T}$ to be introduced
below in (\ref{au40}). The crucial feature of the propagators $S$ and $S^{\ho}$ is their semi-group
property
\begin{equation}\label{au39}
S_{t_0\rightarrow t_2}=S_{t_1\rightarrow t_2}S_{t_0\rightarrow t_1}\quad\mbox{and}\quad
S^{\ho}_{t_0\rightarrow t_2}=S^{\ho}_{t_1\rightarrow t_2}S^{\ho}_{t_0\rightarrow t_1}\quad\mbox{for}\;\;t_0\le t_1\;\le t_2,
\end{equation}
cf Lemma~\ref{Laux2}.
We also note that, for all finite $t$, $S_{t\rightarrow\infty}$ and $S^{\ho}_{t\rightarrow\infty}$ are the Leray projections
(ie the $L^2$-orthogonal projections onto divergence-free fields) wrt to the metric $a$ and $a_{\ho}$
respectively.
In this notation we clearly have 
\begin{align}\label{au59}
q(t)=S_{0\rightarrow t}(ae),
\end{align}
so that it follows from (\ref{au39}) that our fluxes are also propagated by $S$:
\begin{align}\label{au41}
q(T)=S_{t\rightarrow T}q(t).
\end{align}

\medskip

We introduce the following norm that measures the stretched exponential moments of order $s\le 2$
of a random variable $F$
\begin{equation}\label{au16}
\|F\|:=\inf\{M>0|\log\langle\exp(|\frac{F}{M}|^s)\rangle\le 1\},
\end{equation}
the associated semi-norm
\begin{equation}\label{au16b}
\|F\|_*:=\inf\{M>0|\log\langle\exp(|\frac{F-\langle F\rangle}{M}|^s)\rangle\le 1\},
\end{equation}
and a norm on stationary random fields with CLT-scaling built in:
\begin{equation}\label{e.normCLT}
\norm F \norm \,:=\,\sup_{R\ge 1} R^\frac{d}{2} \|F_R\|_*.
\end{equation}
The main result is the following propagator estimate:
\begin{theorem}\label{Tau}
Suppose $\langle\cdot\rangle$ is an ensemble of $\lambda$-uniformly elliptic coefficient
fields which is stationary and of unity range of dependence. 
There exists an exponent $p<\infty$ only depending on $d$, $\lambda>0$, and $s<2$, such that
for all $T\ge t\ge 1$,
\begin{eqnarray}
\norm q(T)-S^{\ho}_{t\rightarrow T}q(t) \norm
&\lesssim& (\frac{1}{\sqrt{t}})^\frac{1}{p},\label{Tau.1}\\
\norm a_\ho \nabla \phi(T)-(a_\ho \nabla \phi(t)+(S^{\ho}_{t\rightarrow T}-1)q(t))\norm 
&\lesssim& (\frac{1}{\sqrt{t}})^\frac{1}{p},\label{Tau.2}
\end{eqnarray}
where $\lesssim$ means $\le C$ with $C<\infty$ a generic constant only depending on $d$, $\lambda$, and~$s$.
(Note that this  estimate also holds for the choice $T=\infty$, in which case $S_{t\rightarrow \infty}^\ho$ is the Leray 
 and $1-S_{t\rightarrow \infty}^\ho$ the Helmholtz projections with respect to $-\nabla \cdot a_\ho \nabla$.) 
\end{theorem}
\begin{remark}
Estimate~\eqref{Tau.2} is written in terms of the propagator. A more explicit form in terms of fields is as follows:
$$
\norm \nabla \phi(T)-(\nabla \phi(t)+\int_t^T \nabla v_{\ho,t} d\tau)\norm 
\,\lesssim\, (\frac{1}{\sqrt{t}})^\frac{1}{p}
$$
where $v_{\ho,t}$ solves
$$
\partial_\tau v_{\ho,t}-\nabla \cdot a_\ho \nabla v_{\ho,t}\,=\,0 \text{ for }\tau>t, \quad v_{\ho,t}(\tau=t)=u(t).
$$
\end{remark}
From this key result, and the flow of lemmas needed in its proof, we may deduce an entire family of 
quantitative results of interest in stochastic homogenization. 

\subsection{Fluctuations of the homogenization commutator}

As a first direct consequence of Theorem~\ref{Tau},  we have the following stochastic cancellations in spatial averages of the triplet $(\nabla\phi,\nabla\sigma,q)$
formed by the gradient of scalar potential (that is, the field), 
the gradient of vector potential, and the flux.
\begin{corollary}\label{cor:CLTdecay}
Suppose $\langle\cdot\rangle$ is an ensemble of $\lambda$-uniformly elliptic coefficient
fields which is stationary and of unity range of dependence. 
When spatially averaged,  $(\nabla\phi,\nabla\sigma,q)$
displays CLT-cancellations in terms of sub-Gaussian integrability: 
For all $d$,  $\lambda>0$, and $s<2$, 
\begin{equation}\label{theo8}
\norm (\nabla\phi,\nabla\sigma,q) \norm \lesssim 1,
\end{equation}
where here and in the sequel $(\nabla\phi,\nabla\sigma,q-\langle q\rangle)$ stands for one
of the components, and where the constant depends on $s$ next to $\lambda$ and $d$.

Likewise, for all $T>0$,
\begin{equation}\label{theo8b}
\norm (\nabla\phi(T),q(T)) \norm,  \norm (\nabla\phi_T,\nabla\sigma_T,q_T) \norm \lesssim 1
\end{equation}
(the constant does not depend on $T$), where $\phi_T,\sigma_T$ and $q_T$
are defined in  \eqref{c62} and~\eqref{o17}.
\end{corollary}

Note that because of $\langle\nabla\phi\rangle=\langle\nabla\sigma\rangle=0$, these two fields do have
not to be re-centered.
While from the point of view of applications (cf quantitative two-scale expansion of Corollary~\ref{2-scale}), we are interested in the stochastic cancellations of the 
augmented corrector $(\nabla\phi,\nabla\sigma)$,
the proof rather passes via the flux $q$; in particular the flux $q$ 
plays the crucial role. However, the proof heavily relies on the version $\sigma_T$ of the
vector potential modified by a massive term. 
\begin{remark}\label{rem:nearly}
Next to Corollary~\ref{cor:CLTdecay}, which yields optimal scaling with nearly-optimal stochastic integrability, we also have the following nearly-optimal scaling
with optimal stochastic integrability:
For all $d$, $\lambda>0$, and $\alpha<\frac d2$,  
\begin{equation}\label{theo8c}
\sup_{R\ge 1} R^\alpha \|(\nabla\phi,\nabla\sigma,q)_R\|_* \lesssim 1,
\end{equation}
for $s=2$, where the constant acquires a dependence on $\alpha$, next to $\lambda$ and $d$.

Likewise, for all $T>0$,
\begin{equation}\label{theo8d}
\sup_{R\ge 1} R^\alpha \| (\nabla\phi(T),q(T))_R\|_*, \sup_{R\ge 1} R^\alpha \| (\nabla\phi_T,\nabla\sigma_T,q_T)_R\|_*  \lesssim 1
\end{equation}
(the constant does not depend on $T$), where $\phi_T,\sigma_T$ and $q_T$
are defined in  \eqref{c62} and~\eqref{o17}.

\end{remark}

\medskip

The canonical nature of the propagator estimates in Theorem \ref{Tau} is revealed by their
close connection to a key object in stochastic homogenization, the ``homogenization commutator'',
which Duerinckx introduced with the two present authors \cite{DGO}. The homogenization commutator is a centered and stationary
random tensor field defined via the corrector, cf (\ref{o56}), and the homogenized coefficient 
$a_{\rm hom}e=\langle a(\nabla\phi+e)\rangle$:
\begin{align}\label{gr4}
\Xi e=a(\nabla\phi+e)-a_{\rm hom}(\nabla\phi+e),
\end{align}
where we recall that $e\in\mathbb{R}^d$ is a fixed direction on which the stationary and centered $\nabla\phi$ depends linearly.
It is obvious that the decay of its spatial averages $\Xi_Re=q_R-a_{\rm hom}(\nabla\phi+e)_R$ for $R\gg 1$ express homogenization
on the level of the corrector in the spirit of H-convergence: 
As $R\uparrow\infty$, the spatial average $q_R$ of the flux $q=a(\nabla\phi+e)$ converges
to the homogenized coefficient $a_{\rm hom}$ applied to the spatial averages of the field $\nabla\phi+e$. The main
insight of the next theorem is that for $R\gg 1$, $\Xi_R\approx\Gamma_R$, where $\Gamma$ is a centered and stationary white noise,
which thus is Gaussian and hence characterized by its covariance tensor ${\mathcal Q}\in\mathbb{R}^{(d\times d)\times(d\times d)}$, 
a four-tensor with the obvious (partial) symmetry, through
\begin{equation}\label{gr7}
\langle(\int\zeta:\Gamma)^2\rangle=\int\zeta:{\mathcal Q}\zeta\quad\mbox{for}\;\zeta\in C_0^\infty(\mathbb{R}^d)^{d\times d}.
\end{equation}
We note that the random field $\Gamma$ 
has the same law as its rescaled version $\mathbb{R}^d\ni\hat x\mapsto R^{\frac{d}{2}}\Gamma(R\hat x)$, so that
we may reformulate the result as saying that for $R\gg 1$, (the law of) 
$\mathbb{R}^d\ni\hat x\mapsto R^{\frac{d}{2}}\Xi(R\hat x)$ and $\Gamma$ are close weakly in $\mathbb{R}^d$, 
which is the point of view taken in Theorem \ref{Tau2}. Since a generic realization of $\Gamma$ is only
a Schwartz distribution in $\mathbb{R}^d$
of order $\frac{d}{2}-$, weak closeness in $\mathbb{R}^d$ will have to be closeness in the sense of the topology of Schwartz distributions. 
In view of the finite order $\frac{d}{2}-$, it turns out that we obtain closeness in the stronger topology defined as the dual
topology to the one on test functions $\zeta$ given by the norm 
\begin{align}\label{gr9}
\sum_{k=0}^{2m}\sup_{\mathbb{R}^d}(1+|x|)^{4m}|\nabla^k\zeta|\quad
\mbox{for fixed}\;m\in\mathbb{N}\;\mbox{with}\;2m>\frac{d}{2}.
\end{align}
We use the word distribution for a field on $\mathbb{R}^d$ that cannot be given a pointwise interpretation (like the Dirac function
or a generic realization of white noise $\Gamma$)
but only an interpretation
as a distribution in the sense of Schwartz, that is, as a continuous linear form on the topological space of Schwartz functions
or rather on the more restricted Banach space defined through (\ref{gr9}).
In Theorem \ref{Tau2}, by a slight abuse of language, we also call
$\mathbb{R}^d\ni\hat x\mapsto R^{\frac{d}{2}}\Xi(R\hat x)$ stationary, 
by which we mean that this random field and any of its translates have
the same law.

\medskip

Theorem \ref{Tau2} reveals that it is the covariance structure of the homogenization commutator $\Xi$ that is 
asymptotically simple, whereas the covariance structures of the flux $q$ and the field $-\nabla\phi$ are implied ones,
coming from the fact that $(q-a_{\rm hom}e,-a_{\rm hom}\nabla\phi)$ are the constant-coefficient
Leray projection $S_{0\rightarrow\infty}^{\rm hom}$ and Helmholtz projections $1-S_{0\rightarrow\infty}^{\rm hom}$ of $\Xi e$,
respectively. 
Equipped with the two propagator estimates, the main idea of the proof of Theorem \ref{Tau2} is clear:
By the triangle inequality, (\ref{Tau.1}) and (\ref{Tau.2}) imply
\begin{align}\label{gr5}
\norm\Xi-\Xi(t)\norm\lesssim(\frac{1}{\sqrt{t}})^\frac{1}{p},
\end{align}
where in analogy to (\ref{gr4}) we set
\begin{align}\label{gr8}
\Xi(t) e=a(\nabla\phi(t)+e)-\bar a_{\ho,t}(\nabla\phi(t)+e),
\end{align}
where $\bar a_{\ho,t}=\langle a(\nabla\phi(t)+e)\rangle$.
One the one hand, according to Lemma \ref{Laux4}, $\Xi(t)$ like its two ingredients $q(t)$ and $\nabla\phi(t)$
is (approximately) local on scale $\sqrt{t}$, so that on scales $R\gg\sqrt{t}$,
it behaves like a tensor-valued stationary white noise, and thus is necessarily Gaussian. 
On the other hand, according to Theorem \ref{Tau} in form of (\ref{gr5}),
the stationary random fields $\Xi$ and $\Xi(t)$ are close provided $\sqrt{t}\gg 1$.

\medskip

Although it might not be clear a priori, the homogenization commutator is also the driving quantity for the propagator estimates in Theorem~\ref{Tau}: 
There is no loss in passing from \eqref{Tau.1}~\&~\eqref{Tau.2} to~\eqref{gr5}.
By Steps~9~\&~10 of the proof of Theorem~\ref{Tau2}, we have for all $t\ge 0$,
$$
\norm q-S^\ho_{0\rightarrow \infty} q(t)\norm+\norm a_\ho \nabla \phi-(a_\ho \nabla \phi(t)+(S^{\ho}_{0\rightarrow \infty}-1)q(t))\norm \,\lesssim\,\norm\Xi-\Xi(t)\norm,
$$
so that \eqref{gr5} implies both  \eqref{Tau.1}~\&~\eqref{Tau.2} for $T=\infty$.

\medskip

There is an easy way to see why the commutator $\Xi$ is a better-behaved quantity than
flux $q=a(\nabla\phi+e)$ or field $\nabla\phi$ individually, in the sense that 
it is a more {\it local} quantity as a function of $a$. To uncover this locality, 
we compare the fields $\Xi=\Xi(a)$ and
$\tilde\Xi=\Xi(\tilde a)$ for two different realizations $a$ and $\tilde a$. 
The better locality of $\Xi$ is uncovered with help of the corrector $\phi'$ and
corresponding vector potential $\sigma'$ of the dual equation (for some direction
$e'\in\mathbb{R}^d$) which according to (\ref{o5800}) satisfies
\begin{align}\nonumber
a'(\nabla\phi'+e')=a_{\rm hom}'e'+\nabla\cdot\sigma'. 
\end{align}
Indeed, a short and elementary 
calculation based on the identities $e'\cdot \Xi e$ 
$=(a'e'-a'_{\rm hom}e')\cdot (\nabla\phi+e)$ 
and $a'e'-a_{\rm hom}'e'$ $=-a'\nabla\phi'+\nabla\cdot\sigma'$ in conjunction
with the skew symmetry of $\sigma'$ in form of 
$(\nabla\cdot\sigma')\cdot\nabla(\phi-\tilde\phi)$
$=-\nabla\cdot(\sigma'\nabla(\phi-\tilde\phi))$ yields 
\begin{align}\label{gr40}
e'\cdot(\Xi-\tilde\Xi)e&=(\nabla\phi'+e')\cdot(a-\tilde a)(\nabla\tilde\phi+e)\nonumber\\
&-\nabla\cdot(\phi'(a-\tilde a)(\nabla\tilde\phi+e))-\nabla\cdot((\phi'a+\sigma')
\nabla(\phi-\tilde\phi)).
\end{align}
The first rhs term in (\ref{gr40}) reveals a completely local dependence:
When one passes from $a$ to $\tilde a$, this term affects the random field 
$\Xi$ only in those regions in space where $a$ and $\tilde a$ differ.
Furthermore, according to Corollary \ref{cor:CLTdecay}, the fields $\nabla\phi'+e'$
and $\nabla\tilde\phi+e$ are bounded (at least after a square average
over a unit ball) with overwhelming probability, so that we see that
this term is not only local, but also begnin in magnitude.
Also the second rhs term of (\ref{gr40}) is completely local. The fact that
it is in {\it divergence form} makes that it is of higher order wrt to 
CLT scaling: When mollified with a Gaussian of scale $R$,
its contribution will be increasingly smaller as $R\uparrow\infty$.
This leverage of the divergence form relies on the fact that,
according to Corollary \ref{cor:corrector}, the fields $\phi'$ and $\sigma'$
are similarly bounded with overwhelming probability (at least for $d>2$, 
and with a logarithmic divergence for $d=2$).
The only non-local effect comes from
the last rhs term in (\ref{gr40}); it is here that the non-locality of the field
enters via $\nabla(\phi-\tilde\phi)=(\nabla\phi+e)-(\nabla\tilde\phi+e)$.
However, it enters only inside the divergence so that also this term is of higher order
wrt to the CLT scaling. Moreover, because of the relation
$\nabla\cdot(\Xi e-\tilde\Xi e)$ $=-\nabla\cdot a_{\rm hom}\nabla(\phi-\tilde\phi)$
and thus
\begin{align*}
a_{\rm hom}\nabla(\phi-\tilde\phi)=
(1-S_{0\rightarrow\infty}^{\rm hom})(\Xi e-\tilde\Xi e)
\end{align*}
in terms of the constant-coefficient Helmholtz projection 
$1-S_{0\rightarrow\infty}^{\rm hom}$, the relationship (\ref{gr40})
allows for a buckling.

\begin{theorem}\label{Tau2}
Suppose $\langle\cdot\rangle$ is an ensemble of $\lambda$-uniformly elliptic coefficient fields
which is stationary and of unit range of dependence.
As $R\uparrow\infty$, the centered and stationary random tensor field
\begin{align}\nonumber
\hat x\mapsto R^\frac{d}{2}\Xi(R\hat x)
\end{align}
converges in law up to almost Gaussian moments, and wrt to the dual topology to (\ref{gr9}), 
to a centered and stationary tensor-valued Gaussian white noise $\Gamma$
characterized by its covariance ${\mathcal Q}\in\mathbb{R}^{(d\times d)\times(d\times d)}$ through (\ref{gr7}).
This covariance tensor ${\mathcal Q}$ may be recovered from (\ref{gr8}) via
\begin{equation}\nonumber
\mathcal Q=\lim_{t\uparrow\infty}\int\langle\Xi(t,z)\Xi(t,0)\rangle dz.
\end{equation}
Moreover, the centered and stationary random tensor fields
\begin{align}\label{Tau2.1}
\hat x\mapsto R^\frac{d}{2}(q(R\hat x)-a_{\rm hom})
\quad\mbox{and}\quad\hat x\mapsto -R^\frac{d}{2}a_{\rm hom}\nabla\phi(R\hat x),
\end{align}
where here $q:=(q_1,\cdots,q_d)$ and $\nabla\phi:=(\nabla\phi_1,\cdots\nabla\phi_d)$,
jointly converge in the above sense to the pair of centered and stationary Gaussian random tensor field 
$S_{0\rightarrow\infty}^{\rm hom}\Gamma$ and $(1-S_{0\rightarrow\infty}^{\rm hom})\Gamma$,
respectively.
We note that $S_{0\rightarrow\infty}^{\rm hom}\Gamma$ and $(1-S_{0\rightarrow\infty}^{\rm hom})\Gamma$ are
independent so that this result implies the one stated first.
\end{theorem}

\subsection{Quantitative stochastic homogenization}

We conclude the statement of main results with quantitative estimates in stochastic homogenization, and
start with the growth of the extended corrector.
\begin{corollary}\label{cor:corrector}
Suppose $\langle\cdot\rangle$ is an ensemble of $\lambda$-uniformly elliptic coefficient
fields which is stationary and of unity range of dependence. 
Then the extended corrector satisfies:
For all $x\in \R^d$,
\begin{equation}
\| (|(\phi,\sigma)-(\phi,\sigma)_1(0)|^2)_1^\frac12(x) \| \,\lesssim \, \mu_d(x)^\frac 12,\label{e.cor-d.1} 
\end{equation}
with stochastic integrability $s=2$ for $d>2$ and any $0\le s<2$ for $d=2$ (in which case the 
constant depends on $s$ next to $\lambda$), 
where $\mu_d(x)=1$ for $d>2$ and $\mu_d(x)=\log (2+|x|)$ for $d=2$.
 \end{corollary}
\begin{remark}\label{rem:modcorr}
A similar argument as for Corollary~\ref{cor:corrector} (based on estimates for the modified corrector instead of the corrector, cf \eqref{theo8b} \& \eqref{theo8d})
yields the following estimate of the growth of the modified corrector in line with \eqref{e.cor-d.1}: For all $T\gg 1$,
\begin{equation*}
\| (|(\phi_T,\sigma_T)|^2)_1^\frac12\| \,\lesssim \, 
\mu_d(T)^\frac12,
\end{equation*}
with the same stochastic integrability  as in Corollary~\ref{cor:corrector}.
\end{remark}
Estimate~\eqref{e.cor-d.1} is optimal both in scaling and in stochastic integrability (uniformly in $\lambda$).
The proof of the critical estimate (that is, in dimension $d=2$) makes crucial use of the semi-group.
This corollary extends \cite[Theorem~3]{GloriaNeukammOtto} to the case of finite range of dependence, and improves the stochastic integrability.
It yields the existence of stationary correctors for $d>2$ by soft arguments (cf proof of \cite[Theorem~3]{GloriaNeukammOtto}).
Since we control $\sigma$, and not only $\phi$, this yields the following quantitative two-scale expansion.
\begin{corollary}[Quantitative two-scale expansion]\label{2-scale}
Suppose $\langle\cdot\rangle$ is an ensemble of $\lambda$-uniformly elliptic coefficient
fields which is stationary and of unity range of dependence, and let $\phi_i$ denote the correctors.
Let $f\in L^2(\R^d)$ be a compactly supported function (with mean-value zero if $d=2$).
For all $\e>0$, let $u_\e\in H^1(\R^d)$ denote the unique weak solution of 
$$
-\nabla \cdot a(\frac{\cdot}{\e})\nabla u_\e\,=\,f,
$$
and let $u_\ho\in H^1(\R^d)$ denote the unique weak solution of the homogenized problem
$$
-\nabla \cdot a_\ho\nabla u_\ho\,=\,f.
$$
Then we have the following quantitative two-scale expansion (summation convention on repeated indices): 
\begin{multline*}\label{eq:2-scale}
\| \Big(\int |\nabla u_\e -\nabla u_\ho-\partial_i u_\ho \nabla \phi_i(\frac{\cdot}{\e})|^2\Big)^\frac{1}{2}\| \\
\lesssim\,
\Big(\int \mu_d(x)|\nabla \nabla u_\ho|^2(x)dx\Big)^\frac{1}{2} \e |\mu_d(\e)|^\frac{1}{2},
\end{multline*}
with the same stochastic integrability as in Corollary~\ref{cor:corrector}.
\end{corollary}
For a proof of Corollary~\ref{2-scale} based on the results of Corollary~\ref{cor:corrector}, we refer 
the reader to~\cite{GloriaNeukammOtto}.

\medskip

Next, we draw consequences of Theorem~\ref{Tau} in relation to the results of \cite{GloriaNeukammOttoInventiones}.
The parabolic equation, acting as here on stationary fields, can be lifted to
probability space and then is the Fokker-Planck equation of the process of the
``environment as seen from the particle''. It has been the starting point of
a quantitative approach to homogenization in \cite{Mourrat,GloriaMourrat,GloriaNeukammOttoInventiones},
where optimal decay rates for algebraic moments have been established in 
\cite[Theorem~7]{GloriaMourrat} for low dimensions and in \cite[Theorem~1]{GloriaNeukammOttoInventiones} for all dimensions.
The main advantage of our parabolic approach over the elliptic approach is that we obtain optimal results in all dimensions at once.
The proof of Theorem~\ref{Tau} yields the following by-product that quantifies the decay in time of the semi-group in the spirit of  \cite[Theorem~1]{GloriaNeukammOttoInventiones}.
\begin{corollary}\label{coro:decay-semi-group}
For all $T\ge 1$,
\begin{eqnarray}\label{e.decay-semi-1}
\text{for all }0\le s<2&:& \|\Big(\int\eta_{\sqrt{2T}}|\nabla u(T)|^2\Big)^\frac{1}{2}\|\,\lesssim_s \,   {T}^{-1-\frac{d}{4}},
\\
\text{for }s=2 \text{ and all } \alpha<\frac d2&:&\|\Big(\int\eta_{\sqrt{2T}}|\nabla u(T)|^2\Big)^\frac{1}{2}\|\,\lesssim_\alpha \,   {T}^{-1-\frac{\alpha}{2}},\label{e.decay-semi-1bis}
\end{eqnarray}
where the constant acquires either a dependence on $s$ or $\alpha$.

In particular, for all $T\ge 1$,
\begin{equation}\label{e.decay-semi-2}
\langle |\nabla u(T)|^2 \rangle^\frac12\,
\lesssim\, T^{-1-\frac{d}{4}} ,
\end{equation}
and for $s=2$  and for all $R\ge 1$,
\begin{equation}\label{e.decay-semi-3}
\|\Big (\int\eta_{R}|\nabla \phi|^2\Big)^\frac{1}{2}\|\,\lesssim \,   1.
\end{equation}
\end{corollary}
This corollary yields the optimal control of the so-called systematic errors, extending the bounds obtained in 
\cite[Corollary~1 and Lemma~8]{GloriaNeukammOttoInventiones} in the case of discrete elliptic equations with i.i.d. conductances
to the continuum setting of nonsymmetric elliptic systems and finite range of dependence, cf also \cite[Proposition~2 and Corollary~2]{GloriaOttoCorrectorEquation} for similar results up to dimension $d=4$
(albeit for scalar equations and under a spectral gap assumption).
\begin{theorem}\label{Tsys}
For all $T>0$ we define the massive approximation $\phi_{Ti}$ of the corrector $\phi_i$ in direction $e_i$ as the unique stationary solution with finite second moment of \eqref{c62}.
Likewise, we denote by $\phi_{Tj}'$ the massive approximation of the adjoint corrector $\phi'_j$, associated with the pointwise transpose coefficient field $a'$ of $a$,
in  direction $e_j$.
For all $\kappa\in \mathbb N$ we define the Richardson extrapolation of $\phi_{Ti}$ wrt $T$ by
$$
\phi_{Ti}^1:=\phi_{Ti},\quad \phi_{Ti}^{\kappa+1}:=\frac{1}{2^\kappa-1}(2^\kappa\phi_{2Ti}^\kappa-\phi_{Ti}^\kappa),
$$
and likewise for ${\phi_{Tj}'^\kappa}$, 
and we define the approximation $a^\kappa_{hT}$ of the homogenized coefficients $a_\ho$ by
\begin{equation}\label{e.sys-def-ahT}
e_j\cdot  a_{hT}^\kappa e_i\,:=\, \langle (\nabla \phi_{Tj}'^\kappa+e_j)\cdot a (\nabla \phi_{Ti}^\kappa+e_i)\rangle .
\end{equation}
If $\langle \cdot \rangle$ has range of dependence unity, then the following estimates of the systematic errors hold true: For all $d\ge 2$ and $\kappa>\frac{d}{4}$, 
\begin{eqnarray}
\langle |\nabla \phi_{Ti}^\kappa-\nabla \phi_i|^2\rangle^\frac12&\lesssim & T^{-\frac{d}{4}},\label{e.sys1} \\
|a_{hT}^\kappa-a_\ho|&\lesssim & T^{-\frac{d}{2}},\label{e.sys2} 
\end{eqnarray}
where the multiplicative constant depends on $\kappa$ next to $d$ and $\lambda$.
\end{theorem}
Recall that stationarity allows one to define a differential calculus in probability through the correspondance for stationary fields $\psi$ (cf \cite[Section~2]{PapanicolaouVaradhan}):
$$
D_i\psi(0) = \lim_{h\downarrow 0} \frac{\psi(a(\cdot+he_i),0)-\psi(a,0)}{h}\,=\, \lim_{h\downarrow 0}\frac{\psi(a,he_i)-\psi(a,0)}{h}\,=\,\nabla_i \psi(a,0).
$$
This defines a Hilbert space: $\mathcal H^1=\{\psi \in L^2(\langle \cdot \rangle)\,|\, \langle |D\psi|^2\rangle<\infty\}$.
In the case when the coefficients $a$ are symmetric, the operator $\mathcal L=-D \cdot a(0) D$ defines a quadratic form on $\mathcal H^1$. We denote by $\mathcal L$ its Friedrichs extension on $L^2(\langle\cdot \rangle)$.
Since $\mathcal L$ is a self-adjoint non-negative operator, by the spectral theorem, it admits the spectral resolution
\begin{equation}\label{eq:spectral-reso}
\mathcal L\,=\,\int_0^\infty \mu P(d\mu).
\end{equation}
We obtain as a corollary of Theorem~\ref{Tsys} the following bounds on the bottom of the spectrum of  $\mathcal L$ projected on $D \cdot a(0)e \in (\mathcal H^1)'$:
\begin{corollary}\label{cor:sys}
Let $\langle\cdot\rangle$ be
a stationary ensemble with range of dependence unity that takes values into the set of symmetric coefficient fields, and $\mathfrak{d}:=D \cdot a(0)e_i$.
Then the spectral resolution $P$ of $\mathcal L$ satisfies: For all $\tilde \mu>0$,
\begin{equation}
\langle \mathfrak{d}P(d\mu)\mathfrak{d} \rangle([0,\tilde \mu])\,\lesssim \,
\tilde \mu^{\frac{d}{2}+1},
\end{equation}
where the multiplicative constant depends on $d$ and $\lambda$.
\end{corollary}
This corollary is a direct extension of \cite[Corollary~1]{GloriaNeukammOttoInventiones} to the continuum setting with finite range of dependence.

\medskip

To conclude, we argue that Theorem~\ref{Tau} yields the validity of the $C^{1,1-}$-regularity theory in the large in the form
introduced in \cite{GloriaNeukammOtto}. Note that we do not use this higher-regularity in the large (and in particular \eqref{e.Lip-reg} below) in the proof of
the results of this paper.
This improved regularity theory over classical elliptic regularity theory holds 
at scales larger than a minimal radius $r_*$, which is an almost-surely finite stationary field defined in \cite[Corallary~2]{GloriaNeukammOtto} for ergodic coefficients.
More precisely, for some $0<\delta\ll 1$ arbitrary yet fixed, the associated minimal radius
is defined as:
\begin{eqnarray}\label{r**}
\lefteqn{r_*:=\inf\Big\{r\ge 1,\;\mbox{dyadic}\Big|}\nonumber\\
&&\forall R\ge r,\;\mbox{dyadic}\quad
\frac{1}{R^2}\fint_{B_R}|(\phi,\sigma)-\fint_{B_R}(\phi,\sigma)|^2\le\delta\Big\}.
\end{eqnarray}
The main interest of $r_*$ (for some well-chosen $\delta\ll 1$) is the following Lipschitz regularity property: Let $R\ge r_*$, if $v\in H^1(B_R)$ solves
$$
-\nabla \cdot a \nabla v \,=\,0 \text{ on }B_R,
$$
then  
\begin{equation}\label{e.Lip-reg}
\fint_{B_R} |\nabla v|^2 \,\lesssim \, \fint_{B_{r_*}} |\nabla v|^2.
\end{equation}
In particular, the smaller $r_*$, the sharper the result.
The following corollary gives the optimal stochastic integrability of $r_*$ under the assumption of finite range of dependence.
\begin{corollary}\label{cor:r*}
Suppose $\langle\cdot\rangle$ is an ensemble of $\lambda$-uniformly elliptic coefficient
fields which is stationary and of unity range of dependence. 
Then the minimal radius $r_*$ defined in \eqref{r**} satisfies 
for some $0<C<\infty$
\begin{equation}\label{r*sharp}
\langle \exp(\frac1C r_*^d) \rangle \,<\, \infty.
\end{equation}
\end{corollary}
This result extends \cite[Theorem~1]{GloriaNeukammOtto} (which is based on functional inequalities) to this class of coefficients.
In particular, it improves the nearly-optimal stochastic integrability of \cite[Theorem~1.2]{ArmstrongSmart} for the validity of \eqref{e.Lip-reg} to optimal stochastic integrability.
The difference between \cite[Theorem~1.2]{ArmstrongSmart} and
Corollary~\ref{cor:r*} stems from the very definition of the minimal radii introduced in \cite{ArmstrongSmart}  and in \cite{GloriaNeukammOtto}, and therefore
the whole regularity theory.
Whereas  \cite{ArmstrongSmart} essentially requires quantitative sublinearity of the corrector at infinity in the form of an algebraic decay rate, 
the intrinsic version of large-scale regularity introduced in  \cite{GloriaNeukammOtto} only requires sublinearity of the extended corrector at infinity in form of a smallness condition.
This weaker condition, which already allowed us to show in \cite[Corollary~1]{GloriaNeukammOtto} that the minimal radius $r_*$ is, unlike the one of \cite{ArmstrongSmart}, finite almost-surely under the mere qualitative assumption of ergodicity, also allows us to capture the optimal stochastic integrability of the minimal radius under the strongest quantitative ergodic assumptions possible (namely in form of a logarithmic-Sobolev inequality in \cite[Theorem~1]{GloriaNeukammOtto}
and in form of a finite range of dependence assumption in Corollary~\ref{cor:r*}).

\section{Structure of the proof of the propagator estimate}

Our proof of Theorem~\ref{Tau} combines two different types of arguments: deterministic and stochastic arguments.
We start the presentation of the structure of the proof by the deterministic part, which culminates in the control of the homogenization error,
ie the difference
between the homogenized semi-group and the heterogeneous semi-group, in relative terms on dyadic time scales, cf Proposition~\ref{Paux1}.
We then turn to the stochastic part, and the crucial notion of relative approximate locality, which culminates in the stochastic version of
Proposition~\ref{Paux1} which encodes CLT cancellations on the rhs, cf Proposition~\ref{Paux2}.


\subsection{Deterministic arguments}

The main deterministic result is the following proposition, on which we comment now. Essentially, it is an estimate of
the homogenization error $F:=q(T)-S^{\rm hom}_{\frac{T}{2}\rightarrow T}q(\frac{T}{2})$
$=(S_{\frac{T}{2}\rightarrow T}-S^{\rm hom}_{\frac{T}{2}\rightarrow T})q(\frac{T}{2})$
on the level of the semi-group over the dyadic time interval $[\frac{T}{2},T]$, and applied to the corrector
in form of its approximate flux $q(\frac{T}{2})$ itself. Like in H-convergence, it
considers the error $F$ in a metric that captures spatially weak convergence, namely 
$\sup_{R\le\sqrt{T}}(\frac{R}{\sqrt{T}})^\frac{d}{2}|F_R|$, where the cut-off scale $\sqrt{T}$ is the natural length scale
associated to the time $T$ via the parabolic equation. At this stage, the CLT-exponent $\frac{d}{2}$ of the 
length-scale ratio $\frac{R}{\sqrt{T}}$
is somewhat arbitrary, it will be convenient when it comes to the stochastic arguments in
Subsection~\ref{sec:sto-arg}.
Note that by the shift-covariance of all involved fields 
it does not matter that we consider the field $F_R=F_R(0)$ just at the origin. 

\medskip

It is of crucial importance that Proposition \ref{Paux1} provides a {\it relative} error estimate in the sense that the rhs of 
(\ref{Paux1.1})
is a small number times the same type of norm on $q(\frac{T}{2})$, that is, of the mollification
$\{(\frac{R}{\sqrt{T}})^\frac{d}{2}q_R(\frac{T}{2})\}_{R\le\sqrt{T}}$. In fact, this is only almost true, since
we have to include (time-space averages of) the entire family 
$\{(\frac{r}{\sqrt{r}})^\frac{d}{2}q_r(t,x)\}_{r\le\sqrt{t},\frac{T}{8}\le t\le\frac{T}{2}}$ for times $t\sim T$ and spatial points
$x$ within distance $O(\sqrt{T})$ of the origin, as encoded by the exponential averaging function $\eta_{2\sqrt{T}}$:
\begin{equation}\label{def:eta}
\eta_R(x):=\frac{1}{R^d}\eta(\frac{x}{R})\quad\mbox{with}\quad\eta(x):=\exp(-|x|).
\end{equation}
The small prefactor on the rhs of (\ref{Paux1.1}) again depends on the (approximate) corrector in form of a weak norm
of its flux, giving the entire estimate a non-linear buckled nature. The small prefactor on the rhs has two origins: 
a random and a systematic homogenization error in the parlance of \cite{GloriaOtto1}. The first
contribution $\delta$ may be assimilated to the random error and is estimated by the size
of the modified corrector, see (\ref{c62}) for its definition, more specifically its flux $q_{t_0}$. 
The second contribution $|a_{hT}-a_{\rm hom}|$, see (\ref{a66}) for the definition of $a_{hT}=\langle q_T\rangle$, 
is the systematic error. It will be crucial when it comes to the stochastic arguments in Subsection~\ref{sec:sto-arg}
that the contribution $\delta$ may be estimated by a large spatial average (scale $\sqrt{T}$) of cancellations
in the flux $q_{t_0}$ at a much earlier time $t_0\ll T$:  This is beneficial because $(q_{t_0}-\langle q_{t_0}\rangle)_{\delta\sqrt{t_0}}$
is (approximately) local on scale $\sqrt{t_0}$ so that there will be a strong effect of concentration of measure,
to the effect that the second (random) contribution to (\ref{Paux1.2}) is in fact very close to its expectation, the first
(deterministic) contribution to (\ref{Paux1.2}).
\begin{proposition}\label{Paux1}
There exists $p\gg 1$ such that for all $T$ and all $R\le \sqrt{T}$ we have for some $a_{hT}$ (see \eqref{a66} below)
\begin{multline}\label{Paux1.1}
\Big(\frac{R}{\sqrt{T}}\Big)^\frac{d}{2} \Big|\Big(q(T)-S^\ho_{\frac{T}{2}\to T}q(\frac{T}{2})\Big)_R \Big|
\\
\lesssim \, 
(\delta+|a_{hT}-a_\ho|)^\frac{1}{p}   \fint_{\frac{T}{8}}^\frac{T}{2}dt\fint_0^{\sqrt{T}}dr(\frac{r}{\sqrt{t}})^{\frac{d}{2}}\int\eta_{2\sqrt{T}}|(q(t)-\langle q(t)\rangle)_r|
\end{multline}
provided $\delta \ll 1$ is such that
\begin{equation}\label{Paux1.2}
\langle|(q_{t_0}-\langle q_{t_0}\rangle)_{\delta\sqrt{t_0}}|\rangle+\int\eta_{\sqrt{T}}
|(q_{t_0}-\langle q_{t_0}\rangle)_{\delta\sqrt{t_0}}|\,\leq \,\delta^p 
\end{equation}
for some $t_0\le \delta^2 T$, and where the quantity $|a_{hT}-a_\ho|$ is estimated by
\begin{equation*}
|a_{hT}-a_\ho| \,\lesssim\, \int_0^\infty \frac{dt}{t} \min\{\frac{t}{T},1\} \fint_0^{\sqrt{T}} dr \Big(\frac{r}{\sqrt{t}}\Big)^\frac{d}{2} \langle |(q(t)-\langle q(t)\rangle)_r|\rangle.
\end{equation*}
\end{proposition}
The following eleven auxiliary lemmas are needed for the proof of this proposition and of the main results.

\medskip

Recall the approximation via the semi-group, and define $u$ as the unique solution of \eqref{e.1} \& \eqref{e.2}.
Formally, on expects $\int_0^T \nabla u d\tau \stackrel{T\uparrow \infty}{\longrightarrow} \nabla \phi$, and $q(t)=a(\int_0^t \nabla u d\tau+e) \stackrel{T\uparrow \infty}{\longrightarrow}  q$.
Our only building blocks for obtaining PDE estimates are standard: semi-group estimates in Lemma \ref{Laux1} and energy estimates
in Lemma \ref{Laux3}, of which we will need both a parabolic version and an elliptic version. 
The only slight twist is that we use spatially localized versions of these estimates, where for
convenience we use the exponential localization given by the averaging function $\eta_R$ from (\ref{def:eta}).
\begin{lemma}[Localized semi-group estimates]\label{Laux1}
Let $a$ be a constant-in-time coefficient field, and let 
 $v$ solve
\begin{eqnarray*}
\partial_\tau v-\nabla \cdot a \nabla v\,=\,0,&& \tau>0,\\
v\,=\,\nabla \cdot q_0,&& \tau=0.
\end{eqnarray*}
Then, for all $R \ge \sqrt{T}>0$ we have 
\begin{multline}\label{Laux1.1}
\Big( \int \eta_R |(T\nabla v(T),\sqrt{T} v(T))|^2\Big)^\frac12
+\Big(\int\eta_R \Big|\int_0^T (\nabla v,\frac{1}{\sqrt{T}} v)d\tau\Big|^2 \Big)^\frac12
\,\lesssim \, \Big(\int \eta_R |q_0|^2 \Big)^\frac12,
\end{multline}
and
\begin{equation}\label{Laux1.2}
\sqrt{T} \Big(\int \eta_R |\nabla v(T)|^2\Big)^\frac12\,\lesssim \, \Big( \int \eta_R  | v(\frac{T}{2})|^2\Big)^\frac12.
\end{equation}
\end{lemma}
The localized semi-group estimates from Lemma \ref{Laux1} ensure existence and uniqueness (in the class of locally but uniformly
$L^2$-bounded functions) of the initial-value problem with $q_0=ae$, so that the propagator introduced in (\ref{au37})
is indeed well-defined. 
\begin{lemma}[Localized energy estimates] \label{Laux3}
Let $a$ be a coefficient field (possibly depending on time in the parabolic case).
Let $T>0$, and let $v$, $f$, and $g$ be related through the elliptic equation on $\R^d$
$$
\frac{1}{T}v-\nabla \cdot a \nabla v\,=\,f+\nabla \cdot g,
$$
then we have for all $R\gg \sqrt{T}$
\begin{equation}\label{Laux3.1}
\int \eta_R \Big|(\frac{v}{\sqrt{T}},\nabla v) \Big|^2 \,\lesssim \, \int \eta_R \Big|(\sqrt{T} f, g)\Big|^2.
\end{equation}
Let $v$, $f$, $g$, and $v_0$ be related through the parabolic equation on $\R_+\times \R^d$
$$
\partial_\tau v-\nabla \cdot a \nabla v \,=\,f+\nabla \cdot g, \text{ for }\tau>0, \quad v=v_0, \text{ for }\tau=0,
$$
then we have for all $R\gg \sqrt{T}$
\begin{multline}\label{Laux3.2}
\sup_{t<T} \Big(\int \eta_R |v|^2\Big)^\frac12+\Big(\int_0^T\int \eta_R  \Big|(\frac{v}{\sqrt{T}},\nabla v) \Big|^2\Big)^\frac12 \\
\lesssim \, \Big(\int\eta_R |v_0|^2\Big)^\frac12+\Big(\int_0^T\int \eta_R \Big|(\sqrt{T} f, g)\Big|^2\Big)^\frac12.
\end{multline}
\end{lemma}
In the sequel, for notational convenience, we shall make the following slight abuse and assume that \eqref{Laux3.1} and \eqref{Laux3.2} hold
for all $R\ge \sqrt{T}$ (which would normally require us to replace $\eta(x)=\exp(-|x|)$ by some $\eta_c(x)=\exp(-\frac{|x|}{c})$ for some universal constant $c>0$).

\medskip

The following lemma shows so-defined propagators have the semi-group property, as claimed in \eqref{au39}.
In order to treat also an intermediate propagator, we allow $a$ to depend on time for this result.

\begin{lemma}[Semi-group property]\label{Laux2}
Let $a$ be a coefficient field (that possibly depends on time). 
For all $T\ge t \ge 0$, introduce the flux propagator $S_{t\to T}$ defined by
$$
S_{t\to T} q_0=q_0+\int_t^T a\nabla vd\tau
$$
where $v$ solves
$$
\partial_\tau v -\nabla \cdot a \nabla v=0, \text{ for }\tau >t,\quad v=\nabla \cdot q_0 \text{ for }\tau=t.
$$
Then $S$ satisfies the semi-group property for three times $t_2\ge t_1\ge t_0\ge 0$:
$$
S_{t_1\to t_2}S_{t_0\to t_1} = S_{t_0\to t_2}.
$$
\end{lemma}
The semi-group property of Lemma~\ref{Laux2} yields in particular the crucial relation $u(T)=\nabla\cdot q(T)$ for all $T\ge 0$.

\medskip

Our proof to Theorem~\ref{Tau} involves two approximations of the flux $q_i=a(\nabla \phi_i+e_i)$: the semi-group approximation $q_i(T)$ up to time $T>0$,
and a Yoshida approximation $q_{Ti}=a(\nabla \phi_{Ti}+e_i)$ by a massive term (infra-red regularization),
where $\phi_{Ti}$ is the solution of 
\begin{equation}\label{c62}
\frac{1}{T}\phi_{Ti}-\nabla \cdot a (\nabla \phi_{Ti}+e_i)\,=\,0 \quad \text{ in }\R^d,
\end{equation}
the existence and uniqueness of which is ensured by Lemma~\ref{Laux3} under mild growth (so that $\phi_T$ is stationary).
Since this massive approximation $q_{Ti}$ of $q_i$ is stationary, it allows us to consider a massive approximation $a_{hT}$ of the homogenized coefficients $a_\ho$
defined  for all $1\le i \le d$ and all $T>0$ by
\begin{equation}\label{a66}
a_{hT} e_i \,:=\, \langle q_{Ti} \rangle.
\end{equation}
The next element of our strategy consists in a suitable representation of the error in the two-scale expansion
based on the modified corrector $\phi_{Ti}$.
In order to write the associated residuum (almost) in divergence form, it is necessary to also appeal to a suitable choice of 
a vector potential $\sigma_{Ti}$ of the flux $q_{Ti}$, which in the case of general $d$ 
is written as a skew-symmetric tensor field $\{\sigma_{Tijk}\}_{j,k=1,\cdots,d}$,
see our initial discussion of $\sigma_{jk}$. The appropriate gauge is given by
\begin{equation}\label{o17}
\frac{1}{T}\sigma_{Tijk}-\triangle\sigma_{Tijk}=\partial_jq_{Tik}-\partial_kq_{Tij}
\end{equation}
which when compared to (\ref{o58}) also contains a massive term. 
Because of the massive terms in both (\ref{c62}) and (\ref{o17}), 
we no longer have $\nabla\cdot\sigma_{Ti}=q_{Ti}-\langle q_{Ti}\rangle$, where the divergence
$\nabla\cdot\sigma$ of a tensor field is defined via $(\nabla\cdot\sigma)_j=\partial_k\sigma_{jk}$. 
To capture the defect in this relation we introduce another auxiliary vector field:
\begin{equation}\label{a4}
g_{Ti}-T\triangle g_{Ti}=q_{Ti}-\langle q_{Ti}\rangle-\nabla\phi_{Ti}.
\end{equation}
Again, (\ref{a4}) is suitably well-posed in the whole space by Lemma~\ref{Laux3} so that it defines a stationary field $g_T$.
Equipped with these notations, we have the following relations
and the following two-scale expansion.
\begin{lemma}[Formulas]\label{Laux4}
The following representation formulas hold:
\begin{itemize}
\item[(i)]
The semi-group and Yoshida approximations of the flux are related for all $T>0$ via 
\begin{equation}\label{Laux4.1}
q_T\,=\, \int_0^\infty \tfrac{dt}{T} \exp(-\tfrac{t}{T}) q(t).
\end{equation}
\item[(ii)]The auxiliary field $g_T$ defined through (\ref{a4}) satisfies
\begin{equation}\label{Laux4.2}
q_{Ti}=a_{hT}e_i+\nabla\cdot\sigma_{Ti}+g_{Ti}.
\end{equation}
\item[(iii)]The tensor $a_{hT}$ defined through (\ref{a66}) is elliptic in the sense of
\begin{equation}\label{Laux4.3}
\xi\cdot a_{hT}\xi\ge\lambda|\xi|^2\quad\mbox{and}\quad
\xi\cdot a_{hT}\xi\ge|a_{hT}\xi|^2.
\end{equation}
Likewise, 
\begin{equation}\label{Laux4.3b}
\xi\cdot a_{\ho}\xi\ge\lambda|\xi|^2\quad\mbox{and}\quad
\xi\cdot a_{\ho}\xi\ge|a_{\ho}\xi|^2.
\end{equation}
\item[(iv)] 
Let $T>0$, and $v$ and $v_h$ be solutions of 
$$
\partial_t v-\nabla \cdot a \nabla v=0, \quad \partial_t v_h-\nabla \cdot a_{hT}\nabla v_h=0, \quad t>0,
$$
with the same initial condition $v=v_h$ at $t=0$.
We define the two-scale expansion error $w:=v-(1+\phi_{Ti}\partial_i)v_h$, and
have the following representation of the field/flux difference
\begin{eqnarray}
\nabla v-\nabla v_h&=&\nabla w+\nabla(\phi_{Ti}\partial_i v_h),\label{Laux4.4}\\
a\nabla v-a_{hT}\nabla v_h&=&a\nabla w+\nabla\cdot(\partial_iv_h\sigma_{Ti})+\partial_i v_hg_{Ti}
+(\phi_{Ti}a-\sigma_{Ti})\nabla\partial_iv_h,\label{Laux4.5}
\end{eqnarray}
and of the residuum 
\begin{equation}
\partial_\tau w-\nabla\cdot a\nabla w
=\, \nabla\cdot((\phi_{Ti}a-\sigma_{Ti})\nabla\partial_i v_h)+\nabla \cdot(\partial_i v_h g_{Ti})-\phi_{Ti} \partial_\tau \partial_iv_h\label{Laux4.6}.
\end{equation}
\end{itemize}
\end{lemma}
We are in the position to introduce an intermediate flux propagator $S^h$, which we obtain
by replacing the heterogeneous coefficients $a$ by the massive approximation $a_{hT}$ of 
the homogenized tensor $a_\ho$ on the dyadic time intervals of the form $(\frac{T}{2},T)$.
More precisely, $S^h_{t\to T}q_0$ is defined for arbitrary $q_0$ by
\begin{equation}\label{au40}
S^h_{t\to T}q_0\,=\, q_0+\int_t^Ta_h \nabla v_h d\tau,
\end{equation}
where 
\begin{eqnarray*}
\partial_\tau v_h-\nabla \cdot a_h \nabla v_h&=&0, \quad \tau>t,\\
v_h&=&\nabla \cdot q_0, \quad \tau=t,
\end{eqnarray*}
and $a_h(\tau)=a_{hT}$ for $\frac{T}{2}<\tau\le T$, $T=2^k$ dyadic, with $a_{hT}$ the massive approximation
defined in \eqref{a66}.
The following lemma gives an estimate of the intermediate homogenization error on dyadic intervals 
provided the extended corrector is small.
\begin{lemma}\label{Laux5}
We have for all $T>0$ dyadic and all $0<R\le \sqrt{T}$,
\begin{equation}\label{Laux5.1}
\Big(\frac{R}{\sqrt{T}}\Big)^{\frac{d}{2}+1} \Big| \Big( (S_{\frac{T}{2}\to T}-S^h_{\frac{T}{2}\to T})(q(\tfrac{T}{2}))\Big)_R\Big|
\,\lesssim \,
\delta^{\frac{1}{\frac{d}{2}+3}} T\Big(\int \eta_{\sqrt{T}}\Big|\nabla u(\tfrac{T}{2})\Big|^2 \Big)^\frac{1}{2}
\end{equation}
provided $\delta \le 1$ satisfies
\begin{equation}\label{Laux5.2}
\Big( \int\eta_{\sqrt{T}}\Big|(\frac{\phi_T}{\sqrt{T}},\frac{\sigma_T}{\sqrt{T}},g_T) \Big|^2 \Big)^\frac12\,\leq \, \delta,
\end{equation}
where we recall that $u$ solves \eqref{e.1} \& \eqref{e.2}.
\end{lemma}
We shall now appeal to two deterministic inner regularity results for elliptic systems.
First, the lhs of \eqref{Laux5.1} is a weak norm of the flux (due the convolution with the Gaussian), whereas the rhs is a strong norm of the field.
In order to buckle, we will need to control the strong norm of the field by weak norms of the flux. Lemma~\ref{Laux6} yields this 
crucial deterministic ingredient.
Second, we need to upgrade the lhs of \eqref{Laux5.1} to the CLT scaling (that is, replace the exponent $\frac{d}{2}+1$ by the exponent $\frac{d}{2}$).
This will be obtained at the price of reducing the exponent on $\delta$ on the rhs by Meyers' estimate, cf Lemma~\ref{Laux7}.
\begin{lemma}\label{Laux6}
For any homogeneous solution $v$ of the semi-group equation
$\partial_t v-\nabla\cdot a\nabla v=0$ for $t\ge 0$
we have for any exponent $p<\infty$ and time $T$
\begin{multline}\label{Laux6.1}
\sqrt{T}\Big(\int\eta_{\sqrt{T}}|\nabla v(\tfrac{T}{2})|^2\Big)^\frac{1}{2}+\Big(\int\eta_{\sqrt{T}}| v(\tfrac{T}{2})|^2\Big)^\frac{1}{2}
\\
\lesssim \, \fint_{\frac{T}{4}}^\frac{T}{2}dt\fint_0^{\sqrt{T}}dr(\frac{r}{\sqrt{t}})^ {p+1}\int\eta_{2\sqrt{T}}|v_{r}(t)| ,
\end{multline}
where here and in the proof, $\lesssim$ refers to $\le C$,
where $C$ denotes a generic constant that only depends on $d$, $\lambda>0$, and on $p<\infty$.

Applied to $v=u$, the solution of \eqref{e.1} \& \eqref{e.2}, this turns into
\begin{multline}\label{Laux6.2}
\sqrt{T}\Big(\int\eta_{\sqrt{T}}|\nabla u(\tfrac{T}{2})|^2\Big)^\frac{1}{2}+\Big(\int\eta_{\sqrt{T}}| u(\tfrac{T}{2})|^2\Big)^\frac{1}{2}
\\
\lesssim \, \frac{1}{\sqrt{T}}  \fint_{\frac{T}{4}}^\frac{T}{2}dt\fint_0^{\sqrt{t}}dr(\frac{r}{\sqrt{t}})^{p}\int\eta_{2\sqrt{T}}|(q(t)-\langle q(t)\rangle)_r| .
\end{multline}
\end{lemma}
While it is well-known that the lhs of (\ref{Laux6.1}) is estimated by the (suitably localized) $H^{-1}$-norm
of $u(t=0)$, the rhs of (\ref{Laux6.1}) can be assimilated with much weaker norm
(for $p$ large) and moreover is an $L^1$-average in $x$ instead of an $L^2$-average.
Estimate (\ref{Laux6.1}) can be seen as a parabolic version of \cite[Lemma~4]{BellaGiuntiOtto}.
The argument however is rather different since it is more about localization in time than in space.
\begin{lemma}[Meyers' estimate]\label{Laux7}
Let $a$ be a (possibly time-dependent) coefficient field.
There exists $\e=\e(d,\lambda)>0$ such that for a solution of the homogeneous parabolic equation 
$\partial_t v-\nabla\cdot a\nabla v=0$ for $t> 0$
and for all $0<R\le \sqrt{T}$,
\begin{equation}\label{Laux7.1}
\Big(\int_{\frac{T}{2}}^T \int \eta_R |\nabla v|^2 d\tau\Big)^\frac12\,\lesssim \, 
\Big(\frac{\sqrt{T}}{R} \Big)^{\frac{d}{2}-\e} \Big( \int_{\frac{T}{4}}^T \int \eta_{\sqrt{T}} |\nabla v|^2 d\tau\Big)^\frac{1}{2}.
\end{equation}
\end{lemma}
Lemmas~\ref{Laux6} and \ref{Laux7} will help us upgrade \eqref{Laux5.1} to 
\begin{multline}
\Big(\frac{R}{\sqrt{T}}\Big)^{\frac{d}{2}} \Big| \Big( (S_{\frac{T}{2}\to T}-S^h_{\frac{T}{2}\to T})(q(\tfrac{T}{2}))\Big)_R\Big|
\\
\lesssim \,
\delta^\frac{1}{p}   \fint_{\frac{T}{8}}^\frac{T}{2}dt\fint_0^{\sqrt{T}}dr(\frac{r}{\sqrt{t}})^{\frac{d}{2}}\int\eta_{2\sqrt{T}}|(q(t)-\langle q(t)\rangle)_r|
\label{e.3}
\end{multline}
for some $p\gg 1$ (note the CLT scaling on both the lhs and rhs).
This estimate will however only hold conditioned on \eqref{Laux5.2}.
The following lemma rephrases this condition in terms of spatial averages of the flux itself.
\begin{lemma}\label{Laux8}
For all $0<\delta\le 1$  and all $\sqrt{t}\le\delta\sqrt{T}$, we have
\begin{eqnarray}\label{Laux8.1}
\lefteqn{\Big(\int\eta_{\sqrt{T}}|(\frac{\phi_T}{\sqrt{T}},\frac{\sigma_T}{\sqrt{T}},g_T)|^2
\Big)^\frac{1}{2}}\nonumber\\
&\lesssim&\delta+\frac{1}{\delta^{\frac{d}{2}+7}}\Big(\langle|(q_t-\langle q_t\rangle)_{\delta\sqrt{t}}|\rangle+\int\eta_{\sqrt{T}}
|(q_t-\langle q_t\rangle)_{\delta\sqrt{t}}|\Big).
\end{eqnarray}
\end{lemma}

To deduce Corollary~\ref{cor:r*} from Theorem~\ref{Tau}, we shall need a slightly modified version of \cite[Proposition~1]{GloriaNeukammOtto}, the proof of which borrows some arguments from the proof of Lemma~\ref{Laux8}.
\begin{lemma}\label{Laux8bis}
We are given a dyadic radius $r$ and a cut-off scale $\sqrt{t}$ with $\sqrt{t}\le r$;
for every dyadic radii $R\ge r$ we are given a cut-off scale $\sqrt{T}$ satisfying
\begin{equation}\label{Cr*-1}
\sqrt{T}\sim(\frac{R}{r})^\alpha\sqrt{t}
\end{equation}
for some exponent $\alpha\in[0,1]$ (so that in particular $\sqrt{t}\lesssim\sqrt{T}\lesssim R$). 
Provided $r\ge r_*$ (the minimal radius defined in \eqref{r**}) we have for any exponent $\beta>0$
\begin{align}\label{Cr*-2}
\lefteqn{\frac{1}{r^2}\fint_{B_r}|(\phi,\sigma)-\fint_{B_r}(\phi,\sigma)|^2}\nonumber\\
&\lesssim&\sup_{R\ge r,\;\mbox{dyadic}}(\frac{R}{r})^\beta\int\eta_{R}\big|(\frac{\phi_T}{\sqrt{T}},\frac{\sigma_T}{\sqrt{T}},
(q_{T}-\langle q_T\rangle)_{\sqrt{T}})\big|^2.
\end{align}
Here $\lesssim$, $\sim$ refer to relations up to a generic constant only depending on the dimension $d$, the ellipticity ratio $\lambda>0$,
and the exponents $\alpha\in[0,1]$ and $\beta>0$.
\end{lemma}

\medskip

Lemma~\ref{Laux5} compares the heterogeneous semi-group to the massive approximation of the homogenized semi-group.
The following lemma compares the homogenized semi-group to its massive approximation.
\begin{lemma}\label{Laux9}
We have for all $T>0$ dyadic and all $0<R\le \sqrt{T}$,
\begin{equation}\label{Laux9.1}
\Big(\frac{R}{\sqrt{T}}\Big)^{\frac{d}{2}} \Big| \Big( (S^h_{\frac{T}{2}\to T}-S^\ho_{\frac{T}{2}\to T})(q(\tfrac{T}{2}))\Big)_R\Big|
\,\lesssim \,
|a_{hT}-a_\ho|\sqrt{T} \Big(\int \eta_{\sqrt{T}}\Big| u(\tfrac{T}{2})\Big|^2 \Big)^\frac{1}{2},
\end{equation}
where we recall that $u$ solves \eqref{e.1} \& \eqref{e.2}.
\end{lemma}
Lemma~\ref{Laux9} allows us to upgrade \eqref{Laux5.1} in Lemma~\ref{Laux5}, in its refined form \eqref{e.3}, to
\begin{multline}\label{e.4}
\Big(\frac{R}{\sqrt{T}}\Big)^{\frac{d}{2}} \Big| \Big( (S_{\frac{T}{2}\to T}-S^\ho_{\frac{T}{2}\to T})(q(\tfrac{T}{2}))\Big)_R\Big|
\\
\lesssim \,
(\delta+|a_{hT}-a_\ho|)^\frac{1}{p}   \fint_{\frac{T}{8}}^\frac{T}{2}dt\fint_0^{\sqrt{T}}dr(\frac{r}{\sqrt{t}})^{\frac{d}{2}}\int\eta_{2\sqrt{T}}|(q(t)-\langle q(t)\rangle)_r|.
\end{multline}
The desired Proposition~\ref{Paux1} will then follow from \eqref{e.4} and the following estimate of the systematic error $|a_{hT}-a_\ho|$
based on the connection
(\ref{Laux4.1}) between $\{q(t)\}_t$ and $\{q_T\}_T$.
\begin{lemma}\label{Laux10}
The systematic error is estimated for all $T>0$ by
\begin{equation}\label{Laux10.1}
|a_{hT}-a_\ho| \,\lesssim\, \int_0^\infty \frac{dt}{t} \min\{\frac{t}{T},1\} \fint_0^{\sqrt{T}} dr \Big(\frac{r}{\sqrt{t}}\Big)^\frac{d}{2} \langle |(q(t)-\langle q(t)\rangle)_r|\rangle.
\end{equation}
\end{lemma}
%


\subsection{Stochastic arguments}\label{sec:sto-arg}

The aim of this second main part of our arguments is to upgrade Proposition~\ref{Paux1} 
by taking the norm $\norm\cdot \norm$ of \eqref{Paux1.1} defined in \eqref{e.normCLT}.
Recall that the stochastic norm $\|\cdot \|$  and semi-norm $\|\cdot \|_*$, defined in \eqref{au16} and \eqref{au16b},
are characterized by the stochastic integrability $0<s\le2$.
We shall iterate in time the following result in the proof of Theorem~\ref{Tau}.
\begin{proposition}\label{Paux2}
There exists $p\gg 1$ such that for all $s<2$  and all dyadic $T$,
\begin{multline}\label{Paux2.1}
\norm q(T)-S_{\frac{T}{2}\to T}^\ho q(\tfrac{T}{2})\norm \,
\lesssim \,(\delta+|a_{hT}-a_\ho|)^\frac1p \sup_{t\le T}  \norm q(t)\norm
\\
+\frac{1}{\delta^p} \sqrt{t_0}^{\frac{d}{2}\frac{1}{s}} \Big(\frac{1}{\sqrt{T}}\Big)^{\frac{d}{2}(\frac{1}{s}-\frac12)}
\end{multline}
provided $\delta$ and $t_0$ satisfy $0<\delta \ll 1$, $t_0\le \delta^2 T$, and 
\begin{equation}\label{Paux2.2}
\langle |(q_{t_0}-\langle q_{t_0})_{\delta \sqrt{t_0}}| \rangle \leq \delta^p.
\end{equation}
(If we replace $S^\ho$ by $S^h$ in the lhs of \eqref{Paux2.1}, we may drop the term $|a_{hT}-a_\ho|$ in the rhs.)
\end{proposition}
The main cause of a CLT behavior of a random variable $F$, as encoded by control of the norm 
$\norm F\norm=\sup_{R\ge 1}R^\frac{d}{2}\|F_R\|_*$,
is its locality in conjunction with the finite range of $\langle\cdot\rangle$, where locality
means that $F=F(a)$ depends on $a$ only through its restriction $a_{|B_r}$ to some ball $B_r$,
see the elementary but crucial Step~3 in the proof of Lemma \ref{Laux12}. 
In Proposition \ref{Paux2}, we are interested in the homogenization error
$F=q(T)-S^{\rm hom}_{\frac{T}{2}\rightarrow T}q(\frac{T}{2})$. By the properties
of the parabolic initial value problem (\ref{e.1})\&(\ref{e.2}), it is clear that $q(T)$ and thus $F$ depends (at least) exponentially little on
$a_{|B_R}$ for $R\gg\sqrt{T}$ as stated in (\ref{Laux11.2}), so that (approximate) locality comes effortless within the semi-group approach.
However, though exponential, these tails of dependence would affect CLT scaling by a logarithm since they are measured in {\it absolute}
terms, whereas the random variable $F$ itself might be already small in the sense of $\|F_{\sqrt{T}}\|_*\ll 1$. 
In order to have an estimate of $\sup_{R\ge\sqrt{T}}R^\frac{d}{2}\|F_R\|$ in terms of $\|F_{\sqrt{T}}\|$,
we need the tails of dependence to be small {\it relative} to $F$ itself. 
By a suitable PDE estimate, this turns out to be the case for $F=q(T)-\langle q(T)\rangle$,
as Lemma \ref{Laux11} shows, provided one does not just take $(q(T)-\langle q(T)\rangle)_{\sqrt{T}}$ into account when measuring the
tails of dependence, but also (spatial averages of)
the entire family $\{(q(t)-\langle q(t)\rangle)_r\}_{t\le T,r\le\sqrt{T}}$, cf (\ref{Laux11.0}), albeit with a strong
concentration on $t\sim T$ and $r\sim \sqrt{T}$. Equipped with this notion of relative approximate locality
(\ref{Laux12.1}), even if it is just polynomial of a degree $p<\frac{d}{2}$, 
we obtain exact CLT scaling (\ref{Laux12.3}) in Lemma \ref{Laux12}.
\begin{lemma}[Relative approximate locality of $q(T)$]\label{Laux11}
The stationary random fields $F=\nabla \phi(T),q(T),S_{\frac{T}{2}\to T}^\ho q(\frac{T}{2}),S_{\frac{T}{2}\to T}^h q(\frac{T}{2})$
are approximately local on scale $\sqrt{T}$ relative to 
\begin{equation}\label{Laux11.0}
\bar F\,:=\,\fint_0^T dt\Big(\frac{\sqrt{t}}{\sqrt{T}}\Big)^\frac{d}{2}\fint_0^{\sqrt{t}}dr \Big(\frac{r}{\sqrt{t}}\Big)^\frac{d}{2} |(q(t)-\langle q(t) \rangle)_r|
\end{equation}
in the sense of 
\begin{equation}\label{Laux11.1}
\Big( \fint_{B_{\sqrt{T}}} |F(a)-F(\tilde a)|^2\Big)^\frac12 \,\lesssim \, \Big(\frac{\sqrt{T}}{R}\Big)^p \int \eta_R (\bar F(a)+\bar F(\tilde a))
\end{equation}
for all $a,\tilde a \in \Omega$ such that $a=\tilde a$ on $B_{2R}:=\{|x|<2R\}$ and $R\ge \sqrt{T}$, for any $p<\infty$.
Note that $F$ is also approximately exponentially local (exp. in short) with respect to~$1$ in the sense of 
\begin{equation}\label{Laux11.2}
\Big( \fint_{B_{\sqrt{T}}} |F(a)-F(\tilde a)|^2\Big)^\frac12 \,\lesssim \,\exp\Big(-\frac{1}{C}\frac{R}{\sqrt{T}}\Big),
\end{equation}
for all $a,\tilde a \in \Omega$ such that $a=\tilde a$ on $B_{2R}$ and $R\ge \sqrt{T}$. 
\end{lemma}
The locality statement \eqref{Laux11.2} is easier than \eqref{Laux11.1}, but not sufficient for our purposes.
The following lemma shows the interest of relative approximate locality in terms of CLT cancellations.
\begin{lemma}[CLT cancellations]\label{Laux12}
Let $F,\bar F$ be stationary random fields (i.~e. $F(a(\cdot+z,x)=F(a,x+z)$ for all $x,z$) such
that $F$ is approximately local on scale $\sqrt{T}\ge 1$ relative to $\bar F$ in the sense of
\begin{equation}\label{Laux12.1}
\Big(\fint_{B_{\sqrt{T}}} |F(a)-F(\tilde a)|^2 \Big)^\frac12 \,\lesssim \, \Big(\frac{\sqrt{T}}{R}\Big)^p \int \eta_R (\bar F(a)+\bar F(\tilde a))
\end{equation}
provided $a=\tilde a$ on $B_{2R}$ and $R\ge \sqrt{T}$ for some $p>\frac{d}{2}$.
Then, for all $0<s\le 2$ 
\begin{equation}\label{Laux12.3}
\sup_{R\ge \sqrt{T}} \Big(\frac{R}{\sqrt{T}}\Big)^{\frac d2} \|F_{R}\|_*\,\lesssim\,  \Big(\sup_{r\le \sqrt{T}} \big(\frac{r}{\sqrt{T}}\big)^\frac{d}{2} \|F_r\|_*\Big)^{1-\frac{d}{2p}} \|\bar F\|^\frac{d}{2p}+\sup_{r\le \sqrt{T}} \big(\frac{r}{\sqrt{T}}\big)^\frac{d}{2} \|F_r\|_*.
\end{equation}
\end{lemma}
\begin{remark}\label{rem:Laux12}
The following more precise version of \eqref{Laux12.3} in Lemma~\ref{Laux12} is the one established in the proof:
Under the same assumptions as in Lemma~\ref{Laux12}, we have for all $0<s\le 2$, $T\ge 1$, and any averaging kernel $\bar G$
\begin{multline}\label{Laux12.2}
 \|\bar G*F_{\sqrt{T}}\|_*\,\lesssim\, \sqrt{T}^\frac{d}{2}\Big(\int\bar G^{2}\Big)^\frac1{2} 
 \\
 \times \bigg( \Big(\sup_{r\le \sqrt{T}} \big(\frac{r}{\sqrt{T}}\big)^\frac{d}{2} \|F_r\|_*\Big)^{1-\frac{d}{2p}} \|\bar F\|^\frac{d}{2p}+\sup_{r\le \sqrt{T}} \big(\frac{r}{\sqrt{T}}\big)^\frac{d}{2} \|F_r\|_*\bigg).
\end{multline}
This more precise version will be needed in the proofs of Corollary~\ref{cor:corrector} and Theorem~\ref{Tau2}.
\end{remark}
The difficulty to pass from Proposition~\ref{Paux1} to Proposition~\ref{Paux2} is that whereas $\delta$ should be seen as a random variable
in \eqref{Paux1.1} in view of condition \eqref{Paux1.2}, it is deterministic in \eqref{Paux2.1} in view of condition  \eqref{Paux2.2}.
We have to take care of this nonlinearity without compromising too much on the stochastic integrability.
The following two lemmas allow us to make use of Lemma~\ref{Laux12} to 
go from Proposition~\ref{Paux1} to Proposition~\ref{Paux2}.
\begin{lemma}\label{Laux13}
Suppose that we have for some $p\gg 1$ and random variables $F,F_0,F_1$,
and some $h>0$
\begin{equation}\label{Laux13.1}
\langle F_0 \rangle + F_0 \leq 3\delta^p \quad \implies \quad |F|\leq (\delta+h)^\frac1p F_1,
\end{equation}
then for all $s<2$ we have
\begin{equation}\label{Laux13.2}
\langle F_0\rangle \le \delta \quad \implies \quad \|F\|_* \,\leq\, (\delta+h)^\frac1p \|F_1\|_*+\Big(\frac{\|F_0\|_2}{\delta^p} \Big)^\frac{2}{s}\|F\|_{\infty},
\end{equation}
where $\|\cdot \|_2$ stands for the norm \eqref{au16} for $s=2$,  $\|F\|_\infty:=\mathrm{ess} \sup_{a\in \Omega} |F|$.
\end{lemma}
We will apply Lemma \ref{Laux13} to $F_0=|(q_{t_0}-\langle q_{t_0}\rangle)_{\delta\sqrt{t_0}}|$
and thus need to capture stochastic cancellations in form of a concentration of measure for
this quantity up to Gaussian moments, ie for $s=2$. Hence we need locality also for this quantity,
where now we don't need relative but just absolute locality (in view of $s=2$, we would not
be able to leverage relative locality and we don't care for optimality here). Since we will apply
Lemma \ref{Laux13} to $F=(\frac{r}{\sqrt{T}})^\frac{d}{2}
((S_{\frac{T}{2}\rightarrow T}-S_{\frac{T}{2}\rightarrow T}^{\rm hom})q(\frac{T}{2}))_r$, we  need
a uniform bound on that quantity, ie a bound of $\|F\|$ for $s=\infty$. Both, approximate locality and uniform bound,
are provided by the following lemma,
which in addition collects exponential approximate locality (in absolute terms) and uniform bounds
for $\nabla\phi(T)$ and $q(T)$ needed in the proof of Theorem \ref{Tau2}.

\begin{lemma}[Uniform bounds]\label{Laux14}
There is $C<\infty$ such that for all $T>0$, $0< r\le \sqrt{T}$, $t_0>0$, and $0<\delta\le 1$,
\begin{eqnarray}
&&\Big(\frac{r}{\sqrt{T}}\Big)^\frac{d}2 |((S_{\frac{T}{2}\to T}-S_{\frac T2 \to T}^\ho)q(\tfrac T2))_r|\,\,\lesssim \,1,\label{Laux14.1}
\\
&&|(q_{t_0}-\langle q_{t_0}\rangle)_{\delta \sqrt{t_0}} |\,\lesssim \, {\delta^{-\frac d2}} ,\label{Laux14.2}
\\
&&|(q_{t_0}-\langle q_{t_0}\rangle)_{\delta \sqrt{t_0}}| \text{ is approximately local on scale }\sqrt{t_0} \text{ relative to } {\delta^{-\frac d2}},\label{Laux14.3}
\\
&& 
 |(\nabla \phi(T),q(T))_{\sqrt{T}}| \,\lesssim\, 1, \label{Laux14.0}
\\
&& \text{$(\nabla \phi(T),q(T))_{\sqrt T}$  is pointwise approximately local on scale $\sqrt{T}$}\nonumber \\
&&\qquad \text{relative to 1 in the sense that if $a=\tilde a$ in $B_R$, $R\ge \sqrt{T}$, then}\nonumber \\
&& \qquad  |(\nabla \phi(T;a)-\nabla \phi(T;\tilde a),q(T;a)-q(T;\tilde a))_{\sqrt{T}}|\,\lesssim \, \exp(-\frac R{C\sqrt T}). \label{Laux14.0b}
\end{eqnarray}
\end{lemma}
The following lemma shows that the norm $\norm\cdot \norm$ is not increased by applying the homogenized propagator $S^\ho$ or its massive approximation $S^h$.
\begin{lemma}\label{Laux15}
For any $t\le T$ and any stationary random field $F$ we have
\begin{equation}\label{Laux15.1}
\norm S_{t\to T}^{\ho} F\norm, \norm S_{t\to T}^h F\norm \,\lesssim\, \norm F\norm.
\end{equation}
\end{lemma}
In particular, this allows for a (probabilistic) definition of the constant-coefficient Leray projection $S_{0\rightarrow \infty}^\ho$
on fields $F$ with $\norm F\norm<\infty$, see Step~10 in the proof of Theorem~\ref{Tau2}.

\medskip

The following  last two lemmas are needed for the  proof of Theorem~\ref{Tau}.
The first lemma yields a suboptimal bound on $\norm q(T) \norm$ and $\norm \nabla \phi(T) \norm$:
\begin{lemma}\label{Laux16}
For any $T>0$  and
$s= 2$,
\begin{eqnarray}
  \norm q(T)\norm,\norm \nabla \phi(T) \norm &\lesssim & \max\{\sqrt{T}^d,1\}.\label{Laux16.2}
\end{eqnarray}
\end{lemma}
The second lemma is an anchoring result,
which is a consequence of \emph{qualitative} homogenization, and for which 
the topology of H-convergence is very convenient.
\begin{lemma}[Anchoring lemma]\label{Laux17}
For all $\delta>0$, 
\begin{equation}\label{Laux17.1}
\lim_{T\uparrow\infty}\langle| (q(T)-\langle q(T)\rangle)_{\delta\sqrt{T}}|\rangle=0,
\end{equation}
where the limit is uniform wrt the ensemble $\langle \cdot \rangle$, up to a dependence on the dimension $d$ and the ellipticity ratio $\lambda>0$.
\end{lemma}
%


\section{Proofs of the main results}

\subsection{Proof of Theorem~\ref{Tau}: propagator estimate}

We start with the proof of \eqref{Tau.1} for the flux, which is the crucial estimate, and
then turn to the proof of \eqref{Tau.2} for the field.

\medskip

\step1 Proof of \eqref{Tau.1}.

Wlog we may assume that the times $t_1\le T$ are dyadic. 
By telescoping and the semi-group property applied to $S$ and $S^{\rm hom}$, cf Lemma \ref{Laux2},
we have operator identity
\begin{eqnarray*}
S_{t_1\rightarrow T}-S_{t_1\rightarrow T}^{\rm hom}
&=&\sum_{t_1<t\le T}(S_{t\rightarrow T}^{\rm hom}S_{t_1\rightarrow t}
-S_{\frac{t}{2}\rightarrow T}^{\rm hom}S_{t_1\rightarrow\frac{t}{2}})\nonumber\\
&=&
\sum_{t_1<t\le T}S_{t\rightarrow T}^{\rm hom}(S_{\frac{t}{2}\rightarrow t}-S_{\frac{t}{2}\rightarrow t}^{\rm hom})
S_{t_1\rightarrow\frac{t}{2}}.
\end{eqnarray*}
We apply this identity to $q(t_1)$, which by (\ref{au41}) may be reformulated as
\begin{align*}
q(T)-S_{t_1\rightarrow T}^{\rm hom}q(t_1)
=
\sum_{t_1<t\le T}S_{t\rightarrow T}^{\rm hom}(q(t)-S_{\frac{t}{2}\rightarrow t}^{\rm hom}q(\frac{t}{2})).
\end{align*}
By the triangle inequality for $\norm \cdot \norm$ this yields
\begin{equation}\nonumber
\norm q(T)-S^{\ho}_{t_1\rightarrow T}q(t_1)\norm
\leq \sum_{t_1< t\le T}\norm S^{\ho}_{t\rightarrow T}\big(q(t)-S^{\ho}_{\frac{t}{2}\rightarrow t}
q({\textstyle\frac{t}{2}})\big)\norm.
\end{equation}
Lemma \ref{Laux15} now allows us to get rid of $S^{\ho}$ in order to estimate the homogenization error
on the large time interval $(t_1,T)$ by homogenization errors on the dyadic intervals $(\frac{t}{2},t)$:
\begin{equation}\label{Tau1-1}
\norm q(T)-S^{\ho}_{t_1\rightarrow T}q(t_1)\norm
\lesssim \sum_{t_1< t\le T}\norm q(t)-S^{\ho}_{\frac{t}{2}\rightarrow t}q({\textstyle\frac{t}{2}})\norm.
\end{equation}
We now may appeal to Proposition~\ref{Paux2} (with $T$ replaced by $t$) for the summands
\begin{multline}\label{Tau1-2}
\norm q(t)-S_{\frac{t}{2}\to t}^\ho q(\tfrac{t}{2})\norm \,
\lesssim \,(\delta+|a_{ht}-a_\ho|)^\frac1p \sup_{t'\le t}  \norm q(t')\norm
+\frac{1}{\delta^p} \sqrt{t_0}^{\frac{d}{2}\frac{1}{s}} \Big(\frac{1}{\sqrt{t}}\Big)^{\frac{d}{2}(\frac{1}{s}-\frac12)}
\end{multline}
provided $\delta$ and $t_0$ satisfy $\delta \ll 1$ and 
\begin{equation}\label{Tau1-3}
t_0\le \delta^2 t \quad \text{ and }\quad \langle |(q_{t_0}-\langle q_{t_0}\rangle)_{\delta \sqrt{t_0}}| \rangle \leq \delta^p,
\end{equation}
cf \eqref{Paux2.2}.
Here and in the sequel, $p$ denotes a generic large exponent that only depends on $d$, $\lambda>0$ and $s<2$.

\medskip

Let us assume that we already knew that for any $\alpha<\frac{d}{2}$
\begin{equation}\label{Tau1-4}
\norm q(t)\norm \lesssim(\frac{1}{\sqrt{t}})^{\alpha-\frac d2},
\end{equation}
where for this paragraph, $\lesssim$ acquires a dependence on $\alpha<\frac{d}{2}$ which will be chosen later.
We will give the self-contained argument for (\ref{Tau1-4}) in Step~2. We note that the stretched exponential bound (\ref{Tau1-4}) 
implies in particular the first-moment bound
\begin{equation}\label{Tau1-5}
\sup_{r}(\frac{r}{\sqrt{t}})^\frac{d}{2}
\langle|(q(t)-\langle q(t)\rangle)_{r}|\rangle\le(\frac{1}{\sqrt{t}})^{\alpha},
\end{equation}
which we insert into \eqref{Laux10.1} in Lemma~\ref{Laux10} to obtain, provided we momentarily restrict to $0<\alpha<1$,
\begin{equation}\label{Tau1-6}
|a_{ht}-a_{\ho}|\lesssim(\frac{1}{\sqrt{t}})^\alpha.
\end{equation}
Rewriting (\ref{Tau1-5}) as
\begin{equation}\nonumber
r^\frac{d}{2}
\langle|(q(t)-\langle q(t)\rangle)_{r}|\rangle\le\sqrt{t}^{\frac{d}{2}-\alpha}
\quad\mbox{for all}\;r, t
\end{equation}
and appealing to \eqref{Laux4.1} in Lemma~\ref{Laux4} in form of $q_{t_0}=\int_0^\infty\frac{1}{t_0}\exp(-\frac{t}{t_0})q(t)dt$, we see that it yields
\begin{equation}\nonumber
r^\frac{d}{2}
\langle|(q_{t_0}-\langle q_{t_0}\rangle)_{r}|\rangle\le\sqrt{t_0}^{\frac{d}{2}-\alpha}\quad
\mbox{for all}\;r, t_0,
\end{equation}
which we specify to $r=\delta\sqrt{t_0}$:
\begin{equation}\nonumber
\langle|(q_{t_0}-\langle q_{t_0}\rangle)_{\delta\sqrt{t_0}}|
\rangle\le(\frac{1}{\delta})^\frac{d}{2}(\frac{1}{\sqrt{t_0}})^\alpha.
\end{equation}
Hence the proviso (\ref{Tau1-3}) holds once we have $\frac{1}{\delta^{\frac1\alpha(p+\frac d 2)}}\leq\sqrt{t_0}\le\delta\sqrt{t}$; we choose $t_0$ to be equal to the lower bound. Hence (up to redefining $p$)
we obtain from (\ref{Tau1-2}) and (\ref{Tau1-6}) that for all $\delta\gg(\frac{1}{\sqrt{t}})^\frac{1}{p}$
\begin{align}\label{Tau1-8}
\norm q(t)-S^{\ho}_{\frac{t}{2}\rightarrow t}q({\textstyle\frac{t}{2}})\norm 
\, \lesssim \, \delta^{\frac{1}{p}}
\sup_{t'\le t}\norm q(t')\norm
+\frac{1}{\delta^p}(\frac{1}{\sqrt{t}})^{\frac{d}{2}(\frac{1}{s}-\frac12)}.
\end{align}
We now make use of the triangle inequality, once again Lemma~\ref{Laux15}, and our assumption (\ref{Tau1-4}) to obtain
for $t'\ge t_1$
\begin{align*}
\norm q(t')\norm 
&\le\, \norm q(t')-S^{\ho}_{t_1\rightarrow t'}q(t_1)\norm
+\norm S^{\ho}_{t_1\rightarrow t'}q(t_1)\norm \\
&\lesssim  \, \norm q(t')-S^{\ho}_{t_1\rightarrow t'}q(t_1)\norm
+\norm q(t_1)\norm 
\\
&\lesssim\, \norm q(t')-S^{\ho}_{t_1\rightarrow t'}q(t_1)\norm
+\sqrt{t_1}^{\frac{d}{2}-\alpha},
\end{align*}
whereas for $t'\le t_1$, $\norm q(t')\norm\,\lesssim\, \sqrt{t'}^{\frac{d}{2}-\alpha}
\le\sqrt{t_1}^{\frac{d}{2}-\alpha}$. These two estimates combine to
\begin{align*}
\sup_{t'\le t} \norm q(t')\norm
&\lesssim\, \sup_{t_1\le t'\le t}\norm q(t')-S^{\ho}_{t_1\rightarrow t'}q(t_1)\norm+\sqrt{t_1}^{\frac{d}{2}-\alpha},
\end{align*}
which we introduce into (\ref{Tau1-8})
\begin{align}\nonumber
{\norm q(t)-S^{\ho}_{\frac{t}{2}\rightarrow t}q({\textstyle\frac{t}{2}})\norm}\lesssim \delta^{\frac{1}{p}}\big(
\sup_{t_1\le t'\le t}\norm q(t')-S^{\ho}_{t_1\rightarrow t'}q(t_1)\norm
+\sqrt{t_1}^{\frac{d}{2}-\alpha}\big)
+\frac{1}{\delta^p}(\frac{1}{\sqrt{t}})^{\frac{d}{2}(\frac{1}{s}-\frac12)}.\nonumber
\end{align}
We now make the Ansatz $\delta=(\frac{1}{\sqrt{t}})^\frac{1}{p_0}$ with an exponent $p_0$ to be chosen below
(and which is admissible provided $p_0\gg 1$) so that the above turns into 
\begin{align}\nonumber
\lefteqn{\norm q(t)-S^{\ho}_{\frac{t}{2}\rightarrow t}q({\textstyle\frac{t}{2}})\norm}\nonumber\\
&\lesssim(\frac{1}{\sqrt{t}})^{\frac{1}{pp_0}}\big(
\sup_{t_1\le t'\le t}\norm q(t')-S^{\ho}_{t_1\rightarrow t'}q(t_1)\norm
+\sqrt{t_1}^{\frac{d}{2}-\alpha}\big)
+(\frac{1}{\sqrt{t}})^{\frac{d}{2}(\frac{1}{s}-\frac12)-\frac{p}{p_0}}.\nonumber
\end{align}
Inserting this into (\ref{Tau1-1}) we obtain, provided $p_0$ is chosen so large that $\frac{d}{2}(\frac{1}{s}-\frac12)-\frac{p}{p_0}>0$
(which is possible since $s<2$)
\begin{align}\nonumber
\lefteqn{\norm q(T)-S^{\ho}_{t_1\rightarrow T}q(t_1)\norm }\nonumber\\
&\lesssim (\frac{1}{\sqrt{t_1}})^{\frac{1}{pp_0}}\big(
\sup_{t_1\le t\le T}\norm q(t)-S^{\ho}_{t_1\rightarrow t}q(t_1)\norm
+\sqrt{t_1}^{\frac{d}{2}-\alpha}\big)
+(\frac{1}{\sqrt{t_1}})^{\frac{d}{2}(\frac{1}{s}-\frac12)-\frac{p}{p_0}}.\nonumber
\end{align}
Provided $t_1\gg 1$, we may absorb the first rhs term to the effect of
\begin{align}\nonumber
\norm q(T)-S^{\ho}_{t_1\rightarrow T}q(t_1)\norm 
\, \lesssim\, (\frac{1}{\sqrt{t_1}})^{\alpha-\frac{d}{2}-\frac{1}{pp_0}}
+(\frac{1}{\sqrt{t_1}})^{\frac{d}{2}(\frac{1}{s}-\frac12)-\frac{p}{p_0}}.\nonumber
\end{align}
Having fixed $p_0$ above, it remains to choose $\alpha<\frac{d}{2}$ so close to $\frac{d}{2}$ that
$\alpha-\frac{d}{2}-\frac{1}{pp_0}>0$.

\medskip

\step2 Proof of  (\ref{Tau1-4}).

We give the argument for (\ref{Tau1-4}) for which we may even assume $s=2$ (cf Remark~\ref{rem:nearly}). For the proof of (\ref{Tau1-4}), we use
(\ref{Tau1-1}) (in addition to a further application of the triangle inequality and Lemma~\ref{Laux15})
with the semi-group $S^{\ho}$ replaced by the propagator $S^h$:
\begin{equation}\nonumber
\norm q(T)\norm
\, \lesssim \, \norm q(t_1)\norm
+\sum_{t_1< t\le T}\norm q(t)-S^{h}_{\frac{t}{2}\rightarrow t}q({\textstyle\frac{t}{2}})\norm.
\end{equation}
Likewise, by Proposition~\ref{Paux2} we have (\ref{Tau1-2}) in form of
\begin{align}\nonumber
\norm q(t)-S^{h}_{\frac{t}{2}\rightarrow t}q({\textstyle\frac{t}{2}})\norm\,
\lesssim\,  \delta^{\frac{1}{p}}
\sup_{t'\le t}\norm q(t')\norm
+\frac{\sqrt{t_0}^p}{\delta^p},\nonumber
\end{align}
provided $\delta$ and $t_0$ are chosen such that (\ref{Tau1-3}) holds for all $t\ge t_1$.
With help of Lemma~\ref{Laux16}, these two estimates turn into
\begin{equation}\nonumber
\norm q(T)\norm 
\lesssim \sqrt{t_1}^\frac{d}{2}
+\sum_{t_1< t\le T}\norm q(t)-S^{h}_{\frac{t}{2}\rightarrow t}q({\textstyle\frac{t}{2}})\norm
\end{equation}
and
\begin{align}\nonumber
\norm q(t)-S^{h}_{\frac{t}{2}\rightarrow t}q({\textstyle\frac{t}{2}})\norm
\lesssim \delta^{\frac{1}{p}}\big(
\sup_{t_1\le t'\le t}\norm q(t')\norm +\sqrt{t_1}^\frac{d}{2}\big)
+\frac{\sqrt{t_0}^p}{\delta^p},\nonumber
\end{align}
respectively. Introducing the sub-CLT exponent $\alpha<\frac{d}{2}$, we rewrite these estimates as
\begin{equation}\nonumber
\sqrt{T}^{\alpha-\frac d2}\norm q(T)\norm 
\lesssim \sqrt{t_1}^{\alpha}
+\sum_{t_1< t\le T}(\frac{\sqrt{t}}{\sqrt{T}})^{\frac{d}{2}-\alpha} \sqrt{t}^{\alpha-\frac d2}
\norm q(t)-S^{h}_{\frac{t}{2}\rightarrow t}q({\textstyle\frac{t}{2}})\norm
\end{equation}
and (using $t\ge t_1\ge 1$)
\begin{align}\nonumber
\sqrt{t}^{\alpha-\frac d2}\norm q(t)-S^{h}_{\frac{t}{2}\rightarrow t}q({\textstyle\frac{t}{2}})\norm\,
\lesssim\, \delta^{\frac{1}{p}}\big(
\sup_{t_1\le t'\le t}\sqrt{t'}^{\alpha-\frac d2}\norm q(t')\norm+\sqrt{t_1}^{\alpha}\big)
+\frac{\sqrt{t_0}^p}{\delta^p}.\nonumber
\end{align}
Inserting the second estimate into the first, we obtain thanks to $\alpha<\frac{d}{2}$
\begin{equation}\label{Tau1-7}
\sqrt{T}^{\alpha-\frac d2}\norm q(T)\norm
\lesssim\delta^{\frac{1}{p}}\sup_{t_1\le t\le T}\sqrt{t}^{\alpha-\frac d2}\norm q(t)\norm+\sqrt{t_1}^{\alpha}
+\frac{\sqrt{t_0}^p}{\delta^p}.\nonumber
\end{equation}
Clearly, by monotonicity of the expression $\sup_{t_1\le t\le T}\sqrt{t}^{\alpha-\frac{d}{2}}\norm q(t)\norm$
in $T$, the lhs of this estimate holding for all $T\ge t_1$ may be replaced by 
$\sup_{t_1\le t\le T}\sqrt{t}^{\alpha-\frac{d}{2}}\norm q(t)\norm$.
We thus may first fix $\delta$ so small that the first rhs term may be absorbed into the lhs.
By the uniform qualitative homogenization established in Lemma \ref{Laux17}, we then may fix a $t_0$ 
so large that the second condition in (\ref{Tau1-3}) is satisfied. Finally, we fix $t_1$ such that
$t_1=\frac{1}{\delta^2}t_0$ so that also the first condition in (\ref{Tau1-3}) is satisfied
for all $t\ge t_1$.

\medskip

\step3 Proof of \eqref{Tau.2}.

In preparation for treating fields like the fluxes, we recall that the main ingredient in 
treating the fluxes was the estimate of $q(t)-S^{\rm hom}_{\frac{t}{2}\rightarrow t}q(\frac{t}{2})=\int_\frac{t}{2}^ta\nabla ud\tau
-\int_\frac{t}{2}^ta_{\rm hom}\nabla v_{\rm hom}d\tau $ provided by Proposition~\ref{Paux2}. An inspection of
the proof of that proposition (leading back to Proposition~\ref{Paux1}) reveals that the difference between the
variable-coefficient solution $u$ and the constant-coefficient solution $v_{\rm hom}$ is not only estimated
on the level of fluxes, but also on the level of fields $\int_\frac{t}{2}^t\nabla ud\tau
-\int_\frac{t}{2}^t\nabla v_{\rm hom}d\tau$, which we may multiply with $a_{\rm hom}$ to obtain
control of $\int_\frac{t}{2}^ta_{hom}\nabla ud\tau-\int_\frac{t}{2}^ta_{hom}\nabla v_{\rm hom}d\tau=
a_{\rm hom}\nabla\phi(t)-a_{\rm hom}\nabla\phi(\frac{t}{2})+(1-S^{\rm hom}_{\frac{t}{2}\rightarrow t})q(\frac{t}{2})$.
Hence the analogue of (\ref{Paux2.1}) on the level of fields reads
\begin{equation}\label{Tau1-9}
\norm a_{\rm hom}\nabla\phi(t)-a_{\rm hom}\nabla\phi(\frac{t}{2})+(1-S^{\rm hom}_{\frac{t}{2}\rightarrow t})q(\frac{t}{2})\norm
\, \lesssim \,\delta^\frac1p \sup_{t\le T}  \norm q(t)\norm
\\
+\frac{1}{\delta^p} \sqrt{t_0}^{\frac{d}{2}\frac{1}{s}} \Big(\frac{1}{\sqrt{T}}\Big)^{\frac{d}{2}(\frac{1}{s}-\frac12)}
\end{equation}
provided $\delta$ and $t_0$ satisfy $0<\delta \ll 1$, $t_0\le \delta^2 T$, and \eqref{Paux2.2}.

\medskip

Starting from telescoping in form of 
\begin{align*}
1-S^{\rm hom}_{t_1\rightarrow T}=\sum_{t_1<t\le T}\big(
(1-S_{\frac{t}{2}\rightarrow T}^{\rm hom})S_{t_1\rightarrow\frac{t}{2}}
-(1-S_{t\rightarrow T}^{\rm hom})S_{t_1\rightarrow t}\big),
\end{align*}
we obtain from the semi-group properties for $S$ and $S^{\rm hom}$ 
in form of
\begin{align*}
\lefteqn{(1-S_{\frac{t}{2}\rightarrow T}^{\rm hom})S_{t_1\rightarrow\frac{t}{2}}
-(1-S_{t\rightarrow T}^{\rm hom})S_{t_1\rightarrow t}}\nonumber\\
&=\big((1-S_{t\rightarrow T}^{\rm hom}S_{\frac{t}{2}\rightarrow t}^{\rm hom})
-(1-S_{t\rightarrow T}^{\rm hom})S_{\frac{t}{2}\rightarrow t}\big)S_{t_1\rightarrow\frac{t}{2}}
\end{align*}
the following operator identity 
\begin{align*}
\lefteqn{1-S^{\rm hom}_{t_1\rightarrow T}
}\nonumber\\
&=-\sum_{t_1<t\le T}(1-S_{t\rightarrow T}^{\rm hom})(S_{\frac{t}{2}\rightarrow t}-
S_{\frac{t}{2}\rightarrow t}^{\rm hom})S_{t_1\rightarrow \frac{t}{2}}
+\sum_{t_1<t\le T}(1-S_{\frac{t}{2}\rightarrow t}^{\rm hom})S_{t_1\rightarrow \frac{t}{2}},
\end{align*}
which we apply to $q(t_1)$ to the effect of
\begin{equation}\nonumber
(1-S^{\rm hom}_{t_1\rightarrow T})q(t_1)=-\sum_{t_1<t\le T}(1-S_{t\rightarrow T}^{\rm hom})(S_{\frac{t}{2}\rightarrow t}-
S_{\frac{t}{2}\rightarrow t}^{\rm hom})q(\frac{t}{2})
+\sum_{t_1<t\le T}(1-S_{\frac{t}{2}\rightarrow t}^{\rm hom})q(\frac{t}{2}).
\end{equation}
Together with the trivial telescoping $a_{\rm hom}\nabla\phi(T)-a_{\rm hom}\nabla\phi(t_1)
=\sum_{t_1< t\le T}(a_{\rm hom}\nabla\phi(t)-a_{\rm hom}\nabla\phi(\frac{t}{2}))$, this yields
\begin{align}\nonumber
\lefteqn{a_{\rm hom}\nabla\phi(T)-a_{\rm hom}\nabla\phi(t_1)+(1-S^{\rm hom}_{t_1\rightarrow T})q(t_1)}\nonumber\\
&=-\sum_{t_1<t\le T}(1-S_{t\rightarrow T}^{\rm hom})(S_{\frac{t}{2}\rightarrow t}-
S_{\frac{t}{2}\rightarrow t}^{\rm hom})q(\frac{t}{2})\nonumber\\
&+\sum_{t_1<t\le T}\big(a_{\rm hom}\nabla\phi(t)-a_{\rm hom}\nabla\phi(\frac{t}{2})+
(1-S_{\frac{t}{2}\rightarrow t}^{\rm hom})q(\frac{t}{2})).\nonumber
\end{align}
The estimate of the first rhs term proceed as in the first step, with the only change that at the beginning,
we use the boundedness of the operator $1-S_{t\rightarrow T}^{\rm hom}$ instead of the one of $S_{t\rightarrow T}^{\rm hom}$.
The estimate of the second rhs term also proceeds as in the first step, now based on (\ref{Tau1-9}).


\subsection{Proof of Corollary~\ref{cor:CLTdecay}: CLT scaling}

For $(q,\nabla \phi)$, this is a direct consequence of Theorem~\ref{Tau}.
Indeed, by the triangle inequality and \eqref{Tau.1} \& \eqref{Tau.2} for $T=\infty$,
\begin{eqnarray*}
\norm q \norm & \lesssim & \norm S^{\ho}_{1\rightarrow \infty }q(1) \norm +1 ,
\\
\norm  \nabla \phi \norm &\lesssim & \norm \nabla \phi(1)\norm +\norm S^{\ho}_{1\rightarrow T}q(1)\norm +\norm q(1)\norm +1,
\end{eqnarray*}
so that the claim follows from Lemmas~\ref{Laux15} and~\ref{Laux16}.
Likewise, this yields for all $t\ge 1$
\begin{equation}\label{e.cltdecay-1}
\norm q(t) \norm, \norm \nabla \phi(t)\norm \,\lesssim \,1 .
\end{equation}

\medskip

We now turn to the gradient $\nabla\sigma$ of the vector potential $\sigma$. Based on (\ref{o58}) we will show
\begin{align}\label{gr3}
\norm \nabla\sigma\norm\lesssim\norm q\norm.
\end{align}
As in the beginning of the proof of Lemma \ref{Laux15} below, this can be reduced to
\begin{align*}
\|(\nabla\sigma)_R\|\lesssim\int_{\frac{R}{2}}^\infty\frac{dr}{r}\| q_r\|_*.
\end{align*}
In order to proceed as in the proof of Lemma \ref{Laux15}, we disintegrate $\nabla\sigma$ by (formally)
writing $\nabla\sigma=\int_0^\infty\nabla vd\tau$, where $v$ solves the constant-coefficient initial value problem
with divergence-form initial data
\begin{align}\label{gr2}
\partial_\tau v_{jk}-\triangle v_{jk}=0,\quad v_{jk}(\tau=0)=\partial_jq_k-\partial_kq_j=\nabla\cdot(q_ke_j-q_je_k),
\end{align}
so that by the triangle inequality, it suffices to show
\begin{align*}
\int_0^\infty\|(\nabla v_{ij})_R(\tau)\|d\tau\lesssim\int_{\frac{R}{2}}^\infty\frac{dr}{r}\|(q_ke_j-q_je_k)_r\|_*
\le\int_{\frac{R}{2}}^\infty\frac{dr}{r}\|q_r\|_*.
\end{align*}
This is exactly the estimate \eqref{L15-10} established in the proo of Lemma \ref{Laux15}
with the only difference that the divergence-form initial data $q$ are replaced by $q_ke_j-q_je_k$, 
the constant coefficient $a_{\rm hom}$ is replaced by ${\rm id}$, and 
the initial time $\tau=t$ is replaced by $\tau=0$. 

\medskip

We finally turn to the estimate of the modified augmented corrector $(\nabla\phi_T,\nabla\sigma_T,q_T)$.
For $q_T$ we appeal to the representation (\ref{Laux4.1}) in terms of an average of $\{q(t)\}_t$,
so that $\norm q_T\norm\le\sup_{t}\norm q(t)\norm$. For $t\ge 1$, we estimated $\norm q(t)\norm$ 
above, cf (\ref{e.cltdecay-1}). For $t\le 1$, we directly appeal to Lemma \ref{Laux16}.
The field $\nabla\phi_T$ can be handled along identical lines. For the gradient $\nabla\sigma$ of the vector potential
of the flux, we have the representation
$\nabla\sigma_T=\int_0^\infty\exp(-\frac{\tau}{T})\nabla v(\tau)d\tau$, where $v$ is defined as in (\ref{gr2})
with the divergence-form initial data $q_ke_j-q_je_k$ replaced by $(q_T)_ke_j-(q_T)_je_k$. Hence we obtain
$\norm \nabla\sigma_T\norm\lesssim\norm q_T\norm$ along the same lines the above estimate (\ref{gr3}).


\subsection{Proof of Theorem~\ref{Tau2}: fluctuations of the homogenization commutator}

We split this self-contained proof into eleven steps.

\medskip 

\step1 Random fields with bounded triple norm are almost-surely Schwartz distributions of order
$2[\frac{d}{4}]+1$. 

More precisely we claim for $m\in\mathbb{N}$ with $2m>\frac{d}{2}$ and any 
centered and stationary random field $F$ that
\begin{align}\label{gr10}
\Big\| \sup\Big\{\int\zeta F\Big|
\sum_{k=1}^{2m}\sup_{x\in\mathbb{R}^d}(|x|+1)^{4m}|\nabla^k\zeta|\le 1\Big\}\Big\|
\lesssim\norm F\norm.
\end{align}
Note that the norm on $\zeta$ above is the one appearing in (\ref{gr9}).
The argument for (\ref{gr10}) is mostly deterministic; by definition of mollification with
the heat kernel we have $\partial_t^k\zeta_{\sqrt{t}}=\triangle^k\zeta_{\sqrt{t}}$ so that
by symmetry of the convolution we get
\begin{align*}
(\frac{d}{dt})^k\int\zeta F_{\sqrt{t}}=\int \triangle^k\zeta F_{\sqrt{t}}
\end{align*}
and thus obtain the estimate
\begin{align}\label{gr11}
|(\frac{d}{dt})^k\int\zeta F_{\sqrt{t}}|\le
\sup_{x}\big((|x|+1)^{4m}|\triangle^k\zeta|\big)\,\int(|x|+1)^{-4m}|F_{\sqrt{t}}|dx.
\end{align}
We now appeal to the weighted Sobolev estimate (in the single variable $t$)
\begin{align*}
\sup_{t\in[0,1]}|f(t)|\lesssim\sum_{k=0}^m\int_0^1 t^\alpha|\frac{d^kf}{dt^k}|dt\quad
\mbox{provided}\;\alpha+1<m,
\end{align*}
which we apply to $f(t)=\int\zeta F_{\sqrt{t}}$ and into which we insert (\ref{gr11})
to the effect of
\begin{align*}
\sup\big\{|\int\zeta F|\big|\sum_{k=0}^m\sup_{x}(|x|+1)^{4m}|\triangle^k\zeta|\le 1\big\}
\lesssim \int_0^1 t^\alpha\int(|x|+1)^{-4m}|F_{\sqrt{t}}|dxdt.
\end{align*}
By stationarity and definition of $\norm\cdot\norm$ this yields
\begin{align*}
\lefteqn{\big\|\sup\big\{|\int\zeta F|\big|
\sum_{k=0}^m\sup_{x}(|x|+1)^{4m}|\triangle^k\zeta|\le 1\big\}\big\|}\nonumber\\
&\lesssim \int_0^1 t^{\alpha-\frac{d}{4}} dt\norm F\norm\quad\mbox{provided}\;4m>d.
\end{align*}
Hence provided $\alpha+1>\frac{d}{4}$ this turns into (\ref{gr10}); we thus recover the condition
$m>\alpha+1>\frac{d}{4}$.

\medskip

\step2 Existence of ``blow-downs'' by weak compactness. 

Let $F$ be a centered and stationary random field
with $\norm F\norm<\infty$. For any radius $R$, we consider its rescaled version $F^R$ given
through $F^R(\hat x):=R^\frac{d}{2}F(R\hat x)$, still a centered and stationary (in the broad sense that
translations of the field have the same finite-dimensional law)
random field. We claim that for any sequence $R\uparrow\infty$, there exists a subsequence
(still denoted by $R\uparrow\infty$) and a centered and stationary random Schwartz distribution $F^\infty$
(defined on an extended probability space) with
\begin{align}\label{gr18}
\|\sup_\zeta\big\{\int\zeta F^\infty\big|\sum_{k=0}^{2m}\sup_{x}(|x|+1)^{4m}|\nabla^k\zeta_n|\le 1\}\|
\lesssim\norm F\norm
\end{align}
(where we use with a slight abuse of notation $\int\zeta F^\infty$ as the symbol for 
the Schwartz distribution $F^\infty$ applied to the Schwartz function $\zeta$)
and which thus (almost-surely) is of order $2[\frac{d}{4}]+1$, such that
\begin{align}\label{gr17}
\lefteqn{\{\zeta\mapsto\int\zeta F^R\}\stackrel{R\uparrow\infty}{\rightarrow}
\{\zeta\mapsto\int\zeta F^\infty\}}\nonumber\\
&\mbox{in law up to near-Gaussian moments},
\end{align}
where $\zeta$ runs over all functions that may be approximated by Schwartz functions
in the norm (\ref{gr9}). An inspection of Step 1 shows that in (\ref{gr10}), we may replace the norm 
$\sum_{k=1}^{2m}\sup_{x}(|x|+1)^{4m}|\nabla^k\zeta|$ by the norm with weaker weight
$\sum_{k=1}^{2m}\sup_{x}(|x|+1)^{4m-}|\nabla^k\zeta|$, so that the closure of Schwartz functions
in this weaker norm contains all $2m$-times continuously differentiable functions for which
the stronger norm is finite. Hence we may assume that $\zeta$ in (\ref{gr17}) runs over this
larger space.

\medskip

We first note that for any radius $R$, the rescaling 
$F^R(\hat x)=R^\frac{d}{2}F(R\hat x)$ interacts with our mollification through
convolution via $\rho^\frac{d}{2}(F^R)_\rho(\hat x)=(\rho R)^\frac{d}{2}F_{\rho R}(R\hat x)$ 
and thus by stationarity does not increase the triple norm $\norm F^R\norm\le\norm F\norm$.
Let $\mathcal S_m$ denote the closure of the space $\mathcal S$ of Schwartz functions under the norm (\ref{gr9}).
This is a Banach space, the norm of which we denote by $\|\cdot\|_{\mathcal S_m}$.  
Let $(\mathcal T_m, \|\cdot\|_{\mathcal T_m})$ be a separable Banach space that compactly embeds into $\mathcal S_m$ (say a similar space with larger weight and more derivatives).
Step~1 yields
\begin{equation}\label{Tau2-tight}
\|\sup_{\zeta\in \mathcal S, \|\zeta\|_{\mathcal T_{m}} \le 1} \int\zeta F^R\| \lesssim \|\sup_{\zeta\in \mathcal S, \|\zeta\|_{\mathcal S_{m}} \le 1} \int\zeta F^R\| \lesssim\norm F\norm.
\end{equation}
%
Let $\mathcal S'_{m} \subset \mathcal T'_{m}$ be the dual spaces of the Banach spaces $\mathcal T_m\subset \mathcal S_{m}$.
We consider the random linear form $\zeta\mapsto \int\zeta F^R$ as a random element on $\mathcal T'_{m}$.
Since $\mathcal S_{m}'$ compactly embeds into $\mathcal T_{m}'$, \eqref{Tau2-tight} yields the tightness of the laws of $\zeta\mapsto \int\zeta F^R$
in $\mathcal T_{m}'$.
By tightness (see e.g. \cite[Theorem~5.1 p.~59]{Billingsley99}), we may select a subsequence
of $R\uparrow\infty$ in such a way that these laws on $\mathcal T_{m}'$ converge
to some law on $\mathcal T_{m}'$ with near-Gaussian moments, where convergence
means testing against functions of near-Gaussian growth.
We denote by $F^\infty$ a random variable on $\mathcal T_m'$ distributed according to that law.
By continuity, the limit $F^\infty$ still satisfies the bound
$$
 \|\sup_{\zeta\in \mathcal S, \|\zeta\|_{\mathcal S_{m}} \le 1} \int\zeta F^\infty\| \lesssim\norm F\norm,
$$
so that 
$\zeta\mapsto\int\zeta F^\infty$ is in particular a random Schwartz distribution of order $2[\frac{d}{4}]+1$; 
it inherits the centeredness and stationarity (in the large sense) from $F^R$. The convergence
(\ref{gr17}) follows by construction.

\medskip

\step3 Relation between blow-downs. 

We make two claims
\begin{itemize}
\item[a)] For any cut-off time $t$, the blow-downs of $\Xi(t)$ and of its mollification $\Xi_{\sqrt{t}}(t)$
are identical. This will follow from (\ref{theo8b}) in Corollary \ref{cor:CLTdecay} 
in form of $\norm\Xi(t)\norm\lesssim 1$.
\item[b)] As $t\uparrow\infty$, the blow-downs of $\Xi(t)$ converge to the blow-downs of $\Xi$ in the
sense of (\ref{gr17}).
This will follow from Theorem \ref{Tau} in form of 
$\norm \Xi(t)-\Xi\norm\stackrel{t\uparrow\infty}{\rightarrow}0$, cf (\ref{gr5}).
\end{itemize}
The first claim follows from the formula $\int\zeta(\Xi_{\sqrt{t}})^R(t)=\int\zeta_{\frac{\sqrt{t}}{R}}\Xi^R(t)$,
the fact that  
$\zeta_{\frac{\sqrt{t}}{R}}\stackrel{R\uparrow\infty}{\rightarrow}\zeta$ in
the norm (\ref{gr9}) for any Schwartz function $\zeta$, 
which together with Step~1 yields $\|\int\zeta_\frac{\sqrt{t}}{R}\Xi^R(t)
-\int\zeta\Xi^R(t)\|$ $\le o(1)\norm\Xi^R(t)\norm$ $\le o(1)\norm\Xi(t)\norm$ $\le o(1)$
as $R\uparrow\infty$. We now turn to the second claim, which reads in detail: Suppose
that there exists a sequence $t\uparrow\infty$ and a (single) sequence $R\uparrow\infty$ along which
every $\Xi(t)$ admits a blow-down $\Xi^\infty(t)$. Then we claim that the blow-downs $\Xi^\infty(t)$ converge
in the sense of (\ref{gr17}) to a centered and stationary random Schwartz distribution $\Xi^\infty$,
which at the same time is the blow-down of $\Xi$ along the sequence $R\uparrow\infty$.
Since by Step 1 we have 
$\|\sup_\zeta\int\zeta(\Xi^R(t)-\Xi^R(t'))\|$ $\lesssim\norm\Xi^R(t)-\Xi^R(t')\norm$ 
$=\norm\Xi(t)-\Xi(t')\norm$ for any two cut-off times $t$ and $t'$
(where here, the sup runs over all $\zeta$ for which the norm (\ref{gr9}) is $\le 1$), we obtain from
(\ref{gr5}) the Cauchy sequence property
\begin{align}\label{gr19}
\lim_{t,t'\uparrow\infty}\|\sup_\zeta\int\zeta(\Xi^R(t)-\Xi^R(t'))\|=0\quad\mbox{uniformly in}\;R.
\end{align}
Because of the uniformity, it is preserved in the limit $R\uparrow\infty$ under our notion of
convergence
(\ref{gr17}): $\lim_{t,t'\uparrow\infty}\|\sup_\zeta\int\zeta(\Xi^\infty(t)-\Xi^\infty(t'))\|=0$, 
where for convenience we pass to a large probability space on which the random fields $\{\Xi^\infty(t)\}_t$ live.
The compactness result argument in Step 2 can also be used as a completeness argument yielding
the existence of a centered and stationary random Schwartz distribution $\Xi^\infty$ with
$\lim_{t\uparrow\infty}\|\sup_\zeta\int\zeta(\Xi^\infty(t)-\Xi^\infty)\|=0$. Completely analogously
to (\ref{gr19}) we have $\lim_{t\uparrow\infty}\|\sup_\zeta\int\zeta(\Xi^R(t)-\Xi^R)\|=0$
uniformly in $R$, so that the convergence of $\Xi^R(t)$ to $\Xi^\infty(t)$ 
in the sense of (\ref{gr17}) along the subsequence
$R\uparrow\infty$ extends as desired to the convergence of $\Xi^R$ to $\Xi^\infty$.

\medskip

\step4 From approximate to exact locality. 

We claim that $\Xi_{\sqrt{t}}(t)$, which is approximately
local on scale $\sqrt{t}$,
may be approximated by $\langle\Xi_{\sqrt{t}}(t)|B_r\rangle$, where $\langle\cdot|D\rangle$ denotes
the expectation conditioned on the restriction $a_{|D}$ of $a$ on some open set $D\subset\mathbb{R}^d$
(see the discussion in Step 6 of the proof of Lemma \ref{Laux12} on conditional expectations in our
context) so that $\langle\Xi_{\sqrt{t}}|B_r\rangle$ is exactly local on scale $r$, as $r\uparrow\infty$,
in the triple norm (but non-uniformly in $t$):
\begin{align}\label{gr25}
\norm\Xi_{\sqrt{t}}(t)-\langle\Xi_{\sqrt{t}}(t)|B_r\rangle\norm
\lesssim\sqrt{t}^\frac{d}{2}\exp(-\frac{1}{C}\frac{r}{\sqrt{t}}).
\end{align}
The main ingredient for (\ref{gr25}) is \eqref{Laux14.0b} in Lemma \ref{Laux14}, which states uniform bounds and approximate locality
(in absolute terms) for $q_{\sqrt{t}}(t)$ and $\nabla\phi_{\sqrt{t}}(t)$, 
which translates to $\Xi_{\sqrt{t}}(t)$ by the triangle inequality:
\begin{align*}
|\Xi_{\sqrt{t}}(t)|&\lesssim 1,\\
|\Xi_{\sqrt{t}}(a,t)-\Xi_{\sqrt{t}}(\tilde a,t)|&\lesssim\exp(-\frac{1}{C}\frac{r}{\sqrt{t}})
\quad\mbox{provided}\;a=\tilde a\;\mbox{on}\;B_r\;\mbox{for}\;r\ge\sqrt{t}.
\end{align*}
As in Step 6 of the proof of Lemma \ref{Laux12}, we reformulate this in terms of conditional expectations:
\begin{align}
|\Xi_{\sqrt{t}}(t)|,|\langle\Xi_{\sqrt{t}}(t)|B_r\rangle|&\lesssim 1,\label{gr30}\\
|\Xi_{\sqrt{t}}(t)-\langle\Xi_{\sqrt{t}}(t)|B_r\rangle|&\lesssim\exp(-\frac{1}{C}\frac{r}{\sqrt{t}}),\label{gr29}
\end{align}
where the second estimate at first only holds for $r\ge\sqrt{t}$ but then by the first
estimate trivially holds for all $r$. As in Step 6 of the proof of Lemma \ref{Laux12}
we consider dyadic differences:
\begin{align}\label{gr23bis}
|\langle\Xi_{\sqrt{t}}(t)|B_{2r}\rangle-\langle\Xi_{\sqrt{t}}(t)|B_r\rangle|
\lesssim\exp(-\frac{1}{C}\frac{r}{\sqrt{t}}).
\end{align}
By the exact locality of
$\langle\Xi_{\sqrt{t}}(t)|B_{2r}\rangle-\langle\Xi_{\sqrt{t}}(t)|B_r\rangle$
on scale $2r$, we obtain like in Step~3 of the proof of Lemma \ref{Laux12} for $R\ge r$ (recall $r\ge \sqrt{t}\gg1$ wlog) that
\begin{align}\label{gr24}
\|(\langle\Xi_{\sqrt{t}}(t)|B_{2r}\rangle-\langle\Xi_{\sqrt{t}}(t)|B_r\rangle)_R\|
\lesssim(\frac{\sqrt{r}}{R})^\frac{d}{2}
\|\langle\Xi_{\sqrt{t}}(t)|B_{2r}\rangle-\langle\Xi_{\sqrt{t}}(t)|B_r\rangle\|;
\end{align}
however for $R\le r$ this estimate is trivial because of the monotonicity
of mollification seen under $\|\cdot\|$. 
Together (\ref{gr23bis}) and (\ref{gr24}) yield
$\norm \langle\Xi_{\sqrt{t}}(t)|B_{2r}\rangle-\langle\Xi_{\sqrt{t}}(t)|B_r\rangle\norm
\lesssim\sqrt{t}^\frac{d}{2}
\exp(-\frac{1}{C}\frac{r}{\sqrt{t}})$ with a change in the value of $C$ by which we absorb the
algebraic factor $(\frac{r}{\sqrt{t}})^\frac{d}{2}$ into the exponential.
The exponential expression on the rhs allows dyadic summation in $r$, which yields (\ref{gr25}).

\medskip

\step5 White-noise property of blow-downs of $\Xi(t)$. 

We claim that a blow-down 
$\Xi^\infty(t)$ of $\Xi(t)$ along any sequence $R\uparrow\infty$ necessarily has the
white-noise property
\begin{align}\label{gr26}
\int\zeta: \Xi^\infty(t)\;\mbox{and}\;\int\zeta':\Xi^\infty(t)
\quad\mbox{are independent for}\quad{\rm supp}\zeta\cap{\rm supp}\zeta'=\emptyset
\end{align}
for a pair of $\mathbb{R}^{d\times d}$-valued Schwartz functions $\zeta$ and $\zeta'$ with
compact support. 
For notational simplicity, we pass to scalar instead of tensorial notation and language.
Let us momentarily introduce the notation $F(r):=\langle\Xi_{\sqrt{t}}(t)|B_r\rangle$;
we note that by Jensen $\norm F(r)\norm\le\norm\Xi_{\sqrt{t}}(t)\norm\le\norm\Xi(t)\norm\lesssim 1$.
We fix a (discrete) sequence $r\uparrow\infty$; by Step 2 and Cantor, we may select a subsequence 
(without changing notation) of our given
sequence $R\uparrow\infty$ along which $F(r)$ admits a blow-down $F^\infty(r)$ 
for {\it all} $r$ in our fixed sequence. After rescaling, (\ref{gr25}) in Step 4 turns into
$\norm\Xi^R_{\sqrt{t}}(t)-F^R(r)\norm$ $\lesssim\sqrt{t}^\frac{d}{2}\exp(-\frac{r}{\sqrt{t}})$,
so that by Step 1 we have $\|\sup_\zeta\int\zeta(\Xi^R_{\sqrt{t}}(t)-F^R(r))\|$
$\lesssim\sqrt{t}^\frac{d}{2}\exp(-\frac{r}{\sqrt{t}})$, where the sup runs over all $\zeta$'s
for which the norm (\ref{gr9}) is $\le 1$. This estimate is preserved under the convergence
(\ref{gr17}) so that 
$\|\sup_\zeta\int\zeta(\Xi^\infty(t)-F^\infty(r))\|$
$\lesssim\sqrt{t}^\frac{d}{2}\exp(-\frac{r}{\sqrt{t}})$,
where we used that according to Step 3 a), $\Xi_{\sqrt{t}}(t)$ and $\Xi(t)$ have the same blow-downs.
In particular, we have analogously to (\ref{gr17})
\begin{align*}
\lefteqn{\{\zeta\mapsto\int\zeta F^\infty(r)\}\stackrel{r\uparrow\infty}{\rightarrow}
\{\zeta\mapsto\int\zeta\Xi^\infty(t)\}}\nonumber\\
&\mbox{in law up to near-Gaussian moments}.\nonumber
\end{align*}
Clearly, this convergence preserves the white-noise property (\ref{gr26}). 
It thus remains to show that for any $r<\infty$ (in our fixed sequence), 
the blow-down $F^\infty(r)$ has the white-noise property, 
using that $F(r)$ is exactly local on scale $r$. 
The latter in conjunction with the finite range of $\langle\cdot\rangle$
implies that $F(r)_{|D}$ and $F(r)|_{|D'}$ are independent
provided the distance between the open sets $D$ and $D'$ is larger than $2r+1$. After rescaling,
this implies that $F^R(r)_{|\hat D}$ and $F^R(r)|_{|\hat D'}$ are independent
provided the distance between the open sets $\hat D$ and $\hat D'$ is larger than $\frac{2r+1}{R}$.
Under the convergence (\ref{gr17}) this turns as desired into:
$F^\infty(r)_{|\hat D}$ and $F^\infty(r)|_{|\hat D'}$ are independent
provided the open sets $\hat D$ and $\hat D'$ are disjoint.

\medskip

\step6 Regularity property of blow-downs of $\Xi(t)$.

We claim that a blow-down            
$\Xi^\infty(t)$ of $\Xi(t)$ along any sequence $R\uparrow\infty$ 
satisfies the following regularity property
\begin{align}\label{gr27}
\|\int\zeta:\Xi^\infty(t)\|\lesssim\big(\int|\zeta|^2\big)^\frac{1}{2}.
\end{align}
For notational simplicity, we again pass to scalar instead of tensorial notation and language.
According to Step 3 a) it is sufficient to consider a blow-down of $\Xi_{\sqrt{t}}(t)$,
hence a limit of $\Xi_{\sqrt{t}}^R(t)$. By our notion of convergence (\ref{gr17}), it is enough
to establish $\|\int\zeta\Xi^R_{\sqrt{t}}(t)d\hat x\|\lesssim\big(\int\zeta^2d\hat x\big)^\frac{1}{2}$,
the form of which is preserved under un-rescaling:
$\|\int\zeta\Xi_{\sqrt{t}}(t)dx\|\lesssim\big(\int\zeta^2dx\big)^\frac{1}{2}$. Such an estimate
holds for $q(t)$ and $\nabla\phi(t)$ separately by Remark \ref{rem:Laux12} right after Lemma \ref{Laux12}:
According to Lemma \ref{Laux11}, these centered and stationary random fields are approximately
local relative to space-time averages of $\{q(t')\}_{t'\le t}$, so that Remark \ref{rem:Laux12}
applies; the rhs of (\ref{Laux12.2}) can then be estimated with help of (\ref{theo8b}) in Corollary
\ref{cor:CLTdecay} by $\frac{1}{\sqrt{t}^\frac{d}{2}}$. This yields the desired estimate.

\medskip

\step7 White-noise random fields are necessarily Gaussian. 

Let the centered and stationary
random tensor-valued Schwartz distribution $F$ have the white-noise property
\begin{align}\label{gr13}
\int\zeta: F\;\mbox{and}\;\int\zeta': F
\quad\mbox{are independent for}\quad{\rm supp}\zeta\cap{\rm supp}\zeta'=\emptyset
\end{align}
and the regularity property
\begin{align}\label{gr12}
\norm\int\zeta: F\norm\lesssim\big(\int\zeta^2\big)^\frac{1}{2}
\end{align}
for all Schwartz functions $\zeta$ with values in $\mathbb{R}^{d\times d}$.
Then $F$ is Gaussian in the sense that there exists a constant tensor
${\mathcal Q}\in\mathbb{R}^{(d\times d)\times(d\times d)}$ such that for all Schwartz functions
\begin{align}\label{gr17bis}
\lefteqn{\int\zeta: F}\nonumber\\
&\mbox{is distributed as a centered Gaussian with variance}
\;\int\zeta:{\mathcal Q}\zeta.
\end{align}
This result is standard, and we include a short proof for completeness.
In the following argument we use scalar instead of tensor notation and language. 
For notational convenience we only consider a single Schwartz function $\zeta$ instead
of an arbitrary finite family $(\zeta_1,\dots,\zeta_n), n\in \N$ of such functions (as needed to characterize finite-dimensional laws), this does not change the argument.
We also decompose $\mathbb{R}^d$ into
cubes $C$ of side length $r$ (which we will send to zero)
and in this step denote by the subscript $r$ on $\zeta_r$ the 
$L^2(\mathbb{R}^d)$-orthogonal projection of $\zeta$ on the piecewise constant functions on $\{C\}$,
so that on each $C$, $\zeta_r$ coincides with $\zeta_C:=\fint_C\zeta$.
We note that by (\ref{gr12}) we may (almost-surely) extend $\int\zeta F$ to
$\zeta\in L^2(\mathbb{R}^d)$ so that (\ref{gr13}) \& (\ref{gr12}) hold also for such $\zeta$.
Since $\big(\int|\zeta_r-\zeta|^2\big)^\frac{1}{2}\lesssim r\big(\int|\nabla\zeta|^2\big)^\frac{1}{2}$ 
we have by (\ref{gr12}) that
\begin{align}\label{gr15}
\norm\int\zeta_r F-\int\zeta F\norm\stackrel{r\downarrow 0}{\rightarrow} 0,
\end{align}
so that we first turn to $\int\zeta_r F=\sum_{C}\zeta_C\fint_CF$. 
By the above-mentioned extension of (\ref{gr13}),
$\{\fint_CF\}_C$ are mutually independent so that we obtain for the characteristic function
with parameter $\sigma\in\mathbb{R}$
\begin{align}\label{gr14}
\langle\exp(i\sigma\int\zeta_r F)\rangle=\prod_{C}\langle\exp(i\sigma\zeta_C\int_CF)\rangle.
\end{align}
By stationarity of $F$ we have $\langle\exp(i\sigma\zeta_C\int_CF)\rangle
=\langle\exp(i\sigma\zeta_C\int_{C_0}F)\rangle$ for some fixed cube $C_0$. Since by
(\ref{gr12}), the centered random variable $r^{-\frac{d}{2}}\int_{C_0}F$ 
has all (algebraic) moments bounded uniformly in $r$, we have
\begin{align*}
\langle\exp(i\sigma r^{-\frac{d}{2}}\int_{C_0}F)\rangle=
\exp\big(-\frac{1}{2}\sigma^2\langle(r^{-\frac{d}{2}}\int_{C_0}F)^2\rangle+O(\sigma^3)\big).
\end{align*}
We insert this into (\ref{gr14})
\begin{align}\label{gr16}
\langle\exp(is\int\zeta_r F)\rangle=
\exp\big(-\frac{1}{2}\sigma^2\int\zeta_r^2\langle(r^{-\frac{d}{2}}\int_{C_0}F)^2\rangle
+O(\sigma^3r^\frac{d}{2}\int|\zeta_r|^3)\big).
\end{align}
We select a subsequence $r\downarrow 0$ such that 
${\mathcal Q}=\lim_{r\downarrow 0}\langle(r^{-\frac{d}{2}}\int_{C_0}F)^2\rangle$ exists; in view of
(\ref{gr15}) we may thus pass to the limit in (\ref{gr16}) along this subsequence to obtain
\begin{align}\nonumber
\langle\exp(i\sigma\int\zeta F)\rangle=
\exp(-\frac{1}{2}\sigma^2\int\zeta^2 {\mathcal Q}),
\end{align}
which by the invariance of Gaussians under the Fourier transform amounts to (\ref{gr17bis}).

\medskip

\step8 Characterization of covariance tensor ${\mathcal Q}$, uniqueness of blow-downs, and conclusion
on $\Xi$. 

The main claim is that if $\Gamma(t)$ is a blow-down of $\Xi(t)$ along some
sequence $R\uparrow\infty$ (and thus by Steps 5-7 a centered Gaussian white noise)
we necessarily have for its covariance tensor ${\mathcal Q}(t)$
\begin{align}\label{gr21}
{\mathcal Q}(t)&=\int\langle\Xi(t,z)\otimes\Xi(t,0)\rangle dz\nonumber\\
&:=\lim_{R\uparrow\infty}\int\exp(-|\frac{z}{R}|^2)
\langle\Xi(t,z)\otimes\Xi(t,0)\rangle dz,
\end{align}
where here, $\lim_{R\uparrow\infty}$ is meant as a {\it continuum} limit (ie for {\it all}
sequences), and we will argue that it exists.
Since the law of a centered Gaussian white noise is determined
by its covariance tensor, we learn from (\ref{gr21}) that
all blow-downs of $\Xi(t)$ are identical in law. Evoking the
compactness result of Step 2, we thus learn that $\Xi^R(t)$ converges for the continuum $R\uparrow\infty$
in the sense of (\ref{gr17}) to the centered Gaussian white noise 
$\Gamma(t)$ characterized by its covariance tensor given by (\ref{gr21}).
We now may conclude: By Step 3 b) we know that in this situation also $\Xi^R$ converges for the 
continuum $R\uparrow\infty$ in the sense of (\ref{gr17}) to a $\Gamma$,
which is the continuum limit $t\uparrow\infty$ of $\Gamma(t)$ and therefore itself a
centered and stationary Gaussian white noise and thus characterized by its covariance tensor
${\mathcal Q}$. In view of (\ref{gr17bis}) and our notion of convergence (\ref{gr17}), 
the covariance tensors converge:
\begin{align*}
\lim_{t\uparrow\infty}{\mathcal Q}(t)={\mathcal Q}.
\end{align*}

\medskip

It remains to establish the identity (\ref{gr21}). 
We note that by construction of ${\mathcal Q}(t)$, cf (\ref{gr17bis}),
we have for any $\mathbb{R}^{d\times d}$-valued Schwartz function $\zeta$
\begin{align}\nonumber
\int\zeta:{\mathcal Q}(t)\zeta=\langle(\int\zeta:\Gamma(t))^2\rangle.
\end{align}
By assumption, $\Gamma(t)$ is the blow-down of $\Xi(t)$ along some discrete sequence $R\uparrow\infty$.
According to Step 3 a), it is also the blow-down of $\Xi_{\sqrt{t}}(t)$ along the same sequence;
the latter being more convenient because we have firm statements on its uniform boundedness
and approximate locality.
Hence we have by our definition of convergence (\ref{gr17}) in particular
\begin{align}\label{gr20}
\int\zeta:{\mathcal Q}(t)\zeta
=\lim_{R\uparrow\infty}\langle(\int\zeta:\Xi^R_{\sqrt{t}}(t))^2\rangle.
\end{align}
We recall $(\Xi_{\sqrt{t}}^R(t))_1(\hat x)$=$R^\frac{d}{2}(\Xi_{\sqrt{t}})_R(R\hat x)$ (cf Step 2)
so that $\int G_1(\hat x)\Xi_{\sqrt{t}}^R(t,\hat x)d\hat x=R^\frac{d}{2}(\Xi_{\sqrt{t}})_R(t)$.
We now choose $\zeta$ to be (spatial) Gaussian $G_1$ times fixed tensor $\zeta_0\in\mathbb{R}^{d\times d}$,
so that we have for the rhs of (\ref{gr20})
\begin{align}\label{gr22}
\langle(\int\zeta:\Xi_{\sqrt{t}}^R(t))^2\rangle=R^d\langle(\zeta_0:(\Xi_{\sqrt{t}})_R(t))^2\rangle.
\end{align}
By the semi-group property and stationarity we obtain
\begin{align}\label{gr23}
R^d\langle(\zeta_0:(\Xi_{\sqrt{t}})_R(t))^2\rangle
=\zeta_0:\int R^dG_{\sqrt{2}R}(z)\langle\Xi_{\sqrt{t}}(t,z)\otimes\Xi_{\sqrt{t}}(t,0)\rangle \zeta_0 dz .
\end{align}
We recall (\ref{gr29}) in form of $|\Xi_{\sqrt{t}}(t)-\langle \Xi_{\sqrt{t}}(t)|B_R\rangle|$
$\lesssim\exp(-\frac{1}{C}\frac{R}{\sqrt{t}})$. By definition of the conditional expectation
in conjunction with the finite range of $\langle\cdot\rangle$,
$\langle\Xi_{\sqrt{t}}(t,z)|B_{R}(z)\rangle$ 
and $\langle\Xi_{\sqrt{t}}(t,0)|B_{R}(0)\rangle$ are independent provided $2R+1\le|z|$.
Therefore in this case $\big\langle\langle\Xi_{\sqrt{t}}(t,z)|B_{R}(z)\rangle
\otimes \langle\Xi_{\sqrt{t}}(t,0)|B_{R}(0)\rangle\big\rangle$
$=\langle\Xi_{\sqrt{t}}(t,z)\rangle\otimes\langle\Xi_{\sqrt{t}}(t,0)\rangle$ $=0$, 
where we used that $\Xi(t)$ is centered. We recall (\ref{gr30}) in form of
$|\Xi_{\sqrt{t}}(t)|,|\langle \Xi_{\sqrt{t}}(t)|B_R\rangle|\lesssim 1$. All three relations
combine to the estimate
$|\langle\Xi_{\sqrt{t}}(t,z)\otimes \Xi_{\sqrt{t}}(t,0)\rangle|$ $\lesssim\exp(-\frac{1}{C}\frac{|z|}{\sqrt{t}})$.
In view of $R^dG_{\sqrt{2}R}(z)$ $=G_{\sqrt{2}}(0)$ $\exp(-|\frac{z}{R}|^2)$,
this is amply enough to ensure the existence of the {\it continuum} limit
$\lim_{R\uparrow\infty}$ $\int dz R^dG_{\sqrt{2}R}(z)$ $\langle\Xi_{\sqrt{t}}(t,z)\otimes \Xi_{\sqrt{t}}(t,0)\rangle$.
We now revert from $\Xi_{\sqrt{t}}(t)$ to $\Xi(t)$ itself:
We note that $\Xi(t)$ has finite second moments (because $\nabla\phi(t)$ has), 
so that $\langle\Xi(t,z)\otimes \Xi(t,0)\rangle$ has a sense; by stationarity and the semi-group property we have 
$\langle\Xi_{\sqrt{t}}(t,z)\otimes \Xi_{\sqrt{t}}(t,0)\rangle$
$=\int G_{\sqrt{2t}}(z-x)\langle\Xi(t,x)\otimes \Xi(t,0)\rangle dx$ and thus
$\int R^dG_{\sqrt{2}R}(z)\langle\Xi_{\sqrt{t}}(t,z)\otimes\Xi_{\sqrt{t}}(t,0)\rangle dz
=\int R^dG_{\sqrt{2(R^2+t)}}(z)\langle\Xi(t,z)\otimes\Xi(t,0)\rangle dz$. Hence also the continuum limit
$\lim_{R\uparrow\infty}$ $\int dz R^dG_{\sqrt{2}R}(z)$ $\langle\Xi(t,z)\Xi(t,0)\rangle$
exists and coincides with the analogous limit for $\Xi_{\sqrt{t}}$. Recalling
$R^dG_{\sqrt{2}R}(z)=G_{\sqrt{2}}(0)\exp(-|\frac{z}{R}|^2)$ this yields
\begin{align}\nonumber
\lim_{R\uparrow\infty}\int R^dG_{2R}(z)\langle\Xi_{\sqrt{t}}(t,z)\otimes\Xi_{\sqrt{t}}(t,0)\rangle dz
=G_{\sqrt{2}}(0)\int \langle\Xi(t,z)\otimes \Xi(t,0)\rangle dz
\end{align}
with the definition of $\int \langle\Xi(t,z)\otimes \Xi(t,0)\rangle dz$ from (\ref{gr21}).
The combination of this with (\ref{gr22}) and (\ref{gr23}) yields for rhs of (\ref{gr20})
\begin{align}\nonumber
\lim_{R\uparrow\infty}\langle(\int\zeta:\Xi^R(t))^2\rangle
=G_{\sqrt{2}}(0)\int \zeta_0:\langle\Xi(t,z)\otimes \Xi(t,0)\rangle \zeta_0dz
\end{align}
for our choice of $\zeta(\hat x)=G_1(\hat x)\zeta_0$. The effect of this choice on the 
lhs of (\ref{gr20}) is
$\int\zeta:{\mathcal Q}\zeta d\hat x=\int G_1^2\zeta_0:{\mathcal Q}\zeta_0$. Noting
that by the semi-group property $\int G_1^2=G_{\sqrt{2}}(0)$, we see that (\ref{gr21}) holds.

\medskip

\step9 Construction of 
$S^{\rm hom}_{0\rightarrow\infty}:=\lim_{T\uparrow\infty}S^{\rm hom}_{0\rightarrow T}$. 

We claim that 
$\lim_{T\uparrow\infty}S_{0\rightarrow T}^{\rm hom}$ exists in the sense of operator-norm convergence
wrt $\norm\cdot\norm$: There exists an exponent $\alpha=\alpha(d)>0$ such that for any centered and stationary
vector field $F$ with bounded second moments we have
\begin{align}\label{gr31}
\norm(S_{0\rightarrow T}^{\rm hom}-S_{0\rightarrow\infty}^{\rm hom})F\norm\lesssim T^{-\alpha}\norm F\norm.
\end{align}
Indeed, for a stationary random vector field $F$ with $\langle|F|^2\rangle<\infty$,
almost-surely, the realizations of $F$ have the property $\int\eta|F|^2dx<\infty$. For such realizations,
$S_{0\rightarrow T}^{\rm hom}F$ is deterministically defined as $F+\int_0^Ta_{\rm hom}\nabla vd\tau$, where $v$ is the 
(unique by Lemma~\ref{Laux1}) solution of the constant-coefficient initial value problem
\begin{align}\label{gr33}
\partial_\tau v-\nabla\cdot a_{\rm hom}\nabla v=0\;\mbox{for}\;\tau>0,\quad
v=\nabla\cdot F\;\mbox{for}\;\tau=0.
\end{align}
By \eqref{e.gr} in the proof of Lemma \ref{Laux15}, we have
$T\|\nabla v_R(2T)\|\lesssim\fint_0^{\sqrt{T}}\|F_{\sqrt{R^2+r^2}}\|dr$
for all times $T$ and averaging radii $R$.
By definition of $\norm\cdot\norm$ this yields $T\|\nabla v_R(2T)\|$
$\lesssim\fint_0^{\sqrt{T}}\frac{1}{\sqrt{R^2+r^2}^\frac{d}{2}}dr\norm F\norm$,
which by rescaling $r=R\hat r$ and taking the supremum over $R$ gives $T\norm\nabla v(2T)\norm$
$\lesssim\fint_0^{\sqrt{T}}\frac{1}{\sqrt{\hat r^2+1}^\frac{d}{2}}d\hat r\norm F\norm$. This implies
for some $\alpha=\alpha(d)>0$ (in fact $\alpha=\frac{1}{2}$ for $d>2$ and any $\alpha<\frac{1}{2}$ for $d=2$)
\begin{align}\label{gr35}
\norm\nabla v(T)\norm \lesssim T^{-1-\alpha}\norm F\norm.
\end{align}
We learn from this that
$S_{0\rightarrow\infty}^{\rm hom}:=\lim_{T\uparrow}S_{0\rightarrow T}^{\rm hom}$ 
is well-defined as operator on the space of a stationary random vector fields
$F$ with $\norm F\norm<\infty$, and that we have operator convergence as $T\uparrow\infty$ in form of (\ref{gr31}).

\medskip

\step{10} Identification of $S^{\rm hom}_{0\rightarrow\infty}$ with the Leray projection. 

We claim for any centered and stationary
random vector field $F$ with finite triple norm and finite second moments
\begin{align}
(1-S^{\rm hom}_{0\rightarrow\infty})F&=0\quad\mbox{provided}\;\nabla\cdot F=0\;\mbox{distributionally},\label{gr32}\\
S^{\rm hom}_{0\rightarrow\infty}a_{\rm hom}F&=0\quad\mbox{provided}\;\nabla\times F=0\;\mbox{distributionally},\label{gr34}
\end{align}
where we use $(d=3)$-notation. Identity (\ref{gr32}) is obvious since $\nabla\cdot F=0$ implies $v=0$ for
the unique solution of (\ref{gr33}) and thus $S^{\rm hom}_{0\rightarrow T}F=F$. By Step~9 this yields
as desired $S^{\rm hom}_{0\rightarrow\infty}F=F$. We now turn to (\ref{gr34}) and to this purpose rewrite
(\ref{gr33}) as
$\partial_t\int_0^tvd\tau-\nabla\cdot a_{\rm hom}(\int_0^t\nabla vd\tau+F)=0$
which we mollify on scale $R>0$ and of which we take a partial derivative in direction $i=1,\cdots,d$ to the effect of
$\partial_t\int_0^t\partial_i v_Rd\tau-\nabla\cdot a_{\rm hom}(\nabla\int_0^t\partial_iv_Rd\tau+\partial_i F_R)=0$.
Making use of $\nabla\times F=0$ in form of $\partial_i F_R=\nabla (F_i)_R$, we may rewrite this again
as homogeneous equation
\begin{align}\nonumber
\partial_t(\int_0^t\partial_i vd\tau+F_i)_R
-\nabla\cdot a_{\rm hom}\nabla(\int_0^t\partial_i vd\tau+F_i)_R=0.
\end{align}
By the semi-group estimate (\ref{Laux1.2}) in conjunction with stationarity we obtain
\begin{equation}\nonumber
\langle|\nabla(\int_0^t\partial_i vd\tau+F_i)_R|^2\rangle\lesssim\frac{1}{t}\langle((F_i)_R)^2\rangle
\le\frac{1}{t}\langle(F_i)^2\rangle.
\end{equation}
Thanks to (\ref{gr35}) we may pass to the limit
$\langle|\nabla(\int_0^\infty\partial_i vd\tau+F_i)_R|^2\rangle=0$, which by ergodicity yields
for the centered random variable $(\int_0^\infty\partial_i vd\tau+F_i)_R$ that
$\langle((\int_0^\infty\partial_i vd\tau+F_i)_R)^2\rangle=0$, and thus
$\norm \int_0^\infty\partial_i vd\tau+F_i\norm=0$. Since by definition, $S^{\rm hom}_{0\rightarrow\infty}a_{\rm hom}F
=a_{\rm hom}(\int_0^\infty\nabla vd\tau+F)$, 
we obtain as desired $\norm S^{\rm hom}_{0\rightarrow\infty}F \norm=0$.

\medskip

\step{11} Conclusion for $q$ and $a_{\rm hom}\nabla\phi$. 

For simplicity, we use
vectorial notation (ie fixing the vector $e$) instead of tensorial notation.
To $\Xi=(q-a_{\rm hom}e)-a_{\rm hom}\nabla\phi$ or rather its rescaled version
$\Xi^R=(q-a_{\rm hom}e)^R-(a_{\rm hom}\nabla\phi)^R$ we apply Step~10 to the effect of
\begin{align*}
S^{\rm hom}_{0\rightarrow\infty}\Xi^R=(q- a_{\rm hom}e)^R\quad\mbox{and}\quad 
(1-S^{\rm hom}_{0\rightarrow\infty})\Xi^R=-(a_{\rm hom}\nabla\phi)^R.
\end{align*}
In order to pass to blow-down limits, we replace $S^{\rm hom}_{0\rightarrow\infty}\Xi^R$
by $S^{\rm hom}_{0\rightarrow T}\Xi^R$ with help of Step~9, cf (\ref{gr31}):
\begin{align*}
\norm S^{\rm hom}_{0\rightarrow T}\Xi^R-(q- a_{\rm hom}e)^R\norm
&+\norm (1-S^{\rm hom}_{0\rightarrow T})\Xi^R+(a_{\rm hom}\nabla\phi)^R\norm\nonumber\\
&\lesssim T^{-\alpha}\norm\Xi^R\norm=T^{-\alpha}\norm\Xi\norm\lesssim T^{-\alpha}.
\end{align*}
We now appeal to Step~1 to gather
\begin{align*}
\lefteqn{\|\sup_{\zeta,\zeta'}\Big(\int\zeta\cdot\big( S^{\rm hom}_{0\rightarrow T}\Xi^R-(q- a_{\rm hom}e)^R\big)}\nonumber\\
&+\int\zeta'\cdot\big( (1-S^{\rm hom}_{0\rightarrow T})\Xi^R+(a_{\rm hom}\nabla\phi)^R\big)\Big)\|\lesssim T^{-\alpha},
\end{align*}
where the supremum is taken over all pairs of $\mathbb{R}^d$-valued Schwartz functions $(\zeta,\zeta')$
with norm (\ref{gr9}) $\le 1$. The advantage of $T<\infty$ is that $S^{\rm hom}_{0\rightarrow T}$ maps
Schwartz functions into Schwartz functions so that appealing to the dual operator, we may bring
the above into a form that is amenable to the limit $R\uparrow\infty$:
\begin{align*}
\lefteqn{\|\sup_{\zeta,\zeta'}\Big(\big(\int(S^{\rm hom}_{0\rightarrow T})^*\zeta\cdot\Xi^R
-\int\zeta\cdot(q- a_{\rm hom}e)^R\big)}\nonumber\\
&+\big(\int(1-S^{\rm hom}_{0\rightarrow T})^*\zeta'\cdot\Xi^R+\int\zeta'\cdot(a_{\rm hom}\nabla\phi)^R\big)\Big)\|
\lesssim T^{-\alpha}.
\end{align*}
We thus obtain by Step~8 for a blow-down $((q-a_{\rm hom}e)^\infty,(a_{\rm hom}\nabla\phi)^\infty)$ 
along some sequence $R\uparrow\infty$
\begin{align*}
&(\zeta,\zeta')\mapsto\int\big(\zeta\cdot(q- a_{\rm hom}e)^\infty-\zeta'\cdot(a_{\rm hom}\nabla\phi)^\infty\big)\nonumber\\
&\mbox{is $T^{-\alpha}$-close in law up to near-Gaussian moments to}\quad\nonumber\\
&(\zeta,\zeta')\mapsto\int\big(\big(S^{\rm hom}_{0\rightarrow T})^*\zeta\cdot\Gamma
+(1-S^{\rm hom}_{0\rightarrow T})^*\zeta'\cdot\Gamma\big).
\end{align*}
Since the triple norm of $\Gamma$ is finite (as the limit of centered and stationary random fields with
uniformly bounded triple norm or by a direct computation based on Gaussianity), we may apply (\ref{gr31}) 
also to $\Gamma$; by the same arguments we used for finite $R$ we have
\begin{align*}
&(\zeta,\zeta')\mapsto\int\big(\zeta\cdot(q- a_{\rm hom}e)^\infty-\zeta'\cdot(a_{\rm hom}\nabla\phi)^\infty\big)\nonumber\\
&\mbox{has the same law as}\quad\nonumber\\
&(\zeta,\zeta')\mapsto\int\big(\zeta\cdot S^{\rm hom}_{0\rightarrow\infty}\Gamma
+\zeta'\cdot(1-S^{\rm hom}_{0\rightarrow\infty})\Gamma\big).
\end{align*}
Since in particular all blow-downs coincide in law, and by the compactness argument
of Step 2 applied to $(q-a_{\rm hom}e,a_{\rm hom}\nabla\phi)$, the blow-down converges for
the continuum $R\uparrow\infty$.


\subsection{Proof of Corollary~\ref{cor:corrector}: growth of the extended corrector}

We only display the proof for the corrector $\phi$. The proof for $\sigma$ is similar starting from the flux, as in the proof of Corollary~\ref{cor:CLTdecay}.

\medskip

We split the proof into two main steps and control $\|(\phi)_1(x)-(\phi)_R(x)\|$ and $\|(\phi)_R(x)-(\phi)_R(0)\|$ separately.

\medskip

\step{1} Proof that for all $R\gg 1$, the stationary field $x\mapsto (\phi)_1(x)-(\phi)_R(x)$ satisfies
\begin{equation}\label{Caux-cor-1}
\|\phi_1-\phi_R\| \,\lesssim \, \mu_d(R)^{\frac 12},
\end{equation}
for the following stochastic integrability: $s=2$ for $d>2$, and any $0<s<2$ for $d=2$.
(Note that in this proof, and in this proof only, $\phi_R$ --- and similarly $\phi_{\sqrt{T}}$ later on --- denotes the convolution of the corrector $\phi$ with the Gaussian $G_R$,
and not the modified corrector.)
We display the proof of this estimate in the critical dimension $d=2$ only (the proof for $d>2$ is simpler, and relies on Remark~\ref{rem:nearly} rather than on Corollary~\ref{cor:CLTdecay}), in which case $s<2$.

\medskip

\substep{1.1} Representation formula.

For any two dyadic times $t\ll T$, we claim that 
\begin{align}\label{Caux-cor-2}
\lefteqn{(\phi_{\sqrt{T}}-\phi_{\sqrt{2t}})(0)}\nonumber\\
&=\int\nabla\zeta_{\sqrt{t}\rightarrow\sqrt{T-t}}\cdot\big(a_{\ho}\nabla\phi(t)-q(t)\big)_{\sqrt{t}}
\nonumber\\
&+\sum_{t<\tau<T}\int\nabla\zeta_{0\rightarrow\sqrt{\tau}}
\cdot\big(a_{\ho}\nabla\phi-(a_{\ho}\nabla\phi(t)+(S_{t\rightarrow\infty}^{\ho}-1)q(t))\big)_{\sqrt{\tau}},
\end{align}
where also $\tau$ runs over dyadic times, and 
$\zeta_{r\rightarrow R}$ is characterized for two radii $0<r<R<\infty$ by
\begin{equation}\label{Caux-cor-3}
-\nabla\cdot a_{\ho}'\nabla\zeta_{r\rightarrow R}=G_R-G_r.
\end{equation}
(Recall that $\nabla \phi=\int_0^\infty \nabla ud\tau$, $\nabla \phi(t)=\int_0^t\nabla ud\tau$, and $q(t)=a\int_0^t \nabla ud\tau+ae$.)
Indeed, the potential $\zeta_{r\rightarrow R}$ (the integrability of $\nabla\zeta_{r\rightarrow R}$ will be implicitly established below) is defined such that
\begin{equation}\label{Caux-cor-4}
(\phi_R-\phi_r)(0)=\int\nabla\zeta_{r\rightarrow R}\cdot a_{\ho}\nabla\phi.
\end{equation}
Recall that $1-S_{t\rightarrow\infty}^{\ho}$ is the Helmholtz projection wrt $a_{\ho}$, so
that for any test function $\zeta$
\begin{equation}\label{Caux-cor-5}
\int\nabla\zeta\cdot(a_{\ho}\nabla\phi(t)+(S_{t\rightarrow\infty}^{\ho}-1)q(t))
=\int\nabla\zeta\cdot(a_{\ho}\nabla\phi(t)-q(t)).
\end{equation}
By an approximation argument, (\ref{Caux-cor-5}) extends to smooth decaying functions with integrable gradient.
By (\ref{Caux-cor-3}) we have $\zeta_{\sqrt{2t}\rightarrow\sqrt{T}}=\sum_{t<\tau<T}
\zeta_{\sqrt{\tau}\rightarrow \sqrt{2\tau}}$, where also $\tau$ runs over dyadic times. 
Since (\ref{Caux-cor-3}) is a  constant coefficient equation
we obtain by the semi-group property of convolution with Gaussians
that $\zeta_{\sqrt{2t}\rightarrow\sqrt{T}}=G_{\sqrt{t}}*\zeta_{\sqrt{t}\rightarrow\sqrt{T-t}}$
and that $\zeta_{\sqrt{\tau}\rightarrow\sqrt{2\tau}}=G_{\sqrt{\tau}}*\zeta_{0\rightarrow\sqrt{\tau}}$.
Hence the representation formula \eqref{Caux-cor-2} from (\ref{Caux-cor-4}) and (\ref{Caux-cor-5}).

\medskip

We estimate the first rhs term and the rhs sum of  \eqref{Caux-cor-2} (which is higher order) in the following two substeps.

\medskip

\substep{1.2} Proof of 
\begin{align}\label{Caux-cor-6}
&\|\sum_{t<\tau<T}\int\nabla\zeta_{0\rightarrow\sqrt{\tau}}
\cdot\big(a_{\ho}\nabla\phi-(a_{\ho}\nabla\phi(t)+(S_{t\rightarrow\infty}^{\ho}-1)q(t))\big)_{\sqrt{\tau}}\|_*\nonumber\\
&\lesssim(\frac{1}{\sqrt{t}})^\frac{1}{p}(\log\frac{T}{t}).
\end{align}
We
momentarily set for abbreviation
\begin{equation}\nonumber
F_t:=a_{\ho}\nabla\phi-(a_{\ho}\nabla\phi(t)+(S_{t\rightarrow\infty}^{\ho}-1)q(t)).
\end{equation}
By the triangle inequality 
for $\|\cdot\|_*$ in conjunction with stationarity of $F_t$ we have
\begin{align}\nonumber
\lefteqn{\|\sum_{t<\tau<T}\int\nabla\zeta_{0\rightarrow\sqrt{\tau}}
\cdot (F_t)_{\sqrt{\tau}}\|_*}\nonumber\\
&\le\sum_{t<\tau<T}\int|\nabla\zeta_{0\rightarrow\sqrt{\tau}}|
\|(F_t)_{\sqrt{\tau}}\|_*a_{\ho}'
\le \sum_{t<\tau<T}\big(\int|\nabla\zeta_{0\rightarrow\sqrt{\tau}}|\big)\frac{1}{\sqrt{\tau}}\norm F_t\norm,\nonumber
\end{align}
where we have inserted the definition of $\norm\cdot\norm$ (for $d=2$) in the second step.
By scaling for (\ref{Caux-cor-3}) we have $\zeta_{0\rightarrow\sqrt{\tau}}(x)=\zeta_{0\rightarrow 1}(\frac{x}{\sqrt{\tau}})$ (for $d=2$),
so that we obtain
\begin{align}\label{Caux-cor-7}
\|\sum_{t<\tau<T}\int\nabla\zeta_{0\rightarrow\sqrt{\tau}}
\cdot (F_t)_{\sqrt{\tau}}\|_*
\le (\log\frac{T}{t})\big(\int|\nabla\zeta_{0\rightarrow 1}|\big)\norm F_t\norm.
\end{align}
We note that $\zeta_{0\rightarrow 1}=K_1-K$, where $K$ is the fundamental solution of the 
constant-coefficient elliptic operator
$-\nabla\cdot a_{\ho}'\nabla$ and recall that we allow for systems so that there is no easy formula.
For $d=2$, we know however $\nabla K$ is smooth outside the origin, homogeneous of degree $-1$, and odd.
Because of the latter, $\nabla\zeta_{0\rightarrow 1}$ decays like $\frac{1}{|x|^3}$ for $|x|\gg 1$,
whereas it blows up like $\frac{1}{|x|}$ for $|x|\ll 1$, so that $\int|\nabla\zeta_{0\rightarrow 1}|$
is finite. Hence \eqref{Caux-cor-6} follows from (\ref{Caux-cor-7}) and (\ref{Tau.2}) for $T=\infty$ in Theorem~\ref{Tau} in form of 
\begin{equation*}
\norm a_{\ho}\nabla\phi-(a_{\ho}\nabla\phi(t)+(S_{t\rightarrow\infty}^{\ho}-1)q(t))\norm \lesssim(\frac{1}{\sqrt{t}})^\frac{1}{p}.
\end{equation*}

\medskip

\substep{1.3} Proof of 
\begin{align}\label{Caux-cor-8}
&\|\int\nabla\zeta_{\sqrt{t}\rightarrow\sqrt{T-t}}\cdot \big(a_{\ho}\nabla\phi(t)-q(t)\big)_{\sqrt{t}}\|_*\lesssim
(\log\frac{T}{t})^\frac{1}{2}.
\end{align}
We'd like to apply Lemma~\ref{Laux12} in the form of \eqref{Laux12.2} to $F=a_{\ho}\nabla\phi(t)-q(t)$ and the averaging
function $\bar G:=\nabla\zeta_{\sqrt{t}\rightarrow\sqrt{T-t}}$. We recall for later use that
by \eqref{theo8b} in Corollary~\ref{cor:CLTdecay}
and the definition of $\norm \cdot\norm$
\begin{equation}\label{Caux-cor-9}
\sup_{r\le\sqrt{t}}(\frac{r}{\sqrt{t}})^\frac{d}{2}\|F\|_*\le (\frac{1}{\sqrt{t}})^\frac{d}{2}
\norm F\norm\lesssim (\frac{1}{\sqrt{t}})^\frac{d}{2}.
\end{equation}
We also note that because of the representation $\bar G=K_{\sqrt{T-t}}-K_{\sqrt{t}}$ 
in terms of the fundamental solution $K$ of the constant-coefficient elliptic operator $-\nabla\cdot a_{\ho}'\nabla$,
we see as above that in case of $d=2$ we have
$|\nabla\bar G(x)|\lesssim\frac{1}{|x|}$ in the intermediate range $\sqrt{t}\le|x|\le\sqrt{T}$,
whereas we have $|\nabla\bar G(x)|\lesssim\frac{1}{\sqrt{t}}$ for the near-field $|x|\le\sqrt{t}$, 
and $|\nabla\bar G(x)|\lesssim\frac{T}{|x|^3}$ for the far-field $|x|\ge\sqrt{T}$. This yields
\begin{equation}\label{Caux-cor-10}
\int|\bar G|^2\lesssim\log\frac{T}{t}.
\end{equation}
We note that according to Lemma \ref{Laux11}, $F$
is indeed approximately local on scale $\sqrt{t}$ relative to flux fluctuations, that is, the stationary random field
\begin{equation}\nonumber
\bar F:=\fint_0^t d\tau(\frac{\sqrt{\tau}}{\sqrt{t}})^\frac{d}{2}
\fint_0^{\sqrt{\tau}}dr(\frac{r}{\sqrt{\tau}})^\frac{d}{2}|(q(\tau)-\langle q(\tau)\rangle)_r|.
\end{equation}
We record for later use that by definition of $\norm\cdot\norm$ and Corollary~\ref{cor:CLTdecay}
\begin{equation}\label{Caux-cor-11}
\|\bar F\|\le(\frac{1}{\sqrt{t}})^\frac{d}{2}\fint_0^td\tau\norm q(\tau)\norm\lesssim(\frac{1}{\sqrt{t}})^\frac{d}{2}.
\end{equation}
Hence we obtain by (\ref{Laux12.2}) with $T$ replaced by $t$
\begin{align*}
&\|\bar G*F_{\sqrt{t}}\|_*\lesssim\sqrt{t}^\frac{d}{2}\big(\int|\bar G|^2\big)^\frac{1}{2}\nonumber\\
&\times\Big(\big(\sup_{r\le\sqrt{t}}(\frac{r}{\sqrt{t}})^\frac{d}{2}\|F_r\|_*\big)^{1-\frac{d}{2p}}
\|\bar F\|^\frac{d}{2p}
+            \sup_{r\le\sqrt{t}}(\frac{r}{\sqrt{t}})^\frac{d}{2}\|F_r\|_*\Big).\nonumber
\end{align*}
By definition of $F$, $\bar G$, and by the estimates (\ref{Caux-cor-10}), (\ref{Caux-cor-9}), and (\ref{Caux-cor-11}),
this yields the desired estimate \eqref{Caux-cor-8}.

\medskip

\substep{1.4} Proof of \eqref{Caux-cor-1}.

Inserting \eqref{Caux-cor-8} and (\ref{Caux-cor-6}) into (\ref{Caux-cor-2}) we obtain
\begin{equation}\label{Caux-cor-12}
\|(\phi_{\sqrt{T}}-\phi_{\sqrt{t}})(0)\|_*\lesssim (\log\frac{T}{t})^\frac{1}{2}
+(\frac{1}{\sqrt{t}})^\frac{1}{p}(\log\frac{T}{t}),
\end{equation}
where we may pass from the semi-norm $\|\cdot\|_*$ to the norm $\|\cdot\|$, since
$(\phi_{\sqrt{T}}-\phi_{\sqrt{t}})(0)$ is a deterministic linear functional of 
the stationary $\nabla\phi$, cf (\ref{Caux-cor-4}), so that $\langle\nabla\phi\rangle=0$ translates into
$\langle(\phi_{\sqrt{T}}-\phi_{\sqrt{t}})(0)\rangle=0$.

\medskip

Estimate (\ref{Caux-cor-12}) is easily upgraded: Applying (\ref{Caux-cor-12}) to the pairs $(T,t)$ and $(t,1)$ 
for some intermediate time $1\ll t\ll T$ to be optimized, we obtain by the triangle inequality
\begin{align*}
\|(\phi_{\sqrt{T}}-\phi_{1})(0)\|&\le
\|(\phi_{\sqrt{T}}-\phi_{\sqrt{t}})(0)\|+\|(\phi_{\sqrt{t}}-\phi_{1})(0)\|\\
&\lesssim (\log\frac{T}{t})^\frac{1}{2}
+(\frac{1}{\sqrt{t}})^\frac{1}{p}(\log\frac{T}{t})+\log t.
\end{align*}
The choice of $t=\log^{2p}T$, which is made such that $(\frac{1}{\sqrt{t}})^\frac{1}{p}(\log\frac{T}{t})
\le(\frac{1}{\sqrt{t}})^\frac{1}{p}(\log T)\le 1$, yields
\begin{equation*}
\|(\phi_{\sqrt{T}}-\phi_{1})(0)\|\lesssim (\log T)^\frac{1}{2}+\log\log T\lesssim (\log T)^\frac{1}{2}.
\end{equation*}
This yields \eqref{Caux-cor-1} for the choice $\sqrt{T}=R$.

\medskip

\step2 Conclusion.

On the one hand, by the triangle inequality and by Poincar\'e's inequality (for $dx$ replaced by $Gdx$), we have for all $x\in \R^d$
$$
\|   (|\phi-\phi_1(0)|^2)_1^\frac12 (x) \| \,\lesssim \, \| \phi_1(x)-\phi_1(0)\|+\| (|\nabla \phi|^2)_1^\frac12 \|,
$$
which, by \eqref{e.decay-semi-3} in Corollary~\ref{coro:decay-semi-group} and the dominance \eqref{dominance} of the Gaussian average by
the exponential average, turns into
$$
\|   (|\phi-\phi_1(0)|^2)_1^\frac12 (x) \| \,\lesssim \, \| \phi_1(x)-\phi_1(0)\|+1.
$$
On the other hand, by the triangle inequality,
$$
\|\phi_1(x)-\phi_1(0)\| \,\leq \, \|\phi_R(x)-\phi_1(x)\|+\|\phi_R(x)-\phi_R(0)\|+\|\phi_R(0)-\phi_1(0)\|.
$$
The first and last rhs terms are controlled by  \eqref{Caux-cor-1}  by stationarity.
It remains to argue that the second rhs term is higher order.
This follows from
\begin{equation}\label{Caux-cor-13}
\|\phi_R(x)-\phi_R(0) \| \,\lesssim \, |x|R^{-1},
\end{equation}
which holds with the following stochastic integrability: $s=2$ for $d>2$, and any $0<s<2$ for $d=2$.
The starting point for \eqref{Caux-cor-13} is the fundamental theorem of calculus in the form
$$
\int (\phi(y+x)-\phi(y))G_R(y)\,=\, |x|\int_0^1 d\tau \frac{x}{|x|}\cdot \Big(\int \nabla \phi(y+\tau x)G_R(y)dy\Big).
$$
By Jensen's inequality and stationarity of $\nabla \phi$,  \eqref{Caux-cor-13}
follows from \eqref{theo8} in Corollary~\ref{cor:CLTdecay} for $d=2$ and $s<2$ and from \eqref{theo8c} in Remark~\ref{cor:CLTdecay}
for $d>2$ and $s=2$ (for the choice $\alpha=1$).


\subsection{Proof of Corollary~\ref{coro:decay-semi-group}: decay of the semi-group}

The starting point for the proof is \eqref{Laux6.2} with $p=\frac{d}2+1$ in Lemma~\ref{Laux6}, which yields
\begin{equation*}
\Big(\int\eta_{\sqrt{2T}}|T\nabla u(T)|^2\Big)^\frac{1}{2}\,\lesssim \,  \fint_{\frac{T}{2}}^Tdt\fint_0^{\sqrt{t}}dr(\frac{r}{\sqrt{t}})^{\frac{d}{2}+1}\int\eta_{2\sqrt{2T}}|(q(t)-\langle q(t)\rangle)_r| ,
\end{equation*}
so that by the triangle inequality and  \eqref{theo8b} in Corollary~\ref{cor:CLTdecay}, we have for all $s<2$
\begin{eqnarray*}
\|\Big(\int\eta_{\sqrt{2T}}|T\nabla u(T)|^2\Big)^\frac{1}{2}\|&\lesssim & \fint_{\frac{T}{2}}^Tdt\fint_0^{\sqrt{t}}dr(\frac{r}{\sqrt{t}})^{\frac{d}{2}+1}\|q(t)_r\|_*
\\
&\lesssim &  \fint_{\frac{T}{2}}^T(\frac{1}{\sqrt{t}})^{\frac{d}{2}+1}dt\fint_0^{\sqrt{t}}r dr\,\lesssim\, \sqrt{T}^{-\frac{d}{2}},
\end{eqnarray*}
that is,  \eqref{e.decay-semi-1}.

\medskip

Likewise, this argument yields \eqref{e.decay-semi-1bis} replacing \eqref{theo8b} by \eqref{theo8d} in Remark~\ref{rem:nearly} with Gaussian integrability and nearly-optimal scaling.

\medskip

Estimate~\eqref{e.decay-semi-2} is an immediate consequence of \eqref{e.decay-semi-1} combined with the stationarity of $|\nabla u(T,\cdot)|^2$:
For all  $T\gg 1$,
$$
\langle |\nabla u(T)|^2 \rangle\,=\, \langle \int\eta_{\sqrt{2T}}|\nabla u(T)|^2 \rangle\,\stackrel{\eqref{e.decay-semi-1}}{\lesssim} \, T^{-1-\frac{d}{4}}.
$$

\medskip

We conclude with the proof of \eqref{e.decay-semi-3}.
Starting point is the triangle inequality in the form of
\begin{eqnarray*}
\Big(\int \eta_R |\nabla \phi|^2\Big)^\frac12 &=& \Big(\int \eta_R \big| \int_0^\infty \nabla u\big|^2d\tau \Big)^\frac12
\\
&\leq & \Big(\int \eta_R \big| \int_0^{R^2} \nabla u\big|^2d\tau \Big)^\frac12  +\sum_{i=1}^\infty \Big(\int_{2^iR^2}^{2^{i+1}R^2} \int \eta_R |\nabla u|^2d\tau\Big)^\frac12 .
\end{eqnarray*}
For the first rhs term, we appeal to the localized semi-group estimate~\eqref{Laux1.1} in Lemma~\ref{Laux1}, which yields the deterministic estimate
$$
\Big(\int \eta_R \big| \int_0^{R^2} \nabla u\big|^2d\tau \Big)^\frac12 \,\lesssim \, 1.
$$
For the rhs sum, we appeal to the Meyers estimate of Lemma~\ref{Laux7}, the triangle inequality, and  \eqref{e.decay-semi-1bis},
to the effect that for $s=2$
\begin{eqnarray*}
\| \Big(\int_{2^iR^2}^{2^{i+1}R^2} \int \eta_R |\nabla u|^2 d\tau\Big)^\frac12\|&\stackrel{\eqref{Laux7.1}}{\lesssim} &
\Big(\frac{\sqrt{2^{i+1}R^2}}{R} \Big)^{\frac{d}{2}-\e} \| \Big( \int_{2^{i-1}R^2}^{2^{i+1}R^2} \int \eta_{\sqrt{2^{i+1}R^2}} |\nabla u|^2 d\tau\Big)^\frac{1}{2} \|
\\
&\leq & \Big(\frac{\sqrt{2^{i+1}R^2}}{R} \Big)^{\frac{d}{2}-\e} \int_{2^{i-1}R^2}^{2^{i+1}R^2}  \| \Big( \int \eta_{\sqrt{2^{i+1}R^2}} |\nabla u|^2 \Big)^\frac{1}{2} \| d\tau
\\
&\stackrel{\eqref{e.decay-semi-1bis}}{\lesssim} & (2^{i})^{\frac d4 -\frac \e 2} (2^i R^2) (2^i R^2)^{-1-\frac d4+\frac\e4}\,=\,R^{-\frac d2+\frac\e2} (2^{\frac \e4})^{-i},
\end{eqnarray*}
which sums to $\sim R^{-\frac d2+\frac\e2}$, and proves \eqref{e.decay-semi-3}.


\subsection{Proof of Theorem~\ref{Tsys}: systematic errors}

We split the proof into two steps.

\medskip

\step{1} Proof of \eqref{e.sys1}.

\noindent The starting point is the representation of $\nabla \phi_T^\kappa$ and $\nabla \phi$ as
\begin{equation*}
\nabla \phi\,=\,\int_0^\infty \nabla u(\tau)d\tau,\qquad
\nabla \phi_T^\kappa\,=\,\int_0^\infty \exp_\kappa(\tau,T) \nabla u(\tau)d\tau,
\end{equation*}
where $\exp_\kappa(\cdot,T)$ is defined as the Richardson extrapolation of $\exp_1(\tau,T):=\exp(-\frac{\tau}{T})$ wrt $T$.
The extrapolation has the effect that
\begin{equation}\label{e.bds-expkapp}
|1-\exp_\kappa(\tau,T)|\,\lesssim\, \min\{(\frac{\tau}{T})^\kappa,1\},\qquad 
|\frac{\partial}{\partial_\tau}\exp_\kappa(\tau,T)|\,\lesssim \, \frac{1}{T}(\frac{\tau}{T})^{\kappa-1} \quad \text{ for all }\tau\ge 0.
\end{equation}
We then split the integral over $(0,+\infty)$ into three contributions.
We start with the contribution on the interval $(0,1)$.
By integration by parts,
\begin{multline*}
\int_0^1 (1-\exp_\kappa(\tau,T))\nabla u(\tau)d\tau \,=\,\int_0^1 \frac{\partial}{\partial \tau}\exp_\kappa(\tau,T) \int_0^\tau \nabla u(t)dt d\tau
\\
+(1-\exp_\kappa(1,T))\int_0^1 \nabla u(t)dt .
\end{multline*}
By Cauchy-Schwarz' inequality in the first time integral and in probability, and by stationarity of $\int_0^1 \nabla u(t)dt$,
\begin{multline*}
\big\langle \big(\int_0^1 (1-\exp_\kappa(\tau,T))\nabla u(\tau)d\tau\big)^2 \big\rangle 
\\
\lesssim\,
\int_0^1  |\frac{\partial}{\partial \tau}\exp_\kappa(\tau,T)|^2d\tau
 \times  \int_0^1 \big \langle \int \eta_{\sqrt{\tau}}  |\int_0^\tau \nabla u(t)dt|^2 \big \rangle d\tau
 \\
+ |1-\exp_\kappa(1,T)|^2 \langle \int \eta |\int_0^1 \nabla u(t)dt|^2 \rangle.
\end{multline*}
Hence, \eqref{Laux1.1} in Lemma~\ref{Laux1} and \eqref{e.bds-expkapp} combine to 
$$
\big\langle \big(\int_0^1 (1-\exp_\kappa(\tau,T))\nabla u(\tau)d\tau\big)^2 \big\rangle \,\lesssim \, T^{-2\kappa},
$$
which is of higher order than the rhs of \eqref{e.sys1}.
We turn now to the contributions on the intervals $(1,T)$ and $(T,\infty)$, for which 
the estimate of the decay of the semi-group \eqref{e.decay-semi-2} in Lemma~\ref{coro:decay-semi-group} 
and \eqref{e.bds-expkapp}  for $\kappa>\frac{d}{4}$ yield 
\begin{eqnarray*}
\big\langle \big(\int_1^\infty (1-\exp_\kappa(\tau,T))\nabla u(\tau)d\tau\big)^2 \big\rangle^\frac12
& \lesssim & \int_1^T (\frac{\tau}{T})^\kappa \tau^{-1-\frac{d}{4}}d\tau+\int_T^\infty  \tau^{-1-\frac{d}{4}}d\tau
\\
&\lesssim & T^{-\frac{d}{4}}.
 \end{eqnarray*}
This proves \eqref{e.sys1}.

\medskip

\step{2} Proof of \eqref{e.sys2}.

\noindent This estimate is a direct consequence of  \eqref{e.sys1}.
By definitions \eqref{e.sys-def-ahT} of $a^\kappa_{hT}$, and the following definition of $a_\ho$, 
$$
e_j\cdot a_\ho e_i\,=\,e_j\cdot\langle  a (\nabla \phi_i+e_i)\rangle\,=\,\langle (\nabla \phi_j'+e_j)\cdot a (\nabla \phi_i+e_i)\rangle,
$$
where $\phi_i$ is solution of  \eqref{o56} for $e=e_i$ (and $\phi_j$ the  solution of \eqref{o56} for $e=e_j$ and $a$ replaced by its pointwise transpose $a'$)
we have
$$
e_j\cdot (a^\kappa_{hT}-a_{\ho})e_i\,=\,\langle(\nabla {\phi_{Tj}'^\kappa}- \nabla \phi'_j)\cdot a(\nabla \phi_{Ti}^\kappa+e_i)\rangle -\langle(\nabla \phi'_j+e_j)\cdot a(\nabla \phi_i-\nabla {\phi_{Ti}^\kappa})\rangle.
$$
Since
$$
\langle (\nabla \phi_j'+e_j)\cdot a(\nabla \phi_i-\nabla \phi_{Ti}^\kappa) \rangle=\langle(\nabla \phi_i-\nabla \phi_{Ti}^\kappa)\cdot a'(\nabla \phi_j'+e_j)\rangle,
$$
the weak form of the corrector equation for $\phi_j'$ in probability and for $\phi_i$ then yields (whence the choice of the adjoint corrector in the definition of $a^\kappa_{hT}$)
$$
\langle (\nabla \phi_i-\nabla \phi_{Ti}^\kappa)\cdot a'(\nabla \phi_j'+e_j)\rangle\,=\,0\,=\langle(\nabla \phi_j'-\nabla \phi_{Tj}'^\kappa)\cdot a(\nabla \phi_i+e_i)\rangle,
$$
and we conclude 
$$
e_j\cdot (a^\kappa_{hT}-a_{\ho})e_i\,=\,\langle (\nabla {\phi_{Tj}'^\kappa}- \nabla \phi_j')\cdot a(\nabla \phi_{Ti}^\kappa- \nabla \phi_i)\rangle,
$$
so that the claim follows from \eqref{e.sys1} (used for both $a$ and $a'$).


\subsection{Proof of Corollary~\ref{cor:sys}: spectral exponents}

For $\tilde \mu \ge 1$, there is nothing to prove and we assume $T=\frac{1}{\tilde \mu}\ge 1$.
The starting point is the spectral theorem which allows one to  rewrite the definition of $\phi_{Ti}^\kappa$ in the form
$$
\phi_{Ti}^\kappa\,=\, g_\kappa(\mathcal L,T) \mathfrak d,
$$
where $g_1(\mu,T)=(T^{-1}+\mu)^{-1}$, and $g_\kappa$ is the Richardson extrapolation of $g_1$ wrt $T$.
Set $g_0(\mu)=\frac{1}{\mu}$.
Then, by the spectral theorem,
\begin{eqnarray*}
\langle \nabla (\phi_{Ti}^\kappa-\phi_i)\cdot a \nabla (\phi_{Ti}^\kappa-\phi_i)\rangle
&=&\langle (\phi_{Ti}^\kappa-\phi_i) \mathcal L (\phi_{Ti}^\kappa-\phi_i)\rangle
\\
&=&
 \int_0^\infty \mu (g_\kappa(\mu,T)-g_0(\mu))^2 \langle \mathfrak d P(d\mu)\mathfrak d\rangle.
\end{eqnarray*}
On the one hand, for $\kappa>\frac{d}{4}$, Theorem~\ref{Tsys} yields 
$$
\langle \nabla (\phi_{Ti}^\kappa-\phi_i)\cdot a \nabla (\phi_{Ti}^\kappa-\phi_i)\rangle \,\le \, \langle |\nabla (\phi_{Ti}^\kappa-\phi_i)|^2\rangle\,\lesssim \,T^{-\frac{d}{2}}.
$$
On the other hand, by induction on $\kappa$ (see for instance \cite[Proof of Lemma~2.5, Step~2]{GloriaNolen}) we have 
$$
|g_\kappa(\mu,T)-g_0(\mu)|\,\gtrsim \,\frac{T^{-\kappa}}{\mu (T^{-1}+\mu)^\kappa},
$$
which we use in the form:  For all $\mu\le \frac{1}{T}$,
$$
\mu (g_\kappa(\mu,T)-g_0(\mu))^2 \, \gtrsim \,\frac{T^{-2\kappa}}{\mu(T^{-1}+\mu)^{2\kappa}}.
$$
The combination of these two estimates directly yields the claim, recalling that $T^{-1}=\tilde \mu$,
\begin{eqnarray*}
\int_0^\frac{1}{T}  \langle \mathfrak d P(d\mu)\mathfrak d\rangle  &\lesssim &
T^{-1} \int_0^\frac{1}{T} \frac{T^{-2\kappa}}{\mu(T^{-1}+\mu)^{2\kappa}} \langle \mathfrak d P(d\mu)\mathfrak d\rangle
\\
&\lesssim &T^{-1} \int_0^\frac{1}{T} \mu (g_\kappa(\mu,T)-g_0(\mu))^2  \langle \mathfrak d P(d\mu)\mathfrak d\rangle
\\
&\lesssim & T^{-1-\frac{d}{2}}.
\end{eqnarray*}
%


\subsection{Proof of Corollary~\ref{cor:r*}: stochastic integrability of the minimal radius}
Starting point is Lemma~\ref{Laux8bis}, and 
we split the proof of the corollary into three steps.

\medskip

\step1 Deterministic estimates on the quantities appearing on the rhs
of (\ref{Cr*-2}). 

For the random stationary field
\begin{equation}\label{Cr*-3}
F_T:=|(\frac{\phi_T}{\sqrt{T}},\frac{\sigma_T}{\sqrt{T}},(q_T-\langle q_T\rangle)_{\sqrt{T}})|^2
\end{equation}
we claim the following uniform bound on the spatially averaged version $(F_T)_{\sqrt{T}}$
\begin{equation}\label{Cr*-4}
(F_T)_{\sqrt{T}}\lesssim 1,
\end{equation}
and the following approximate locality on the same quantity in an absolute (and exponential) sense:
\begin{align}\label{Cr*-5}
\lefteqn{\big(\fint_{B_{\sqrt{T}}}|(F_T)_{\sqrt{T}}(a)-(F_T)_{\sqrt{T}}(\tilde a)|^2\Big)^\frac{1}{2}\lesssim\exp(-\frac{1}{C}\frac{R}{\sqrt{T}})}\nonumber\\
&\quad\mbox{provided}\;a=\tilde a\;\mbox{outside of}\;B_R\;\mbox{and}\;R\ge\sqrt{T}.
\end{align}

\medskip

We first address the uniform bound (\ref{Cr*-4})
which by the estimate on the Gaussian/exponential kernels $G_L \lesssim \eta_L$ in form of 
\begin{equation}\label{dominance}
F\ge 0 \,\implies \, F_L \,\lesssim\, \int \eta_L F,
\end{equation}
which we refer to as the 
dominance of the Gaussian average by
the exponential average, follows from a uniform bound on
$\int\eta_{\sqrt{T}}|(\frac{\phi_T}{\sqrt{T}},\frac{\sigma_T}{\sqrt{T}},q_T)|^2$ (since $\langle q_T\rangle =a_{hT}$
is bounded by \eqref{Laux4.3}),
which in turn is a consequence of the localized elliptic energy estimates \eqref{Laux3.1}
for equations \eqref{c62} and \eqref{o17}
\begin{equation}\label{Cr*-6}
\int\eta_{\sqrt{T}}|(\frac{\phi_T}{\sqrt{T}},\nabla\phi_T)|^2\lesssim 1\quad\mbox{and}
\int\eta_{\sqrt{T}}|(\frac{\sigma_T}{\sqrt{T}},\nabla\sigma_T)|^2\lesssim 1.
\end{equation}

\medskip

We now address the approximate locality (\ref{Cr*-5}). 
By Jensen's inequality and the dominance
of the Gaussian by the exponential average we have
\begin{align*}\nonumber
\lefteqn{\fint_{B_{\sqrt{T}}}|(F_T)_{\sqrt{T}}(a)-(F_T)_{\sqrt{T}}(\tilde a)|}\nonumber\\
\lesssim&\int\eta_{\sqrt{T}}\big|(|\frac{\phi_T}{\sqrt{T}}|^2-|\frac{\tilde\phi_T}{\sqrt{T}}|^2,
|\frac{\sigma_T}{\sqrt{T}}|^2-|\frac{\tilde\sigma_T}{\sqrt{T}}|^2,
|q_T-\langle q_T\rangle|^2-|\tilde q_T-\langle q_T\rangle|^2)\big|,
\end{align*}
where $(\phi_T,\tilde\phi_T)$ stands for $(\phi_T(a),\phi_T(\tilde a))$ and the analogue
notational convention for $\sigma_T$ and $q_T$.
By the Cauchy-Schwarz inequality in conjunction with (\ref{Cr*-6}), the lhs of the above is estimated by
the square root of
\begin{equation}\nonumber
\lesssim\int\eta_{\sqrt{T}}\big|(\frac{\phi_T-\tilde\phi_T}{\sqrt{T}},
\frac{\sigma_T-\tilde\sigma_T}{\sqrt{T}},
q_T-\tilde q_T)\big|^2.
\end{equation}
In order to establish that this term is estimated by the rhs of (\ref{Cr*-5}), we note that
$\phi_T-\tilde\phi_T$ satisfies the equation
\begin{equation}\nonumber
\frac{1}{T}(\phi_T-\tilde\phi_T)-\nabla\cdot \tilde a\nabla(\phi_T-\tilde\phi_T)=\nabla\cdot(a-\tilde a)(\nabla\phi_T+e),
\end{equation}
so that we obtain by the localized elliptic energy estimate \eqref{Laux3.1}
\begin{equation}\nonumber
\int\eta_{\sqrt{T}}\big|(\frac{\phi_T-\tilde\phi_T}{\sqrt{T}},
\nabla(\phi_T-\tilde\phi_T))\big|^2\lesssim
\int\eta_{\sqrt{T}}|(a-\tilde a)(\nabla\phi_T+e)|^2.
\end{equation}
Because of $q_T-\tilde q_T=\tilde a\nabla(\phi_T-\tilde\phi_T)+(a-\tilde a)(\nabla\phi_T+e)$, this yields
\begin{equation}\nonumber
\int\eta_{\sqrt{T}}\big|(\frac{\phi_T-\tilde\phi_T}{\sqrt{T}},
q_T-\tilde q_T)\big|^2\lesssim
\int\eta_{\sqrt{T}}|(a-\tilde a)(\nabla\phi_T+e)|^2.
\end{equation}
Giving up a bit of the exponential cut-off on the rhs in form of the identity $\eta_{\sqrt{T}}(x)$ $=\exp(-\frac{|x|}{2\sqrt{T}})\eta_{2\sqrt{T}}(x)$, 
using that $\exp(-\frac{|x|}{2\sqrt{T}})|a-\tilde a|^2\le\exp(-\frac{R}{2\sqrt{T}})$ under our assumption
on the relation of $a$ and $\tilde a$, and appealing to the uniform bound (\ref{Cr*-6}) in the slightly averaged form of
$\int\eta_{2\sqrt{T}}|\nabla\phi_T+e|^2\lesssim 1$ yields
\begin{equation}\label{Cr*-7}
\int\eta_{\sqrt{T}}\big|(\frac{\phi_T-\tilde\phi_T}{\sqrt{T}},
q_T-\tilde q_T)\big|^2\lesssim
\exp(-\frac{R}{2\sqrt{T}}).
\end{equation}
We note that $\sigma_T-\tilde\sigma_T$ satisfies, in convenient 3-d notation,
\begin{equation}\nonumber
\frac{1}{T}(\sigma_T-\tilde\sigma_T)-\triangle(\sigma_T-\tilde\sigma_T)=\nabla\times(q_T-\tilde q_T),
\end{equation}
so that by the localized energy estimate in form of
\begin{equation}\nonumber
\int\eta_{\sqrt{T}}|\frac{\sigma_T-\tilde\sigma_T}{\sqrt{T}}|^2
\lesssim\int\eta_{\sqrt{T}}|q_T-\tilde q_T|^2,
\end{equation}
we see that (\ref{Cr*-7}) may be upgraded to (\ref{Cr*-5}).

\medskip

\step2 Stochastic estimates on the quantities appearing on the rhs of (\ref{Cr*-2}).

We claim that the random variable defined by the dyadic sum
\begin{equation}\label{Cr*-8}
F:=\sum_{R\ge r}(\frac{R}{r})^\beta \int\eta_RF_T
\end{equation}
satisfies
\begin{equation}\label{Cr*-9}
\langle I(F>\delta)\rangle\lesssim\exp(-\frac{1}{C}\delta^{d+2}r^d)\quad\mbox{for all}\;\delta\gg \frac1{\sqrt{t}}.
\end{equation}

\medskip

Equipped with the deterministic estimates on the building blocks $F_T$ from Step 1, we now may embark on their stochastic estimate.
We start with the fluctuations and to this purpose
apply Lemma \ref{Laux12} with the generic $F$ replaced by $(F_T)_{\sqrt{T}}$. 
Because of the absolute approximate locality (\ref{Cr*-5}), we may choose $\bar F=1$.
Because of the uniform bound (\ref{Cr*-4}), the outcome (\ref{Laux12.3}) of Lemma \ref{Laux12} simplifies to
\begin{equation}\nonumber
\|((F_T)_{\sqrt{T}})_R\|_*\lesssim (\frac{\sqrt{T}}{R})^\frac{d}{2}\quad\mbox{for}\;R\ge\sqrt{T}.
\end{equation}
By the semi-group property of convolution with Gaussians and the monotonicity of
$R\mapsto \|F_R\|_*$, this yields
\begin{equation}\nonumber
\|(F_T)_R\|_*\lesssim (\frac{\sqrt{T}}{R})^\frac{d}{2}\quad\mbox{for}\;R\ge\sqrt{T}.
\end{equation}
By the dominance of Gaussian averages by exponential averages,  stationarity,  and Jensen's inequality, this implies
\begin{equation}\label{Cr*-10}
\|\int\eta_RF_T\|_*\lesssim (\frac{\sqrt{T}}{R})^\frac{d}{2}\quad\mbox{for}\;R\gtrsim\sqrt{T}.
\end{equation}

\medskip

We now turn to the expectation of $F_T$ and claim the following (suboptimal) bound
\begin{equation}\label{Cr*-11}
\langle \int\eta_R F_T\rangle\lesssim\frac{1}{\sqrt{T}}.
\end{equation}
We recall our bounds on the augmented corrector $(\phi,\sigma)$ and its flux $q$, in their optimal
form in terms of scaling, which here we just need on the level of their second stochastic moments
and of the modified corrector $(\phi_T,\sigma_T)$ and its flux $q_T$ (actually, we don't even
need the optimal scaling), cf Remark~\ref{rem:modcorr} and \eqref{theo8b} in Corollary~\ref{cor:CLTdecay}. They yield for $T\gg 1$
\begin{align*}
\langle|(\phi_T,\sigma_T)|^2\rangle\lesssim\mu_d(T)
\quad\mbox{and}\quad\langle|(q_T-\langle q_T\rangle)_{\sqrt{T}}|^2\rangle\lesssim\frac{1}{\sqrt{T}^d}.
\end{align*}
The first estimate limits the final result, which we use in the suboptimal form of
\begin{equation*}
\langle F_T\rangle\lesssim\frac{1}{\sqrt{T}},
\end{equation*}
from which (\ref{Cr*-11}) follows by stationarity.

\medskip

We now consider the ``aggregated'' random variable $F$ defined in (\ref{Cr*-8}). 
We recall that the scales $(r,\sqrt{t})$ and $(R,\sqrt{T})$ are subject to the constraints
of Lemma~\ref{Laux8bis}, so that in particular $R\gtrsim\sqrt{T}$.
By (\ref{Cr*-11})  we obtain for its expectation
\begin{equation}\nonumber
\langle F\rangle\lesssim\sum_{R\ge r}(\frac{R}{r})^\beta \frac{1}{\sqrt{T}};
\end{equation}
by (\ref{Cr*-10}) in conjunction with the triangle inequality we obtain for its fluctuations
\begin{equation}\nonumber
\| F\|_*\lesssim\sum_{R\ge r}(\frac{R}{r})^\beta(\frac{\sqrt{T}}{R})^\frac{d}{2},
\end{equation}
since (\ref{Cr*-1}) ensures that $r\ge\sqrt{t}$ propagates to $R\gtrsim\sqrt{T}$.
We now substitute $\sqrt{T}$ according to $\sqrt{T}\stackrel{(\ref{Cr*-1})}{\sim}(\frac{R}{r})^\alpha\sqrt{t}$. 
Provided we choose the exponents $\alpha\in[0,1]$ and $\beta>0$ such that
\begin{equation}\nonumber
\beta<\alpha\quad\mbox{and}\quad\beta<(1-\alpha)\frac{d}{2}
\end{equation}
--- which can be easily done: $\alpha=\frac{1}{2}$ and $\beta=\frac{1}{4}$ is such a choice ---
the two estimates above turn into
\begin{equation}\nonumber
\langle F\rangle\lesssim\frac{1}{\sqrt{t}}\quad\mbox{and}\quad\| F\|_*\lesssim(\frac{\sqrt{t}}{r})^\frac{d}{2}.
\end{equation}
By definition of $\|\cdot\|$ for $s=2$, the second estimate implies by Chebychef's inequality
for all $\delta>0$
$\langle I(F-\langle F\rangle>\delta)\rangle\lesssim\exp(-\frac{1}{C}(\frac{r}{\sqrt{t}})^d\delta^2)$,
which in conjunction with the first estimate can be rewritten as
$\langle I(F>\delta)\rangle\lesssim\exp(-\frac{1}{C}(\frac{r}{\sqrt{t}})^d\delta^2)$
provided $\delta\gg\frac{1}{\sqrt{t}}$, which in turn yields (\ref{Cr*-9}).

\medskip

\step3 Conclusion.

Since $r_*$ is the minimal {\it dyadic} value greater or equal to one for which
\begin{equation}\label{Cr*-12}
\frac{1}{r^2}\fint_{B_r}|(\phi,\sigma)-\fint_{B_r}(\phi,\sigma)|^2\le\delta\quad\mbox{for all}\;r\ge r_*,
\end{equation}
it is enough to show
\begin{equation}\label{Cr*-13}
\langle I(r_*=r)\rangle\lesssim\exp(-\frac{1}{C}r^d)\quad\mbox{for any dyadic}\;r\gg 1,
\end{equation}
Indeed, specifying (\ref{Cr*-13}) to
\begin{equation}\label{Cr*-14}
\langle I(r_*=r)\rangle\le C_0\exp(-\frac{1}{C_0}r^d)\quad\mbox{for any dyadic}\;r\ge C_0^\frac{1}{d}
\end{equation}
for some specific constant $C_0$ with the above dependence, we recover (\ref{r*sharp}) in form of
\begin{eqnarray*}\nonumber
\langle\exp(\frac{1}{2C_0}r_*^d)\rangle&=&
\sum_{r\ge 1,\mbox{dyadic}}\exp(\frac{1}{2C_0}r^d)\langle I(r_*=r)\rangle\\
&\stackrel{(\ref{Cr*-14})}{\le}&\sum_{r\ge C_0^\frac{1}{d},\mbox{dyadic}}C_0\exp(-\frac{1}{2C_0}r^d)
+\sum_{r< C_0^\frac{1}{d},\mbox{dyadic}}\exp(\frac{1}{2}).
\end{eqnarray*}

\medskip

Let us fix such a dyadic $r\ge 1$ for (\ref{Cr*-13}).
By definition (\ref{Cr*-12}) of the minimal radius $r_*$, the event $r_*=r$ entails the event
\begin{equation}\nonumber
r_*\le r \quad\mbox{and}\quad 
\frac{1}{(r/2)^2}\fint_{B_{r/2}}|(\phi,\sigma)-\fint_{B_{r/2}}(\phi,\sigma)|^2>\delta,
\end{equation}
which because of $\fint_{B_{r/2}}\cdot\le 2^d\fint_{B_R}\cdot$ implies
\begin{equation}\nonumber
r_*\le r \quad\mbox{and}\quad 
\frac{1}{r^2}\fint_{B_{r}}|(\phi,\sigma)-\fint_{B_{r}}(\phi,\sigma)|^2\ge 2^{-(d+2)}\delta.
\end{equation}
In view of Lemma~\ref{Laux8bis}, this event entails the event
\begin{equation}\nonumber
\sup_{R\ge r,\;\mbox{dyadic}}(\frac{R}{r})^\beta\int\eta_R F_{T}\gtrsim\delta,
\end{equation}
where we use the abbreviation (\ref{Cr*-3}). This in turn implies that
\begin{equation}\nonumber
F\gtrsim \delta
\end{equation}
in the notation introduced in (\ref{Cr*-8}). Therefore (\ref{Cr*-9}) yields (\ref{Cr*-14}).


\section{Proofs of the deterministic auxiliary results}

\subsection{Proof of Proposition~\ref{Paux1}: deterministic propagator estimate}

By the triangle inequality,
\begin{multline*}
(\frac{R}{\sqrt{T}})^{\frac{d}{2}} |(q(T)-S_{\frac T2 \to T}^\ho q(\tfrac T2))_R|
\,\leq \, (\frac{R}{\sqrt{T}})^{\frac{d}{2}} |(S_{\frac T2\to T}-S_{\frac T2 \to T}^h) q(\tfrac T2))_R|
\\+
(\frac{R}{\sqrt{T}})^{\frac{d}{2}} |(S_{\frac T2\to T}^h-S_{\frac T2 \to T}^\ho) q(\tfrac T2))_R|.
\end{multline*}
By Lemmas~\ref{Laux9} and \ref{Laux6}, we may bound the second rhs term by
\begin{multline*}
(\frac{R}{\sqrt{T}})^{\frac{d}{2}} |(S_{\frac T2\to T}^h-S_{\frac T2 \to T}^\ho) q(\tfrac T2))_R|\\
\lesssim \, |a_{hT}-a_\ho| \fint_{\frac{T}{4}}^\frac{T}{2}dt\fint_0^{\sqrt{T}}dr(\frac{r}{\sqrt{t}})^{\frac{d}{2}}\int\eta_{2\sqrt{T}}|(q(t)-\langle q(t)\rangle)_r|,
\end{multline*}
and $|a_{hT}-a_\ho|$ by Lemma~\ref{Laux10}.

Hence, it remains to prove that for some $p\gg 1$ we have
\begin{multline}\label{a-prop1.1}
(\frac{R}{\sqrt{T}})^{\frac{d}{2}} |(S_{\frac T2\to T}-S_{\frac T2 \to T}^h) q(\tfrac T2))_R|
\\
\lesssim \delta^\frac1p \fint_{\frac{T}{4}}^\frac{T}{2}dt\fint_0^{\sqrt{T}}dr(\frac{r}{\sqrt{t}})^{\frac{d}{2}}\int\eta_{2\sqrt{T}}|(q(t)-\langle q(t)\rangle)_r|
\end{multline}
provided $\delta \ll 1$ is such that
\begin{equation}\label{a-prop1.2}
\Big( \int\eta_{\sqrt{T}}\Big|(\frac{\phi_T}{\sqrt{T}},\frac{\sigma_T}{\sqrt{T}},g_T) \Big|^2 \Big)^\frac12\,\leq \, 2\delta,
\end{equation}
which, by Lemma~\ref{Laux8}, holds if we have for some $t_0\le \delta^2 T$
$$
\Big(\langle|(q_{t_0}-\langle q_{t_0}\rangle)_{\delta\sqrt{t_0}}|\rangle+\int\eta_{\sqrt{T}}
|(q_{t_0}-\langle q_{t_0}\rangle)_{\delta\sqrt{t_0}}|\Big) \leq \delta^p
$$
for (some generic) $p\gg 1$.
On the one hand,  \eqref{Laux5.1} in Lemma~\ref{Laux5}, combined with \eqref{Laux6.2} in Lemma~\ref{Laux6}, takes the form 
\begin{multline}\label{a-prop1.3}
\Big(\frac{R}{\sqrt{T}}\Big)^{\frac{d}{2}} \Big| \Big( (S_{\frac{T}{2}\to T}-S^h_{\frac{T}{2}\to T})(q(\frac{T}{2}))\Big)_R\Big|
\\
\lesssim \,
\delta^{\frac{1}{p} }\frac{\sqrt{T}}{R}
\fint_{\frac{T}{4}}^\frac{T}{2}dt\fint_0^{\sqrt{T}}dr(\frac{r}{\sqrt{t}})^{\frac{d}{2}}\int\eta_{2\sqrt{T}}|(q(t)-\langle q(t)\rangle)_r|
\end{multline}
provided \eqref{a-prop1.2} holds.
In order to establish \eqref{a-prop1.1}, we combine this inequality with the following alternative estimate
\begin{multline}\label{a-prop1.4}
\Big(\frac{R}{\sqrt{T}}\Big)^{\frac{d}{2}} \Big| \Big( (S_{\frac{T}{2}\to T}-S^h_{\frac{T}{2}\to T})(q(\frac{T}{2}))\Big)_R\Big|
\\
\lesssim \,
\Big( \frac{R}{\sqrt{T}}\Big)^\frac1p
\fint_{\frac{T}{8}}^\frac{T}{2}dt\fint_0^{\sqrt{T}}dr(\frac{r}{\sqrt{t}})^{\frac{d}{2}}\int\eta_{2\sqrt{T}}|(q(t)-\langle q(t)\rangle)_r|,
\end{multline}
for some $p\gg 1$ depending only on $\lambda$ and $d$, and which follows from the Meyers estimate of Lemma~\ref{Laux7} as we argue below.
Let us directly draw the conclusion.
Since 
$$
\min\{\delta^{\frac{1}{p}} \frac{\sqrt{T}}{R}, \Big( \frac{R}{\sqrt{T}} \Big)^\frac1p
\} \le \delta^{\frac{1}{p(p+1)}}
$$
the combination of \eqref{a-prop1.3} and \eqref{a-prop1.4} yields \eqref{a-prop1.1} for some generic $p\gg 1$.

\medskip

We conclude with the proof of \eqref{a-prop1.4}.
As noted at the beginning of the proof of Lemma \ref{Laux2},
we have the representation
\begin{align*}
q(T)-S^h_{\frac{T}{2}\rightarrow T}q({\textstyle\frac{T}{2}})
=a\int_{\frac{T}{2}}^T\nabla ud\tau-a_h\int_{\frac{T}{2}}^T\nabla u_h d\tau,
\end{align*}
where $u_h$ is the solution of the initial value problem
\begin{align*}
\partial_\tau u_h-\nabla\cdot a_{hT}\nabla u_h=0\;\;\mbox{for}\;\tau>\tfrac{T}{2},\quad
u_h= u\;\;\mbox{for}\;\tau=\tfrac{T}{2}.
\end{align*}
Hence we have by the triangle inequality
\begin{align*}
|(q(T)-S^h_{\frac{T}{2}\rightarrow T}q({\textstyle\frac{T}{2}}))_R|
\le\int_{\frac{T}{2}}^T|(a\nabla u)_R|d\tau+\int_{\frac{T}{2}}^T|(\nabla\tilde u)_R|d\tau,
\end{align*}
where $\tilde u$ solves the homogeneous equation $\partial_\tau\tilde u-\nabla\cdot \tilde a\nabla\tilde u=0$
with the space-time varying coefficient field $\tilde a$ given by
\begin{align}\label{au46}
\tilde a(\tau)=a_{hT}\;\;\mbox{for}\;\tau>\tfrac{T}{2},\quad
\tilde a(\tau)=a\;\;\mbox{for}\;\tau<\tfrac{T}{2}.
\end{align}
By the dominance of the Gaussian kernel by the exponential kernel (and Cauchy-Schwarz' inequality, the Meyers estimate in form of Lemma~\ref{Laux7} applies to both $u$ and $\tilde u$, each a solution of a homogeneous
parabolic equation with uniformly elliptic (space-time-varying) coefficient fields, and entails
the equi-integrability of the square gradient in form of
\begin{align*}
\lefteqn{\int_{\frac{T}{2}}^T|(a\nabla\tilde u)_R|d\tau}\\
&\lesssim\sqrt{T}\big(\int_{\frac{T}{2}}^T\int\eta_R|\nabla\tilde u|^2d\tau\big)^\frac{1}{2}
\lesssim\sqrt{T}(\frac{\sqrt{T}}{R})^{\frac{d}{2}-\frac{1}{p}}
\big(\int_{\frac{T}{4}}^T\int\eta_{\sqrt{T}}|\nabla\tilde u|^2d\tau\big)^\frac{1}{2},
\end{align*}
where $p$ denotes a large generic exponent that now depends also on $\lambda>0$ next to $d$;
and the analogue estimate for $u$. Hence we obtain by the localized energy estimate \eqref{Laux3.2} in Lemma~\ref{Laux3}
 for (\ref{au46})
\begin{align}\label{au49}
(\frac{R}{\sqrt{T}})^{\frac{d}{2}-\frac{1}{p}}|(q(T)-S^h_{\frac{T}{2}\rightarrow T}q({\textstyle\frac{T}{2}}))_R|
\lesssim
\sqrt{T}\big(\int\eta_{\sqrt{T}}|u({\textstyle\frac{T}{4}})|^2\big)^\frac{1}{2}.
\end{align}
Combined with \eqref{Laux6.2} in Lemma~\ref{Laux6}, this yields the desired estimate~\eqref{a-prop1.4}.

\medskip

For future reference, note that this argument also yields 
\begin{multline}\label{Paux1-50}
R^{\frac d2} \int_{\frac T2}^T (|(a\nabla u)_R|+|(\nabla u)_R|)d\tau \\
\lesssim \, \Big( \frac{R}{\sqrt{T}}\Big)^\frac1p
\fint_{\frac{T}{8}}^\frac{T}{2}dt\fint_0^{\sqrt{T}}dr \int\eta_{2\sqrt{T}}r^{\frac{d}{2}} |(q(t)-\langle q(t)\rangle)_r|.
\end{multline}
%


\subsection{Proof of Lemma~\ref{Laux3}: localized elliptic and parabolic energy estimates}

W.l.o.g. we take $T=1$.

\medskip

\step 1 Proof of \eqref{Laux3.1}.

Consider $v$ solution of 
$v-\nabla\cdot a\nabla v=f+\nabla\cdot g$. 
We test this equation with $\eta_R v$
with scale $R$ to be specified below and obtain
\begin{equation}\nonumber
\int(\eta_Rv^2+\nabla(\eta_Rv)\cdot a\nabla v)=\int(\eta_Rvf-\nabla(\eta_Rv)\cdot g).
\end{equation}
Using Leibniz' rule and the bounds (\ref{1.3}) on $a$, this yields
\begin{equation}\label{Laux3-1}
\int\eta_R(v^2+|\nabla v|^2)\lesssim\int|\nabla\eta_R||v|(|\nabla v|+|g|)+\eta_R(|v||f|+|\nabla v||g|).
\end{equation}
Because of $|\nabla\eta_R|\le\frac{1}{R}\eta_R$ for our exponential averaging function,
we may absorb the first rhs term into the lhs for $R\gg 1$
to the effect of
\begin{equation}\nonumber
\int\eta_R(v^2+|\nabla v|^2)\lesssim\int\eta_R\big(|v|(|f|+|g|)+|\nabla v||g|\big),
\end{equation}
which with help of Young's inequality turns into
\begin{equation*}
\int\eta_R(v^2+|\nabla v|^2)\lesssim\int\eta_R(|f|^2+|g|^2)\quad\mbox{for}\;R\ge R_0
\end{equation*}
for some $R_0\le C$, that is, \eqref{Laux3.1}.

\medskip

\step 2 Proof of \eqref{Laux3.2}.

Proceeding as in Step~1, we obtain
\begin{eqnarray*}
\frac{d}{dt}\int \eta_R \frac12 v^2 &=&\int \eta_R v\partial_t v
\\
&\leq & \int \eta_R (-\lambda |\nabla v|^2+\frac{1}{R}|v||\nabla v|+|v||f|+|\nabla v||g|+\frac1R|v||g|).
\end{eqnarray*}
For $R\gg 1$, this yields for some $C(\lambda)\ge 2$,
\begin{equation}\label{a-aux3.1}
\frac{d}{dt}\int \eta_R \frac12 v^2 \,\leq \, \int \eta_R (-\frac{1}{C}|\nabla v|^2+Cv^2 + C(f^2+|g|^2)).
\end{equation}
We then multiply by $\exp(-Ct)$, absorb the second rhs term in the derivative of the lhs, and integrate in time from $0$ to $t$ for all $t\le 1$ to get
\begin{multline*}
\sup_{t \le 1} \exp(-C t)\int \eta_R \frac12 v(t)^2+\frac1C \int_0^1\exp(-C\tau)\int \eta_R |\nabla v|^2d\tau  \\
\leq \, C\int_0^1 \exp(-C\tau)\int (f^2+|g|^2)d\tau+ \int \eta_R \frac12 |v_0|^2,
\end{multline*}
which yields the claim.


\subsection{Proof of Lemma~\ref{Laux1}: localized semi-group estimates}

We split the proof into three steps.

\medskip

\step{1} Reductions.

By scaling, it is again enough to consider $\sqrt{T}=1$ and $R\ge 1$.
Let $U(t):=\int_0^t vd\tau$.
We split estimate \eqref{Laux1.1} into four parts:
\begin{eqnarray}
\sup_{t\le 1} \int \eta_R |U|^2 &\lesssim & \int\eta_R|q_0|^2,\label{a-aux1.1}\\
\int_0^1\int\eta_R |\nabla U|^2dt&\lesssim&\int\eta_R|q_0|^2,\label{a-aux1.2}\\
\sup_{t\in[\frac{1}{2},1]}\int\eta_R|\nabla v|^2&\lesssim&\int_0^1 t\int\eta_R v^2dt,\label{a-aux1.3}\\
\int_0^1 t\int\eta_R v^2dt&\lesssim&\int \eta_R |q_0|^2,\label{a-aux1.4}
\end{eqnarray}
and shall obtain \eqref{Laux1.2} as a by-product of the proof of \eqref{a-aux1.3}.
In order to deduce \eqref{Laux1.1} from \eqref{a-aux1.1}--\eqref{a-aux1.4}, we need to upgrade
the time-averaged estimates \eqref{a-aux1.2} \& \eqref{a-aux1.4} into pointwise-in-time estimates.
For \eqref{a-aux1.4} we appeal to \eqref{a-aux3.1} in the proof of Lemma~\ref{Laux3} which implies that 
$t\mapsto \exp(-{C}t) \int\eta_R v^2(t)$ 
is non-increasing, and therefore allows one to upgrade \eqref{a-aux1.4} to $\int \eta_R v^2(T)\lesssim \int \eta_R |q_0|^2$ for all $R\ge 1$.
For \eqref{a-aux1.2} we first argue that 
\begin{equation}\label{a-aux1.5}
\int\eta_R|\partial_t \nabla U|^2{\lesssim} \,
\int \eta_R |q_0|^2
\end{equation}
which follows from the identity $\int\eta_R|\partial_t \nabla U|^2=\int\eta_R|\nabla v|^2$ combined with \eqref{a-aux1.3} \& \eqref{a-aux1.4}.
The desired estimate $\int\eta_R |\nabla U|^2(T)\,\lesssim\,\int\eta_R|q_0|^2$ then follows from this Lispchitz-type estimate \eqref{a-aux1.5} and the 
time-averaged estimate \eqref{a-aux1.2}.

\medskip

The proof of \eqref{a-aux1.1} \& \eqref{a-aux1.2} is elementary and deferred to the end of this step.
The proof of \eqref{a-aux1.3} \& \eqref{a-aux1.4} is more subtle. For future reference, we state and prove a slightly more general
version of \eqref{a-aux1.3}  \& \eqref{Laux1.2} in Step~2, from which we deduce \eqref{a-aux1.4} in Step~3 by duality.
Let us now prove \eqref{a-aux1.1} \& \eqref{a-aux1.2}.
Since $U$ is a solution of 
\begin{equation}\nonumber
\partial_t U-\nabla\cdot a\nabla U=\nabla\cdot q_0,
\end{equation}
the estimate \eqref{a-aux3.1} with $f=0$, cf~proof of Lemma~\ref{Laux3}, takes the form for $R\ge 1$
\begin{equation}\label{b86}
 \frac{d}{dt}\exp(-Ct)\int\eta_R |U|^2+\frac{1}{C}\exp(-Ct)\int\eta_R|\nabla U|^2\le C\exp(-Ct)\int\eta_R |q_0|^2.
\end{equation}
In view of $U(t=0)=0$, this yields 
\begin{equation}\label{b85}
\sup_{t\le 1}\int \eta_R |U|^2 +\int_0^1\int\eta_R|\nabla U|^2dt\lesssim \int_0^1\int\eta_R|q_0|^2dt=\int\eta_R|q_0|^2,
\end{equation}
that is,  \eqref{a-aux1.1} \& \eqref{a-aux1.2}.

\medskip

\step{2} Localized parabolic higher-order energy estimates: For all $\tau \in (0,1)$, 
$\beta\ge 0$, and $R\ge 1$,
\begin{eqnarray}\label{p1}
\int\eta_R|\nabla w|^2(\tau) &\lesssim& \frac{1}{\tau} \int\eta_R w^2(\frac{\tau}{2}),\\
\int\eta_R|\nabla w|^2(1)& \lesssim& \int_0^\frac{1}{2}(2t)^\beta \int \eta_R w^2(t)dt,  \label{p1bis}
\end{eqnarray}
for any solution
of  the homogeneous equation $\partial_t w-\nabla \cdot a \nabla w=0$.
Estimate \eqref{p1} implies \eqref{Laux1.2}, whereas estimate \eqref{p1bis} implies \eqref{a-aux1.3} by rescaling of time.

The argument for \eqref{p1} \& \eqref{p1bis} would be straightforward for symmetric $a$
in which case one computes $\frac{d}{dt}\frac{1}{2}\int\eta\nabla w\cdot a\nabla w$
to find that this energy is non-increasing up to lower-order terms, ie terms that involve
$\nabla\eta_R$. In the general case, one has to pass via $\frac{d}{dt}\frac{1}{2}\int\eta_R w^2$
and $\frac{d}{dt}\frac{1}{2}\int\eta_R(\nabla\cdot a\nabla w)^2$, and their interpolation, instead. 
Indeed, we have from
the equation $\partial_tw-\nabla\cdot a\nabla w=0$, integration by parts, and Leibniz' rule
\begin{equation}\label{p2.bis}
\frac{d}{dt}\frac{1}{2}\int\eta_R w^2
=-\int(\eta_R\nabla w\cdot a\nabla w+w\nabla\eta_R\cdot a\nabla w),
\end{equation}
and thus by the bounds (\ref{1.3}) on $a$, the property $|\nabla\eta_R|\le\eta_R$ for our
exponential localization function with $R\ge 1$, and Young's inequality
\begin{equation}\label{p2}
\frac{d}{dt}\int\eta_R w^2
\le-\int\eta_R(\lambda|\nabla w|^2-\frac{1}{\lambda}w^2).
\end{equation}
Since $a$ does not depend on time, $\partial_t$ commutes with the parabolic operator $\partial_t-\nabla\cdot a
\nabla$
so that $V:=\nabla\cdot a\nabla w=\partial_tw$ is also a solution of the homogeneous equation to
which we may apply (\ref{p2}). For a parameter $\delta>0$ still to be chosen, this yields
in combination
\begin{align}\label{p4}
\lefteqn{\frac{d}{dt}\big(\delta t^2\int\eta_R V^2+\int\eta_R w^2\big)}\\
&\le-\lambda\big(\delta t^2\int\eta_R|\nabla V|^2+\int\eta_R|\nabla w|^2\big)
+2\delta t\int\eta_R V^2+\frac{1}{\lambda}\big(\delta t^2\int\eta_R V^2+\int\eta_R w^2\big).\nonumber
\end{align}

\medskip

In order to absorb $\int\eta_R V^2$ into $\int\eta_R|\nabla V|^2$ and $\int\eta_R|\nabla w|^2$,
we appeal to the following interpolation estimate
\begin{equation}\label{p3}
\int\eta_R V^2\le2\big(\int\eta_R|\nabla V|^2\int\eta_R|\nabla w|^2\big)^\frac{1}{2}+\int\eta_R|\nabla w|^2,
\end{equation}
which follows from integration by parts and Leibniz' rule in form of 
$\int\eta_R V^2=-\int(\eta_R\nabla V\cdot a\nabla w+V\nabla\eta_R\cdot a\nabla w)$, the upper bound (\ref{1.3})
on $a$ and the choice of $\eta_R$, which yield
$\int\eta_R V^2\le\int\eta_R(|\nabla V||\nabla w|+|V||\nabla w|)$, and the inequalities of Cauchy-Schwarz
and Young. By choosing $\delta\ll 1$, we learn from applying once more Young's inequality to
the rhs of (\ref{p3}) that (\ref{p4}) turns into the differential inequality
\begin{align}\label{p5}
\frac{d}{dt}\big(\delta t^2\int\eta_R V^2+\int\eta_R w^2\big)
\le\frac{1}{\lambda}\big(\delta t^2\int\eta_R V^2+\int\eta_R w^2\big).
\end{align}

\medskip

We now integrate (\ref{p5}) rewritten in form of 
$\frac{d}{dt}\big(\exp(-\frac{t}{\lambda})\big(\delta t^2\int\eta_R V^2+\int\eta_R w^2\big)\big)\le0$ over $t\in(0,\tau)$ 
to obtain 
\begin{equation}\label{p7}
\big(\tau^2\int\eta_R V^2(\tau)+\int\eta_R w^2(\tau)\big)\lesssim\int\eta_R w^2(0).
\end{equation}
We use once more an interpolation inequality, this time
\begin{equation}\label{p6}
\lambda\tau \int\eta_R|\nabla w|^2(\tau)\le2\big(\tau^2\int\eta_R V^2(\tau)\int\eta_R w^2(\tau)\big)^\frac{1}{2}+\frac{\tau}{\lambda}\int\eta_R w^2(\tau),
\end{equation}
which follows from the lower bound on $a$ in (\ref{1.3}) in form of $\lambda\int\eta_R|\nabla w|^2
\le\int\eta_R\nabla w\cdot a\nabla w$, integration by parts and Leibniz' rule in form of
$\int\eta_R\nabla w\cdot a\nabla w=-\int\eta_R w V-\int w\nabla\eta_R\cdot a\nabla w$,
the upper bound on $a$ in (\ref{1.3}), $|\nabla\eta_R|\le\eta_R$, and the inequalities of Cauchy-Schwarz
and Young. Applying once more Young's inequality to (\ref{p6}) and inserting the result into (\ref{p7})
yields for all $\tau\in (0,1)$,
\begin{equation}\label{p2.ter}
\tau \int\eta_R|\nabla w|^2(\tau) \lesssim\int\eta_R w^2(0).
\end{equation}
By a rescaling in space and time we may use this with the time interval $(0,\tau)$ replaced
by $(\frac{\tau}{2},\tau)$, which yields (\ref{p1}).
In combination with (\ref{p2}), we obtain  (\ref{p1bis}).
Indeed, (\ref{p2}) in form of $\frac{d}{dt}\big(\exp(-\frac{t}{\lambda}) \int\eta_R w^2\big)\le0$ 
implies that, for all $\beta\ge 0$, $\int\eta_R w^2(\frac{1}{2})
\lesssim\int_0^\frac{1}{2}(2t)^\beta \int \eta_R w^2dt$, from which the claim follows.

\medskip

\step{3} Proof of \eqref{a-aux1.4}.

The starting point is \eqref{p1bis} in the simple form of 
\begin{equation}\label{p32}
\int\eta_R|\nabla w|^2_{|t=1}\lesssim\int_0^1\int\eta_R w^2dt
\end{equation}
for any solution of the homogeneous equation $\partial_tw-\nabla\cdot a\nabla w=0$
and any $R\ge 1$. 
Since the only property used on the cut-off is the inequality $|\nabla \eta_R|\le \eta_R$, 
(\ref{p32}) also holds with $\eta_R$ replaced by $\frac{1}{\eta_R}$ for any $R\ge 1$:
\begin{equation}\nonumber
\int\frac{1}{\eta_R}|\nabla w|^2_{|t=1}\lesssim\int_0^1\int\frac{1}{\eta_R} w^2dt.
\end{equation}
If in addition, the initial data are given by $w(t=0)=w_0$, we obtain from (\ref{p2})
with $\eta_R$ replaced by $\frac{1}{\eta_R}$ that $\int_0^1\int\frac{1}{\eta_R} w^2dt
\lesssim\int\frac{1}{\eta_R}w_0^2$ so that the above turns into
\begin{equation}\nonumber
\int\frac{1}{\eta_R}|\nabla w|^2_{|t=1}\lesssim\int\frac{1}{\eta_R} w_0^2.
\end{equation}
By parabolic rescaling, this yields $t\int\frac{1}{\eta_{\sqrt{t}\rho}}|\nabla w(t)|^2
\lesssim\int\frac{1}{\eta_{\sqrt{t}\rho}} w_0^2$. For $t\le 1$ we may choose $\rho=\frac{R}{\sqrt{t}}\ge 1$
to the effect of
\begin{equation}\label{p33}
t\int\frac{1}{\eta_R}|\nabla w(t)|^2\lesssim\int\frac{1}{\eta_R} w_0^2.
\end{equation}
Clearly, (\ref{p33}) also holds with $a$ replaced by its transpose $a^*$;
note that transposition exactly preserves (\ref{1.3}).
If now we have a solution of the inhomogeneous $\partial_t W-\nabla\cdot a^*\nabla W=f$ with
homogeneous initial data $W_{|t=0}=0$,
we may write $W$ in terms of the semi-group $s_*(t)$ for $-\nabla\cdot a^*\nabla$ as
$W_{|t=1}=\int_0^1s_*(1-t)f(t)dt$, so that (\ref{p33}), which amounts to the semi-group estimate
$t\int\frac{1}{\eta}|\nabla s_*(t) w_0|^2\lesssim\int\frac{1}{\eta} w_0^2$, entails
\begin{equation}\nonumber
\int\frac{1}{\eta}|\nabla W|^2_{|t=1}\lesssim\int\frac{1}{1-t}\int\frac{1}{\eta}|f|^2dt.
\end{equation}
The dualization (by $t\leadsto 1-t$) of this estimate yields for our solution $v$ of the homogeneous 
equation $\partial_t v-\nabla\cdot a\nabla v=0$ with inhomogeneous initial data $u(t=0)=\nabla\cdot q_0$
that
\begin{equation}\nonumber
\int t\int\eta|v|^2dt\lesssim\int\eta|q_0|^2,
\end{equation}
which turns into (\ref{a-aux1.4}).


\subsection{Proof of Lemma~\ref{Laux2}: semi-group property}

Let the (spatial) vector field $q_0$ be given and set
$q_1:=S_{t_0\rightarrow t_1}q_0$. Let  the space-time scalar field $v_i$,
where $i=0,1$, denote the solution of the initial value problem
\begin{align}\label{a-au24}
\partial_\tau v_i-\nabla\cdot a\nabla v_i=0\;\;\mbox{for}\;\;\tau>t_i,\quad
v_i=\nabla\cdot q_i\;\;\mbox{for}\;\;\tau=t_i.
\end{align}
By definition of $S_{t_0\rightarrow t_1}$, we have
\begin{align}\label{a-au25}
q_1=q_0+\int_{t_0}^{t_1}a\nabla v_0d\tau.
\end{align}
Based on this we now argue that
\begin{align}\label{a-au26}
v_1=v_0\;\;\mbox{for}\;\;\tau>t_1.
\end{align}
Indeed, applying $\int_{t_0}^{t_1}d\tau$ to (\ref{a-au24}) for $i=0$, we obtain
$v_0(t_1)-\nabla\cdot q_0-\nabla\cdot\int_{t_0}^{t_1}a\nabla v_0d\tau=0$, which
we rewrite as $v_0(t_1)=\nabla\cdot(q_0+\int_{t_0}^{t_1}a\nabla v_0d\tau)\stackrel{(\ref{a-au25})}{=}\nabla\cdot q_1$.
Hence by uniqueness for the initial value (\ref{a-au24}) for $i=1$, we obtain (\ref{a-au26}). 
Equipped with (\ref{a-au26}), we are in the position to conclude: By definition of $S_{t_1\rightarrow t_2}$ we have
$S_{t_1\rightarrow t_2}q_1=q_1+\int_{t_1}^{t_2}a\nabla v_1d\tau$, into which we plug (\ref{a-au25}):
\begin{align*}
S_{t_1\rightarrow t_2}q_1=q_0+\int_{t_0}^{t_1}a\nabla v_0d\tau+\int_{t_1}^{t_2}a\nabla v_1d\tau
\stackrel{(\ref{a-au26})}{=}q_0+\int_{t_0}^{t_2}a\nabla v_0d\tau=S_{t_0\rightarrow t_2}q_0.
\end{align*}


\subsection{Proof of Lemma~\ref{Laux4}: formulas}

We split the proof into four steps.

\medskip

\step 1 Proof of \eqref{Laux4.1}

We apply $\int_0^\infty d\tau \exp(-\frac{\tau}{T})$ to the parabolic equation \eqref{e.1} \& \eqref{e.2} to the effect of
$$
0\,=\,\frac1T \int_0^\infty d\tau \exp(-\frac\tau T)u(\tau)-\nabla \cdot (a e)-\nabla\cdot a\nabla \int_0^\infty d\tau \exp(-\frac\tau T)u(\tau).
$$
From the uniqueness of the solution of \eqref{c62} (cf estimate \eqref{Laux3.1} in Lemma~\ref{Laux3}), we then learn that
\begin{equation}\label{Laux4-20}
\phi_T =  \int_0^\infty d\tau \exp(-\frac \tau T) u(\tau).
\end{equation}
Applying $a(\nabla +e)$, we obtain
\begin{eqnarray*}
q_T&=&a\Big(\int_0^\infty d\tau \exp(-\frac \tau T)\nabla u(\tau)+e\Big)
\\
&=&a \Big(\int_0^\infty \frac{dt}{T}\exp(-\frac t T) \int_0^t d\tau \nabla u(\tau)+\int_0^\infty \frac{dt}{T}\exp(-\frac t T)e \Big),
\end{eqnarray*}
so that \eqref{Laux4.1} follows from the definition $q(t)=a(\int_0^t \nabla ud\tau+e)$.

\medskip

\step 2 Proof of \eqref{Laux4.2}

By (\ref{a66}), \eqref{Laux4.2} takes the form
\begin{equation}\label{a-aux4.1}
q_{Ti}=\langle q_{Ti}\rangle+\nabla\cdot\sigma_{Ti}+g_{Ti}.
\end{equation}
We fix and drop the indices $Ti$.
Taking the divergence of (\ref{o17}), ie taking the derivative wrt $x_k$ and summing over $k=1,\cdots,d$ yields
\begin{equation}\nonumber
\frac{1}{T}(\nabla\cdot\sigma)_j-\triangle(\nabla\cdot\sigma)_{j}
=\partial_j\nabla\cdot q-\triangle q_{j}.
\end{equation}
Inserting (\ref{c62}) yields the identity of vector fields
\begin{equation}\nonumber
\nabla\cdot\sigma-T\triangle\nabla\cdot\sigma
=\nabla\phi-T\triangle q.
\end{equation}
Adding (\ref{a4}) yields
\begin{equation}\nonumber
(\mathrm{id}-T\triangle)(\nabla\cdot\sigma+g)
=q-\langle q\rangle-T\triangle q,
\end{equation}
which implies (\ref{a-aux4.1}) and thus (\ref{Laux4.2}) by invertibility of ${\rm id}-T\triangle$ on bounded fields.

\medskip

\step 3 Proof of \eqref{Laux4.3} \& \eqref{Laux4.3b}.

For $\xi\in\mathbb{R}^d$ we define
$\phi_{T\xi}:=\xi_i\phi_{Ti}$, so that from (\ref{c62}) we obtain by linearity
$\frac{1}{T}\phi_{T\xi}-\nabla\cdot a(\nabla\phi_{T\xi}+\xi)=0$. Since by uniqueness
the solution $\phi_{T\xi}(a,x)$ is shift-covariant or stationary in the sense of
$\phi_{T\xi}(a,x+z)=\phi_{T\xi}(a(\cdot+z),x)$, we get from the stationarity of $\langle\cdot\rangle$
\begin{equation}\nonumber
\frac{1}{T}\langle\phi_{T\xi}^2\rangle+\langle\nabla\phi_{T\xi}\cdot a(\nabla\phi_{T\xi}+\xi)\rangle=0,
\end{equation}
which implies in particular
\begin{equation}\label{a-aux4.2}
\langle(\nabla\phi_{T\xi}+\xi)\cdot a(\nabla\phi_{T\xi}+\xi)\rangle
\le\xi\cdot\langle a(\nabla\phi_{T\xi}+\xi)\rangle\stackrel{(\ref{a66}),(\ref{c62})}{=}\xi\cdot a_{hT}\xi.
\end{equation}
On the one hand, (\ref{a-aux4.2}) yields the lower bound in (\ref{Laux4.3}). Indeed, using Jensen's inequality for
$\langle\cdot\rangle$ and stationarity of $\phi_\xi$ in form of $\langle\nabla\phi_{T\xi}\rangle
=\nabla\langle\phi_{T\xi}\rangle=0$, we have
\begin{equation}\nonumber
\xi\cdot a_{hT}\xi\stackrel{(\ref{a-aux4.2}),(\ref{1.3})}{\ge}\lambda\langle|\nabla\phi_{T\xi}+\xi|^2\rangle
\ge\lambda|\langle\nabla\phi_{T\xi}\rangle+\xi|^2=\lambda |\xi|^2.
\end{equation}
On the other hand, we obtain from (\ref{a-aux4.2}), followed by (\ref{1.3}) and Jensen's inequality
$$
\xi \cdot a_{hT} \xi \ge \langle(\nabla\phi_{T\xi}+\xi)\cdot a(\nabla\phi_{T\xi}+\xi)\rangle
\ge \langle |a(\nabla\phi_{T\xi}+\xi)|^2\rangle  \ge |\langle a(\nabla\phi_{T\xi}+\xi)\rangle|^2 = |a_{hT}\xi|^2,
$$
which is the upper bound in (\ref{Laux4.3}).

By qualitative stochastic homogenization, $a_\ho =\lim_{T\uparrow \infty} a_{hT}$, so that \eqref{Laux4.3b} follows
from \eqref{Laux4.3}.

\medskip

\step 4 Proof of \eqref{Laux4.4}, \eqref{Laux4.5}, and \eqref{Laux4.6}.

Identity (\ref{Laux4.4}) follows immediately from the definition of $w$.

\medskip

We now give the argument for (\ref{Laux4.5}), and eventually (\ref{Laux4.6}); for notational
simplicity we omit the index $T$. By definition of $w$, we obtain by Leibniz' rule
\begin{equation}\label{a8}
\nabla w=\nabla u-\partial_iv_h(e_i+\nabla\phi_i)-\phi_i\nabla\partial_iv_h.
\end{equation}
We obtain (\ref{Laux4.5}) by applying $a$ to this identity and using (\ref{c62}) \& (\ref{Laux4.2})
\begin{eqnarray}
a\nabla w&=&a\nabla u-\partial_iv_h(a_{hT}e_i+\nabla\cdot\sigma_i+g_i)-\phi_ia\nabla\partial_iv_h\nonumber\\
&=&a\nabla u-a_{hT}\nabla v_h-\partial_iv_h\nabla\cdot\sigma_i-\partial_i v_h g_i-\phi_ia\nabla\partial_iv_h\nonumber\\
&=&a\nabla u-a_{hT}\nabla v_h-\nabla\cdot(\partial_iv_h\sigma_i)-\partial_i v_h g_i
-(\phi_ia-\sigma_i)\nabla\partial_iv_h.\nonumber
\end{eqnarray}
Applying $-\nabla\cdot a$ to (\ref{a8}) yields by definition (\ref{c62}) of $q_i$
\begin{equation}\nonumber
-\nabla\cdot a\nabla w=-\nabla\cdot a\nabla u+(\nabla\partial_iv_h)\cdot q_i
+\partial_iv_h\nabla\cdot q_i
+\nabla\cdot(\phi_ia\nabla\partial_iv_h).
\end{equation}
We now make use of (\ref{Laux4.2}) and (\ref{c62}) to obtain
\begin{equation}\nonumber
-\nabla\cdot a\nabla w=-\nabla\cdot a\nabla u+(\nabla\partial_iv_h)\cdot(a_{hT}e_i+\nabla\cdot\sigma_i+g_i)
+\partial_iv_h\frac{1}{T}\phi_i
+\nabla\cdot(\phi_ia\nabla\partial_iv_h).
\end{equation}
Since $a_{hT}$ is constant, $(\nabla\partial_iv)\cdot a_{hT}e_i=\nabla\cdot a_{hT}\nabla v$.
In addition, we appeal to the general formula
$(\nabla\zeta)\cdot(\nabla\cdot\sigma)=-\nabla\cdot(\sigma\nabla\zeta)$
(which is a consequence of
the skew-symmetry of $\sigma$ in combination with the symmetry of $\nabla^2\zeta$) to obtain
\begin{equation}\nonumber
-\nabla\cdot a\nabla w=-\nabla\cdot a\nabla u+\nabla\cdot a_{hT}\nabla v_h
-\nabla\cdot(\sigma_i\nabla\partial_iv_h)
+g_i\cdot\nabla\partial_iv_h
+\partial_iv_h\frac{1}{T}\phi_i
+\nabla\cdot(\phi_ia\nabla\partial_iv_h).
\end{equation}
Taking the sum of this with $\partial_\tau w{=}\partial_\tau u-\partial_\tau v-\phi_{i}\partial_i\partial_\tau v$
and appealing to the defining equations for $v$ and $v_h$, we obtain (\ref{Laux4.6}).


\subsection{Proof of Lemma~\ref{Laux5}: estimate of intermediate homogenization error}

We split the proof into three steps.

\medskip

\step 1 Preliminaries.

By definition of the semi-groups, we have
$$
(S_{\frac T2\to T}-S^h_{\frac T2\to T})q(\tfrac T2)\,=\,\int_{\frac T2}^T (a\nabla u -a_{hT} \nabla u_h)d\tau,
$$
where
\begin{eqnarray}
\partial_\tau u-\nabla \cdot a \nabla u=0, & \tau>\tfrac T2, & \label{a-aux5.1} \\
\partial_\tau u_h-\nabla \cdot a_{hT} \nabla u_h=0,& \tau>\frac T2,&u_h(\tfrac T2)=u(\tfrac T2). \label{a-aux5.2}
\end{eqnarray}
W.l.o.g. we may assume $T=1$, and drop the index $T$.

\medskip

\step 2 Proof of 
\begin{multline} \label{a-aux5.3}
R^{\frac d2 +1} \Big(\int_{\frac 12}^1 |(a\nabla u-a_{h}\nabla u_h)_R|^2\Big)^\frac12
\,\lesssim \,\delta \Big[\Big(\int_{\frac 12}^1 \sup_x \eta |(\nabla,\nabla^2,\nabla^3)u_h|^2 d\tau \Big)^{\frac12} 
\\
+\Big(\sup_x \eta|\nabla {u_h}(\tfrac12)|^2\Big)^\frac12\Big].
\end{multline}
For future reference, we shall prove this estimate for any $u_h,u$ satisfying \eqref{a-aux5.1} \& \eqref{a-aux5.2} with $T=1$ and $u_h(\frac12)=u(\frac 12)$ (not necessarily given by $\nabla \cdot q(\frac12)$).

By \eqref{Laux4.6} in Lemma~\ref{Laux4}, $w=u-u_h$ satisfies
\begin{eqnarray*}
\tau>\frac12&:&\partial_\tau w-\nabla\cdot a\nabla w
=\, \nabla\cdot \underbrace{((\phi_{i}a-\sigma_{i})\nabla\partial_i v_h+ \partial_i v_h g_{i})}_{\dps =:g} -\underbrace{\phi_{i} \partial_\tau \partial_iv_h}_{\dps :=f},
\\
\tau=\tfrac12&:& w=-\phi_i \partial_i u_h,
\end{eqnarray*}
so that energy estimate \eqref{Laux3.2} in Lemma~\ref{Laux3} combined with $\eta_{\frac12}\sim \eta^2$ 
and $|\nabla \partial_\tau u_h|\lesssim |\nabla^3 u_h|$  yields
\begin{eqnarray}
\Big(\int_{\frac12}^1 \int \eta_{\frac 12} |\nabla w|^2 d\tau\Big)^\frac12 
&\lesssim &\Big( \int \eta_{\frac12}|w_{|\tau=\frac12}|^2\Big)^\frac12+
\Big( \int_{\frac12}^1 \int \eta_{\frac12} |(f,g)|^2d \tau \Big)^\frac12\nonumber
\\
&\lesssim & \Big(\int \eta |(\phi,\sigma,g)|^2\Big)^\frac12 \Big[ 
\Big(\int_{\frac12}^1 \sup_x  \eta |(\nabla,\nabla^2,\nabla^3 )u_h|^2d\tau \Big)^\frac12\nonumber 
\\
&&
\qquad\qquad \qquad\qquad\qquad +\Big(\sup_x \eta |\nabla {u_h}(\tfrac12)|^2\Big)^\frac12\label{a-aux5.4}
\Big].
\end{eqnarray}
We then appeal to \eqref{Laux4.5} in Lemma~\ref{Laux4} in form of  
$$
a\nabla u-a_{h}\nabla u_h\,=\,\underbrace{a\nabla w
+(\phi_{i}a-\sigma_{i})\nabla\partial_iu_h+\partial_i u_hg_{i}}_{\dps=:f}+\nabla\cdot \underbrace{(\partial_iv_h\sigma_{i})}_{\dps =:g}.
$$
On the one hand, since $R\le \sqrt{T}=1$, $G_R\lesssim G_{2 R} \eta_1$ 
so that by Cauchy-Schwarz' inequality
$$
|f_R| \lesssim  \Big(\int G_{2R}^2\Big)^\frac12  \Big(\int \eta^2 |f|^2\Big)^\frac12 \,\sim \,  \frac1{R^{\frac d2}} \Big(\int \eta_{\frac12} |f|^2\Big)^\frac12.
$$
On the other hand, by the semi-group property, $(\nabla \cdot g)_{\sqrt{2}\delta} = (\nabla \cdot g_\delta)_\delta$. 
Since the gradient of a Gaussian
is estimated by the Gaussian with twice the variance, cf $|\nabla G_{\delta}(z)|=|\frac{z}{\delta}G_{R}(z)| \lesssim G_{\sqrt{2}\delta}(z)$, this implies $|(\nabla \cdot g)_{\sqrt{2}\delta}|
\lesssim\frac{1}{\delta}|(g)_\delta|_{\sqrt{2}\delta}$, and thus by
Jensen's inequality $|\nabla \cdot g_{\sqrt{2}\delta}|^2
\lesssim\frac{1}{\delta^2}|(g)_\delta|_{\sqrt{2}\delta}^2$. Hence, 
arguing as above, and using Jensen's inequality in form of 
$ |g_\delta|^2\leq |g|_\delta^2$ followed by 
\eqref{c84} from the proof of Lemma~\ref{Laux6} below, we obtain
\begin{equation}\label{Laux5-1ref}
|(\nabla \cdot g)_{\sqrt{2} \delta}| \,\lesssim \, \frac{1}{\delta}|g_\delta|_{\sqrt{2}\delta} \,\lesssim \, 
\frac1\delta  \Big(\int G_{2\sqrt{2}\delta}^2\Big)^\frac12  \Big(\int \eta^2 |g_\delta|^2\Big)^\frac12
\,\stackrel{\eqref{c84}}{\lesssim}\, \frac{1}{\delta^{\frac d2+1}}   \Big(\int \eta^2 |g|^2\Big)^\frac12.
\end{equation}
Wiith $R=\sqrt{2}{\delta}$, these last two estimates combine to 
\begin{eqnarray*}
\lefteqn{|(a\nabla u-a_h \nabla u_h)_R|}
\\
 &\lesssim & \frac1{R^{\frac d2}} \Big(\int \eta_{\frac12} |f|^2\Big)^\frac12+\frac1{R^{\frac d2+1}} \Big( \int \eta_{\frac12} |g|^2\Big)^\frac12
\\
&\lesssim &\frac1{R^{\frac d2+1}}\Big[ \Big(\eta_{\frac12} |\nabla w|^2 \Big)^\frac12+\Big(\int \eta |(\phi,\sigma,g)|^2\Big)^\frac12 \Big(\sup_x \eta |(\nabla,\nabla^2)u_h|^2 \Big)^\frac12\Big].
\end{eqnarray*}
We integrate in time the square of this inequality, insert \eqref{a-aux5.4}, and recall the definition of $\delta$ to obtain \eqref{a-aux5.3}. 

\medskip

\step 3 Post-processing and proof that for all $\rho\le 1$,
\begin{equation}\label{a-aux5.5}
R^{\frac d2 +1} \Big(\int_{\frac 12}^1 |(a\nabla u-a_{h}\nabla u_h)_R|^2\Big)^\frac12
\,\lesssim \,(\frac1{\rho^{\frac d2+2}}\delta +\rho) \Big(\int \eta|\nabla {u}(\tfrac12)|^2\Big)^\frac12,
\end{equation}
which implies the claim \eqref{Laux5.1} by optimization in $\rho$.

In order to prove \eqref{a-aux5.5} we split the initial data into high and low pass:
$$
u_0:=u(\tfrac12)=u_0^{>}+u_0^{<}, \text{ with }u_0^<=(u_{0})_\rho, u_0^>=u_0-(u_{0})_\rho,
$$
where we recall that $(u_0)_\rho$ is the convolution of $u_0$ with the $\rho$-rescaled Gaussian.
We then solve
\begin{eqnarray*}
\partial_\tau u^{</>}-\nabla \cdot a \nabla u^{</>}\,=\,0\quad\text{for } \tau>\tfrac12, &&\partial_\tau u_h^{</>}-\nabla \cdot a_h \nabla u_h^{</>}\,=\,0\quad \text{for } \tau>\tfrac12,\\
u^{</>}=u_0^{</>} \quad \text{for } \tau=\tfrac12,&&u_h^{</>}=u_0^{</>} \quad \text{for } \tau=\tfrac12,
\end{eqnarray*}
so that (by uniqueness), $u=u^>+u^<$ and $u_h=u_h^>+u_h^<$, and thus 
$a\nabla u-a_h\nabla u_h=(a\nabla u^>-a_h\nabla u_h^>)+(a\nabla u^<-a_h\nabla u_h^<)$.

\medskip

We use a coarse estimate based on Lemma~\ref{Laux3} for the high-pass.
Since we may bound the Gaussian kernel by the exponential cut-off, we have
\begin{eqnarray*}
\lefteqn{R^\frac d2 \Big(\int_{\frac12}^1 |(a\nabla u^>-a_h\nabla u_h^>)_R|^2d\tau\Big)^\frac12}
\\
&\lesssim &\Big(\int_{\frac12}^1 \int \eta |a\nabla u^>-a_h\nabla u_h^>|^2d\tau\Big)^\frac12
\\
&\leq & \Big(\int_{\frac12}^1 \int \eta |\nabla u^>|^2d\tau\Big)^\frac12+\Big(\int_{\frac12}^1 \int \eta |\nabla u_h^>|^2d\tau\Big)^\frac12.
\end{eqnarray*}
By \eqref{Laux3.2} in Lemma~\ref{Laux3} followed by a convolution estimate, this turns into
\begin{equation}\label{a-aux5.6}
R^\frac d2 \Big(\int_{\frac12}^1 |(a\nabla u^>-a_h\nabla u_h^>)_R|^2d\tau\Big)^\frac12\,\stackrel{\eqref{Laux3.2}}{\lesssim} \, \Big(\int \eta | u_0-(u_0)_\rho|^2\Big)^\frac12
\,\lesssim \, \rho \Big(\int \eta |\nabla u_0|^2\Big)^\frac12.
\end{equation}

\medskip

For the low-pass, we use Step~2 in form of
\begin{multline} \label{a-aux5.7}
R^{\frac d2 +1} \Big(\int_{\frac 12}^1 |(a\nabla u^<-a_{h}\nabla u_h^<)_R|^2\Big)^\frac12
\,\lesssim \,\delta \Big[\Big(\int_{\frac 12}^1 \sup_x \eta |(\nabla,\nabla^2,\nabla^3)u_h^<|^2 d\tau \Big)^{\frac12} 
\\
+\Big(\sup_x \eta|\nabla {u_h^<}(\tfrac12)|^2\Big)^\frac12\Big],
\end{multline}
together with regularity theory for the constant-coefficients equation.
On the one hand, by Sobolev embedding applied to $\eta \nabla u_h^<$,
$$
\Big( \sup_x \eta |(\nabla,\nabla^2,\nabla^3)u_h^<|^2  \Big)^{\frac12} \,\lesssim \, \sum_{m=0}^{[\frac d2]+2}\Big(\int \eta |\nabla^m \nabla u_h^<|^2\Big)^\frac12.
$$
On the other hand, since $\nabla^m$ commutes with the constant-coefficient operator, we obtain by the localized energy estimate \eqref{Laux3.2} in Lemma~\ref{Laux3}
$$
\Big(\int_{\frac12}^1 \sup_x \eta|(\nabla,\nabla^2,\nabla^3)u_h^<|^2d\tau\Big)^\frac12
+\Big(\sup_x \eta|\nabla u_h^<(\tfrac12)|^2\Big)^\frac12
\,\lesssim \, \sum_{m=0}^{[\frac d2]+2}\Big(\int \eta |\nabla^m \nabla u_0^<|^2\Big)^\frac12.
$$
Recalling that $u_0^<=(u_0)_\rho$ for some $\rho\le 1$ (the convolution with the $\rho$-rescaled Gaussian), we have
by an elementary inverse estimate
$$
\Big(\int \eta |\nabla^m \nabla u_0^<|^2\Big)^\frac12\,\lesssim \, \frac1{\rho^m} \Big(\int \eta |\nabla u_0|^2\Big)^\frac12.
$$
Since $\rho\le 1$ and by definition of $u_0$, we may therefore upgrade \eqref{a-aux5.7} to 
\begin{equation*} \label{a-aux5.8}
R^{\frac d2 +1} \Big(\int_{\frac 12}^1 |(a\nabla u^<-a_{h}\nabla u_h^<)_R|^2\Big)^\frac12
\,\lesssim \,\delta  \frac1{\rho^{\frac d2+2}} \Big(\int \eta |\nabla u(\tfrac12)|^2\Big)^\frac12.
\end{equation*}
Combined with \eqref{a-aux5.6}, this yields the desired estimate \eqref{a-aux5.5}.


\subsection{Proof of Lemma~\ref{Laux6}: control of strong norms by averages}

By scaling we may wlog assume $T=1$.
The proof of \eqref{Laux6.1} is based on the localized parabolic estimates of Lemma~\ref{Laux1},
an interpolation estimate (cf Step~1), and an ODE argument (cf Step~2).
We deduce \eqref{Laux6.2} from \eqref{Laux6.1} in Step~3.
We shall also make use of some elementary statements to relate Gaussian and exponential kernels, 
that we defer to Step~4.

\medskip

\step{1} Interpolation estimate. For an arbitrary function $v$ we have
\begin{align}\label{b59}
\lefteqn{\Big(\int\eta v^2\Big)^\frac{1}{2}}\nonumber\\
&\lesssim\;
\Big[\Big(\int\eta |\nabla v|^2\Big)^\frac{1}{2}\Big]^{1-\frac{1}{\frac{d}{2}+p+2}}
\Big[\int_0^{1}drr^p\int\eta_2|v_{r}| \Big]^{\frac{1}{\frac{d}{2}+p+2}}\nonumber\\
&\;+\int_0^{1}drr^p\int\eta_2|v_{r}|.
\end{align}
Here comes the argument: For $r\ll 1$ we split $u$ into a low pass and a high pass
according to the convolution with the $r$-rescaled Gaussian. This yields
\begin{equation}\label{b57}
\int\eta v^2\lesssim\int\eta(v-v_{\sqrt{2}r})^2+\int\eta(v_{\sqrt{2}r})^2.
\end{equation}
The first rhs term is a regularization error. 
Since the exponential localization function $\eta$ dominates the Gaussian convolution kernel, it satisfies
\begin{equation}\label{b55bis}
\int\eta (v-v_{\sqrt{2}r})^2 \lesssim r^2\int\eta|\nabla v|^2.
\end{equation}
Indeed, if $G_\delta$ denotes the Gaussian kernel of scale $\delta$, we have
\begin{equation}\nonumber
(v_\delta-v)(x)=\int_0^1\int\nabla v(x+sz)\cdot zG_\delta(z)dzds;
\end{equation}
since $|zG_\delta(z)|\lesssim\delta G_{2\delta}(z)$, this yields
$|v_\delta-v|\lesssim\delta\int_0^1|\nabla v|_{2\delta s}ds$,
and thus by Jensen's inequality
$|v_\delta-v|^2\lesssim\delta^2\int_0^1|\nabla v|_{2\delta s}^2ds$ and therefore by the
symmetry of the convolution operator (and estimate \eqref{c84} from Step~4)
\begin{equation}\nonumber
\int\eta|v-v_\delta|^2\lesssim\delta^2\int_0^1\int\eta|\nabla v|_{2\delta s}^2ds
\stackrel{(\ref{c84})}{\lesssim}\delta^2\int\eta|\nabla v|^2,
\end{equation}
which is (\ref{b55bis}) for $\delta =\sqrt{2}r$.

We turn to the second rhs term and seek to replace the $L^2$-norm by an $L^1$-norm
(at the expense of decreasing the convolution scale and increasing the averaging scale)
in the sense of
\begin{equation}\label{b56bis}
\Big(\int\eta(v_{\sqrt{2}r})^2\Big)^\frac{1}{2}\lesssim\frac{1}{r^\frac{d}{2}}\int\eta_2|v_{r}|.
\end{equation}
Starting from the semi-group property and Jensen's inequality
in form of $|v_{\sqrt{2}r}|=|(v_{r})_{r}|\le|u_r|_r$ we obtain
since the Gaussian convolution kernel is dominated by the exponential localization function
of the same scale
\begin{equation}\nonumber
|v_{\sqrt{2}r}|\lesssim\int\eta_r|v_r|\lesssim\frac{1}{r^d}\int\eta_2|v_r|.
\end{equation}
Evoking translation invariance to replace the origin by a general point $y$, this implies
by the specific properties of the exponential localization function
\begin{equation}\nonumber
\eta_2(y)|(v_{\sqrt{2}r})(y)|\lesssim\frac{1}{r^d}\int\eta_2(y)\eta_2(x-y)|v_r(x)|dx
\lesssim\frac{1}{r^d}\int\eta_2(x)|v_r(x)|dx,
\end{equation}
so that using $\eta\lesssim\eta_2^2$ we obtain
\begin{equation}\nonumber
\Big(\int\eta(v_{\sqrt{2}r})^2\Big)^\frac{1}{2}\lesssim
\big(\int\eta_2|v_{\sqrt{2}r}|\big)^\frac{1}{2}\big(\frac{1}{r^d}\int\eta_2|v_{r}|\big)^\frac{1}{2},
\end{equation}
which turns into (\ref{b56bis}) by appealing once more to $|v_{\sqrt{2}r}|\le|v_r|_r$ and to (\ref{c84}) (cf Step~4).

\medskip

Inserting (\ref{b55bis}) and (\ref{b56bis}) into (\ref{b57}) yields for all $r\ll 1$
\begin{equation}\label{b58}
\Big(\int\eta v^2\Big)^\frac{1}{2}
\lesssim
r\Big(\int\eta|\nabla v|^2\Big)^\frac{1}{2}+\frac{1}{r^\frac{d}{2}}\int\eta_2|v_{r}|.
\end{equation}
By the monotonicity property (\ref{c86+}) and estimate \eqref{c84} (see Step~4) we have $\int\eta_2|v_{r}|\lesssim\int\eta_2|v_{r'}|$
for any $r'\le r\lesssim 1$ and thus
\begin{equation}\nonumber
\int\eta_2|v_{r}|\lesssim\frac{1}{r^{p+1}}\int_0^rdr' {r'}^p\int\eta_2|v_{r'}|.
\end{equation}
Inserting this into (\ref{b58}) yields the interpolation inequality in its additive form
\begin{equation}\nonumber
\Big(\int\eta v^2\Big)^\frac{1}{2}
\lesssim
r\Big(\int\eta|\nabla v|^2\Big)^\frac{1}{2}+\frac{1}{r^{\frac{d}{2}+p+1}}
\int_0^1dr'{r'}^p\int\eta_2|v_{r'}|,
\end{equation}
from which we obtain the desired form (\ref{b59}) by optimization in $r\ll 1$.

\medskip

\step{2} Proof of \eqref{Laux6.1} by an ODE argument.

\noindent
Let us momentarily introduce the abbreviations
\begin{equation}\nonumber
E_1(t):=\int\eta\nabla v\cdot a\nabla v,\quad
E_0(t):=\int\eta v^2,\quad
E_{-1}(t):=(\int_0^1r^p\int\eta|v_r|dr)^2,
\end{equation}
and note that $E_1(t)\sim\int\eta|\nabla v|^2$ by the bounds (\ref{1.3}) on $a$,
so that the results of the previous step turn into the differential relations
\begin{equation}\label{b60}
\frac{d}{dt}[\exp(-Ct)E_0]+\frac{1}{C}\exp(-Ct)E_1\le0,\quad
\frac{d}{dt}[\exp(-Ct)E_1]\le0
\end{equation}
and the ``algebraic'' relation
\begin{equation}\label{b61}
E_0\le C(E_1^{1-\theta} E_{-1}^{\theta}+E_{-1}),
\end{equation}
where we have set $\theta:=\frac{1}{\frac{d}{2}+p+2}\in(0,1)$. From (\ref{b60}), we infer
by integration by parts for any exponent $\alpha>0$
\begin{equation}\label{b62}
E_0(t=1)+\int_0^1 t^{1+\alpha}E_1dt\lesssim\int_0^1 t^\alpha E_0dt\quad\mbox{and}\quad
E_1(t=1)\lesssim\int_0^1 t^{1+\alpha}E_1dt;
\end{equation}
note that for $t\in[0,1]$, $\exp(-Ct)\sim 1$ so that we may safely neglect this factor
in (\ref{b60}).
From (\ref{b61}) we infer by H\"older's inequality for $\alpha=\frac{1-\theta}{\theta}$
\begin{equation}\label{b63}
\int_0^1 t^\alpha E_0dt \lesssim \big(\int_0^1t^{1+\alpha}E_1dt\big)^{1-\theta}
\big(\int_0^1E_{-1}dt\big)^{\theta}+\int_0^1E_{-1}dt.
\end{equation}
Inserting (\ref{b62}) into (\ref{b63}), we in a first stage obtain thanks to $\theta>0$ by Young's inequality
$\int_0^1 t^\alpha E_0dt \lesssim\int_0^1E_{-1}dt$ which then may be upgraded by using
both estimates in (\ref{b62}) to the desired result \eqref{Laux6.1} in form of
$E_0(t=1)+E_1(t=1)\lesssim\int_0^1E_{-1}dt$.

\medskip

\step{3} Proof of \eqref{Laux6.2}.

We apply \eqref{Laux6.1} to $u$, and need to post-process the rhs
\begin{equation}\nonumber
\fint_{\frac T4}^{\frac T2} dt\fint_0^{\sqrt{t}}dr (\frac{r}{\sqrt{t}})^{p+1}\int\eta_{2\sqrt{T}}|u_r(t)|.
\end{equation}
Since 
$t \mapsto \int\eta_{\sqrt{t}}|u_r|$ is decreasing with increasing convolution scale,
as a consequence of Jensen's inequality and (\ref{c84}) below, we have
\begin{equation}\label{L11-8}
\fint_{\frac T4}^{\frac T2} dt\fint_0^{\sqrt{t}}dr (\frac{r}{\sqrt{t}})^{p+1}\int\eta_{2\sqrt{T}}|u_r(\tau)|
\lesssim\fint_\frac{T}{4}^\frac{T}{2} dt\fint_0^{\sqrt{t}}dr (\frac{r}{\sqrt{t}})^{p+1}\int\eta_{\sqrt{t}}|u_r(t)|.
\end{equation}
Finally, by the semi-group property of Lemma~\ref{Laux2}  in form of $u(t)=\nabla\cdot(q(t)-\langle q(t)\rangle)$ (note that $\nabla \cdot \langle q(t)\rangle=0$
by stationarity of $q(t)$) and using $|(\nabla f)_r| $ $=|(\nabla f_{\frac{1}{\sqrt{2}}r})_{\frac{1}{\sqrt{2}}r}|$
$\lesssim\frac{1}{r}|f_{\frac{1}{\sqrt{2}}r}|_r$  as well as
$\int\eta_{\sqrt{t}}f_r$ $\stackrel{(\ref{c84})}{\lesssim}\int\eta_{\sqrt{t}}f$
for $r\lesssim\sqrt{t}$, we obtain
\begin{equation}\nonumber
\int\eta_{\sqrt{t}}|u_r(t)|
\lesssim \frac{1}{r}\int\eta_{\sqrt{t}}|(q(t)-\langle q(t)\rangle)_r|.
\end{equation}
This yields \eqref{Laux6.2} by \eqref{Laux6.1} and (\ref{L11-8}).

\medskip

\step{4} Some auxiliary statements.

In the course of our proofs, we  need some statements (to relate 
different Gaussian and exponential kernels) which hold for a function $f=f(x)\ge 0$ and two scales $r\le\frac{1}{2}R$:
\begin{align}
\mbox{(holds even for $r\lesssim R$)}\quad\int\eta_Rf_r\sim&\int\eta_R f,\label{c84}\\
(\frac{R}{r})^\frac{d}{2}\int\eta_{2R}(y)\big(\int\eta_r(x-y)f^2(x)dx\big)^\frac{1}{2}dy
\gtrsim&\big(\int\eta_R f^2\big)^\frac{1}{2},\label{c85}\\
\int\eta_R(y)\int\eta_r(x-y)f(x)dxdy
\sim&\int\eta_R f.\label{c86}
\end{align}
We shall also make use of the following monotonicity in the convolution scale,
based on the semi-group property $f_{\sqrt{t'}}=(f_{\sqrt{t}})_{\sqrt{t'-t}}$
for $t'\ge t$
of convolution with Gaussian yielding by Jensen's inequality 
\begin{equation}\label{c86+}
|f_{\sqrt{t'}}|^2\le|f_{\sqrt{t}}|^2_{\sqrt{t'-t}}.
\end{equation}
To prove \eqref{c84}--\eqref{c86}, we may assume $R=1$. Clearly, (\ref{c84}) is equivalent to $(\eta_1)_r\sim\eta_1$,
that is,  
$$\frac{1}{r^d}\int\exp(-\frac{1}{2}|\frac{y}{r}|^2)\exp(|x|-|x-y|)dy\sim 1,$$
by the triangle inequality we have $\exp(|x|-|x-y|)$ $\in[\exp(-|y|),\exp(|y|)]$, so that
the statement follows from $\int\exp(-\frac{1}{2}|\frac{y}{r}|^2\pm|y|)dy$ $\sim r^d$, which
can be seen from applying Young's inequality on $-\frac{1}{2}|\frac{y}{r}|^2\pm|y|$.
The approximate identity (\ref{c86}) is equivalent to
$\int\eta_1(y)\eta_r(x-y)dy=\int\eta_r(y)\eta_1(x-y)dy\sim\eta_1(x)$. Since
$\frac{\eta_1(x-y)}{\eta_1(x)}$
$=\exp(-|x-y|+|x|)$ $\in[\exp(-|y|),\exp(|y|)]$ this follows from $r\le\frac{1}{2}$.
For (\ref{c85}) we start with
$\eta_{1}(y)\int\eta_r(x-y)f^2(x)dx$ $\le(\frac{1}{r})^d\int\eta_{1}(x)f^2(x)dx$, a consequence
of 
$$\eta_1(y)\eta_r(x-y)=\frac{1}{r^d}\exp(-|y|-|\frac{x-y}{r}|)\le\frac{1}{r^d}\exp(-|x|)=\frac{1}{r^d}\eta_1(x),$$
 of which we take the square root:
$2^d\eta_{2}(y)\big(\int\eta_r(x-y)f^2(x)dx\big)^\frac{1}{2}$ $\le(\frac{1}{r})^\frac{d}{2}
\big(\int\eta_{1}f^2\big)^\frac{1}{2}$. We combine this with
(\ref{c86}) in form of
\begin{multline*}
\int\eta_1f^2 \lesssim\int\eta_1(y)\int\eta_r(x-y)f^2(x)dxdy
=\int2^d\eta_2(y)\eta_2(y)\int\eta_r(x-y)f^2(x)dxdy
\\
\le(\frac{1}{r})^\frac{d}{2}
\big(\int\eta_{1}f^2\big)^\frac{1}{2}\int\eta_2(y)\big(\int\eta_r(x-y)f^2(x)dx\big)^\frac{1}{2}dy.
\end{multline*}
Dividing by $\big(\int\eta_1f^2\big)^\frac{1}{2}$ yields (\ref{c85}).


\subsection{Proof of Lemma~\ref{Laux7}: quantified equi-integrability}

We proceed into steps.

\medskip

\step{1} Reduction.

\noindent By scale and translation invariance, it is enough to show for a solution of
\begin{equation}\label{c53}
\partial_tv-\nabla\cdot a\nabla v=0\quad\mbox{for}\quad t<0
\end{equation}
that
\begin{equation}\label{c60}
\big(\int_{-\frac{1}{4}}^0\int\eta_{R}|\nabla v|^2dxdt\big)^\frac{1}{2}
\lesssim(\frac{1}{R})^{\frac{d}{2}-\varepsilon}
\big(\int_{-1}^0\int\eta_1|\nabla v|^2dxdt\big)^\frac{1}{2}\quad\mbox{for}\;R\le 1.
\end{equation}

\medskip

\step{2} Meyers estimate.

\noindent 
Estimate (\ref{c60}) relies on the parabolic Meyers estimate, which states that for
the inhomogeneous equation
\begin{equation}\nonumber
\partial_tw-\nabla\cdot a\nabla w=\nabla\cdot g
\end{equation}
on the entire time-space $\mathbb{R}\times\mathbb{R}^d$ with compactly supported
$v$ and $g$, there exists an exponent $p=p(d,\lambda)>2$ (possible very close to $2$) such that we have
the Calder\'on-Zygmund type estimate
\begin{equation}\label{c47}
\big(\int\int|\nabla w|^pdxdt\big)^\frac{1}{p}\lesssim\big(\int\int|g|^pdxdt\big)^\frac{1}{p}.
\end{equation}
We need to post-process this estimate in the following way: For
\begin{equation}\label{c46}
\partial_tw-\nabla\cdot a\nabla w=f+\nabla\cdot g
\end{equation}
with compactly supported $w$, $f$ and $g$, we claim
\begin{equation}\label{c51}
\big(\int\int|\nabla w|^pdxdt\big)^\frac{1}{p}\lesssim
\big(\int\big(\int|f|^qdx\big)^\frac{p}{q}dt\big)^\frac{1}{q}+\big(\int\int|g|^pdxdt\big)^\frac{1}{p},
\end{equation}
where the exponent $q$ is determined by scaling $q:=\frac{dp}{d+p}$.
In fact, we will use (\ref{c51}) in the non-scale invariant form of
\begin{eqnarray}\label{c55}
\big(\int\int|\nabla w|^pdxdt\big)^\frac{1}{p}&\lesssim&
\sup_t\Big(\big(\int f^2dx\big)^\frac{1}{2}+\big(\int|g|^pdx\big)^\frac{1}{p}\Big),\nonumber\\
&&\mbox{provided}\;{\rm supp}(w,f,g)\subset[-1,0]\times\{|x|\le 1\},
\end{eqnarray}
where we used that wlog we may assume $p\le\frac{2d}{d-2}$ so that $q\le 2$.
Here comes the argument for (\ref{c51}): By the scale and translation invariance of the estimate, we may wlog
assume that all functions involved in (\ref{c46})
have compact support in $[-1,0]\times\{|x|\le 1\}$. Hence $f$ must have
vanishing space-time average. In other words, introducing the spatial average $\bar f(t):=\int fdx$
we have $\int \bar f dt=0$. On the one hand, therefore there exists $\bar F=\bar F(t)$ supported in $t\in[-1,0]$
and such that $\frac{d\bar F}{dt}=\bar f$. By ${\rm supp}f\subset[-1,0]\times\{|x|\le 1\}$, we clearly have
\begin{equation}\label{c50}
\big(\int|\bar F|^pdt\big)^\frac{1}{p}
\le\sup_{t}|\bar F|\lesssim\int\int|f|dxdt
\lesssim\big(\int\big(\int|f|^qdx\big)^\frac{p}{q}dt\big)^\frac{1}{p}.
\end{equation}
On the one hand, $\bar F$ allows us to rewrite (\ref{c46})
as $\partial_t(w-\chi\bar F)-\nabla\cdot a\nabla(w-\chi\bar F)=f-\chi\bar f+\nabla\cdot(g+\bar Fa\nabla\chi)$,
where $\chi=\chi(x)$ is a cut-off function for $\{|x|\le 1\}$ in $\{|x|\le 2\}$ with $\int\chi dx=1$ (that is, $\chi\ge0$, $\chi \gtrsim 1$ on $\{|x|\le 1\}$
and $\chi=0$ on $\{|x|>2\}$).
Fixing the time $t$, this allows us to solve the Neumann problem
\begin{equation}\nonumber
-\triangle W=f-\bar f\chi\;\mbox{for}\;|x|<2,\quad -x\cdot\nabla W=0\;\mbox{for}\;|x|=2
\end{equation}
on the ball $\{|x|\le 2\}$. By maximal regularity for this Neumann problem 
and Sobolev's inequality (using $p=\frac{dq}{d-q}$) we have
\begin{equation}\label{c48}
\big(\int_{|x|\le 2}|\nabla W|^p\big)^\frac{1}{p}\lesssim
\big(\int_{|x|\le 2}|f-\bar f\chi|^qdx\big)^\frac{1}{q}
\lesssim\big(\int |f|^qdx\big)^\frac{1}{q}.
\end{equation}
Note that by the Neumann boundary
condition, when extending $-\nabla W$ trivially outside of $\{|x|\le 2\}$ to a vector field $\delta g$
it still solves $\nabla\cdot\delta g=f-\bar f\chi$, now in all space. Also the estimate (\ref{c48})
turns into
\begin{equation}\label{c49}
\big(\int\int|\delta g|^pdxdt\big)^\frac{1}{p}\lesssim
\big(\int\big(\int|f|^qdx\big)^\frac{p}{q}dt\big)^\frac{1}{p}.
\end{equation}
Hence with the abbreviations $\tilde w:=w-\bar F\chi$ and $\tilde g:=g+\bar Fa\nabla\chi+\delta g$,
we obtain $\partial_t\tilde w-\nabla\cdot a\nabla\tilde w=\nabla\cdot\tilde g$, so that an application
of (\ref{c47}) yields by the triangle inequality, by the smoothness and compact support of $\chi$,
and by the bound (\ref{1.3}) on $a$,
\begin{equation}\nonumber
\big(\int\int|\nabla w|^pdxdt\big)^\frac{1}{p}
\lesssim \big(\int|\bar F|^pdt\big)^\frac{1}{p}
+\big(\int\int|g|^pdxdt\big)^\frac{1}{p}
+\big(\int\int|\delta g|^pdxdt\big)^\frac{1}{p}.
\end{equation}
Inserting (\ref{c50}) and (\ref{c49}) into this yields (\ref{c51}).

\medskip

\step 3 Localized semi-group estimates.

We now turn to a solution of (\ref{c53}) and claim that there exists a constant $c$ such that
\begin{equation}\label{c52}
\sup_{t\in(-\frac{1}{2},0)}\Big(\big(\int_{|x|\le\frac{1}{2}}|\nabla v|^2dx\big)^\frac{1}{2}
+\big(\int_{|x|\le\frac{1}{2}}|v-c|^pdx\big)^\frac{1}{p}\Big)\lesssim
\big(\int_{-1}^0\int \eta_1 |\nabla v|^2\big)^\frac{1}{2}.
\end{equation}
We start with the first estimate in (\ref{c52}).
Since equation \eqref{p2.bis} in the proof of Lemma~\ref{Laux1} holds with $w$ replaced by $w-\int \eta_R w$, estimates 
\eqref{p2}, \eqref{p2.ter}, and ultimately \eqref{p1bis} hold with $w$ replaced by $w-\int \eta_R w$ as well.
In the case of $v$, this takes the form
$$
\sup_{t\in(-\frac{1}{2},0)}\big(\int \eta_1 |\nabla v|^2dx\big)^\frac{1}{2} \,\lesssim\,\big(\int_{-1}^0\int \eta_1 (v-\int \eta_1 v)^2\big)^\frac{1}{2}.
$$
By the Poincar\'e inequality for $dx$ replaced by $\eta_1 dx$ (cf \cite[Lemma~2.1]{BobkovLedoux}), this yields
$$
\sup_{t\in(-\frac{1}{2},0)}\big(\int_{|x|\le \frac12}  |\nabla v|^2dx\big)^\frac{1}{2} \,\lesssim\,\big(\int_{-1}^0\int \eta_1 |\nabla v|^2\big)^\frac{1}{2}.
$$

\medskip

We now turn to the second estimate in (\ref{c52}). We select a smooth function $\chi=\chi(x)$ supported
in $\{|x|\le\frac{1}{2}\}$ with $\int\chi dx=1$. Denoting by $\bar u=\bar u(t)$ the average of
$u$ as defined by $\bar v=\int\chi vdx$ we have by the Poincar\'e-Sobolev inequality
\begin{equation}\nonumber
\big(\int_{|x|\le\frac{1}{2}}|v-\bar v|^pdx\big)^\frac{1}{p}
\lesssim\big(\int_{|x|\le\frac{1}{2}}|\nabla v|^2dx\big)^\frac{1}{2},
\end{equation}
where we use again that $p\le\frac{2d}{d-2}$.
Hence in view of the first estimate in (\ref{c52}), it remains to show
\begin{equation}\nonumber
\inf_{c}\sup_{t\in(-\frac{1}{2},0)}|\bar v-c|\lesssim
\big(\int_{-\frac{1}{2}}^0|\frac{d\bar v}{dt}|^2dt\big)^\frac{1}{2}
\lesssim\big(\int_{-\frac{1}{2}}^0\int_{|x|\le\frac{1}{2}}|\nabla v|^2dxdt\big)^\frac{1}{2}.
\end{equation}
The latter follows from the equation which gives $\frac{d\bar v}{dt}=\int\nabla\chi\cdot a\nabla vdx$
via the Cauchy-Schwarz inequality.

\medskip

\step4 Conclusion.

We now may conclude and to this purpose consider $w=\chi(v-c)$ where $\chi=\chi(t,x)$ is a cut-off
function for $[-\frac{1}{4},0]\times\{|x|\le\frac{1}{4}\}$ in
$[-\frac{1}{2},0]\times\{|x|\le\frac{1}{2}\}$ and $c$ is the constant in (\ref{c52}).
We note that by (\ref{c53}),
$w$ satisfies (\ref{c46}) with $f=(v-c)\partial_t\chi-\nabla\chi\cdot a\nabla v$
and $g=-(v-c)a\nabla\chi$, and that $w$, $f$, and $g$ have support in
$(-\frac{1}{2},0)\times\{|x|\le\frac{1}{2}\}$.
Hence by an application of (\ref{c55}) we obtain
\begin{equation}\nonumber
\big(\int_{-\frac{1}{4}}^0\int_{|x|\le\frac{1}{4}}|\nabla v|^pdxdt\big)^\frac{1}{p}\lesssim
\sup_{t\in(-\frac{1}{2},0)}
\Big(\big(\int_{|x|\le \frac 12}|\nabla v|^2dx\big)^\frac{1}{2}+\big(\int|v-c|^pdx\big)^\frac{1}{p}\Big).
\end{equation}
Inserting (\ref{c52}) now yields the ``reverse H\"older'' inequality
\begin{equation}\nonumber
\big(\int_{-\frac{1}{4}}^0\int_{|x|\le\frac{1}{4}}|\nabla v|^pdxdt\big)^\frac{1}{p}\lesssim
\big(\int_{-1}^0\int \eta_1 |\nabla v|^2dxdt\big)^\frac{1}{2}.
\end{equation}
In combination with H\"older's inequality in form of
$$ \big(\int_{-\frac{1}{4}}^0\int_{|x|\le r}|\nabla v|^2dxdt\big)^\frac{1}{2}
\lesssim r^{d(\frac{1}{2}-\frac{1}{p})}\big(\int_{-\frac{1}{4}}^0\int_{|x|\le r}|\nabla v|^pdxdt\big)^\frac{1}{p},
$$
we obtain
\begin{equation}\label{c57}
\big(\int_{-\frac{1}{4}}^0\int_{|x|\le r}|\nabla v|^2dxdt\big)^\frac{1}{2}
\lesssim r^\varepsilon\big(\int_{-1}^0\int \eta_1 |\nabla v|^2dxdt\big)^\frac{1}{2}
\end{equation}
with $\varepsilon=\varepsilon(d,\lambda):=d(\frac{1}{2}-\frac{1}{p})>0$ for $r\le\frac{1}{4}$
and also for $r\le 1$ (since the estimate is trivial for $r\in[\frac{1}{4},1]$).

\medskip

It remains to introduce the exponential averaging functions $\eta_R$.
For given $R\le\frac{1}{2}$ we appeal to the co-area formula in form of
$\int\eta_R|\nabla v|^2dx$
$=\frac{1}{R^d}\int\exp(-\frac{|x|}{R})|\nabla v|^2dx$
$=\frac{1}{R^d}\int_0^\infty\exp(-\frac{r}{R})\int_{|x|=r}|\nabla v|^2dxdr$
$=\frac{1}{R^{d+1}}\int_0^\infty\exp(-\frac{r}{R})\int_{|x|\le r}|\nabla v|^2dxdr$.
From splitting the last integral into $r\le 1$ and $r\ge 1$ and using (\ref{c57})
on the first part we obtain
\begin{align}\label{c58}
\int_{-\frac{1}{4}}^0\int\eta_R|\nabla v|^2dxdt
&\lesssim\frac{1}{R^{d+1}}\int_0^1\exp(-\frac{r}{R})r^{2\varepsilon}dr\int_{-1}^0\int \eta_1 |\nabla v|^2dxdt\nonumber\\
&+\frac{1}{R^{d+1}}
\int_1^\infty\exp(-\frac{r}{R})\int_{-1}^0\int_{|x|\le r}|\nabla v|^2dxdtdr.
\end{align}
For the first rhs term we use $\frac{1}{R^{d+1}}\int_0^1\exp(-\frac{r}{R})r^{2\varepsilon}dr
\le\frac{1}{R^{d+1}}\int_0^\infty\exp(-\frac{r}{R})r^{2\varepsilon}dr\sim\frac{1}{R^{d-2\varepsilon}}$
so that
\begin{equation}\label{c59}
\frac{1}{R^{d+1}}\int_0^1\exp(-\frac{r}{R})r^{2\varepsilon}dr\int_{-1}^0\int \eta_1 |\nabla v|^2dxdt
\lesssim
\frac{1}{R^{d-2\varepsilon}}\int_{-1}^0\int\eta_1|\nabla v|^2dxdt.
\end{equation}
For the second rhs term in (\ref{c58}) we note that
for $R\le\frac{1}{2}$ and $r\ge 1$ we have $\exp(-\frac{r}{R})$ $\lesssim R^{1+2\varepsilon}\exp(-r)$
(since in this range $\exp(-\frac{r}{R})\exp(r)$ $=\exp(-(\frac{1}{R}-1)r)\le\exp(-\frac{1}{2R})$)
so that we have for the second term
\begin{multline*}
\frac{1}{R^{d+1}}\int_1^\infty\exp(-\frac{r}{R})\int_{-1}^0\int_{|x|\le r}|\nabla v|^2dxdtdr
\\
\lesssim\frac{1}{R^{d-2\varepsilon}}\int_0^\infty\exp(-r)\int_{-1}^0\int_{|x|\le r}|\nabla v|^2dxdtdr,
\end{multline*}
which by the co-area formula yields
\begin{equation}\nonumber
\frac{1}{R^{d+1}}
\int_1^\infty\exp(-\frac{r}{R})\int_{-1}^0\int_{|x|\le r}|\nabla v|^2dxdtdr
\lesssim\frac{1}{R^{d-2\varepsilon}}\int_{-1}^0\int\eta_1|\nabla v|^2dxdt.
\end{equation}
Inserting this and (\ref{c59}) into (\ref{c58}) yields (\ref{c60}) for $R\le \frac12$. The range $\frac 12 \le R \le 1$ is trivial.


\subsection{Proof of Lemma~\ref{Laux8}: strong smallness of the extended corrector}

In order to avoid indices, we use $(d=3)$-notation when it comes
to the vector potential $\sigma_T$, that is, we write $\frac{1}{T}\sigma_T-\triangle\sigma_T=\nabla\times q_T$
instead of (\ref{o17}).

\medskip

\step{1} We claim that for all $t\le T$ and $0<\delta\le 1$
\begin{align}
\Big(\int\eta_{\sqrt{T}}|\frac{\phi_t}{\sqrt{t}}|^2\Big)^\frac{1}{2}
&\lesssim\frac{1}{\delta}\Big(\int\eta_{\sqrt{T}}
|(q_t-\langle q_t\rangle)_{\delta\sqrt{t}}|^2\Big)^\frac{1}{2}+\delta,\label{b40}\\
\Big(\int\eta_{\sqrt{T}}|\frac{\sigma_t}{\sqrt{t}}|^2\Big)^\frac{1}{2}
&\lesssim\Big(\int\eta_{\sqrt{T}}
|(q_t-\langle q_t\rangle)_{\delta\sqrt{t}}|^2\Big)^\frac{1}{2}+\delta,\label{b41}\\
\Big(\int\eta_{\sqrt{T}}|g_t|^2\Big)^\frac{1}{2}
&\lesssim\Big(\int\eta_{\sqrt{T}}
|(q_t-\langle q_t\rangle-\nabla\phi_t)_{\delta\sqrt{t}}|^2\Big)^\frac{1}{2}+\delta.\label{b42}
\end{align}
It is enough to prove these estimates for $\sqrt{T}=\sqrt{t}$: We
replace the origin by $y$ and then apply $\int\eta_{\sqrt{T}}(y)\cdot dy$, using
$\eta_{\sqrt{T}}*\eta_{\sqrt{t}}\lesssim\eta_{\sqrt{T}}$
(for $\sqrt{t}\ll\sqrt{T}$, cf~(\ref{c86}))
for our exponential localization function.
By scale invariance, it is enough to establish them for $t=1$,
we shall thus write $(\phi,\sigma,g)$ and $\eta$
instead of $(\phi_t,\sigma_t,g_t)$ and $\eta_{\sqrt{t}}$, respectively. 
For all functions $(\phi,\sigma,g)$ we start from the splitting
\begin{equation}\label{b43}
\int\eta|(\phi,\sigma,g)|^2\lesssim\int\eta|(\phi,\sigma,g)_\delta|^2+\int\eta|(\phi,\sigma,g)-
(\phi,\sigma,g)_\delta|^2.
\end{equation}
By \eqref{b55bis} in the proof of Lemma~\ref{Laux6}, we have for the second term
\begin{equation}\label{b55}
\int\eta|(\phi,\sigma,g)-(\phi,\sigma,g)_\delta|^2\lesssim\delta^2\int\eta|\nabla(\phi,\sigma,g)|^2.
\end{equation}
Thanks to the massive term in the equations (\ref{c62}), (\ref{o17}), (\ref{a4}) for $\phi$, $\sigma$, and $g$,
we may localize the elliptic energy estimates (see (\ref{Laux3.1}) in Lemma~\ref{Laux3} for ${t}=1$)
\begin{equation}\label{L4.1}
\int\eta|(\nabla\phi,q)|^2\lesssim\int\eta|ae|^2\lesssim 1,
\end{equation}
and
\begin{equation}\nonumber
\int\eta|\nabla\sigma|^2\lesssim\int\eta|q|^2,\quad
\int\eta|\nabla g|^2\lesssim\int\eta|(\nabla\phi,q-\langle q\rangle)|^2,
\end{equation}
which together with $|\langle q\rangle|\le\langle|q|^2\rangle\lesssim 1$, which follows
from (\ref{L4.1}) by stationarity, yield
\begin{equation}\nonumber
\int\eta|(\nabla\phi,\nabla\sigma,\nabla g,q-\langle q\rangle)|^2\lesssim1.
\end{equation}

\medskip

It remains to estimate the first term in (\ref{b43}). This is easy for $(\sigma,g)$: Since these
satisfy constant-coefficient equations we have $\sigma_\delta-\triangle\sigma_\delta
=\nabla\times(q-\langle q\rangle)_\delta$ and
$g_\delta-\triangle g_\delta=(q-\langle q\rangle-\nabla\phi)_\delta$, so that by the same localized
elliptic energy estimate
\begin{equation}\nonumber
\int\eta|\sigma_\delta|^2\lesssim\int\eta|(q-\langle q\rangle)_\delta|^2,\quad
\int\eta|g_\delta|^2\lesssim\int\eta|(q-\langle q\rangle-\nabla\phi)_\delta|^2.
\end{equation}
For $\phi$ we note that $\phi=\nabla\cdot(q-\langle q\rangle)$, cf~(\ref{c62}), and thus also
$\phi_{\sqrt{2}\delta}=(\nabla\cdot(q-\langle q\rangle)_\delta)_\delta$ by the semi-group
property of convolution with the Gaussian. 
Since the gradient of a Gaussian
is estimated by the Gaussian with twice the variance, cf $|\nabla G_{\kappa}(z)|=|\frac{z}{\kappa}G_{\kappa}(z)| \lesssim G_{\sqrt{2}\kappa}(z)$, this implies $|\phi_{\sqrt{2}\delta}|
\lesssim\frac{1}{\delta}|(q-\langle q\rangle)_\delta|_{\sqrt{2}\delta}$, and thus by
Jensen's inequality $|\phi_{\sqrt{2}\delta}|^2
\lesssim\frac{1}{\delta^2}|(q-\langle q\rangle)_\delta|_{\sqrt{2}\delta}^2$. Hence we have
\begin{equation}\nonumber
\int\eta|\phi_{\sqrt{2}\delta}|^2
\lesssim\frac{1}{\delta^2}\int\eta_{\sqrt{2}\delta}|(q-\langle q\rangle)_\delta|^2
\stackrel{(\ref{c84})}{\lesssim}\frac{1}{\delta^2}\int\eta|(q-\langle q\rangle)_\delta|^2.
\end{equation}
It remains to replace $\delta$ by $\sqrt{2}\delta$ in (\ref{b43}) when it comes to $\phi$.

\medskip

\step{2} We claim that for $t\le T$
\begin{eqnarray}\label{a37}
\lefteqn{\Big(\int\eta_{\sqrt{T}}|(\frac{\phi_T-\phi_t}{\sqrt{T}},\frac{\sigma_T-\sigma_t}{\sqrt{T}},
\nabla\phi_T-\nabla\phi_t,q_T-q_t)|^2\Big)^\frac{1}{2}}\nonumber\\
&\lesssim&\Big(\int\eta_{\sqrt{T}}|(\frac{\phi_t}{\sqrt{t}},\frac{\sigma_t}{\sqrt{t}},
(q_{t}-\langle q_{t}\rangle)_{\sqrt{t}})|^2\Big)^\frac{1}{2}.
\end{eqnarray}
We note that by equation (\ref{c62}), $u:=\phi_T-\phi_t$ satisfies
$\frac{1}{T}u-\nabla\cdot a\nabla u=(\frac{1}{t}-\frac{1}{T})\phi_t$.
Using (\ref{c62}) in form of $\frac{1}{t}\phi_t=\nabla\cdot(q_t-\langle q_t\rangle)$,
we split the rhs according to $\phi_t=(\phi_{t}-\phi_{t,\sqrt{t}})
+t\nabla\cdot(q_{t}-\langle q_{t}\rangle)_{\sqrt{t}}$. We now perform a localized elliptic energy estimate
which means testing with $\eta_{\sqrt{T}}u$; thanks to the massive term we obtain
as in the argument for (\ref{Laux3.1})
\begin{equation}\label{b56}
\int\eta_{\sqrt{T}}|(\frac{u}{\sqrt{T}},\nabla u)|^2
\lesssim (1-\frac{t}{T})
\int\eta_{\sqrt{T}}u\frac{1}{t}(\phi_{t}-\phi_{t,\sqrt{t}})
+\int\eta_{\sqrt{T}}|(q_t-\langle q_t\rangle)_{\sqrt{t}}|^2.
\end{equation}
For the contribution of the first rhs term we use symmetry of the convolution operator
\begin{equation}\nonumber
(1-\frac{t}{T})\int\eta_{\sqrt{T}}u\frac{1}{t}(\phi_{t}-\phi_{t,\sqrt{t}})
=(1-\frac{t}{T})\int\big(\eta_{\sqrt{T}}u-(\eta_{\sqrt{T}}u)_{\sqrt{t}}\big)\frac{1}{t}\phi_t.
\end{equation}
Denoting by $G_{\sqrt{t}}$ the Gaussian convolution kernel on scale $\sqrt{t}$, we have
\begin{eqnarray}\nonumber
\lefteqn{\big(\eta_{\sqrt{T}}u-(\eta_{\sqrt{T}}u)_{\sqrt{t}}\big)(x)}\nonumber\\
&=&\int_0^1\int\big(\eta_{\sqrt{T}}\nabla u
+u\nabla\eta_{\sqrt{T}}\big)(x+sz)\cdot zG_{\sqrt{t}}(z)dzds\nonumber
\end{eqnarray}
and thus by $|zG_{\sqrt{t}}(z)|\lesssim \sqrt{t}G_{2\sqrt{t}}(z)$,
$|\nabla\eta_{\sqrt{T}}|\lesssim\frac{1}{\sqrt{T}}\eta_{\sqrt{T}}$,
and Jensen's inequality
\begin{eqnarray}\nonumber
\lefteqn{\big|\eta_{\sqrt{T}}u-(\eta_{\sqrt{T}}u)_{2\sqrt{t}}\big|(x)}\nonumber\\
&\lesssim&\sqrt{t}\int_0^1\int\big(\eta_{\sqrt{T}}|(\frac{u}{\sqrt{T}},\nabla u)|\big)(x+sz)G_{2\sqrt{t}}(z)dzds
\nonumber\\
&=&\sqrt{t}\int_0^1\big(\eta_{\sqrt{T}}|(\frac{u}{\sqrt{T}},\nabla u)|\big)_{2\sqrt{t}s}(x)ds\nonumber
\end{eqnarray}
and thus by Cauchy-Schwarz' inequality, the symmetry of the convolution operator $\cdot_{2\sqrt{t}s}$ (that we use twice),
and Jensen's inequality
\begin{eqnarray}\nonumber
\lefteqn{\big|(1-\frac{t}{T})\int\eta_{\sqrt{T}}u\frac{1}{t}(\phi_{t}-\phi_{t,\sqrt{t}})\big|}\nonumber\\
&\lesssim&\int_0^1\int\big(\eta_{\sqrt{T}}|(\frac{u}{\sqrt{T}},\nabla u)|\big)_{2\sqrt{t}s}
|\frac{\phi_t}{\sqrt{t}}|ds\nonumber\\
&\lesssim&\Big(\int\eta_{\sqrt{T}}|(\frac{u}{\sqrt{T}},\nabla u)|^2\Big)^\frac{1}{2}
\Big(\int_0^1\int(\eta_{\sqrt{T}})_{2\sqrt{t}s}(\frac{\phi_t}{\sqrt{t}})^2ds\Big)^\frac{1}{2}\nonumber\\
&\stackrel{(\ref{c84})}{\lesssim}&\Big(\int\eta_{\sqrt{T}}|(\frac{u}{\sqrt{T}},\nabla u)|^2\Big)^\frac{1}{2}
\Big(\int\eta_{\sqrt{T}}(\frac{\phi_t}{\sqrt{t}})^2\Big)^\frac{1}{2}.\nonumber
\end{eqnarray}
Equipped with this estimate we return to (\ref{b56}) and
obtain for $\phi_T-\phi_t$ and $|q_T-q_t|\le|\nabla(\phi_T-\phi_t)|$ as desired
\begin{align}\label{a30}
\lefteqn{\Big(\int\eta_{\sqrt{T}}|(\frac{\phi_T-\phi_t}{\sqrt{T}},
\nabla\phi_T-\nabla\phi_t,q_T-q_t)|^2\Big)^\frac{1}{2}}\nonumber\\
&\lesssim\Big(\int\eta_{\sqrt{T}}(\frac{\phi_t}{\sqrt{t}})^2\Big)^\frac{1}{2}
+\Big(\int\eta_{\sqrt{T}}|(q_t-\langle q_t\rangle)_{\sqrt{t}}|^2\Big)^\frac{1}{2}.
\end{align}

\medskip

We now turn to the contribution of $\sigma$. As for $\phi$, we note that
by (\ref{o17}), $u:=\sigma_T-\sigma_t$ satisfies $\frac{1}{T}u-\triangle u
=(1-\frac{t}{T})\frac{1}{t}\sigma_t+\nabla\times(q_T-q_t)$
and split the first rhs term into $\frac{1}{t}\sigma_t=\frac{1}{t}(\sigma_t-\sigma_{t,\sqrt{t}})
+\frac{1}{t}\sigma_{t,\sqrt{t}}$. Applying the convolution $\cdot_{\sqrt{t}}$ to
$\frac{1}{t}\sigma_t-\triangle\sigma_t=\nabla\times q_t$ leads to
\begin{equation}\label{b45}
\frac{1}{t}\sigma_{t,\sqrt{t}}-\triangle\sigma_{t,\sqrt{t}}=\nabla\times (q_t-\langle q_t\rangle)_{\sqrt{t}}.
\end{equation}
Hence we write
\begin{eqnarray}\label{b46}
\lefteqn{\frac{1}{T}u-\triangle u}\\
&=&(1-\frac{t}{T})\big(\frac{1}{t}(\sigma_t-\sigma_{t,\sqrt{t}})+\triangle\sigma_{t,\sqrt{t}}
+\nabla\times (q_t-\langle q_t\rangle)_{\sqrt{t}}\big)
+\nabla\times(q_T-q_t).\nonumber
\end{eqnarray}
By  (\ref{Laux3.1}) in Lemma~\ref{Laux3} we may exponentially
localize the energy estimate for (\ref{b45})
on scale $\sqrt{t}$, and a fortiori on scale $\sqrt{T}$
\begin{equation}\label{a36}
\Big(\int\eta_{\sqrt{T}}|(\frac{\sigma_{t,\sqrt{t}}}{\sqrt{t}},\nabla\sigma_{t,\sqrt{t}})|^2\Big)^\frac{1}{2}
\lesssim\Big(\int\eta_{\sqrt{T}}|(q_t-\langle q_t\rangle)_{\sqrt{t}}|^2\Big)^\frac{1}{2}.
\end{equation}
Likewise, the localized elliptic energy estimate \eqref{Laux3.1} applied to (\ref{b46}) takes the form
(without using Cauchy-Schwarz' inequality on the first rhs term):
\begin{eqnarray*}\nonumber
\lefteqn{\int\eta_{\sqrt{T}}|(\frac{u}{\sqrt{T}},\nabla u)|^2}\nonumber\\
&\lesssim&
(1-\frac{t}{T})\int\eta_{\sqrt{T}}u\frac{1}{t}(\sigma_{t}-\sigma_{t,\sqrt{t}})
+\int\eta_{\sqrt{T}}|(\nabla\sigma_{t,\sqrt{t}},(q_t-\langle q_t\rangle)_{\sqrt{t}},q_T-q_t)|^2.
\end{eqnarray*}
Treating the first rhs term as in case of $\phi$, we obtain
\begin{eqnarray}\label{a35}
\lefteqn{\Big(\int\eta_{\sqrt{T}}|\frac{\sigma_T-\sigma_t}{\sqrt{T}}|^2\Big)^\frac{1}{2}}\nonumber\\
&\lesssim&\Big(\int\eta_{\sqrt{T}}|(\frac{\sigma_t}{\sqrt{t}},
\nabla\sigma_{t,\sqrt{t}},(q_t-\langle q_t\rangle)_{\sqrt{t}},q_T-q_t)|^2\Big)^\frac{1}{2}.
\end{eqnarray}
The combination of (\ref{a30}), and (\ref{a36}) inserted into (\ref{a35}) yields (\ref{a37}).

\medskip


\step{3} We claim that provided $t\le T$ and $0<\delta<1$
\begin{eqnarray}\label{b48}
\lefteqn{\Big(\int\eta_{\sqrt{T}}|(\frac{\phi_T}{\sqrt{T}},\frac{\sigma_T}{\sqrt{T}},g_T)|^2
\Big)^\frac{1}{2}}\\
&\lesssim&\frac{1}{\delta}\Big(\Big(\int\eta_{\sqrt{T}}
|(\nabla\phi_t,q_t-\langle q_t\rangle)_{\delta\sqrt{t}}|^2\Big)^\frac{1}{2}
+\langle|(q_t-\langle q_t\rangle)_{\delta\sqrt{t}}|^2\rangle^\frac{1}{2}\Big)+\delta.\nonumber
\end{eqnarray}
For the estimate on $(\phi_T,\sigma_T)$,
we use (\ref{a37}) in conjunction with the triangle inequality
\begin{eqnarray*}\nonumber
\lefteqn{\Big(\int\eta_{\sqrt{T}}|(\frac{\phi_T}{\sqrt{T}},\frac{\sigma_T}{\sqrt{T}},
\nabla\phi_T-\nabla\phi_t,q_T-q_t)|^2\Big)^\frac{1}{2}}\nonumber\\
&\lesssim&\Big(\int\eta_{\sqrt{T}}
|(\frac{\phi_t}{\sqrt{t}},\frac{\sigma_t}{\sqrt{t}},(q_t-\langle q_t\rangle)_{\sqrt{t}}|^2\Big)^\frac{1}{2},
\nonumber
\end{eqnarray*}
and insert (\ref{b40}) \& (\ref{b41}):
\begin{eqnarray}\label{b44}
\lefteqn{\Big(\int\eta_{\sqrt{T}}|(\frac{\phi_T}{\sqrt{T}},\frac{\sigma_T}{\sqrt{T}},
\nabla\phi_T-\nabla\phi_t,q_T-q_t)|^2\Big)^\frac{1}{2}}\nonumber\\
&\lesssim&\frac{1}{\delta}\Big(\int\eta_{\sqrt{T}}
|(q_t-\langle q_t\rangle)_{\delta\sqrt{t}}|^2\Big)^\frac{1}{2}+\delta,
\end{eqnarray}
where we made use of the monotonicity \eqref{c86+} in the convolution scale,
based on the semi-group property $q_{\sqrt{t'}}=(q_{\sqrt{t}})_{\sqrt{t'-t}}$
for $t'\ge t$
of convolution with Gaussian yielding by Jensen's inequality $|q_{\sqrt{t'}}|^2\le|q_{\sqrt{t}}|^2_{\sqrt{t'-t}}$
and thus by (\ref{c84})
\begin{equation}\label{b47}
\int\eta_{\sqrt{T}}|q_{\sqrt{t'}}|^2
\lesssim\int\eta_{\sqrt{T}}|q_{\sqrt{t}}|^2.
\end{equation}
For later use we note that because of stationarity of the fields $q_T-q_t$
and $(q_t-\langle q_t\rangle)_{\delta\sqrt{t}}$, taking the square expectation
of (\ref{b44}) implies in particular
\begin{eqnarray}\label{b49}
|\langle q_T\rangle-\langle q_t\rangle|\le\langle|q_T-q_t|^2\rangle^\frac{1}{2}
&\lesssim&\frac{1}{\delta}\langle
|(q_t-\langle q_t\rangle)_{\delta\sqrt{t}}|^2\rangle^\frac{1}{2}+\delta.
\end{eqnarray}

\medskip

We now turn to the estimate of $g_T$. We use (\ref{b42}) for $t=T$ and obtain by the
triangle inequality
\begin{eqnarray*}\nonumber
\lefteqn{\Big(\int\eta_{\sqrt{T}}|g_T|^2\Big)^\frac{1}{2}}\nonumber\\
&\lesssim&
\Big(\int\eta_{\sqrt{T}}|(q_T-q_t,\nabla\phi_T-\nabla\phi_t,
q_t-\langle q_t\rangle,\nabla\phi_t)_{\delta\sqrt{T}}|^2\Big)^\frac{1}{2}\nonumber\\
&&+|\langle q_T\rangle-\langle q_t\rangle|+\delta.\nonumber
\end{eqnarray*}
By the monotonicity (\ref{b47}) in the convolution scale, we may get rid of the convolution on the first two arguments of the first rhs term, which  yields
\begin{eqnarray*}\nonumber
\lefteqn{\Big(\int\eta_{\sqrt{T}}|g_T|^2\Big)^\frac{1}{2}}\nonumber\\
&\lesssim&
\Big(\int\eta_{\sqrt{T}}|(q_T-q_t,\nabla\phi_T-\nabla\phi_t)|^2\Big)^\frac{1}{2}\nonumber\\
&&+\Big(\int\eta_{\sqrt{T}}|(q_t-\langle q_t\rangle,\nabla\phi_t)_{\delta\sqrt{t}}|^2\Big)^\frac{1}{2}
+|\langle q_T\rangle-\langle q_t\rangle|+\delta.\nonumber
\end{eqnarray*}
Inserting (\ref{b44}) and (\ref{b49}) into this yields (\ref{b48}).

\medskip

\step{4} Post-processing.

\noindent We first note that we may get rid of the $\nabla\phi_t$-term on the rhs of (\ref{b48}):
Indeed, by (\ref{c62}) we have $\nabla\phi_t=t\nabla\nabla\cdot (q_t-\langle q_t\rangle)$ and thus
$(\nabla\phi_t)_{\sqrt{2}\delta \sqrt{t}}
=t(\nabla\nabla\cdot (q_t-\langle q_t\rangle)_{\delta \sqrt{t}})_{\delta \sqrt{t}}$, yielding
$|(\nabla\phi_t)_{\sqrt{2}\delta \sqrt{t}}|
\lesssim\frac{1}{\delta^2}|(q_t-\langle q_t\rangle)_{\delta \sqrt{t}}|_{2\delta \sqrt{t}}$
and therefore $|(\nabla\phi_t)_{\sqrt{2}\delta \sqrt{t}}|^2
\lesssim\frac{1}{\delta^4}|(q_t-\langle q_t\rangle)_{\delta \sqrt{t}})|_{2\delta t}^2$ (cf  argument for \eqref{Laux5-1ref} 
in the proof of Lemma~\ref{Laux5}).
Since $(\eta_{\sqrt{T}})_{2\delta \sqrt{t}}\stackrel{(\ref{c84})}{\lesssim}
\eta_{\sqrt{T}}$ this implies as desired
\begin{equation}\label{b51}
\Big(\int\eta_{\sqrt{T}}|(\nabla\phi_t)_{\sqrt{2}\delta \sqrt{t}}|^2\Big)^\frac{1}{2}
\lesssim\frac{1}{\delta^2}
\Big(\int\eta_{\sqrt{T}}|(q_t-\langle q_t\rangle)_{\delta \sqrt{t}}|^2\Big)^\frac{1}{2}.
\end{equation}
\medskip

Next we argue that
\begin{equation}\label{b52}
\int\eta_{\sqrt{T}}
|(q_t-\langle q_t\rangle)_{\delta\sqrt{t}}|^2
\lesssim\frac{1}{\delta^{\frac{d}{2}}}\int\eta_{\sqrt{T}}
|(q_t-\langle q_t\rangle)_{\delta\sqrt{t}}|,
\end{equation}
which amounts to showing
\begin{equation}\label{b50}
\sup|(q_t-\langle q_t\rangle)_{r}|
\lesssim(\frac{\sqrt{t}}{r})^\frac{d}{2}
\end{equation}
for $r\le\sqrt{t}$, and which entails by stationarity
\begin{equation}\label{b53}
\langle|(q_t-\langle q_t\rangle)_{\delta\sqrt{t}}|^2\rangle
\lesssim\frac{1}{\delta^{\frac{d}{2}}}\langle|(q_t-\langle q_t\rangle)_{\delta\sqrt{t}}|\rangle.
\end{equation}
In order to establish (\ref{b50}), we note that by H\"older's inequality and since
the exponential localization function dominates the Gaussian kernel of the same scale
\begin{equation}\label{b54}
|(q_t)_r|\lesssim\Big(\int\eta_{r}|q_t|^2\Big)^\frac{1}{2}
\lesssim(\frac{\sqrt{t}}{r})^\frac{d}{2}\Big(\int\eta_{\sqrt{t}}|q_t|^2\Big)^\frac{1}{2}
\lesssim(\frac{\sqrt{t}}{r})^\frac{d}{2},
\end{equation}
where the last estimate follows from the uniform bound \eqref{Laux3.1} 
from Lemma~\ref{Laux3} applied to \eqref{c62} in form of $(\int\eta_{\sqrt{t}}|q_t|^2)^\frac{1}{2}\lesssim 1$.
Taking the expectation of (\ref{b54}), we obtain by stationarity $|\langle q_t\rangle|=|\langle(q_t)_r\rangle|
\le\langle|(q_t)_r|\rangle\lesssim 
(\frac{\sqrt{t}}{r})^\frac{d}{2}$ so that also
$|(q_t-\langle q_t\rangle)_r|\lesssim 
(\frac{\sqrt{t}}{r})^\frac{d}{2}$.
Since this holds with the origin replaced by any point, we obtain (\ref{b50}).

\medskip

The estimates (\ref{b51}), (\ref{b52}) \& (\ref{b53}) combine to
the following estimate on the rhs of (\ref{b48})
\begin{eqnarray*}
\lefteqn{\frac{1}{\delta}
\Big(\int\eta_{\sqrt{T}}|(\nabla\phi_t,q_t-\langle q_t\rangle)_{\sqrt{2}\delta t}|^2
+\langle|(q_t-\langle q_t\rangle)_{\delta\sqrt{t}}|^2\rangle\Big)^\frac{1}{2}}\nonumber\\
&\lesssim&\frac{1}{\delta^{\frac{d}{4}+3}}
\Big(\int\eta_{\sqrt{T}}|(q_t-\langle q_t\rangle)_{\delta t}|+\langle|(q_t-\langle q_t\rangle)_{\delta\sqrt{t}}|\rangle\Big)^\frac{1}{2}\\
&\lesssim&\frac{1}{\delta^{\frac{d}{2}+7}}\Big(\int\eta_{\sqrt{T}}|(q_t-\langle q_t\rangle)_{\delta t}|
+\langle|(q_t-\langle q_t\rangle)_{\delta\sqrt{t}}|\rangle\Big)
+\delta,
\end{eqnarray*}
where we used Young's inequality in the last step.


\subsection{Proof of Lemma~\ref{Laux8bis}: control of the corrector by the modified corrector}

We split the proof into three steps and start with a reduction argument.

\medskip

\step{1} Reduction.

We shall prove the following stronger version of (\ref{Cr*-2}),
which amounts to an estimate of the systematic error $(\phi-\phi_T,\sigma-\sigma_T)$:
\begin{eqnarray}\label{o30}
\lefteqn{\sup_{R\ge r,\;\mbox{dyadic}}(\frac{R}{r})^\gamma
\fint_{B_R} |\nabla(\phi-\phi_T,\sigma-\sigma_T)|^2}\nonumber\\
&\lesssim&\sup_{R\ge r,\;\mbox{dyadic}}(\frac{R}{r})^\gamma\int \eta_R \Big(\frac{1}{T}|(\phi_T,\sigma_T)|^2
+|(q_{T})_{\sqrt{T}}-\langle q_T \rangle |^2\Big).
\end{eqnarray}
Estimate (\ref{o30}) implies (\ref{Cr*-2}), since by the triangle
inequality, Poincar\'e's inequality,
 $\eta_R \gtrsim R^{-d} 1_{B_R}$, and $\sqrt{T}\lesssim R$,
\begin{eqnarray*}
\lefteqn{\frac{1}{R^2}\fint  \Big|(\phi,\sigma)-\fint_{B_R}(\phi,\sigma)\Big|^2}\\
&\lesssim&\frac{1}{R^2}\fint_{B_R} \Big|(\phi-\phi_T,\sigma-\sigma_T)-\fint_{B_R}(\phi-\phi_T,\sigma-\sigma_T)\Big|^2
+\frac{1}{R^2} \int \eta_R |(\phi_T,\sigma_T)|^2\\
&\lesssim&\fint_{B_R}  |\nabla(\phi-\phi_T,\sigma-\sigma_T)|^2
+\int \eta_R  \frac{1}{T}|(\phi_T,\sigma_T)|^2.
\end{eqnarray*}
We establish (\ref{o30}) by a Campanato iteration based on $\gamma>0$, reducing it to the
one-step but iterable estimate
\begin{eqnarray}\label{o60}
\lefteqn{\fint_{B_r} \Big( |\nabla(\phi-\phi_t)|^2+(\frac{r}{R})^d|\nabla(\sigma-\sigma_t)|^2\Big)}\nonumber\\
&\lesssim&\fint_{B_R}  \Big( |\nabla(\phi-\phi_T)|^2+(\frac{r}{R})^d|\nabla(\sigma-\sigma_T)|^2\Big)\nonumber\\
&+&\int \eta_R \Big(
\frac{1}{T}|(\phi_T,\sigma_T)|^2+(\frac{T}{t}+(\frac{R}{r})^d)\frac{1}{t}|(\phi_t,\sigma_t)|^2
+(\frac{R}{r})^d| (q_t)_{\sqrt{t}}-\langle q_{t}\rangle |^2 \Big),
\end{eqnarray}
which amounts to passing from the pair of radius/cut-off $(r,\sqrt{t})$ to the
pair $(R,\sqrt{T})$, where on the latter we only assume
\begin{equation}\label{o51}
r_*\le r\ll R,\quad\sqrt{t}\lesssim\sqrt{T},\quad\sqrt{t}\lesssim r,\quad\sqrt{T}\lesssim R.
\end{equation}
We split the rest of the proofs into two steps: the proof of  the one-step estimate (\ref{o60}) and the Campanato iteration proper.

\medskip

\step{2} Proof of  (\ref{o60}).

We split the one-step estimate (\ref{o60}) into a part for $\phi$, namely
\begin{multline}\label{o33}
{\fint_{B_r}  |\nabla(\phi-\phi_t)|^2} \lesssim\, \fint_{B_R} |\nabla(\phi-\phi_T)|^2
\\ + \int \eta_R \Big(\frac{1}{T}\phi_T^2+(\frac{T}{t}+(\frac{R}{r})^d)\frac{1}{t}\phi_t^2
+(\frac{R}{r})^d|(q_t)_{\sqrt{t}}-\langle q(t)\rangle|^2\Big)
\end{multline}
and a subordinate part for $\sigma$:
\begin{multline}\label{o50}
{\fint_{B_r} |\nabla(\sigma-\sigma_t)|^2}\, \lesssim \, \fint_{B_R}|\nabla(\sigma-\sigma_T)|^2+(\frac{R}{r})^d|\nabla(\phi-\phi_T)|^2
 \\
+\int \eta_R \Big( \frac{1}{T}|(\phi_T,\sigma_T)|^2+(\frac{T}{t}+(\frac{R}{r})^d)\frac{1}{t}|(\phi_t,\sigma_t)|^2+(\frac{R}{r})^d|(q_t)_{\sqrt{t}}-\langle q(t) \rangle|^2\Big).
\end{multline}
Clearly (\ref{o60}) follows from adding $(\frac{r}{R})^d\times(\ref{o50})$ to (\ref{o33}).

\medskip

\substep{2.1} Proof of (\ref{o33}).

We note that according to (\ref{o56}) and (\ref{c62}),
$\phi-\phi_t$ satisfies
$-\nabla\cdot a\nabla(\phi-\phi_t)=\frac{1}{t}\phi_t$.
In preparation to the energy estimate, 
Hence we may split $\phi-\phi_t=u+w$  on $B_R$
by solving the following two auxiliary boundary value problems on that ball:
\begin{eqnarray}
-\nabla\cdot a \nabla u=0\quad\mbox{in}\;B_R,&&
u=\phi-\phi_t\quad\mbox{on}\;\partial B_R,\label{o34}\\
-\nabla\cdot a\nabla w=\frac{1}{t}\phi_t\quad\mbox{in}\;B_R,&&
w=0\quad\mbox{on}\;\partial B_R.\label{o35}
\end{eqnarray}
By the energy estimate for (\ref{o34}) (i.~e.~testing with $u-(\phi-\phi_t)$ and using the uniform ellipticity)
we have
\begin{equation}\label{o38}
\int_{B_R}|\nabla u|^2\lesssim\int_{B_R}|\nabla(\phi-\phi_t)|^2.
\end{equation}
In preparation for the energy estimate for (\ref{o35}), we split the r.~h.~s. according to
low and high pass.
\begin{equation}\nonumber
\frac{1}{t}\phi_t=\frac{1}{t}(\phi_t)_{\sqrt{t}}+\frac{1}{t}(\phi_{t}-(\phi_t)_{\sqrt{t}})
\stackrel{(\ref{c62})}{=}\nabla\cdot ( (q_t)_{\sqrt t}-\langle q_t\rangle )+\frac{1}{t}(\phi_{t}-(\phi_t)_{\sqrt t}).
\end{equation}
Hence from testing
(\ref{o35}) with $w$ we obtain
\begin{equation}\label{o35bis}
\lambda\int_{B_R}|\nabla w|^2\le -\int_{B_R}\nabla w\cdot( (q_t)_{\sqrt t}-\langle q_t\rangle )
+\int_{B_R}w\frac{1}{t}(\phi_{t}-(\phi_t)_{\sqrt t}).
\end{equation}
Extending $w$ trivially (and continuously)
onto the entire $\mathbb{R}^d$, we see that the second term may be reformulated as
\begin{equation*}
|\int_{B_R}w\frac{1}{t}(\phi_{t}-(\phi_t)_{\sqrt t})|\,=\,|\int(w-(w)_{\sqrt{t}})\frac{1}{t}\phi_{t}|.
\end{equation*}
We insert the weight function $x\mapsto \exp(\frac{|x|}{R})$ to the effect
\begin{multline*}
|\int_{B_R}w\frac{1}{t}(\phi_{t}-(\phi_t)_{\sqrt t})|\,\leq \, \big(\int \exp(\frac{|x|}{R}) |w-(w)_{\sqrt t}|^2(x)dx\big)^\frac12
\\
\times \big( \int\exp(-\frac{|x|}{R}) (\frac{1}{t}\phi_t)^2dx\big)^\frac{1}{2}.
\end{multline*}
As we quickly argue below, since $w$ is supported in $B_R$,  $\sqrt{t}\lesssim r\ll R$ so that $\sqrt{t}\le R$, and the convolution with the Gaussian has scale $\sqrt{t}$, we may bound
the first rhs term by
$$
\int \exp(\frac{|x|}{R}) |w-(w)_{\sqrt t}|^2(x)dx \, \lesssim t \int |\nabla w|^2.
$$
We prove this estimate in its rescaled form $\sqrt{t}\le R=1$ for a function $v$. 
By Jensen's inequality (used twice in a row) and the fundamental theorem of calculus,
\begin{eqnarray*}
\lefteqn{\int \exp(|x|)|v(x)-v_{\sqrt{t}}(x)|^2 dx}
\\
&\leq & \int dz G(z)\int \exp(|x|)|v(x)-v(x-\sqrt t z)|^2 dx 
\\
&\le & t\int dz |z|^2 G(z) \int_0^1 ds \int \exp(|x|)|\nabla v(x-s\sqrt{t} z)|^2dx
\\
&=& t \int \int_0^1 \int |z|^2 G(z)\exp(|y+s\sqrt{t}z|)dzds |\nabla v(y)|^2dy
\\
&\le & t \int \int |z|^2 G(z)\exp(|z|)dz \exp(|y|) |\nabla v(y)|^2 dy
\\
&\lesssim &  t \int |\nabla v(y)|^2 dy,
\end{eqnarray*}
where we used that $\exp(|y|)\le e$ on the support of $\nabla v$ and that $\exp(s\sqrt t|z|)\le \exp(|z|)$.

The combination of these last two estimates yields
$$
|\int_{B_R}w\frac{1}{t}(\phi_{t}-(\phi_t)_{\sqrt t})| \,\lesssim \, \big( t \int |\nabla w|^2  \int\exp(-\frac{|x|}{R}) (\frac{1}{t}\phi_t)^2dx\big)^\frac12.
$$
Inserting this into (\ref{o35bis}), we see that
\begin{equation}\label{o39}
\int_{B_R}|\nabla w|^2\lesssim R^d \int  \Big(\eta_{R} |(q_t)_{\sqrt t}-\langle q_t\rangle |^2
+\eta_{R} \frac{1}{t}\phi_t^2\Big).
\end{equation}

\medskip

We now use the mean-value property on the $a$-harmonic function $u$ in $B_R$
(cf \cite[Lemma~2]{GloriaNeukammOtto}), which is admissible because of $r_*\le r\le R$) and a trivial
estimate on $w$:
\begin{equation}\nonumber
\fint_{B_r}|\nabla u|^2\lesssim\fint_{B_R}|\nabla u|^2\quad\mbox{and}\quad
\fint_{B_r}|\nabla w|^2\le(\frac{R}{r})^d\fint_{B_R}|\nabla w|^2.
\end{equation}
The combination of this with (\ref{o38}) \& (\ref{o39}) yields
\begin{equation}\label{o55}
\fint_{B_r}|\nabla(\phi-\phi_t)|^2\lesssim
\fint_{B_R}|\nabla(\phi-\phi_t)|^2
+\int \eta_R (\frac{R}{r})^d|(q_t)_{\sqrt t}-\langle q_t \rangle |^2
+\int \eta_{R} (\frac{R}{r})^d\frac{1}{t}\phi_t^2.
\end{equation}
We now turn to the estimate of $\phi_T-\phi_t$, which 
by (\ref{c62}) satisfies
$\frac{1}{T}(\phi_T-\phi_t)-\nabla\cdot a\nabla(\phi_T-\phi_t)=(\frac{1}{t}-\frac{1}{T})\phi_t$.
Hence by the Caccioppoli estimate, that is, by testing with $\chi^2(\phi_T-\phi_t)$
where $\chi$ is a cut-off function for $B_R$ in $B_{2R}$, we obtain at first
\begin{eqnarray*}
\lefteqn{\int\chi^2\frac{1}{T}(\phi_T-\phi_t)^2}\\
&+&\nabla(\chi^2(\phi_T-\phi_t))\cdot a\nabla(\phi_T-\phi_t)
=(1-\frac{t}{T})\int\chi^2(\phi_T-\phi_t)\frac{1}{t}\phi_t
\end{eqnarray*}
and then, by Young's inequality to absorb the r.~h.~s. into the massive l.~h.~s. term 
and the usual argument for the elliptic term:
\begin{eqnarray}\label{o40}
\fint_{B_R}|\nabla(\phi_T-\phi_t)|^2
&\lesssim&\fint_{B_{2R}}\frac{1}{R^2}(\phi_T-\phi_t)^2+(\frac{T}{t})\frac{1}{t}\phi_t^2\nonumber\\
&\stackrel{\frac{1}{R^2}\lesssim\frac{1}{T}\lesssim\frac{1}{t}}{\lesssim}&
\fint_{B_{2R}}\frac{1}{T}\phi_T^2+(\frac{T}{t})\frac{1}{t}\phi_t^2.
\end{eqnarray}
Together with (\ref{o55}), this yields
by the triangle inequality
\begin{eqnarray*}
\lefteqn{\fint_{B_r}|\nabla(\phi-\phi_t)|^2}\nonumber\\
&\lesssim&\fint_{B_{R}}|\nabla(\phi-\phi_T)|^2
+\int \eta_R (\frac{R}{r})^d|(q_t)_{\sqrt t}-\langle q_t\rangle |^2\\
&&+\int \eta_{2R}\Big( \frac{1}{T}\phi_T^2+(\frac{T}{t}+(\frac{R}{r})^d)\frac{1}{t}\phi_t^2 \Big)\\
&\lesssim&\fint_{B_{2R}}|\nabla(\phi-\phi_T)|^2+
\int \eta_{2R} \Big(\frac{1}{T}\phi_T^2
+(\frac{T}{t}+(\frac{R}{r})^d)\frac{1}{t}\phi_t^2
+(\frac{R}{r})^d|(q_t)_{\sqrt t}-\langle q_t\rangle |^2\Big),
\end{eqnarray*}
which yields (\ref{o33}) by replacing $2R$ with $R$.

\medskip

\substep{2.2} Proof of (\ref{o55}). 

We note that according to (\ref{o58}) and (\ref{o17}),
$\sigma-\sigma_t$ satisfies $-\triangle(\sigma-\sigma_t)=\nabla\times(q-q_t)+\frac{1}{t}\sigma_t$.
On $B_R$, we split $\sigma-\sigma_t=u+w$ according to 
\begin{eqnarray}
-\triangle u=0\quad\mbox{in}\;B_R,&&
u=\sigma-\sigma_t\quad\mbox{on}\;\partial B_R,\label{o45}\\
-\triangle w=\nabla\times(q-q_t)+\frac{1}{t}\sigma_t\quad\mbox{in}\;B_R,&&
w=0\quad\mbox{on}\;\partial B_R.\label{o46}
\end{eqnarray}
By the mean value inequality and the variational characterization
for the harmonic $u$ defined through (\ref{o45}) 
we have (incidentally as clean inequalities)
\begin{equation}\label{o47}
\fint_{B_r}|\nabla u|^2\le\fint_{B_R}|\nabla u|^2\le\fint_{B_R}|\nabla(\sigma-\sigma_t)|^2.
\end{equation}
Also for $w$ defined through (\ref{o46}) we proceed as before: Writing the second part of the r.~h.~s. as
\begin{equation}\nonumber
\frac{1}{t}\sigma_t=\frac{1}{t}(\sigma_t)_{\sqrt t}+\frac{1}{t}(\sigma_t-(\sigma_t)_{\sqrt t})
\stackrel{(\ref{o17})}{=}\triangle(\sigma_t)_{\sqrt t}+\nabla\times((q_t)_{\sqrt t}-\langle q_t\rangle)
+\frac{1}{t}(\sigma_t-(\sigma_t)_{\sqrt t}),
\end{equation}
we obtain from the energy estimate
\begin{equation}\nonumber
\int_{B_R}|\nabla w|^2\lesssim\int_{B_R}|q-q_t|^2+|\nabla (\sigma_t)_{\sqrt t}|^2
+|(q_t)_{\sqrt t}-\langle q_t\rangle|^2+R^d \int \eta_{2R} \frac{1}{t}|\sigma_t|^2,
\end{equation}
where only the term with $\nabla(\sigma_t)_{\sqrt t}$ is of a new type, but can be easily controlled
through the inverse estimate
\begin{equation}\nonumber
\int_{B_R}|\nabla (\sigma_t)_{\sqrt t}|^2\lesssim R^d \int \eta_{2R}\frac{1}{t}|\sigma_t|^2
\end{equation}
and thus is contained in one of the existing terms. Together with $|q-q_t|\le|\nabla(\phi-\phi_t)|$, 
and $\sqrt{t}\ll R$ we obtain
\begin{equation}\nonumber
\int_{B_R}|\nabla w|^2\lesssim\int_{B_R}|\nabla(\phi-\phi_t)|^2
+R^d \int \eta_{2R} \big(|(q_t)_{\sqrt t}-\langle q_t\rangle|^2+ \frac{1}{t}|\sigma_t|^2 \big).
\end{equation}
In combination with the trivial estimate $\fint_{B_r}|\nabla w|^2\le(\frac{R}{r})^d\fint_{B_R}|\nabla w|^2$
and (\ref{o47}) we obtain by the triangle inequality 
\begin{eqnarray}\label{o49}
{\fint_{B_r}|\nabla(\sigma-\sigma_t)|^2}
&\lesssim&
\fint_{B_{R}}|\nabla(\sigma-\sigma_t)|^2+(\frac{R}{r})^d|\nabla(\phi-\phi_t)|^2
\nonumber\\
&&+\int \eta_{2R} (\frac{R}{r})^d\big(\frac{1}{t}|\sigma_t|^2+|(q_t)_{\sqrt t}-\langle q_t\rangle|^2 \big).
\end{eqnarray}
We now consider $\sigma_T-\sigma_t$ which in view of (\ref{o17}) 
satisfies $\frac{1}{T}(\sigma_T-\sigma_t)
-\triangle(\sigma_T-\sigma_t)=\nabla\times a\nabla(\phi_T-\phi_t)+(\frac{1}{t}-\frac{1}{T})\sigma_t$.
We thus obtain from the Caccioppoli estimate 
\begin{eqnarray}
\lefteqn{\fint_{B_R}\frac{1}{T}|\sigma_T-\sigma_t|^2+|\nabla(\sigma_T-\sigma_t)|^2}\nonumber\\
&\lesssim&\fint_{B_{2R}}\frac{1}{R^2}|\sigma_T-\sigma_t|^2+|\nabla(\phi_T-\phi_t)|^2
+\frac{T}{t^2}|\sigma_t|^2\nonumber\\
&\stackrel{\frac{1}{R^2}\lesssim\frac{1}{T}\lesssim\frac{1}{t}}{\lesssim}&
\fint_{B_{2R}}|\nabla(\phi_T-\phi_t)|^2+\frac{1}{T}|\sigma_T|^2+(\frac{T}{t})\frac{1}{t}|\sigma_t|^2\nonumber\\
&\stackrel{(\ref{o40}) \mbox{ with } R\leadsto2R}{\lesssim}&
\fint_{B_{4R}}\frac{1}{T}|(\phi_T,\sigma_T)|^2+(\frac{T}{t})\frac{1}{t}|(\phi_t,\sigma_t)|^2.\nonumber
\end{eqnarray}
Combined (via the triangle inequality) with (\ref{o49}) and once more (\ref{o40}),
this gives (\ref{o50}) in form of
\begin{eqnarray*}
\lefteqn{\fint_{B_r}|\nabla(\sigma-\sigma_t)|^2}\nonumber\\
&\lesssim&\fint_{B_{R}}|\nabla(\sigma-\sigma_T)|^2
+(\frac{R}{r})^d|\nabla(\phi-\phi_T)|^2\nonumber\\
&&+\int \eta_{2R} \Big(\frac{1}{T}|(\phi_T,\sigma_T)|^2+(\frac{T}{t}+(\frac{R}{r})^d)\frac{1}{t}|(\phi_t,\sigma_t)|^2
+(\frac{R}{r})^d |(q_t)_{\sqrt t}-\langle q_t\rangle|^2\Big).
\end{eqnarray*}

\medskip

\step{3} Campanato iteration.

We now address the argument that leads from the one-step estimate
(\ref{o60}) to (\ref{o30}). In line with (\ref{Cr*-1}), we define for some dyadic
$\theta\ll 1$ to be selected later
\begin{equation}\label{o57}
(R_0,T_0)=(r,t)\quad\mbox{and}\quad (R_n,T_n)=(\theta^{-1} R_{n-1},\theta^{-\alpha} T_{n-1}).
\end{equation}
For any $n\in\mathbb{N}$ we now apply (\ref{o60}) to $(r,t)=(R_{n-1},T_{n-1})$ and $(R,T)=(R_n,T_n)$ yielding
\begin{multline}\label{o41}
{\fint_{B_{R_{n-1}}}|\nabla(\phi-\phi_{T_{n-1}})|^2+\theta^d|\nabla(\sigma-\sigma_{T_{n-1}})|^2}
\\
\le\,C_0\Big(\fint_{B_{R_n}}|\nabla(\phi-\phi_{T_n})|^2+\theta^d|\nabla(\sigma-\sigma_T)|^2 \Big)
\\
 +C_0\int \eta_{R_n} \Big( \theta^{-d}|(q_{T_{n-1}})_{\sqrt{T_{n-1}}} - \langle q_{T_{n-1}}\rangle |^2\\
+\frac{1}{T_n}|(\phi_{T_n},\sigma_{T_n})|^2+
\theta^{-d}\frac{1}{T_{n-1}}|(\phi_{T_{n-1}},\sigma_{T_{n-1}})|^2\Big),
\end{multline}
where $C_0<\infty$ is a constant only depending on $d$ and $\lambda$ whose value we'd like
to retain for a moment.
It is convenient to introduce
\begin{equation}\nonumber
\Lambda:=\sup_{n=0,1,\cdots}(\frac{R_n}{r})^\gamma\fint_{B_{R_n}}|\nabla(\phi-\phi_{T_{n}})|^2
+\theta^d|\nabla(\sigma-\sigma_{T_n})|^2,
\end{equation}
and for simplicity we shall assume that it is finite: Clearly, this assumption
requires an approximation argument,
a possible approximation argument being to replace $(\phi,\sigma)$ by $(\phi_{\bar T},\sigma_{\bar T})$
and to stop increasing $T_n$ once $\bar T$ is reached; the estimate of this proposition
will then hold uniformly in $\bar T$, so that the desired statement is obtained in the limit 
$\bar T\uparrow\infty$ almost-surely, 
using the qualitative theory in form of $(\phi_{\bar T},\sigma_{\bar T})\rightarrow(\phi,\sigma)$
following from mere ergodicity and stationarity of the ensemble $\langle\cdot\rangle$.
We may rewrite (\ref{o41}) as
\begin{multline*}
\lefteqn{(\frac{R_{n-1}}{r})^\gamma\fint_{B_{R_{n-1}}}|\nabla(\phi-\phi_{T_{n-1}})|^2
+\theta^d|\nabla(\sigma-\sigma_{T_{n-1}})|^2}
\\
\le\,C_0\theta^\gamma\bigg(\Lambda +(\frac{R_n}{r})^\gamma\int \eta_{R_n} \Big( \frac{1}{T_n}|(\phi_{T_n},\sigma_{T_n})|^2
+\frac{\theta^{-d}}{T_{n-1}}|(\phi_{T_{n-1}},\sigma_{T_{n-1}})|^2
\\
+\theta^{-d}|(q_{T_{n-1}})_{\sqrt{T_{n-1}}}-\langle q_{T_{n-1}} \rangle|^2\Big)\bigg).
\end{multline*}
Now choosing $\theta=\theta(d,\lambda,\gamma)>0$ so small that $C_0\theta^\gamma\le\frac{1}{2}$
and taking the supremum over $n\in\mathbb{N}$ we obtain
\begin{eqnarray*}
\Lambda&\lesssim&\sup_{n\in\mathbb{N}}(\frac{R_n}{r})^\gamma \int \eta_{R_n}\Big(\frac{1}{T_n}|(\phi_{T_n},\sigma_{T_n})|^2
+\frac{1}{T_{n-1}}|(\phi_{T_{n-1}},\sigma_{T_{n-1}})|^2
\\
&& \qquad \qquad  \qquad \qquad  \qquad \qquad  \qquad \qquad +|(q_{T_{n-1}})_{\sqrt{T_{n-1}}}-\langle q_{T_{n-1}} \rangle|^2\Big)\\
&\sim&\sup_{n\in\mathbb{N}}(\frac{R_n}{r})^\gamma\int \eta_{R_n} \Big(
\frac{1}{T_{n-1}}|(\phi_{T_{n-1}},\sigma_{T_{n-1}})|^2
+|(q_{T_{n-1}})_{\sqrt{T_{n-1}}}-\langle q_{T_{n-1}} \rangle|^2 \Big).
\end{eqnarray*}
Since $T_{n-1}\sim T_n=(\frac{R_n}{r})^\alpha t$ by (\ref{o57}), we obtain (\ref{o30}) with scales related by (\ref{Cr*-1}).


\subsection{Proof of Lemma~\ref{Laux9}: approximation of homogenized semi-group}

By definition of the propagators $S^h_{\frac{T}{2}\rightarrow T}$
and $S^{\ho}_{\frac{T}{2}\rightarrow T}$, cf (\ref{au37}) with $S$ replaced by $S^{\ho}$ and (\ref{au40}),
combined with $\nabla\cdot q({\textstyle\frac{T}{2}})=u({\textstyle\frac{T}{2}})$,
we have
\begin{align}\label{au54}
(S^h_{\frac{T}{2}\rightarrow T}-S^{\ho}_{\frac{T}{2}\rightarrow T})q({\textstyle\frac{T}{2}})
=\int_{\frac{T}{2}}^T(a_{hT}\nabla u_hd\tau-a_{\ho}\nabla u_{\ho})d\tau,
\end{align}
where $u_h$ and $u_{\ho}$ are the solutions of the following initial value problem with constant coefficients
\begin{align}
\partial_\tau u_h-\nabla\cdot a_{hT}\nabla u_h=0\;\;\mbox{for}\;\tau>{\textstyle\frac{T}{2}},&\quad
u_h=u\;\;\mbox{for}\;\tau={\textstyle\frac{T}{2}},\nonumber\\
\partial_\tau u_{\ho}-\nabla\cdot a_{\ho}\nabla u_{\ho}=0\;\;\mbox{for}\;\tau>{\textstyle\frac{T}{2}},&\quad
u_{\ho}=u\;\;\mbox{for}\;\tau={\textstyle\frac{T}{2}},\label{au53}
\end{align}
so that their difference $w:=u_h-u_{\ho}$ satisfies
\begin{align}\label{au52}
\partial_\tau w-\nabla\cdot a_{hT}\nabla w=-\nabla\cdot(a_{hT}-a_{\ho})\nabla u_{\ho}&
\quad\mbox{for}\;\tau>{\textstyle\frac{T}{2}},\\
w=0&\quad\mbox{for}\;\tau={\textstyle\frac{T}{2}}.\nonumber
\end{align}
By the localized energy estimates of Lemma~\ref{Laux3} for (\ref{au52}) and (\ref{au53}) we have
\begin{align*}
\big(\int_{\frac{T}{2}}^T\int\eta_{\sqrt{T}}|\nabla w|^2d\tau\big)^\frac{1}{2}
&\lesssim|a_{hT}-a_{\ho}|\big(\int_{\frac{T}{2}}^T\int\eta_{\sqrt{T}}|\nabla u_{\ho}|^2d\tau\big)^\frac{1}{2},\\
\big(\int_{\frac{T}{2}}^T\int\eta_{\sqrt{T}}|\nabla u_{\ho}|^2d\tau\big)^\frac{1}{2}
&\lesssim\big(\int\eta_{\sqrt{T}}|u({\textstyle\frac{T}{2}})|^2\big)^\frac{1}{2},
\end{align*}
so that, inserting the second energy estimate into the first, which we upgrade to
\begin{align*}
\lefteqn{\big(\int_{\frac{T}{2}}^T\int\eta_{\sqrt{T}}|a_h\nabla u_h-a_{\ho}\nabla u_{\ho}|^2d\tau\big)^\frac{1}{2}}\\
&\lesssim|a_{hT}-a_{\ho}|\big(\int_{\frac{T}{2}}^T\int\eta_{\sqrt{T}}|\nabla u_{\ho}|^2d\tau\big)^\frac{1}{2},
\end{align*}
we obtain
\begin{align*}
\big(\int_{\frac{T}{2}}^T\int\eta_{\sqrt{T}}|a_h\nabla u_h-a_{\ho}\nabla u_{\ho}|^2d\tau\big)^\frac{1}{2}
&\lesssim|a_{hT}-a_{\ho}|\big(\int\eta_{\sqrt{T}}|u({\textstyle\frac{T}{2}})|^2\big)^\frac{1}{2}.
\end{align*}
By Cauchy-Schwarz' inequality in time and space (followed by the dominance of the Gaussian kernel by the exponential cut-off) in form of
\begin{align*}
\lefteqn{(\frac{R}{\sqrt{T}})^\frac{d}{2}|(\int_{\frac{T}{2}}^T(a_h\nabla u_h-a_{\ho}\nabla u_{\ho})d\tau)_R|}\nonumber\\
&\lesssim\sqrt{T}
\big(\int_{\frac{T}{2}}^T\int\eta_{\sqrt{T}}|a_h\nabla u_h-a_{\ho}\nabla u_{\ho}|^2d\tau\big)^\frac{1}{2}
\end{align*}
and the representation (\ref{au54}), the claim (\ref{Laux9.1}) follows from Lemma~\ref{Laux6}.


\subsection{Proof of Lemma~\ref{Laux10}: systematic error}

 In view of the relation $q_T=\int_0^\infty q(t)\exp(-\frac{t}{T})\frac{dt}{T}$,
which entails $q_T-q=\int_0^\infty(q(t)-q)\exp(-\frac{t}{T})\frac{dt}{T}$, and the definition
$q(t)=a(\int_0^\infty\nabla u(\tau)d\tau+e)$, we obtain the formula
\begin{align*}
q_T-q=\int_0^\infty(1-\exp(-\frac{\tau}{T}))a\nabla u(\tau) d\tau.
\end{align*}
By the definition $(a_{hT}-a_{\ho})e=\langle q_T-q\rangle$, we arrive at
\begin{align*}
(a_{hT}-a_{\ho})e=\int_0^\infty(1-\exp(-\frac{\tau}{T}))\langle a\nabla u(\tau)\rangle d\tau,
\end{align*}
which we will use in form of the inequality
\begin{align}\label{au42}
|(a_{hT}-a_{\ho})e|\le\int_0^\infty(1-\exp(-\frac{\tau}{T}))\langle|\nabla u(\tau)|\rangle d\tau.
\end{align}
We now appeal to \eqref{Laux6.2} in Lemma~\ref{Laux6} in form 
\begin{align*}
\tau\int\eta_{\sqrt{2\tau}}|\nabla u(\tau)|
\lesssim\fint_{\frac{\tau}{2}}^\tau dt\fint_0^{\sqrt{t}}dr(\frac{r}{\sqrt{t}})^\frac{d}{2}\int\eta_{2\sqrt{2\tau}}|(q(t)-\langle q(t))_r|,
\end{align*}
which by stationarity turns into
\begin{align*}
\tau\langle|\nabla u(\tau)|\rangle
\lesssim\frac{1}{\tau}\int_{\frac{\tau}{2}}^\tau dt
\fint_0^{\sqrt{t}}dr(\frac{r}{\sqrt{t}})^\frac{d}{2}\langle|(q(t)-\langle q(t))_r|\rangle
\end{align*}
We now insert this into (\ref{au42}) to the effect of
\begin{align}\nonumber
\lefteqn{|(a_{hT}-a_{\ho})e|}\nonumber\\
&\lesssim\int_0^\infty dt\int_{\frac t2}^{t}\frac{d\tau}{t^2}(1-\exp(-\frac{2\tau}{T}))
\fint_0^{\sqrt{\tau}}dr(\frac{r}{\sqrt{t}})^\frac{d}{2}\langle|(q(\tau)-\langle q(\tau))_r|\rangle\nonumber\\
&\lesssim\int_0^\infty \frac{dt}{t}\min\{\frac{t}{T},1\}
\fint_0^{\sqrt{t}}dr(\frac{r}{\sqrt{t}})^\frac{d}{2}\langle|(q(t)-\langle q(t))_r|\rangle\nonumber.
\end{align}
%


\section{Proofs of the stochastic auxiliary results}

\subsection{Proof of Proposition~\ref{Paux2}: stochastic propagator estimate}
We split the proof into two steps.

\medskip

\step1 Proof of 
\begin{eqnarray}
\norm (S_{\frac T2 \to T}-S^\ho_{\frac T2\to T})q(\tfrac T2)\norm
&\lesssim&  \sup_{R\le \sqrt{T}} R^\frac{d}{2} \| ((S_{\frac T2 \to T}-S^\ho_{\frac T2\to T})q(\tfrac T2))_R\|_*\nonumber
\\
&&+\Big( \sup_{R\le \sqrt{T}} R^\frac{d}{2} \| ((S_{\frac T2 \to T}-S^\ho_{\frac T2\to T})q(\tfrac T2))_R\|_* \Big)^{\frac12}\nonumber
\\
&& \qquad \times \Big(\sup_{t\le T,R\le \sqrt{t}} R^\frac{d}{2} \|(q(t))_R\|_*\Big)^\frac12 \label{Paux2-1}
\end{eqnarray}
for all $T\ge 1$, $s\le 2$.

By Lemma~\ref{Laux11},  the fields $F=q(T),S_{\frac{T}{2}\to T}^\ho q(\frac{T}{2}),S_{\frac{T}{2}\to T}^h q(\frac{T}{2})$ are local on scale $\sqrt{T}\ge 1$ relative to 
$$
\bar F\,:=\,\fint_0^T dt\Big(\frac{\sqrt{t}}{\sqrt{T}}\Big)^\frac{d}{2}\fint_0^{\sqrt{t}}dr \Big(\frac{r}{\sqrt{t}}\Big)^\frac{d}{2} |(q(t)-\langle q(t) \rangle)_r|,
$$
so that by \eqref{Laux11.1},
$$
\Big( \int_{B_{\sqrt{T}}} |F(a)-F(\tilde a)|^2\Big)^\frac12 \,\lesssim\, \Big(\frac{\sqrt{T}}{R}\Big)^p \int \eta_R (\bar F(a)+\bar F(\tilde a))
$$
for all $a,\tilde a \in \Omega$ such that $a=\tilde a$ on $B_{2R}:=\{|x|<2R\}$ and $R\ge \sqrt{T}$, for any $p<\infty$.
Hence, by \eqref{Laux12.3} in Lemma~\ref{Laux12}, we have for all $s\le 2$,
\begin{eqnarray*}
\lefteqn{\sup_{R\ge \sqrt{T}} \Big(\frac{R}{\sqrt{T}}\Big)^\frac{d}{2} \| ((S_{\frac T2 \to T}-S^\ho_{\frac T2\to T})q(\tfrac T2))_R\|_*}
\\
&\lesssim & \sup_{r\le \sqrt{T}} \Big(\frac{r}{\sqrt{T}}\Big)^\frac{d}{2} \| ((S_{\frac T2 \to T}-S^\ho_{\frac T2\to T})q(\tfrac T2))_r\|_*
\\
&&+\Big( \sup_{r\le \sqrt{T}} \Big(\frac{r}{\sqrt{T}}\Big)^\frac{d}{2} \| ((S_{\frac T2 \to T}-S^\ho_{\frac T2\to T})q(\tfrac T2))_r\|_* \Big)^{\frac12}
\\
&& \qquad \times \Big(\sup_{t\le T} \Big(\frac{\sqrt{t}}{\sqrt{T}}\Big)^\frac{d}{2}\sup_{r\le \sqrt{t}} \Big(\frac{r}{\sqrt{t}}\Big)^\frac{d}{2} \|(q(t))_r\|_*\Big)^\frac12.
\end{eqnarray*}
The corresponding inequality being obvious in the range $R\le \sqrt{T}$ on the lhs, this yields \eqref{Paux2-1} for all $T\ge 1$ and $s\le 2$.

\medskip

If $\delta$ was deterministic in Proposition~\ref{Paux1}, the desired estimate \eqref{Paux2.1} would follow from \eqref{Paux2-1} and from taking the $\|\cdot\|_*$-norm of
\eqref{Paux1.1}. However, $\delta$ is random in view of \eqref{Paux1.2}, and we shall need to appeal to Lemma~\ref{Laux13} --- which we compel us to restrict to $s<2$. 

\medskip

\step2 Proof of \eqref{Paux2.1}.

Fix $T\ge 1$, $R\le\sqrt{T}$, $\delta\le 1$, $\sqrt{t_0}\le\delta\sqrt{T}$ and introduce
the following abbreviations for the random variables under consideration
\begin{eqnarray}\nonumber
F&:=&(\frac{R}{\sqrt{T}})^{\frac{d}{2}}
|(q(T)-S^{\ho}_{\frac{T}{2}\rightarrow T}q({\textstyle\frac{T}{2}}))_R|,\nonumber\\
F_0&:=&\int\eta_{\sqrt{T}}|(q_{t_0}-\langle q_{t_0}\rangle)_{\delta\sqrt{t_0}}|,\nonumber\\
F_1&:=&C\fint_{\frac{T}{8}}^\frac{T}{2}dt\fint_0^{\sqrt{t}}dr
(\frac{r}{\sqrt{t}})^\frac{d}{2}
\int\eta_{\sqrt{T}}|(q(t)-\langle q(t)\rangle)_r|\nonumber.
\end{eqnarray}
According to Proposition~\ref{Paux1} we have for $\delta\ll 1$
(and provided $C\gg 1$)
\begin{align}\label{P2-1}
F_0+\langle F_0\rangle\le 3\delta^p\quad\Longrightarrow\quad |F|\le(\delta+\partial a)^\frac{1}{p}|F_1|,
\end{align}
where $\partial a:=|a_{hT}-a_{\ho}|$ is deterministic and
where we recall that $p$ denotes here a large generic exponent only depending on $d$ and $\lambda>0$.
By Lemma~\ref{Laux13},  (\ref{P2-1}) implies on this abstract level for all $s<2$
\begin{align}\label{P2-2}
\langle F_0\rangle\le\delta^p\quad\Longrightarrow\quad
\|F\|_{s}\leq  (\delta+\partial a)^\frac{1}{p}\|F_1\|_s+\big(\frac{\|F_0\|_{*,2}}{\delta^p}\big)^\frac{2}{s}\|F\|_{\infty}.
\end{align}

\medskip

We are in the position to conclude the proof of the proposition.
By \eqref{Laux14.1} in
Lemma~\ref{Laux14} we have
$\|F\|_\infty\lesssim 1$. Since
$F_0\stackrel{\eqref{c84}}{\sim} \int\eta_{\sqrt{T}}|(q_{t_0}-\langle q_{t_0}\rangle)_{\delta\sqrt{t_0}}|_{\sqrt{T}}$ we have by
stationarity $\|F_0\|_{*,2}\lesssim\||(q_{t_0}-\langle q_{t_0}\rangle)_{\delta\sqrt{t_0}}|_{\sqrt{T}}\|_{*,2}$
so that by the approximate locality \eqref{Laux14.3} and the uniform bound \eqref{Laux14.2} in Lemma~\ref{Laux14}  we obtain by Lemma~\ref{Laux12} 
$\|F_0\|_{*,2}\lesssim (\frac{\sqrt{t_0}}{\delta\sqrt{T}})^\frac{d}{2}$. Finally, we clearly
have 
$$\|F_1\|_s\le C\sup_{\frac{T}{4}\le t\le\frac{T}{2}}\sup_{r\le\sqrt{t}}
(\frac{r}{\sqrt{t}})^\frac{d}{2}\|(q(t)-\langle q(t)\rangle)_r\|_{s}$$
 and obviously
$\|F\|_s=(\frac{R}{\sqrt{T}})^\frac{d}{2}\|(q(T)-S^{\ho}_{\frac{T}{2}\rightarrow T}q({\textstyle\frac{T}{2}}))_R\|_{s}$.
Hence, taking  the supremum over $R\le\sqrt{T}$ yields provided  $\langle F_0\rangle\le\delta^p$
\begin{equation}\label{Paux2-2}
\sup_{R\le \sqrt{T}} R^\frac{d}{2}\|(q(T)-S^{\ho}_{\frac{T}{2}\rightarrow T}q({\textstyle\frac{T}{2}}))_R\|_{s}
\,\lesssim \, (\delta+\partial a)^\frac{1}{p} \sup_{t\le \sqrt{T}} \norm q(t)\norm+ \Big(\frac{1}{\sqrt{T}}\Big)^{\frac{d}{2}(\frac{2}{s}-1)} \Big(\frac{{t_0}^\frac{d}{4}}{\delta^p}\Big)^\frac{2}{s},
\end{equation}
up to redefining $p$.
The desired estimate (\ref{Paux2.1}) then follows from the combination of \eqref{Paux2-1} and \eqref{Paux2-2}.

\subsection{Proof of Lemma~\ref{Laux11}: relative approximate locality}

We shall prove a stronger statement and replace the lhs of \eqref{Laux11.1} and \eqref{Laux11.2}
by $\Big(\int_{B_R} |F(a)-F(\tilde a)|^2\Big)^\frac12$.
We split the proof into two steps.
For the proof of the approximate locality relative to 1, we refer to the proof of Lemma~\ref{Laux14}.

\medskip

\step1 A stronger version of the parabolic Caccioppoli estimate. 

For some radius
$R\ge 1$ consider a solution $w$ of
\begin{align}\label{L11-1}
\left.\begin{array}{rl}
\partial_\tau w-\nabla\cdot\bar a\nabla w=0&\quad\mbox{for}\;\tau>0\\
w=0&\mbox{for}\;\tau=0
\end{array}\right\}\quad\mbox{in}\;(0,1)\times B_{2R},
\end{align}
where $\bar a$ is a $\lambda$-uniform elliptic coefficient field, which may depend on
space and time.

Then we have for any exponent $p<\infty$ 
\begin{equation}\label{L11-2}
\Big(\int_0^1\int_{B_R}|\nabla w|^2 dt\Big)^\frac{1}{2}\lesssim\frac{1}{R^p}\int_0^1t^p\big(\int_{B_{2R}\setminus B_R}w^2\Big)^\frac{1}{2}dt.
\end{equation}
This estimate expresses the fact that a solution $w$ of the homogeneous equation
with homogeneous initial data
in $(0,1)\times B_R$ depends only sub-algebraically on $(0,1)\times(B_{2R}\setminus B_R)$,
and even less so on $(0,t)\times(B_{2R}\setminus B_R)$ for $t\ll 1$.

\medskip

Provided the function $\chi=\chi(x)$ is supported in $B_{2R}$ we obtain from (\ref{L11-1}) that
$\frac{d}{dt}\frac{1}{2}\int\chi^2w^2=-\int\nabla(\chi^2 w)\cdot a\nabla w$. Together with
Leibniz' rule $\nabla(\chi^2 w)\cdot a\nabla w=\chi^2\nabla w\cdot a\nabla w+2\chi w\nabla\chi\cdot a\nabla w$,
the bounds (\ref{1.3}) on $a$ yielding $\nabla(\chi^2 w)\cdot a\nabla w\ge\lambda\chi^2|\nabla w|^2-
2|\chi||w||\nabla\chi||\nabla w|\ge\frac{1}{C}\chi^2|\nabla w|^2-C|\nabla\chi|^2w^2$ we obtain the
differential inequality
\begin{equation}\nonumber
\frac{d}{dt}\int\chi^2w^2+\frac{1}{C}\int\chi^2|\nabla w|^2\le C\int|\nabla\chi|^2w^2.
\end{equation}
Replacing $\chi$ by $\chi^{\tilde p}$ where $\chi$ is a cut-off for $B_R$ in $B_{2R}$, this turns into
\begin{equation}\nonumber
\frac{d}{dt}\int\chi^{2{\tilde p}}w^2+\frac{1}{C}\int\chi^{2\tilde p}|\nabla w|^2\le \frac{C}{R^2}\int_{B_{2R}\setminus B_R}\chi^{2(\tilde p-1)}w^2.
\end{equation}
Applying H\"older's inequality on the rhs term this yields
\begin{equation}\label{L11-3}
\frac{d}{dt}\int\chi^{2\tilde p}w^2+\frac{1}{C}\int\chi^{2\tilde p}|\nabla w|^2\le \frac{C}{R^2}
\big(\int\chi^{2\tilde p}w^2\big)^{1-\frac{1}{\tilde p}}\big(\int_{B_{2R}\setminus B_R}w^2\big)^{\frac{1}{\tilde p}},
\end{equation}
so that neglecting the second lhs term we obtain
\begin{equation}\nonumber
\frac{d}{dt}\big(\int\chi^{2\tilde p}w^2\big)^\frac{1}{\tilde p}\lesssim\frac{1}{R^2}\big(\int_{B_{2R}\setminus B_R}w^2\big)^\frac{1}{\tilde p}.
\end{equation}
Integrating this in time while making use of the homogeneous initial data, cf (\ref{L11-1}), this yields
\begin{equation}\label{L11-4}
\sup_{t\le 1}\int\chi^{2\tilde p}w^2\lesssim\Big(\frac{1}{R^2}\int_0^1\big(\int_{B_{2R}\setminus B_R}w^2\big)^\frac{1}{\tilde p}dt\Big)^{\tilde p}.
\end{equation}
A direct integration of (\ref{L11-3}) in time gives
\begin{equation}\nonumber
\int_0^1\int\chi^{2\tilde p}|\nabla w|^2 dt\lesssim
\big(\sup_{t\le 1}\int\chi^{2\tilde p}w^2\big)^{1-\frac{1}{\tilde p}}\frac{1}{R^2}\int_0^1\big(\int_{B_{2R}\setminus B_R}w^2\big)^{\frac{1}{\tilde p}}dt,
\end{equation}
into which we plug (\ref{L11-4}) to obtain by choice of $\chi$
\begin{equation}\nonumber
\int_0^1\int_{B_R}|\nabla w|^2 dt\le\int_0^1\int\chi^{2\tilde p}|\nabla w|^2 dt\lesssim
\Big(\frac{1}{R^2}\int_0^1\big(\int_{B_{2R}\setminus B_R}w^2\big)^{\frac{1}{\tilde p}}dt\Big)^{\tilde p}.
\end{equation}
Applying H\"older's inequality in time on the rhs, this yields for all $p\ge 0$ in case of $\tilde p>2$
\begin{equation}\nonumber
\big(\int_0^1\int_{B_R}|\nabla w|^2 dt\big)^\frac{1}{2}
\lesssim\frac{1}{R^{\tilde p}}
\Big(\int_0^1(\frac{1}{t^p})^\frac{2}{\tilde p-2}dt\Big)^\frac{\tilde p-2}{2}\int_0^1 t^p\big(\int_{B_{2R}\setminus B_R}w^2\big)^{\frac{1}{2}}dt.
\end{equation}
For $p<\frac{\tilde p-2}{2}$, the first rhs integral is finite. Since in this case we have in particular $p\le \tilde p$,
we may replace $\frac{1}{R^{\tilde p}}$ by $\frac{1}{R^p}$ to obtain (\ref{L11-2}).

\medskip

\step2 Conclusion. 

By scaling, we may assume $T=1$. By the triangle inequality, it is enough to show that
the four components
$\nabla \phi(1), S^{h}_{\frac{1}{2}\rightarrow 1}q({\textstyle\frac{1}{2}})$,
$S^{\ho}_{\frac{1}{2}\rightarrow 1}q({\textstyle\frac{1}{2}})$, and
$q(1)$ are approximately local on scale $1$ wrt to the same $\bar F$.
For simplicity, we focus on
$F:=S^{h}_{\frac{1}{2}\rightarrow 1}q({\textstyle\frac{1}{2}})\stackrel{(\ref{au59})}{=}
S^{h}_{\frac{1}{2}\rightarrow 1}S_{0\rightarrow\frac{1}{2}}ae$
and note that by the same argument used in establishing the semi-group property for a propagator in
the proof of Theorem \ref{Tau}, it may be written as follows
\begin{align*}
F=
ae+\int_0^1\bar a\nabla \bar ud\tau,
\end{align*}
where $\bar u$ solves the initial value problem
\begin{align*}
\partial_\tau\bar u-\nabla\cdot\bar a\nabla \bar u=0\;\;\mbox{for}\;\tau\ge 0,\quad
\bar u=\nabla\cdot(\bar ae)\;\;\mbox{for}\;\tau=0,
\end{align*}
and the space-time-dependent coefficient field $\bar a$ is defined via 
\begin{align*}
\bar a(\tau)=a\;\;\mbox{for}\;\tau<{\textstyle\frac{1}{2}},\quad
\bar a(\tau)=a_{h1}\;\;\mbox{for}\;\tau>{\textstyle\frac{1}{2}}.
\end{align*}
Hence if we are given two (spatial) coefficient fields $a$ and $\tilde a$ which agree on $B_{2R}$,
also $\bar a$ and $\bar{\tilde a}$ agree on $B_{2R}$, so that $w:=\bar u-\bar{\tilde u}$ satisfies
on $(0,1)\times B_{2R}$ the homogeneous equation (\ref{L11-1}). Moreover, we have $F-\tilde F=\int_0^1\bar a\nabla wd\tau$
on $B_{2R}$ and thus $|F-\tilde F|\le\big(\int_0^1|\nabla w|^2d\tau\big)^\frac{1}{2}$.
Hence, (\ref{L11-2}) implies
\begin{align}\label{L11-5}
\lefteqn{\Big(\int_{B_R}|F-\tilde F|^2\Big)^\frac{1}{2}}\nonumber\\
&\lesssim\frac{1}{R^p}\Big(\int_0^1t^p\big(\int_{B_{2R}\setminus B_R}\bar u^2\big)^\frac{1}{2}dt
+\int_0^1t^p\big(\int_{B_{2R}\setminus B_R}\bar{\tilde u}^2\big)^\frac{1}{2}dt\Big),
\end{align}
where we used the triangle inequality on the rhs. In order to obtain (\ref{Laux11.1}) (with a redefined $p$)
it thus remains to compare
\begin{align}\label{L11-6}
\lefteqn{\frac{1}{R^p}\int_0^1t^p\big(\int_{B_{2R}}u^2\Big)^\frac{1}{2}dt}\nonumber\\
&\quad\mbox{vs}\quad\frac{1}{R^p}\int_0^1dt\, t^p\fint_0^{\sqrt{t}}dr\int\eta_R(\frac{r}{\sqrt{t}})^p
|(q(t)-\langle q(t)\rangle)_r|,
\end{align}
for generic $a$, which we now do using arguments from the proof of Lemma~\ref{Laux6}.

\medskip

We first note that for $R\ge 1\ge\sqrt{t}$ we have using that $\eta_{4R}\lesssim \eta_R$
\begin{align}\label{L11-7}
\lefteqn{\big(\int_{B_{2R}}u^2\Big)^\frac{1}{2}}\nonumber\\
&\lesssim R^\frac{d}{2}\big(\int\eta_{2R} u^2\Big)^\frac{1}{2} 
\stackrel{(\ref{c85})}{\lesssim} R^\frac{d}{2}(\frac{R}{\sqrt{t}})^\frac{d}{2}
\int\eta_R(y)\big(\int\eta_{\sqrt{t}}(\cdot-y) u^2\big)^\frac{1}{2}dy,
\end{align}
where the functions are evaluated at time $t$.
We then appeal to Lemma~\ref{Laux7}, that is (\ref{Laux6.2}) with $\frac T2$ replaced by $t$:
\begin{equation}\nonumber
\big(\int\eta_{\sqrt{t}}u^2(t)\big)^\frac{1}{2}
\lesssim\fint_\frac{t}{2}^td\tau\frac{1}{\sqrt{\tau}}\fint_0^{\sqrt{\tau}}dr
\int\eta_{\sqrt{\tau}}(\frac{r}{\sqrt{\tau}})^{p-1}|(q(\tau)-\langle q(\tau)\rangle)_r|.
\end{equation}
We then insert the latter, translated by $y$, into (\ref{L11-7}) and make use of (\ref{c86})
\begin{equation}\nonumber
\big(\int_{B_{2R}}u(t)^2\Big)^\frac{1}{2}
\lesssim R^d\fint_{\frac{t}{2}}^td\tau\,\frac{1}{\sqrt{\tau}^{\frac{d}{2}+1}}
\fint_0^{\sqrt{\tau}}dr\int\eta_R(\frac{r}{\sqrt{\tau}})^{p-1}|(q(\tau)-\langle q(\tau)\rangle)_r|.
\end{equation}
We finally apply $\frac{1}{R^p}\int_0^1dt t^p$ to obtain (\ref{L11-6}) in form
\begin{align}\nonumber
\lefteqn{\frac{1}{R^p}\int_0^1dt\,t^p \big(\int_{B_{2R}}u(t)^2\Big)^\frac{1}{2}}\nonumber\\
&\lesssim\frac{1}{R^{p-d}}\fint_0^1d\tau\,\tau^{p-\frac{1}{2}(\frac{d}{2}+1)}
\fint_0^{\sqrt{\tau}}dr\int\eta_R(\frac{r}{\sqrt{\tau}})^{p-1}|(q(\tau)-\langle q(\tau)\rangle)_r|.
\end{align}
Together with (\ref{L11-5}), this yields (\ref{Laux11.1}) with $p$ replaced by $\min\{p-d,p-\frac{1}{2}(\frac{d}{2}+1),p-1\}$.


\subsection{Proof of Lemma~\ref{Laux12}: CLT cancellations}

We split the proof into six steps.

\medskip

\step 1 A property of centered random variables. 

For a centered $F$ we have for all $s_*\ge 2$
\begin{multline}\label{L12-1b}
\log\langle\exp(\nu F)\rangle\le |\nu|^{s_*}\;\mbox{for all}\;|\nu|\ge 1
\\
\Longrightarrow\quad
\log\langle\exp(\nu F)\rangle\lesssim\nu^2+|\nu|^{s_*}\;\mbox{for all}\;\nu.
\end{multline}
To this purpose we note that
\begin{eqnarray*}\nonumber
\frac{d}{d\nu}\log\langle\exp(\nu F)\rangle
&=&\frac{\langle F\exp(\nu F)\rangle}{\langle\exp(\nu F)\rangle},\nonumber\\
\frac{d^2}{d\nu^2}\log\langle\exp(\nu F)\rangle
&=&\frac{\langle F^2\exp(\nu F)\rangle}{\langle\exp(\nu F)\rangle}
-\big(\frac{\langle F\exp(\nu F)\rangle}{\langle\exp(\nu F)\rangle}\big)^2.
\end{eqnarray*}
Since $F$ is centered, we infer from Jensen's inequality
that $\langle\exp(\nu F)\rangle\ge \exp(\nu \langle F\rangle)=1$
and from
the first formula that
\begin{equation}\nonumber
\log\langle\exp(\nu F)\rangle_{|\nu=0}=0,\quad
\frac{d}{d\nu}_{|\nu=0}\log\langle\exp(\nu F)\rangle=0.
\end{equation}
From the second formula, the lower bound $\langle\exp(\nu F)\rangle\ge 1$, and Jensen's inequality,
we deduce that
\begin{equation}\nonumber
0\le\frac{d^2}{d\nu^2}\log\langle\exp(\nu F)\rangle
\le\langle F^2(\exp(F)+\exp(-F))\rangle\quad\mbox{for}\;|\nu|\le 1.
\end{equation}
By the real variable inequality $q^2(\exp(q)+\exp(-q))\le(\exp(q)+\exp(-q))^2=\exp(2q)+\exp(-2q)+2$,
this implies $\langle F^2(\exp(F)+\exp(-F))\rangle\le\langle\exp(2F)\rangle+\langle\exp(-2F)\rangle+2
\le2(\exp(2^2+2^{s_*})+1)$ by the assumption in (\ref{L12-1b}).
Hence the second derivative of $\log\langle\exp(\nu F)\rangle$
in $\nu\in[-1,1]$ is bounded, while the function vanishes to first order in the origin.
This establishes the conclusion in (\ref{L12-1b}).

\medskip

\step 2 An equivalent norm based on the characteristic function: We claim that for all $1<s\le 2$,
\begin{align*}
\|F\|_*\sim\inf\big\{M>0\,|\,\forall\;\nu\in\mathbb{R}\;\log\langle\exp(\nu\frac{F-\langle F\rangle}{M})\rangle\le
\max\{\frac{1}{2}\nu^2,\frac{1}{s_*}|\nu|^{s_*}\}\big\},
\end{align*}
where $\frac1s+\frac1{s_*}=1$.

We start with showing that the lhs is controlled by the rhs. By homogeneity, we may assume that $\langle F \rangle =0$ and that
$$
\log \langle \exp(\nu F)\rangle \,\leq \,\max\{\frac{1}{2}\nu^2,\frac{1}{s_*}|\nu|^{s_*}\} \qquad \text{for all }\nu \in \R,
$$
which, thanks to $\exp |F|\le \exp(F)+\exp(-F)$, implies
$$
\langle \exp(\nu |F|)\rangle \,\leq\, 2\exp (\max\{\frac{1}{2}\nu^2,\frac{1}{s_*}|\nu|^{s_*}\}) \qquad \text{for all }\nu \in \R.
$$
Thus, by Chebyshev' inequality, we have for any threshold $N\in [0,\infty)$
$$
\langle I(|F|\ge N)\rangle \leq 2 \exp(\max\{\frac{1}{2}\nu^2,\frac{1}{s_*}|\nu|^{s_*}\}- \nu N) \qquad \text{for all }\nu \in \R,
$$
which yields by optimization in $\nu$:
$$
\langle I(|F|\ge N)\rangle \leq 2 \exp(-\min\{\frac1s N^s,\frac12 N^2\}).
$$
From this we obtain for any $M\in (0,\infty)$
\begin{eqnarray*}
\langle \exp( (\tfrac{|F|}{M})^s)\rangle &=& \int_0^\infty \langle I(|F|\ge N)\rangle \frac{d}{dN}\exp((\tfrac{N}{M})^s)dN
\\ 
&\lesssim & \frac{1}{M^{s}} \int_1^\infty N^{s-1} \exp(-(\tfrac1s-\tfrac1{M^s})N^s)dN+\frac{1}{M^s}\exp(\frac1{M^s}).
\end{eqnarray*}
Clearly this expression is smaller than $1$ for $M\gg 1$, which proves the first part of the claim.

We turn now to the proof that the lhs controls the rhs.
By homogeneity we may wlog assume that $\langle F\rangle =0$ and that
$$
\log \langle \exp(|F|^s) \rangle \le 1.
$$
For any $M \in (0,\infty)$ we have by Young's inequality for some $k\ge 1$ to be fixed later
$$
|\nu \tfrac{F}{M}|\,\leq \, \frac1{sk^s} |F|^s+(1-\frac1{sk^s})\frac{k^{s_*}}{s_*(1-\frac1{sk^s})} (\tfrac{|\nu|}{M})^{s_*}
$$
and thus by convexity of $F\mapsto \log \langle \exp F \rangle$
\begin{eqnarray*}
\log \langle \exp|\nu \tfrac{F}{M}|\rangle &\le & \frac1{sk^s} \log \langle \exp|F|^s\rangle +\frac{k^{s_*}}{s_*} (\tfrac{|\nu|}{M})^{s_*} \\
&\leq & \frac1{sk^s} +\frac{k^{s_*}}{s_*} (\tfrac{|\nu|}{M})^{s_*} \qquad \text{for all }\nu \in \R.
\end{eqnarray*}
In particular, for $k\gg 1$ and then $M \gg 1$,
$$
\log \langle \exp |\nu \tfrac{F}{M}|\rangle \,\ll \, \frac1{s_*} |\nu|^{s_*} \qquad \text{for all }|\nu|\ge 1.
$$
It remains to argue that because of $\langle F \rangle =0$, we have 
$$
\sup_{|\nu|\le 1} \frac1{\nu^2} \log \langle \exp(\nu \tfrac FM) \rangle \,\lesssim \, \log \langle \exp|\tfrac FM|\rangle,
$$
which follows from Step~1.
The combination then yields as desired for $M\gg 1$
$$
\log \langle \exp(\nu \tfrac FM)\rangle \,\leq \, \max\{\frac12 \nu^2,\frac1{s_*}|\nu|^{s_*}\} \qquad \text{for all }\nu \in \R.
$$

\medskip

\step3 CLT for exactly local random variables.
For a random variable $F$ that is exactly local on a scale $r\ge 1$ in the sense that
\begin{align}\label{L12-1}
F(\tilde a)=F(a)\quad\mbox{provided}\;\tilde a=a\;\mbox{on}\;B_r
\end{align}
we claim that for all $s\le 2$ and all convolution kernels $\bar G$, 
\begin{align}\label{L12-2}
\| \bar G*G_\frac12*F\|_* \,\lesssim \, r^\frac{d}{2} \Big(\int \bar G^{2}\Big)^\frac1{2} \|F\|_*,
\end{align}
where $G_\frac12$ is our Gaussian kernel (this is only needed for a neat treatment of the case $s<2$), where we identify
$F=F(a)$ with its spatial extension $F(x,a)=F(a(\cdot+x))$.
By scaling, it is enough to assume that $r=1$, $\langle F \rangle =0$,  and $\|F\|=1$, in which case, by Step~2,
we reformulate the claim as: For all $M>0$ and $\nu \in \R$
\begin{multline}\label{L12-3}
\log\langle\exp(\nu \frac{F}{M})\rangle\, \leq\, \max\{\frac{1}{2}\nu^2,\frac{1}{s_*}|\nu|^{s_*}\} 
\\
\implies 
\log\langle\exp(\nu \frac{\bar G*G_\frac12*F}{(\int  \bar G^2)^\frac1{2} M})\rangle\,\lesssim \,\max\{\frac{1}{2}\nu^2,\frac{1}{s_*}|\nu|^{s_*}\} 
\end{multline}
for a random function $F$ that is {\it exactly} local on scale one.
Set $G'(x):=\bar G*G_\frac12(-x)$.
In order to establish (\ref{L12-3}), we write $F'=F'(0)=\int_{\mathbb{R}^d}( G' F)(x)dx$
as convex combination of a {\it sum} of independent
random variables:
\begin{equation}\nonumber
F'=\fint_{[0,3)^d}\sum_{z\in\mathbb{Z}^d}3^d( G'F)(3z+x)dx.
\end{equation}
Note that for fixed $x$, the points in the lattice $\{3z+x\}_{z\in\mathbb{Z}^d}$ have distance at least
3 and thus the numbers $\{F(3z+x)\}_{z\in\mathbb{Z}^d}$
and therefore $\{3^d( G'F)(3z+x)\}_{z\in\mathbb{Z}^d}$
are independent since $F$ is one-local
and $\langle \cdot \rangle$ has range unity. Together with the convexity of $F\mapsto\log\langle\exp F\rangle$ this implies
\begin{eqnarray*}\nonumber
\log\langle\exp(\nu \frac{F'}{M})\rangle
&\le&\fint_{[0,3)^d}\log\big\langle\exp\big(\nu\sum_{z\in\mathbb{Z}^d}3^d\frac{(G'F)(3z+x)}{M}\big)\big\rangle dx
\nonumber\\
&=&\fint_{[0,3)^d}\sum_{z\in\mathbb{Z}^d}\log\langle\exp(\nu3^d\frac{(G'F)(3z+x)}{M})\rangle dx.
\end{eqnarray*}
By the assumption in (\ref{L12-3}) with $\nu$ replaced by $3^d G'(3z+x)\nu$
in conjunction with stationarity of $\langle\cdot\rangle$, this yields
\begin{eqnarray*}\nonumber
\log\langle\exp(\nu \frac{ F'}{M})\rangle
&\lesssim &\fint_{[0,3)^d}\sum_{z\in\mathbb{Z}^d}  G'^2(3z+x)\nu^2dx+\fint_{[0,3)^d}\sum_{z\in\mathbb{Z}^d}  G'^{s_*}(3z+x)\nu^{s_*}dx\\
&\lesssim &\nu^2 \int_{\mathbb{R}^d} G'^2 +\nu^{s_*}  \int_{\mathbb{R}^d} G'^{s_*} \\
&{\lesssim} &\max\{\frac12 (\nu (\int  G'^2)^\frac12)^2,\frac1{s_*} (|\nu| (\int  G'^{s_*})^\frac1{s_*})^{s_*}\}.
\end{eqnarray*}
By Young's inequality on $G'(x)=\bar G *G_\frac12(-x)$ with exponents $\frac12+\frac{s_*+2}{2s_*}=\frac1{s_*}+1$,
$$
 \Big(\int   G'^{s_*}\Big)^\frac1{s_*} \,\leq\,  \Big(\int   \bar G^{2}\Big)^\frac1{2} \Big(\int  G_\frac12^{\frac{2s_*}{s_*+2}}\Big)^\frac{s_*+2}{2s_*}\,\lesssim \,   \Big(\int   \bar G^{2}\Big)^\frac1{2},
$$
so that the above turns into the conclusion in (\ref{L12-3}) in form of
\begin{equation*}\nonumber
\log\langle\exp(\nu \frac{ \bar G*G_\frac12*F}{M})\rangle
\,{\lesssim} \, \max\{\frac12 (\nu (\int  \bar G^2)^\frac1{2})^2,\frac1{s_*} (|\nu| (\int  \bar G^2)^\frac1{2})^{s_*}\}.
\end{equation*}

\medskip

\step4 Change of averaging kernel.

Let $G'$ denote a Schwartz function different from our Gaussian $G$,
let $F_r'$ denote the averaging wrt to this variable, that is, $F_r'(a)=\int dz G_r'(z)F(a(\cdot-z))$
where $G'_r(z)=\frac{1}{r^d}G(\frac{z}{r})$. Then we claim that for all $s \le 2$,
\begin{align}\label{L12-4}
\sup_{r\le 1}r^{\frac d2}\|F'_r\|_*\sim \sup_{r\le 1}r^{\frac d2}\|F_r\|_*.
\end{align}
In fact, the exponent $\frac d2$
could be replaced by any finite exponent $p_0$.
For future reference, we note that by scaling, $r\le 1$ can be replaced by  $r\le R$ for any $R>0$ in \eqref{L12-4}, so that (by taking the supremum over $R>0$),  \eqref{L12-4} also implies 
\begin{align}\label{L12-eq-clt}
\sup_{R>0}R^{\frac d2}\|F'_R\|_*\sim \sup_{R>0}R^{\frac d2}\|F_R\|_*.
\end{align}
Let $\psi$ and $\psi'$ be two Schwartz functions
of unit integral; let $\psi_r:=\frac{1}{r^d}\psi(\frac{\cdot}{r})$ and $\psi_r'$ denote the rescalings
to some scale $r$. 
For a given integrable random variable $F$ with vanishing expectation let us momentarily denote by
$F_r(a)=\int\psi_r(-z)F(a(z+\cdot))dz$ and $F_r'$ the corresponding convolutions. We claim that for any $p<\infty$
and all $\nu \in \R$
\begin{equation}\label{L12-5}
\sup_{r\le 1}\log\langle\exp(\nu r^p F_r)\rangle\le \nu^2+|\nu|^{s_*}
\quad\Longrightarrow\quad
\log\langle\exp(\nu F'_1)\rangle\lesssim \nu^2+|\nu|^{s_*}.
\end{equation}
Let us argue that this yields \eqref{L12-4}. Indeed, by scaling, \eqref{L12-5} takes the form for all $R>0$,
\begin{equation*}
\sup_{r\le R}\log\langle\exp(\nu r^p F_r)\rangle\le \nu^2+|\nu|^{s_*}
\quad\Longrightarrow\quad
\log\langle\exp(\nu R^p F'_R)\rangle\lesssim \nu^2+|\nu|^{s_*},
\end{equation*}
Taking the supremum over $R\le 1$ on both sides yields
\begin{equation*}
\sup_{r\le 1}\log\langle\exp(\nu r^p F_r)\rangle\le \nu^2+|\nu|^{s_*}
\quad\Longrightarrow\quad
\sup_{R\le 1} \log\langle\exp(\nu R^p F'_R)\rangle\lesssim \nu^2+|\nu|^{s_*},
\end{equation*}
from which \eqref{L12-4} follows by exchanging the roles of $G$ and $G'$.

\medskip

In order to prove \eqref{L12-5}, we start by arguing that we may assume in addition that $\psi$ has vanishing
spatial moments of order $1,\cdots,p-1$, that is,
\begin{equation}\label{L12-6}
\int x^\alpha\psi dx=\left\{\begin{array}{cl}1&\mbox{for}\;\alpha=0\\0&\mbox{for}\;0<|\alpha|<p\end{array}\right\},
\end{equation}
where for every multi-index $\alpha\in\{0,1,\cdots\}^d$ we use the standard notation
$x^\alpha:=x_1^{\alpha_1}\cdots x_d^{\alpha_d}$ and $|\alpha|:=\alpha_1+\cdots+\alpha_d$.
In order to establish this, we need to show that for a Schwartz kernel $\psi$, there {\it exists} a
(necessarily un-signed) Schwartz kernel $\psi'$ so that (\ref{L12-6}) holds and such that
for all $\nu \in \R$
\begin{equation}\label{L12-7}
\sup_{r\le 1}\log\langle\exp(\nu r^p F_r)\rangle\le \nu^2+|\nu|^{s_*}
\quad\Longrightarrow\quad
\sup_{r\le 1}\log\langle\exp(\nu r^p F'_r)\rangle\lesssim \nu^2+|\nu|^{s_*}.
\end{equation}
In fact, we shall construct a Schwartz function $\omega$ such that
\begin{equation}\label{L12-8}
\psi'=\omega*\psi
\end{equation}
satisfies (\ref{L12-6}). We first argue that the form (\ref{L12-8}) implies (\ref{L12-7}). Indeed,
(\ref{L12-8}) transfers to the rescaled level $\psi'_r=\omega_r*\psi_r$, which in turn implies in
terms of the random variable $F$ that the two averaged versions are related by 
$F_r'(a)=\int\omega_r(-z)F_r(a(z+\cdot))dz$, which entails
\begin{equation}\nonumber
|F_r'(a)|\le\int\frac{|\omega_r(-z)|}{m} m|F_r(a(z+\cdot))|dz\quad\mbox{where}\;m:=\int|\omega|,
\end{equation}
that is, $|F_r'|$ is a convex combination of translates of $m|F_r|$. Hence by convexity of
$g\mapsto\log\langle\exp(g)\rangle$ and stationarity of $\langle\cdot\rangle$ we obtain
\begin{equation}\nonumber
\log\langle\exp|\nu r^p F_r'|\rangle\le\log\langle\exp|m\nu r^p F_r|\rangle.
\end{equation}
Using $\log\langle\exp|F|\rangle\le\max\{\log\langle\exp(F)\rangle,\log\langle\exp(-F)\rangle\}+\log 2$,
we see that we obtain from the lhs of (\ref{L12-7}) 
\begin{equation}\label{L12-11}
\log\langle\exp|\nu r^p F_r'|\rangle\le(m\nu)^2+(m|\nu|)^{s_*}+\log 2.
\end{equation}
Since $\langle F\rangle=0$ translates by stationarity of $\langle\cdot\rangle$ to
$\langle F_r'\rangle=0$, Step~1 allows to upgrade this to the
rhs of (\ref{L12-7}).

\medskip

We now turn to the construction of $\psi'$ of the form (\ref{L12-8}) with (\ref{L12-6}). The basis 
for the construction is the following identity
\begin{equation}\nonumber
\int x^\beta\partial^\alpha\psi dx
=\left\{\begin{array}{cl}0&\mbox{for}\;|\beta|\le|\alpha|,\beta\not=\alpha\\
(-1)^{|\alpha|}\alpha!&\mbox{for}\;\beta=\alpha\end{array}\right\},
\end{equation}
which follows from integration by parts, and
where we use the standard notation of $\partial^\alpha=\partial_1^{\alpha_1}\cdots\partial_d^{\alpha_d}$
and $\alpha!=\alpha_1!\cdots\alpha_d!$. Because of its triangular structure (in the sense of matrices
indexed by $\alpha,\beta$) we
learn from this identity that the linear map
\begin{equation}\nonumber
(\omega_\beta)_{|\beta|<p}\mapsto
\big(\int x^\alpha\sum_{|\beta|<p}\omega_\beta\partial^\beta\psi dx\big)_{|\alpha|<p}
\end{equation}
is invertible. Since for any $\alpha,\beta$ we have $\lim_{r\downarrow0}\int x^\alpha\partial^\beta(\psi_r*\psi)dx
=\int x^\alpha\partial^\beta\psi dx$, there exists a convolution scale $r>0$ such that
\begin{equation}\nonumber
(\omega_\beta)_{|\beta|<p}\mapsto
\big(\int x^\alpha\sum_{|\beta|<p}\omega_\beta\partial^\beta(\psi_r*\psi) dx\big)_{|\alpha|<p}
\end{equation}
still is invertible. In particular, there exists a set of coefficients $(\omega_\beta)_{|\beta|<p}$ such that
\begin{equation}\nonumber
\int x^\alpha\sum_{|\beta|<p}\omega_\beta\partial^\beta(\psi_r*\psi) dx
=\left\{\begin{array}{cl}1&\mbox{for}\;\alpha=0\\0&\mbox{for}\;0<|\alpha|<p\end{array}\right\}.
\end{equation}
Since $\sum_{|\beta|<p}\omega_\beta\partial^\beta(\psi_r*\psi)=\omega*\psi$ 
with $\omega:=\sum_{|\beta|<p}\omega_\beta\partial^\beta\psi_r$ being a Schwartz function, 
we obtain the desired structure (\ref{L12-8}).
This closes the argument that we may wlog assume (\ref{L12-6}) while proving (\ref{L12-5}).

\medskip

We will now establish (\ref{L12-5}) under the assumption (\ref{L12-6}) in the weaker form of: For all $\nu \in \R$,
\begin{equation}\label{L12-9}
\sup_{r\le 1}\log\langle\exp(\nu r^{\tilde p} F_r)\rangle\le \nu^2+|\nu|^{s_*}
\quad\Longrightarrow\quad
\log\langle\exp(\nu F'_1)\rangle\lesssim \nu^2+|\nu|^{s_*}
\end{equation}
for any $\tilde p<\frac{p}{s_*}$; it is only seemingly weaker since $p$ in (\ref{L12-6}) was arbitrary.
For this purpose we will establish the
following representation of $\psi'$ in terms of $\{\psi_r\}_{r\le 1}$: For any ratio $\theta\ll 1$
(to be fixed in the sequel) there exist weight functions $\{\omega_n\}_{n=0,1,\cdots}$ such that
\begin{equation}\label{L12-10}
\psi'=\sum_{n=0}^\infty\omega_n*\psi_{r_n}\quad\mbox{with}\quad r_n:=\theta^{-n}\;\mbox{and}\;
\int|\omega_n|\lesssim (C_0\theta^p)^{n},
\end{equation}
where $C_0<\infty$ depends only on $\psi$, thus in particular not on $\theta$. 
Let us first argue how (\ref{L12-10}) implies (\ref{L12-9}), very much like we
argued above how (\ref{L12-8}) implies (\ref{L12-11}): Using (\ref{L12-10}) we estimate $|F_1'|$ as a convex
combination of translates of $\{|m F_{r_n}|\}_{n=0,1,\cdots}$:
\begin{equation}\nonumber
|F_1'(a)|\le\frac{1}{m}\sum_{n=0}^\infty\int|\omega_n(-z)||m F_{r_n}(a(z+\cdot))|dz\quad
\mbox{where}\;m:=\sum_{n=0}^\infty\int|\omega_n|.
\end{equation}
Hence we obtain
\begin{equation}\nonumber
\log\langle\exp|\nu F_1'|\rangle\le\frac{1}{m}\sum_{n=0}^\infty
m_n\log\langle\exp|m\nu F_{r_n}|\rangle\quad\mbox{where}\;m_n:=\int|\omega_n|.
\end{equation}
Inserting the lhs of (\ref{L12-9}) this yields for any $\nu\in\mathbb{R}$
\begin{equation}\nonumber
\log\langle\exp|\nu F_1'|\rangle\le\frac{1}{m}\sum_{n=0}^\infty
m_n((\frac{m\nu}{r_n^{\tilde p}})^2+(\frac{m|\nu|}{r_n^{\tilde p}})^{s_*}).
\end{equation}
Inserting the estimate (\ref{L12-10}), which provided $\theta\ll 1$
yields in particular $m\lesssim 1$,  this turns into
\begin{equation}\nonumber
\log\langle\exp|\nu F_1'|\rangle\lesssim\sum_{n=0}^\infty
(C_0\theta^{p})^n((\frac{\nu}{\theta^{\tilde p n}})^2+(\frac{|\nu|}{\theta^{\tilde p n}})^{s_*}).
\end{equation}
Since by assumption $2\tilde p\le s_*\tilde p<p$, we may choose $\theta$ so small that this series converges
to the effect of $\log\langle\exp|\nu F_1'|\rangle\lesssim\nu^2+|\nu|^{s_*}$,
that is the statement on the rhs of (\ref{L12-9}).

\medskip

We finally turn to the argument for the representation (\ref{L12-10}), which relies on the
two estimates for an arbitrary Schwartz function $\omega$
\begin{align}
\int|\omega-\psi_r*\omega|&\le C_0 r^p\int|\nabla^p\omega|,\label{L12-12}\\
\int|\nabla^p(\omega-\psi_r*\omega)|&\le C_0\int|\nabla^p\omega|\label{L12-13}
\end{align}
with some constant $C_0=C_0(\psi)$.
Indeed, equipped with (\ref{L12-12}) \& (\ref{L12-13}) we define the weight functions
$\{\omega_n\}_{n=0,1,\cdots}$ recursively via
\begin{equation}\label{x39}
\omega_0=\psi'\quad\mbox{and}\quad \omega_{n+1}=\omega_n-\psi_{r_{n+1}}*\omega_n,
\end{equation}
where $r_n=\theta^n$ in line with (\ref{L12-10}).
From (\ref{L12-13}) we learn at first 
that $\int|\nabla^p\omega_n|\le C_0^n\int|\nabla\psi'|\lesssim C_0^{n}$ and then from (\ref{L12-12}) that
\begin{equation}\label{x40}
\int|\omega_{n+1}|\lesssim C_0 r_{n+1}^p C_0^{n}=(C_0\theta^p)^{n+1},
\end{equation}
which is the desired estimate in (\ref{L12-10}).
From (\ref{x40}) we learn that provided $\theta\ll 1$ we have in particular 
$\int|\omega_n|\rightarrow 0$ as $n\uparrow\infty$
so that the recursive definition (\ref{x39}) in form of $\psi'=\omega_0=\omega_0*\psi_{r_0}+\omega_1
=\omega_0*\psi_{r_0}+\omega_1*\psi_{r_1}+\omega_2=\cdots$ yields the representation
in (\ref{L12-10}), as an identity in $L^1(\mathbb{R}^d)$.

\medskip

We close this step by giving the argument for (\ref{L12-12}) \& (\ref{L12-13}). Estimate (\ref{L12-13}) is
an immediate consequence of the triangle inequality, $\nabla^p(\psi_r*\omega)=\psi_r*\nabla^p\omega$,
and the convolution estimate $\int|\psi_r*\nabla^p\omega|\le\int|\psi|\int|\nabla^p\omega|$, so
that the constant is given by $C_1=1+\int|\psi|$. For estimate (\ref{L12-12}) we need the moment
condition (\ref{L12-6}) --- it allows to write
$(\psi_r*\omega-\omega)(x)$ $=\int\psi_r(-z)\big(\omega(x+z)-\sum_{|\alpha|<p}(\partial^\alpha\omega)(x)z^\alpha\big)dz$.
Together with Taylor's representation
$u(x+z)-\sum_{|\alpha|<p}(\partial^\alpha\omega)(x)z^\alpha$
$=\int_0^1(1-s)^{p-1}\sum_{|\alpha|=p}(\partial^\alpha\omega)(x+sz)z^\alpha ds$ 
and Cauchy-Schwarz' inequality in form of $|\sum_{|\alpha|=p}(\partial^\alpha\omega)(x+sz)z^\alpha|\le
|\nabla^p\omega(x+sz)||z|^p$ (provided $|\nabla^p\omega|^2$ is defined as the sum of the squares of all
partial derivatives of order $p$) we obtain the pointwise inequality
\begin{equation}\nonumber
|(\omega-\psi_r*\omega)(x)|\le\int_0^1(1-s)^{p-1}\int|z|^p|\psi_r(-z)||\nabla^p\omega(x+sz)|dz ds,
\end{equation}
which yields by Fubini
$\int|\omega-\psi_r*\omega|dx$ $\le\big(\int_0^1(1-s)^{p-1}ds\big)$ $\big(\int|z|^p|\psi_r(-z)|dz\big)$
$\big(\int|\nabla^p\omega|dx\big)$. This gives (\ref{L12-12}) with constant
$C_2=\frac{1}{p!}\int|\hat z|^p|\psi(\hat z)|d\hat z$ given by the $p$-th moment of $|\psi|$.
It remains to choose $C_0=\max\{C_1,C_2\}$.

\medskip

\step5 Reduction of statement from random stationary fields to simple random variables. 

We claim that it is enough to establish the following:
For any random variable $F$ that is local (on scale 1) relative to the family of
random variables $\{\bar F^{R}\}_{R\ge 1}$ in the sense of
\begin{align}\label{L12-14}
|F(\tilde a)-F(a)|\le(\frac{1}{R})^{p_0}(\bar F^R(a)+\bar F^R(\tilde a))\\
\quad\mbox{provided}\;\;\tilde a=a\;\;\mbox{on}\;B_{2R}\;\mbox{and}\;R\ge 1\nonumber
\end{align}
we have for all $s\le 2$ and all convolution kernels $\bar G$,
\begin{align}\label{L12-15}
 \|  \bar G*F_\frac12\|_*= \|  \bar G*G_\frac12*F\|_*\, \lesssim\, \Big(\int  \bar G^{2}\Big)^\frac1{2}\Big( \|F\|_*^{1-\frac{2d}{p_0}}\big(\sup_{R\ge 1}\|\bar F^R\|\big)^\frac{2d}{p_0}+\|F\|_*\Big),
\end{align}
with the usual understanding that we identify $F=F(a)$ with its stationary spatial extension $F(a,x)=F(a(\cdot +x))$.
First of all, by scaling it is enough to establish (\ref{Laux12.2}) for $\sqrt{T}=1$, that is,
\begin{align}\label{L12-16}
\|\bar G *G_1* F\|_*  \lesssim \Big(\int  \bar G^{2}\Big)^\frac1{2} \Big(\big(\sup_{r\le 1}r^\frac{d}{2}\|F_r\|_*\big)^{1-\frac{2d}{p_0}}\|\bar F\|^\frac{2d}{p_0}
+\sup_{r\le 1}r^\frac{d}{2}\|F_r\|_*\Big).
\end{align}
(To obtain estimate \eqref{Laux12.3}, it suffices to replace $G_1$ in \eqref{L12-16} by $G_{\frac 12}$, and use the semi-group property.) 
We have to show that (\ref{L12-15}) for any family of random variables $(F,\bar F^R)$ with (\ref{L12-14})
implies (\ref{L12-16}) for any pair of stationary random fields $(F,\bar F)$ with
\begin{align}\label{L12-17}
\Big(\fint_{B_1}|F(\tilde a)-F(a)|^2\Big)^\frac12  \le(\frac{1}{R})^{p_0}\int\eta_R(\bar F(a)+\bar F(\tilde a))\\
\quad\mbox{provided}\;\;\tilde a=a\;\;\mbox{on}\;B_{2R}\;\mbox{and}\;R\ge 1.\nonumber
\end{align}
To this purpose, we select a Schwartz kernel $G'\ge 0$ supported in $B_1$. For $r\le 1$ we claim that
the family of random variables 
$(r^{\frac d2} F_r',\bar F^R:=\int\eta_R\bar F)$ satisfies (\ref{L12-14}). Indeed,
the passage from (\ref{L12-17}) to (\ref{L12-14}) is obvious on the rhs, and follows because of $R\ge 1\ge r$
from $r^{\frac d2} |F'_r(\tilde a)-F'_r(a)|\le (\int G'^2 )^\frac12 \Big(\int_{B_1}|F(\tilde a)-F(a)|^2\Big)^\frac12$
for the lhs, where we crucially use that $G'$ is supported in $B_1$. For the family $(r^{\frac d2}F_r',\int\eta_R\bar F)$,
(\ref{L12-15}) takes the form 
\begin{align*}
 \| r^{\frac d2} \bar G*G_\frac12 *F_r'\|_*
\lesssim\, \Big(\int  \bar G^{2}\Big)^\frac1{2} \Big(\|r^{\frac d2} F_r'\|_*^{1-\frac{2d}{p_0}}\big(\sup_{R\ge 1}\|\int\eta_R\bar F\|\big)^\frac{2d}{p_0}+\|r^{\frac d2} F_r'\|_*\Big),
\end{align*}
which by the commutation of convolution in form of 
$\bar G*G_\frac12*F_r'=(\bar G*G_\frac12*F)_r'$ and
Jensen's inequality in form of $\|\int\eta_R \bar F\|\le\|\bar F\|$ turns into
\begin{align*}
r^{\frac d2} \|(\bar G*G_\frac12*F)_r'\|_*
\lesssim \Big(\int  \bar G^{2}\Big)^\frac1{2}\Big(\big(r^{\frac d2}\|F_r'\|_*\big)^{1-\frac{2d}{p_0}}\|\bar F\|^\frac{2d}{p_0}+r^{\frac d2}\|F_r'\|_*\Big).
\end{align*}
We take the supremum of this estimate in $r\le 1$ and apply Step~4, in form of (\ref{L12-4}), to $\bar G*G_\frac12*F$ and $F$:
\begin{align*}
\sup_{r\le 1}r^{\frac d2}\|(\bar G*G_\frac12*F)_r\|_*
\lesssim \Big(\int  \bar G^{2}\Big)^\frac1{2} \Big(\big(\sup_{r\le 1}r^{\frac d2}\|F_r\|_*)^{1-\frac{2d}{p_0}}\|\bar F\|^\frac{2d}{p_0}+\sup_{r\le 1}r^{\frac d2}\|F_r\|_*\Big).
\end{align*}
For our Gaussian convolution, we have the semi-group property $G_r*G_\frac12*(\bar G*F)=G_{\sqrt{\frac14+r^2}}*(\bar G*F)$ so that  we have for the lhs 
(by evaluating at $r=\frac{\sqrt{3}}{2}$)
$$
\sup_{r\le 1}r^{\frac d2} \|(\bar G*G_\frac12*F)_r\|_*\gtrsim \|(\bar G*F)_{1}\|_* = \|\bar G*G_1*F\|_*.
$$
This yields the desired estimate (\ref{L12-16}).

\medskip

\step6 Proof of (\ref{L12-15}) by a martingale argument.

We denote by $\langle\cdot|B_r\rangle$ the conditional expectation wrt to the $\sigma$-sub algebra
generated by the measurable functions $F=F(a)$ on $\Omega$ that depend on $a$ only through $a_{|B_r}$;
it can be considered as the $L^2(\langle\cdot\rangle)$-orthogonal projection on these functions,
see for instance \cite[Section 4.1]{Durrett}.
Since $\Omega$ is a compact topological space when endowed with the notion of H-convergence, 
the conditional expectation has the following regularity property \cite[Section 4.1.c]{Durrett}: 
For $\langle\cdot\rangle$-a.e. $a\in\Omega$, $\langle\cdot|B_r\rangle(a)$ is a probability
measure on $\Omega$. It is easy to check from the characterizing property of $\langle\cdot|B_r\rangle$ that
we have $\langle\langle(F(\tilde a)-F(a))^2|B_r\rangle_{\tilde a}\rangle_a=0$ for any continuous function $F=F(a)$ on $\Omega$ 
that depends on $a$ only through $a_{|B_r}$. By the definition of the topology of H-convergence, this implies
\begin{equation}\nonumber
\tilde a=a\quad\mbox{Lebesgue-a.e. in}\;B_r\quad\mbox{for $\langle\cdot|B_r\rangle(a)$-a.e.}\;\tilde a\quad
\mbox{for $\langle\cdot\rangle$-a.e.}\; a.
\end{equation}
Let now $F$ be such that $\|F\|<\infty$ for some $s\le 2$. In terms of conditional expectations, we may rewrite
(\ref{L12-14}) as 
\begin{align}\label{L12-24}
|F-\langle F|B_{2R}\rangle|\le(\frac{1}{R})^{p_0}(\bar F^R+\langle\bar F^R|B_{2R}\rangle)\quad\mbox{for}\;R\ge 1.
\end{align}
For a radius $r\ge 1$ to be optimized at the end, we introduce
\begin{align*}
F^0:=\langle F|B_r\rangle-\langle F\rangle\quad\mbox{and}\quad F^n:=\langle F|B_{2^nr}\rangle-\langle F|B_{2^{n-1}r}\rangle\;\;
\mbox{for}\;n\in\mathbb{N},
\end{align*}
so that we have the representation
\begin{align}\label{L12-25}
F-\langle F\rangle=\sum_{n=0}^\infty F^n,
\end{align}
where the almost-sure convergence of the telescopic sum follows from (\ref{L12-24}) in conjunction with the bounds $\langle F^2\rangle^\frac12 \lesssim \|F\|<\infty$
and $\sup_{R\ge 1} \|\bar F^R\|<\infty$, that we may assume to hold wlog.

On the one hand, we have by (\ref{L12-24}) for $n\in\mathbb{N}$
\begin{align*}
|F^n|\le(\frac{1}{2^{n-1}r})^{p_0}(\bar F^{2^nr}+\langle\bar F^{2^nr}|B_{2^nr}\rangle
+\bar F^{2^{n-1}r}+\langle\bar F^{2^{n-1}r}|B_{2^{n-1}r}\rangle)
\end{align*}
and thus by the triangle inequality
\begin{align}\label{L12-26}
\|F^n\|_*=\|F^n\|\le(\frac{1}{2^{n-1}r})^{p_0}4\sup_{R\ge 1}\|\bar F^R\|,
\end{align}
where the first identity follows from $\langle F^n\rangle=0$. We also note that
\begin{align}\label{L12-27}
\|F^0\|_*\le\|F\|_*.
\end{align}
On the other hand, we have in the sense of (\ref{L12-1}) for all $n\in\mathbb{N}$
\begin{align*}
F^n\;\mbox{is exactly local on scale}\;2^{n+1}r
\end{align*}
and thus by Step~3, cf (\ref{L12-2}),
\begin{align}\label{L12-28}
\|\bar G*G_\frac12*F^n\|_*\lesssim (2^nr)^\frac{d}{2} \Big(\int \bar G^{2}\Big)^\frac1{2} \|F^n\|_*.
\end{align}
We now obtain from the triangle inequality
\begin{eqnarray*}
\|\bar G*G_\frac12*F\|_*
&\stackrel{(\ref{L12-25})}{\le}&
\sum_{n=0}^\infty \|\bar G*G_\frac12*F^n\|_*\\
&\stackrel{(\ref{L12-28})}{\lesssim}&
 \Big(\int \bar G^{2}\Big)^\frac1{2}\sum_{n=0}^\infty(2^nr)^\frac{d}{2}\|F^n\|_*\\
&\stackrel{(\ref{L12-26}),(\ref{L12-27})}{\lesssim}&
 \Big(\int \bar G^{2}\Big)^\frac1{2}\Big(r^\frac{d}{2}\|F\|_*+\sum_{n=1}^\infty(2^nr)^{\frac{d}{2}-p_0}\sup_{R\ge 1}\|\bar F^R\|\Big)
 \\
&\stackrel{p_0>\frac{d}{2}}{\sim}& \Big(\int \bar G^{2}\Big)^\frac1{2}\Big(
r^\frac{d}{2}\|F\|_*+r^{\frac{d}{2}-p_0}\sup_{R\ge 1}\|\bar F^R\|\Big).
\end{eqnarray*}
The optimization in $r\ge 1$ now yields (\ref{L12-15}).

\subsection{Proof of Lemma~\ref{Laux13}: CLT-scaling and stochastic integrability}

Note that $\langle F_0\rangle\le\delta^p$ entails
\begin{align}\nonumber
|F|&=I(F_0-\langle F_0\rangle\le\delta^p)|F|+I(F_0-\langle F_0\rangle>\delta^p)|F|\nonumber\\
&\stackrel{(\ref{Laux13.1})}{\le}(\delta+h)^\frac{1}{p}|F_1|
+\|F\|_\infty I(F_0-\langle F_0\rangle>\delta^p).\nonumber
\end{align}
By the triangle inequality this yields
\begin{align}\nonumber
\|F\|_s\le(\delta+h)^\frac{1}{p}\|F_1\|_s
+\|F\|_\infty\|I(F_0-\langle F_0\rangle>\delta^p)\|_s.
\end{align}
It remains to note that for some event ${\mathcal B}$ we have by definition of the norm 
applied to a characteristic function $\|I({\mathcal B})\|_s=
\inf\{M>0\,|\,\langle I({\mathcal B})\rangle(\exp(\frac{1}{M^s})-1)+1\le e\}$
so that we obtain the exact relationship $\|I({\mathcal B})\|_s=\|I({\mathcal B})\|_2^{\frac{2}{s}}$. Finally, we have
Chebychef's inequality in form of $\|I(F_0-\langle F_0\rangle>\delta^p)\|_2\le\frac{\|F_0\|_{*,2}}{\delta^p}$.


\subsection{Proof of Lemma~\ref{Laux14}: uniform bounds}

We split the proof into three steps.

\medskip

\step1 Proof of uniform bound \eqref{Laux14.0} and locality \eqref{Laux14.0b} for $(\nabla \phi(T),q(T))_{\sqrt{T}}$.

By scaling, we may restrict to $T=1$.
Since for all $r\le 1$, $G_r\leq (\frac{1}{r})^\frac{d}{2} G \lesssim  (\frac{1}{r})^\frac{d}{2} \eta$ (where $\eta$ is the exponential kernel), 
 \eqref{Laux14.0} will follow by Jensen's inequality once we establish  
\begin{equation}\label{Laux14-1}
\big(\int\eta |(\int_0^1\nabla ud\tau,q(1))|^2\big)^\frac{1}{2}\lesssim 1.
\end{equation}

\medskip

We apply the localized energy estimate \eqref{Laux3.2} to $\phi(t):=\int_0^t ud\tau$, which is characterized by the
initial value problem
\begin{equation}\nonumber
\partial_t\phi-\nabla\cdot a\nabla \phi=\nabla\cdot ae,\quad \phi(t=0)=0.
\end{equation}
Applied to $v=\phi$, $f=v_0=0$, $g=ae$, and $R=2\sqrt{T}=2$ 
(recall we assumed with a slight abuse that (\ref{Laux3.2}) holds for all $\eta_r$ with $r\ge 1$), we learn that
\begin{equation}\label{Laux14-2}
\int_0^1\int\eta_2|(\nabla \phi,q=a(\nabla \phi+e))|^2dt\lesssim\int_0^1\int\eta_2|ae|^2dt\sim 1.
\end{equation}
We shall upgrade this estimate to \eqref{Laux14-1} at the end of this step, and turn our attention to
the locality statement  \eqref{Laux14.0b}, which takes the form
\begin{equation}\label{Laux14-5}
|(\nabla \phi)_1-(\nabla\tilde \phi)_1|,|(q-\langle q \rangle)_1-(\tilde q-\langle q \rangle )_1|
\lesssim\exp(-\frac{R}{C})\quad
\mbox{provided}\;a=\tilde a\;\mbox{on}\;B_R,
\end{equation}
with the notation $q=q(a)$, $\tilde q=q(\tilde a)$, $\phi=\phi(a)$, and $\tilde \phi=\phi(\tilde a)$.
We compare $\phi=\phi(a,t,x)$ to $\tilde \phi:=\phi(\tilde a,t,x)$ and note that
the difference satisfies
\begin{equation}\nonumber
\partial_t(\phi-\tilde \phi)-\nabla\cdot \tilde a\nabla(\phi-\tilde \phi)=
\nabla\cdot\big((a-\tilde a)(\nabla \phi+e)\big),\quad (\phi-\tilde \phi)(t=0)=0,
\end{equation}
so that once more by (\ref{Laux3.2})
\begin{equation}\nonumber
\int_0^1\int\eta|(\nabla(\phi-\tilde \phi),q-\tilde q)|^2dt\lesssim\int_0^1\int\eta|a-\tilde a|^2|q|^2dt
\lesssim\sup_{B_R^c}\eta_2\int_0^1\int\eta_2|q|^2dt,
\end{equation}
where in the second estimate, we used $\eta\sim\eta_2^2$ and our assumption in form of
${\rm supp}(a-\tilde a)\subset B_R^c$.
We now insert (\ref{Laux14-2}) and obtain as intermediate result
\begin{equation}\label{Laux14-3}
\int_0^1\int\eta|(\nabla(\phi-\tilde \phi),q-\tilde q)|^2dt\lesssim\exp(-\frac{R}{2}).
\end{equation}

\medskip

In order to upgrade the time-averaged estimates (\ref{Laux14-3}) \& (\ref{Laux14-2}) to the pointwise-in-time estimates, we need the following a priori estimate
\begin{equation}\label{Laux14-4}
\int\eta|\partial_t \nabla \phi|^2=\int\eta|\nabla u|^2\lesssim 1\quad\mbox{for}\;t\in[\frac12, 1],
\end{equation}
which follows from the semi-group estimate \eqref{Laux1.1} in Lemma~\ref{Laux1}.

\medskip

We upgrade (\ref{Laux14-3}): From the latter in form of
\begin{equation}\nonumber
\int_\frac{1}{2}^1\int\eta|\nabla(\phi-\tilde \phi)|^2dt\lesssim\exp(-\frac{R}{2})
\end{equation}
and (\ref{Laux14-4}) in form of
\begin{equation}\nonumber
\sup_{t\in(\frac{1}{2},1)}\int\eta|\partial_t \nabla(\phi-\tilde \phi)|^2\lesssim 1
\end{equation}
we obtain via the elementary interpolation inequality
\begin{equation}\nonumber
\int\eta|g|^2_{|t=1}
\lesssim \frac{1}{\gamma}\int_{1-\gamma}^1\int\eta|g|^2dt
+\gamma^2\sup_{t\in(1-\gamma,1)}\int\eta|\partial_t g|^2\quad\mbox{for}\;\gamma\ll 1
\end{equation}
which we use for $g=\nabla(\phi-\tilde \phi)$, that
\begin{equation}\nonumber
\int\eta|\nabla(\phi-\tilde \phi)|^2_{|t=1}
\lesssim \frac{1}{\gamma}\exp(-\frac{R}{2})+\gamma^2.
\end{equation}
The choice of $\gamma=\exp(-\frac{R}{6})$ yields
\begin{equation}\nonumber
\int\eta|(\nabla \phi-\nabla\tilde \phi,q-\tilde q)|^2_{|t=1}\sim
\int\eta|\nabla(\phi-\tilde \phi)|^2_{|t=1}
\lesssim\exp(-\frac{R}{3}).
\end{equation}
Using again that $G\lesssim \eta$, $|(\nabla \phi)_1-(\nabla\tilde \phi)_1|,$ $|(q-\langle q\rangle)_1-(\tilde q-\langle q\rangle)_1|$
${\lesssim}$ $\big(\int\eta|(\nabla \phi-\nabla\tilde \phi,$ $q-\tilde q)|^2\big)^\frac{1}{2}$,
this yields (\ref{Laux14.0b}) in its rescaled version \eqref{Laux14-5}.
The upgrade of  (\ref{Laux14-2}) to  \eqref{Laux14.0} is similar.

\medskip

\step2 Proof of the uniform bound \eqref{Laux14.1} on the homogenization error.

As in Step~1, we may wlog assume $T=1$, and by the relation between the Gaussian and exponential kernels, it is enough to prove
that
\begin{equation}\label{Laux14-14}
\big(\int\eta |(S_{\frac12\to1}-S^\ho_{\frac12\to1})q(\tfrac12)|^2\big)^\frac{1}{2}\lesssim 1.
\end{equation}
Let $v$ be the solution of the initial value problem
\begin{equation}\nonumber
\partial_t v-\nabla \cdot a_\ho \nabla v=0 \text{ for } t>\tfrac12, \quad v(\tfrac12)=u(\tfrac12)=\nabla \cdot q(\tfrac12), 
\end{equation}
to the effect that $(S_{\frac12 \to 1}-S_{\frac 12 \to 1}^\ho) q(\frac 12)=\int_\frac12^1a\nabla ud\tau-\int_\frac12^1a_\ho \nabla vd\tau$.

We apply \eqref{Laux3.2} to $V(t):=\int_\frac12^t (v-u)d\tau$, which is characterized by the
initial value problem
\begin{equation}\nonumber
\partial_tV-\nabla\cdot a_\ho \nabla V=\nabla\cdot (a-a_\ho) \nabla u \text{ for }t>\tfrac12,\quad V(\tfrac12)=0,
\end{equation}
so that by \eqref{Laux1.1} for $u$,
\begin{equation}\nonumber
\int_\frac12^1\int\eta|\nabla V|^2dt\stackrel{\eqref{Laux3.2}}{\lesssim}\int_\frac12^1\int\eta|\nabla u|^2dt\stackrel{\eqref{Laux1.1}}{\sim} 1.
\end{equation}
As in Step~1, we may upgrade this time-averaged estimate to the pointwise-in-time estimate
\begin{equation}\label{Laux14-13}
\int\eta|\nabla V(1)|^2\,\lesssim\, 1.
\end{equation}
We are in the position to conclude the proof of \eqref{Laux14-14}.
By the triangle inequality,
\begin{multline*}
\big(\int\eta |(S_{\frac12\to1}-S^\ho_{\frac12\to1})q(\tfrac12)|^2\big)^\frac{1}{2}\,\leq \,
\big(\int\eta |(a-a_\ho)\int_\frac12^1\nabla ud\tau|^2\big)^\frac{1}{2}
\\
+\big(\int\eta |a_\ho\int_\frac12^1(\nabla v-\nabla u)d\tau|^2\big)^\frac{1}{2},
\end{multline*}
so that  \eqref{Laux14-14} follows from \eqref{Laux14-1} (with $[0,1]$ replaced by $[\frac12,1]$), \eqref{Laux14-13},
and the upper bounds (\ref{1.3}) \&  (\ref{Laux4.3b}) on $a$ and $a_\ho$.

\medskip

\step3 Proof of uniform bound \eqref{Laux14.2} and locality \eqref{Laux14.3} for $|(q_{t_0}-\langle q_{t_0}\rangle)_{\delta \sqrt{t_0}}|$.

By scaling, we may wlog assume $t_0=1$ and will thus drop the index $t_0=1$ when writing
$q=q_{t_0=1}$, $\eta_\tau=\eta_{\tau \sqrt{t_°}}$, so that the locality in form of \eqref{Laux11.2} follows from the stronger pointwise statement
\begin{equation}\label{Laux14-8}
||q-\langle q\rangle|_{\tau}(a)-|q-\langle q\rangle|_{\tau}(\tilde a)|
\lesssim\,\tau^{-\frac d2} \exp(-\frac{R}{C})\quad
\mbox{provided}\;a=\tilde a\;\mbox{on}\;B_R.
\end{equation}
From the defining equation $\phi-\nabla\cdot a(\nabla\phi+e)=0$ we obtain by the
localized elliptic energy estimate \eqref{Laux3.1} for all $r\ge 1$
\begin{equation}\label{Laux14-6}
\int\eta_r (\phi^2+|\nabla\phi|^2) \lesssim 1
\end{equation}
which, noting that $\eta_\tau\leq \tau^{-d}\eta$, turns into
\begin{eqnarray}
\int\eta_\tau (\phi^2+|\nabla \phi|^2) &\stackrel{\eqref{Laux14-6}}{\lesssim} &\tau^{-d}.\label{Laux14-7}
\end{eqnarray}
From (\ref{Laux14-7}) we thus obtain by Jensen's inequality
\begin{equation*}
|\int G_\tau q| \,\leq \, \Big(\int G_\tau |q|^2\Big)^\frac12\,\lesssim\,  \Big(\int \eta_\tau |q|^2\Big) \lesssim \tau^{-\frac d2},
\end{equation*}
which settles the uniform bound \eqref{Laux14.2} (since the estimate is deterministic and the origin plays no role).

\medskip

As in Step~1, we now compare $\phi=\phi(a)$ with
$\tilde \phi:=\phi(\tilde a)$; from the difference of the equations
\begin{equation}\nonumber
(\phi-\tilde\phi)-\nabla\cdot\tilde a\nabla(\phi-\tilde\phi)=\nabla\cdot(a-\tilde a)(\nabla\phi+e)
\end{equation}
we get from the localized elliptic energy estimate \eqref{Laux3.1} (and the same post-precessing as above
to go from $\eta$ to $\eta_\tau$) that
$\int\eta |\nabla(\phi-\tilde\phi)|^2\lesssim
 \int\eta |(a-\tilde a)(\nabla\phi+e)|^2$. Since
$q-\tilde q=\tilde a\nabla(\phi-\tilde\phi)+(a-\tilde a)(\nabla\phi+e)$ this yields in particular
\begin{equation}\label{Laux14-9}
\int\eta_\tau |q-\tilde q|^2\lesssim
\tau^{-d} \int\eta|q-\tilde q|^2\lesssim
\tau^{-d} \int\eta |(a-\tilde a)(\nabla\phi+e)|^2\lesssim \tau^{-d}\int_{B_R^c}\eta|q|^2,
\end{equation}
where we used the assumption $a=\tilde a$ on $B_R$.

\medskip

We now post-process (\ref{Laux14-9}) in two ways:
On the one hand, we use as in Step~1 that we have  $\eta \sim\eta_{2}^2$ so that by (\ref{Laux14-6})
\begin{equation}\label{Laux14-10}
\int_{B_R^c}\eta |q|^2\lesssim (\sup_{B_R^c}\eta_{2})\int\eta_{2}|q|^2\lesssim\exp(-\frac{R}{C}).
\end{equation}
On the other hand, we have $|f_\tau|^2 {\lesssim}\int\eta_\tau |f|^2$, so that
by the triangle inequality
\begin{eqnarray}\label{Laux14-11}
||q-\langle q\rangle|_{\tau}-|\tilde q-\langle q\rangle|_{\tau}|^2
&\lesssim&\int\eta_\tau ||q-\langle q\rangle|-|\tilde q-\langle q\rangle||^2\nonumber\\
&\le&\int\eta_\tau |q-\tilde q|^2.
\end{eqnarray}
Inserting (\ref{Laux14-10}) and (\ref{Laux14-11}) into (\ref{Laux14-9}) yields (\ref{Laux14-8}).


\subsection{Proof of Lemma~\ref{Laux15}:  homogenized semi-group and CLT-norm $\norm \cdot \norm$}

We start by claiming that it is enough to show
\begin{align}\label{L15-1}
\|(S^{h}_{t\rightarrow T}q)_R\|_*\lesssim\int_{\frac{R}{2}}^\infty\frac{dr}{r}\|q_r\|_*,
\end{align}
and the same estimate with the propagator $S^h$ based on the spatially constant coefficients $a_h$ 
(which however vary in time) replaced by the semi-group $S^{\ho}$ based on the constant coefficients $a_{\ho}$. 
Indeed, multiplying \eqref{L15-1} by $R^\frac{d}2$ yields by definition of $\norm q\norm$
\begin{align*}
R^\frac{d}{2}\|(S^{h}_{t\rightarrow T}q)_R\|_*\lesssim R^\frac{d}{2} \int_{\frac{R}{2}}^\infty\frac{dr}{r^{\frac{d}2+1}} \Big(r^\frac{d}{2}\|q_r\|_*\Big)
\,\lesssim \, \norm q \norm,
\end{align*}
so that the claim follows from taking the supremum over $R$ in the lhs.

\medskip

In order to establish (\ref{L15-1}) we note that by the definition
of the propagator $S^h_{t\rightarrow T}$, we have
\begin{align}\label{L15-2}
S^h_{t\rightarrow T}q=q+\int_t^Ta_h\nabla vd\tau,
\end{align}
where $v$ is the solution of the initial value problem
\begin{align}\label{L15-3}
\partial_\tau v-\nabla\cdot a_h\nabla v=0\;\;\mbox{for}\;\tau\ge t,\quad
v=\nabla\cdot q\;\;\mbox{for}\;\tau=t,
\end{align}
where $q$ is a given field.
By the triangle inequality, the representation (\ref{L15-2}) yields
\begin{align}\label{L15-4}
\|(S^h_{t\rightarrow T}q)_R\|_*\le\|q_R\|_*+\int_t^\infty\|a_h\nabla v\|_*d\tau
\le\|q_R\|_*+2\int_t^\infty\|\nabla v\|d\tau,
\end{align}
where the second inequality follows from the fact that $a_h$ is constant in space (although time-dependent), $|a_h|\le 1$, and
$\|\cdot\|_*\le 2\|\cdot\|$.

\medskip

We first turn to the estimate of $\|\nabla v(\tau)\|$ for $0<\tau-t\le R^2$.
By Sobolev's estimate applied to $\eta_{R}\nabla v(\tau)$ we obtain in particular
\begin{align}\label{L15-5}
|\nabla v(\tau)|\lesssim\sum_{m=0}^{[\frac{d}{2}]+1}R^m\big(\int\eta_{R}|\nabla^m\nabla v(\tau)|^2\big)^\frac{1}{2}.
\end{align}
Since $\nabla \nabla^m$ commutes with the spatially-constant coefficient equation (\ref{L15-3}), we obtain by the localized energy estimate
\eqref{Laux3.2}
\begin{align*}
\big(\int\eta_{R}|\nabla^m\nabla v(\tau)|^2\big)^\frac{1}{2}
\lesssim
\big(\int\eta_{R}|\nabla^m\nabla v(t)|^2\big)^\frac{1}{2}.
\end{align*}
Since also convolution (with our Gaussian) commutes with equation (\ref{L15-3}),
we may apply the two previous estimates to $v_R$ to the effect of
\begin{align*}
|\nabla v_R(\tau)|\lesssim
\sum_{m=0}^{[\frac{d}{2}]+1}R^m\big(\int\eta_{R}|\nabla^m\nabla\nabla\cdot(q-\langle q\rangle)_R|^2\big)^\frac{1}{2}.
\end{align*}
In combination with the inverse estimates
\begin{align*}
\big(\int\eta_{R}|\nabla^m\nabla\nabla\cdot (q-\langle q\rangle)_R|^2\big)^\frac{1}{2}
&\lesssim\frac{1}{R^{m+2}}\big(\int\eta_{R}|(q-\langle q\rangle)_{\frac{7}{8}R}|_{R}^2\big)^\frac{1}{2}\\
&\lesssim\frac{1}{R^{m+2}}\int\eta_{R}|(q-\langle q\rangle)_{\frac{3}{4}R}|
\end{align*}
this yields 
\begin{align*}
|\nabla v_R(\tau)|\lesssim\frac{1}{R^2}\int\eta_{R}|(q-\langle q\rangle)_{\frac{3}{4}R}|
\end{align*}
and thus by the triangle inequality, stationarity, and the monotonicity of $R\mapsto\|F_R\|_*$
\begin{align*}
\|\nabla v_R(\tau)\|\lesssim\frac{1}{R^2}\|q_{\frac{3}{4}R}\|\lesssim\frac{1}{R^2}\int_{\frac{R}{2}}^{\frac{3R}{4}}\frac{dr}{r}\|q_r\|_*,
\end{align*}
so that
\begin{align}\label{L15-6}
\int_t^{t+R^2}\|\nabla v_R(\tau)\|d\tau\lesssim\int_{\frac{R}{2}}^\infty\frac{dr}{r}\|q_r\|_*.
\end{align}

\medskip

For the range $T:=\tau-t\ge R^2$ we have to capture more cancellations: They are provided by \eqref{Laux6.1} in Lemma~\ref{Laux6} which yields
\begin{align}\label{L15-7}
\lefteqn{\big(\int_{t+T}^{t+2T}\int\eta_{\sqrt{T}}|\nabla v(\tau)|^2d\tau\big)^\frac{1}{2}}\nonumber\\
&\lesssim\fint_t^{t+T}d\tau\fint_0^{\sqrt{T}}dr(\frac{r}{\sqrt{T}})^{\frac{d}{2}+1}
\big(\int\eta_{\sqrt{T}}|v_r(\tau)|^2\big)^\frac{1}{2}.
\end{align}
We first turn to the rhs and note that by the localized energy estimate followed by an inverse estimate
for $r\le\sqrt{T}$
\begin{align*}
\big(\int\eta_{\sqrt{T}}|v_r(\tau)|^2\big)^\frac{1}{2}
\lesssim\big(\int\eta_{\sqrt{T}}|\nabla\cdot q_r|^2\big)^\frac{1}{2}
\lesssim\frac{1}{r}(\frac{\sqrt{T}}{r})^\frac{d}{2}\int\eta_{\sqrt{T}}|(q-\langle q\rangle)_{\frac{r}{2}}|,
\end{align*}
so that (\ref{L15-7}) turns into
\begin{align}\label{L15-8}
\big(\int_{t+T}^{t+2T}\int\eta_{\sqrt{T}}|\nabla v(\tau)|^2d\tau\big)^\frac{1}{2}\lesssim
\frac{1}{\sqrt{T}}\fint_0^{\sqrt{T}}dr\int\eta_{\sqrt{T}}|(q-\langle q\rangle)_{r}|.
\end{align}
We now turn to the lhs of this estimate.
By the localized energy estimate for (\ref{L15-3}), which also holds for spatial derivatives,
we have 
\begin{align*}
\sqrt{T}^{m+1}\big(\int\eta_{\sqrt{T}}|\nabla^m\nabla v(t+2T)|^2\big)^\frac{1}{2}
\lesssim \big(\int_{t+T}^{t+2T}\int\eta_{\sqrt{T}}|\nabla v(\tau)|^2d\tau\big)^\frac{1}{2}.
\end{align*}
Together with the Sobolev estimate (\ref{L15-5}) (with $R$ replaced by $\sqrt{T}$) this yields
\begin{align*}
\sqrt{T}|\nabla v(t+2T)|
\lesssim \big(\int_{t+T}^{t+2T}\int\eta_{\sqrt{T}}|\nabla v(\tau)|^2d\tau\big)^\frac{1}{2}.
\end{align*}
Inserting this into (\ref{L15-8}) yields the deterministic estimate
\begin{align*}
T|\nabla v(t+2T)|
\lesssim \fint_0^{\sqrt{T}}dr\int\eta_{\sqrt{T}}|(q-\langle q\rangle)_{r}|,
\end{align*}
which we apply to $v$ replaced by $v_R$, obtaining by the semi-group property of convolution with a Gaussian
\begin{align}\label{e.gr}
T|\nabla v_R(t+2T)|
\lesssim \fint_0^{\sqrt{T}}dr\int\eta_{\sqrt{T}}|(q-\langle q\rangle)_{\sqrt{R^2+r^2}}|.
\end{align}
By the triangle inequality and stationarity, this implies the stochastic estimate
\begin{align*}
T\|\nabla v_R(t+2T)\|
\lesssim \fint_0^{\sqrt{T}}dr \|q_{\sqrt{R^2+r^2}}\|_*
=\int_0^{\sqrt{T}}\frac{dr}{r} \frac{r}{\sqrt{T}}\|q_{\sqrt{R^2+r^2}}\|_*,
\end{align*}
which by the change of variables $\hat r=\sqrt{R^2+r^2}$ (with $\frac{dr}{r}=\frac{d\hat r}{\hat r}$) turns into
\begin{align*}
T\|\nabla v_R(t+2T)\|
\lesssim\int_R^{\sqrt{T+R^2}}\frac{d\hat r}{\hat r} \frac{r}{\sqrt{T}}\|q_{\hat r}\|_*.
\end{align*}
Integrating in $2T\in(R^2,\infty)$ yields
\begin{align}\label{{L15-9}}
\int_{t+R^2}^\infty\|\nabla v_R(\tau)\|d\tau
\lesssim\int_R^\infty\frac{d\hat r}{\hat r}\int_{r^2}^\infty\frac{dT}{T}\frac{r}{\sqrt{T}}\|q_{\hat r}\|_*
\sim\int_R^\infty\frac{d\hat r}{\hat r}\|q_{\hat r}\|_*.
\end{align}

\medskip

The combination of (\ref{L15-6}) and (\ref{{L15-9}}) gives
\begin{align}\label{L15-10}
\int_{t}^\infty\|\nabla v_R(\tau)\|d\tau
\lesssim\int_{\frac{R}{2}}^\infty\frac{dr}{r}\|q_{r}\|_*.
\end{align}
In combination with $\|q_R\|_*\lesssim\int_{\frac{R}{2}}^R\frac{dr}{r}\|q_{r}\|_*$ (once more
by monotonicity in $R$), this shows that (\ref{L15-4}) implies (\ref{L15-1}).

\subsection{Proof of Lemma~\ref{Laux16}: suboptimal estimate of the CLT-norms}
We only treat $q(T)$, and first address the case $T\ge 1$.
Applying $\|\cdot\|_*$ to \eqref{Laux14.0} in Lemma~\ref{Laux14}, we obtain
\begin{equation}\label{L16-1}
(\frac{r}{\sqrt{T}})^\frac{d}{2} \|(q(T))_r\|_* \,\lesssim \, 1 \text{ for } r\leq \sqrt{T}.
\end{equation}

\medskip

In case $T\ge 1$, by \eqref{Laux14.0b}, we may apply Lemma~\ref{Laux12} (with $F=q(T)$ and $\bar F=1$).
Inserting \eqref{L16-1} into \eqref{Laux12.3}, we obtain
$$
(\frac{R}{\sqrt{T}})^\frac{d}{2}\|(q(T))_R\|_*\,\lesssim \, 1 \text{ for }R\ge \sqrt{T}.
$$
The combination of these two estimates yields \eqref{Laux16.2} in case $T\ge 1$.

\medskip

We now turn to the case $T\le 1$. By the monotonicity of $r\mapsto\|(q(T))_r\|_*$, used in the range of
$\sqrt{T}\le r\le 1$, we see that  \eqref{L16-1} implies 
\begin{align*}
r^\frac{d}{2}\|(q(T))_r\|_*\lesssim 1\;\mbox{for}\;r\le 1.
\end{align*}
Since by (\ref{Laux14.0b}), $q(T)$ is in particular approximately local on scale $1$ (relative to $1$),
we may again apply Lemma \ref{Laux12}, this time for scale $1$ (the smallest scale on which we have finite range), 
so that the above upgrades to
\begin{align*}
R^\frac{d}{2}\|(q(T))_R\|_*\lesssim 1\;\mbox{for}\;R\ge 1.
\end{align*}
Again, the combination of the last two estimates yields (\ref{Laux16.2}) in case of $T\le 1$.

\subsection{Proof of Lemma~\ref{Laux17}: anchoring lemma}

We split the proof into four steps.
In this proof, and in this proof only, we denote by $\|\cdot \|$ the $L^2(\R^d)$-norm.

\medskip

\step{1} Rescaling: We have
\begin{equation}\label{p13}
\langle|(q(T)-\langle q(T)\rangle)_{\delta\sqrt{T}}|^2\rangle=
\langle|(q(1)-\langle q(1)\rangle_{\sqrt{T}})_{\delta}|^2\rangle_{\sqrt{T}},
\end{equation}
where the measure $\langle\cdot\rangle_{\sqrt{T}}$ denotes the push-forward of $\langle\cdot\rangle$
under the map $a\mapsto \hat a$, where $\hat a(\hat x)=a(\sqrt{T}\hat x)$ is the coefficient field
rescaled by $\sqrt{T}$. Here comes the argument: This transformation of the coefficient field
is based on the spatial change of variables $x=\sqrt{T}\hat x$, which we complement with the temporal
change of variables $t=T\hat t$, so that we have for the elliptic operator $T(\partial_t-\nabla\cdot a\nabla)
=(\partial_{\hat t}-\hat\nabla\cdot\hat a\hat\nabla)$. In view of the transformation of the initial
condition $\sqrt{T}\nabla\cdot(ae)=\hat\nabla(\hat ae)$, the solution of the parabolic initial value
problem (\ref{e.1})-(\ref{e.2}) transforms according to $\sqrt{T}u=\hat u$, which leads to $T\nabla u=\hat\nabla\hat u$
and thus $\int_0^T\nabla udt=\int_0^1\hat\nabla\hat ud\hat t$ as well as $q(T)=\hat q(1)$. This establishes
(\ref{p13}). 

\medskip

\step{2} Continuity wrt to H-convergence. 

For any $\delta>0$, the expression
\begin{equation}\nonumber
\big(\int_0^1\nabla ud\tau,q(1)=a(\int_0^1\nabla ud\tau+e)\big)_\delta,
\end{equation}
seen as a function of the $\lambda$-uniform coefficient field $a$, is continuous wrt
to H-convergence. We recall that the topology defining H-convergence is the coarsest
topology on the space $\Omega$ of $\lambda$-uniformly elliptic coefficient fields $a=a(x)$
for which the following functionals $F$ are continuous: For any ball $B\subset\mathbb{R}^d$,
any two vector fields $h,\tilde h\in L^2(B)$ we consider $F=\int\tilde h\cdot(\nabla u,a\nabla u)$,
where $u$ is the Lax-Milgram solution of $-\nabla\cdot a\nabla u=\nabla\cdot h$ in $B$
that vanishes outside of $B$.

\medskip

Here comes the argument: We will argue step by step that certain classes of function(al)s are continuous wrt
to H-convergence. We start by replacing the elliptic equation $-\nabla\cdot a\nabla u=\nabla\cdot h$
by the resolvent equation
\begin{equation}\label{p17}
-zu-\nabla\cdot a\nabla u=\nabla\cdot h
\end{equation}
for any $z\in\mathbb{C}$ in the
concave sector $\{({\calR}z)_+<\lambda|z|\}$, where ${\calR}z$ denotes the real part of $z$,
and $({\calR}z)_+$ the positive part of ${\calR}z$. 
We then replace the ball $B$ by the whole space $\mathbb{R}^d$,
localizing the topology with help of our exponential weight function $\eta$, which requires restricting
to a version of the sector shifted to the left. 
We then use the standard complex-variable argument to pass from the resolvent in this shifted sector 
to the semi-group,
which in turn yields the result.

\medskip

Passing to the resolvent equation primarily requires the following estimate for (\ref{p17})
\begin{equation}\label{p18}
\big(\int|\nabla u|^2\big)^\frac{1}{2}
\le\frac{|z|+({\calR}z)_+}{\lambda|z|-({\calR}z)_+}
\big(\int|h|^2\big)^\frac{1}{2},
\end{equation}
an elementary estimate at the origin of the notion of sectorial operators.
Here comes the argument for (\ref{p18}): Testing (\ref{p17}) with the complex conjugate
$\bar u$ yields the identity
\begin{equation}\label{p20}
-z\int|u|^2+\int\nabla\bar u\cdot a\nabla u=\int\nabla\bar u\cdot h.
\end{equation}
Of this equation we take the real part; noting that because $a$ is real the first inequality in
(\ref{1.3}) implies $\nabla\bar u\cdot a\nabla u\ge\lambda|\nabla u|^2$, we obtain with help of Cauchy-Schwarz'
inequality on the rhs of (\ref{p20}) the estimate
\begin{equation}\label{p21}
-({\calR}z) \|u\|^2 + \lambda\|\nabla u\|^2 \le \|\nabla u\| \|h\|,
\end{equation}
where we temporarily introduced the notation $\|\cdot\|$ for the $L^2$-norm of functions
and vector fields (note that it does not matter whether we take $L^2(B)$ or $L^2(\mathbb{R}^d)$
since all functions and fields can be assumed to vanish outside $B$). We also may use (\ref{p20})
to estimate the first lhs term as follows, using the upper bound $|a\xi|\le|\xi|$ provided by
(\ref{1.3}):
\begin{equation}\label{p22}
|z|\|u\|^2\le\|\nabla u\|^2+\|\nabla u\|\|h\|.
\end{equation}
Rewriting (\ref{p21}) as $\lambda\|\nabla u\|^2\le\|\nabla u\|\|h\|+({\calR}z)_+\|u\|^2$,
inserting (\ref{p22}), and dividing by $\|\nabla u\|$ yields
$\lambda\|\nabla u\|\le\|h\|+\frac{({\calR}z)_+}{|z|}(\|\nabla u\|+\|h\|)$,
which entails (\ref{p18}).

\medskip

Let us now argue in favor of H-continuity of $a\mapsto \int\tilde h\cdot (\nabla u,a\nabla u)$,
where $h,\tilde h\in L^2(\mathbb{R}^d)$ and $u$ is the Lax-Milgram solution of (\ref{p17}) in some ball $B$
(always with the understanding of vanishing boundary conditions). 
To probe continuity we give ourselves a sequence $\{a_n\}_{n\uparrow\infty}$
of $\lambda$-uniformly elliptic coefficient fields that H-converge to $a$ on $\mathbb{R}^d$
and thus a fortiori on $B$ and have to show that $(\nabla u_n,a_n\nabla u_n)$ weakly
converges to $(\nabla u,a\nabla u)$ in $L^2$.  
According to (\ref{p18}), the corresponding $\{\nabla u_n\}_n$ are
bounded in $L^2$ and thus, after passing to a subsequence which
we don't indicate in our notation, weakly convergent
to some $\nabla u$. Since the domain $B$ is bounded, we obtain by Rellich's compactness theorem
that $u_n$ strongly converges to $u$ in $L^2$. Consider the Lax-Milgram solution $\tilde u_n$
defined through $-\nabla\cdot a_n\nabla \tilde u_n=\nabla\cdot h+zu$ in $B$. 
On the one hand, rewriting the equation for $u_n$ as $-\nabla\cdot a_n\nabla u_n=\nabla\cdot h+zu_n$,
we learn from the energy estimate in conjunction with Poincar\'e's estimate that
$\nabla\tilde u_n-\nabla u_n$ and thus also $a_n\nabla\tilde u_n-a_n\nabla u_n$
converges strongly (and a fortiori weakly) to zero in $L^2$; in particular the weak limit of $\nabla\tilde u_n$
has to agree with the weak limit $\nabla u$ of $\nabla u_n$. On the other hand,
by definition of H-convergence, $(\nabla \tilde u_n,a_n\nabla\tilde u_n)$ weakly converges to $(\nabla u, a\nabla u)$,
so that $(\nabla\tilde u_n,a_n\nabla\tilde u_n)$ has to weakly converge
to $(\nabla u,a\nabla u)$ too, where $u$ must be the Lax-Milgram solution of
$-\nabla\cdot a\nabla u=\nabla\cdot h+zu$ in $B$. 

\medskip

For $r\ge\frac{1}{\lambda}$, 
we now argue in favor of H-continuity of $a\mapsto\int\eta_r\tilde h\cdot(\nabla u,a\nabla u)$
where the vector field $h$ and the function $f$ satisfy $\int\eta_r|(h,f)|^2<\infty$ and
the vector field $\tilde h$ is in the dual space $\int\frac{1}{\eta_r}|\tilde h|^2<\infty$,
and where $u$ is related to $(h,f)$ by the resolvent equation
\begin{equation}\label{p23}
-zu-\nabla\cdot a\nabla u=\nabla\cdot h+f
\end{equation}
on the whole space $\mathbb{R}^d$ and $z\in\mathbb{C}$ lies in the shifted concave sector
\begin{equation}\label{p26}
{\calR}(z+1)<\frac{\lambda}{4}|z+1|,
\end{equation}
which is slightly less concave than for (\ref{p18}).
As above, this relies essentially on the following a priori estimate for (\ref{p23})
\begin{eqnarray}\label{p24}
\lefteqn{\big(\int\eta_r(|u|^2+|\nabla u|^2)\big)^\frac{1}{2}}\nonumber\\
&\le&
\frac{2|z+1|+{\calR}(z+1)}{\frac{\lambda}{2}|z+1|-2{\calR}(z+1)}
\big(\big(\int\eta_r|f|^2\big)^\frac{1}{2}+\big(\int\eta_r|h|^2\big)^\frac{1}{2}\big),
\end{eqnarray}
which holds for both a ball $B$ and the whole space $\mathbb{R}^d$ and in particular provides
well-posedness in the latter case. We first turn to the argument for (\ref{p24}); to this purpose,
we shift the equation rather than the sector $z$ is in, and thus consider
\begin{equation}\nonumber
-zu+({\rm id}-\nabla\cdot a\nabla)u=\nabla\cdot h+f.
\end{equation}
Testing this equation with $\eta_r\bar u$ and using Leibniz' rule yields the identity
\begin{align}\nonumber
\lefteqn{-z\int\eta_r|u|^2+\int\big(\eta_r(|u|^2+\nabla \bar u\cdot a\nabla u)+\bar u\nabla\eta_r\cdot a\nabla u\big)}\\
&=\int\big(\eta_r(\bar uf-\nabla \bar u\cdot h)-\bar u\nabla\eta_r\cdot h\big).\nonumber
\end{align}
Using next to the bounds (\ref{1.3}) on $a$ 
that for our exponential localization function $|\nabla\eta_r|\le\frac{1}{r}\eta_r$,
we have the following inequalities for $r\ge\frac{1}{\lambda}$:
\begin{eqnarray*}
\int\big(\eta_r(|u|^2+\nabla \bar u\cdot a\nabla u)+\bar u\nabla\eta_r\cdot a\nabla u\big)&\ge&
\frac{\lambda}{2}\int\eta_r(|u|^2+|\nabla u|^2),\\
\big|\int\big(\eta_r(|u|^2+\nabla \bar u\cdot a\nabla u)+\bar u\nabla\eta_r\cdot a\nabla u\big)\big|&\le&
2\int\eta_r(|u|^2+|\nabla u|^2),\\
\big|\int\big(\eta_r(\bar uf-\nabla \bar u\cdot h)-\bar u\nabla\eta_r\cdot h\big)\big|&\le&
\big(\int\eta_r(|u|^2+|\nabla u|^2)\big)^\frac{1}{2}\\
&\;\times&\big(\big(\int\eta_r|f|^2\big)^\frac{1}{2}+\big(\int\eta_r|h|^2\big)^\frac{1}{2}\big).
\end{eqnarray*}
We now may argue as for (\ref{p18}) with $\big(\int|\nabla u|^2\big)^\frac{1}{2}$
replaced by $\big(\int\eta_r(|u|^2+|\nabla u|^2)\big)^\frac{1}{2}$, the lower bound
$\lambda$ in (\ref{p21}) replaced by $\frac{\lambda}{2}$, the upper bound $1$ in (\ref{p22})
replaced by $2$, and the size of the rhs $\big(\int|\nabla u|^2\big)^\frac{1}{2}$
replaced by $\big(\int\eta_r|f|^2\big)^\frac{1}{2}+\big(\int\eta_r|h|^2\big)^\frac{1}{2}$. Hence this yields
\begin{equation}\nonumber
\big(\int\eta_r(|u|^2+|\nabla u|^2)\big)^\frac{1}{2}\le\frac{2|z|+{\calR}z}{\frac{\lambda}{2}|z|-2{\calR}z}
\big(\big(\int\eta_r|f|^2\big)^\frac{1}{2}+\big(\int\eta_r|h|^2\big)^\frac{1}{2}\big),
\end{equation}
which turns into (\ref{p24}) after shifting back.

\medskip

Equipped with the estimate (\ref{p24}), we now may turn to the H-continuity of the new class
of functionals $F$. We establish continuity by showing that such an $F$ is the limit
of the corresponding functionals $\{F_R\}_{R\uparrow\infty}$, 
for which $u$ is replaced by the Lax-Milgram solution $u_R$
for (\ref{p23}) with the whole space replaced by the centered ball $B_R$ of radius $R$, 
and which we thus know to be continuous 
by the previous argument. Since we establish this limit in the uniform topology, the continuity
of $F_R$ transfers to $F$. More precisely, we shall show that
\begin{equation}\label{p25}
\big(\int\eta_r|\nabla u_R-\nabla u|^2\big)^\frac{1}{2}
\lesssim\big(\frac{1}{R}+\exp(-\frac{R}{2r})\big)
\big(\int\eta_r(|f|^2+|h|^2)\big)^\frac{1}{2}.
\end{equation}
To this purpose we consider $w=u_R-u$ which satisfies the homogeneous equation $-zw-\nabla\cdot a\nabla w=0$,
but only in $B_R$ (without vanishing bc). We thus consider $\tilde w=\tilde\eta w$, where $\tilde\eta$ is
a smooth cut-off for $B_\frac{R}{2}$ in $B_R$, and note that it satisfies
$-z\tilde w-\nabla\cdot a\nabla \tilde w=\nabla\cdot\tilde h+\tilde f$, now on all of $\mathbb{R}^d$, 
but with the rhs $\tilde h:=-wa\nabla\tilde\eta$ and $\tilde f:=-\nabla\tilde\eta\cdot a\nabla w$.
Hence an application of (\ref{p24}) yields
\begin{equation}\nonumber
\int\eta_r|\nabla\tilde w|^2\lesssim\int\eta_r(|\tilde f|^2+|\tilde h|^2),
\end{equation}
which in view of the choice of $\tilde\eta$ turns into
\begin{equation}\nonumber
\int_{B_{\frac{R}{2}}}\eta_r|\nabla w|^2\lesssim\frac{1}{R^2}\int\eta_r(|w|^2+|\nabla w|^2).
\end{equation}
Using $\eta_r(x)\lesssim\exp(-\frac{|x|}{2r})\eta_{2r}(x)$ and 
$\sup_{x\in B_{\frac{R}{2}}^c}\exp(-\frac{|x|}{2r})\le\exp(-\frac{R}{4r})$ we have
\begin{equation}\nonumber
\int_{B_{\frac{R}{2}}^c}\eta_r|\nabla w|^2\lesssim\exp(-\frac{R}{4r})\int\eta_r|\nabla w|^2.
\end{equation}
The combination of the two last estimates yields by definition of $w=u_R-u$
\begin{eqnarray}\nonumber
\lefteqn{\big(\int\eta_r|\nabla u_R-\nabla u|^2\big)^\frac{1}{2}}\\
&\lesssim&\big(\frac{1}{R}+\exp(-\frac{R}{2r})\big)
\big(\int\eta_r(|u_R|^2+|u|^2+|\nabla u_R|^2+|\nabla u|^2)\big)^\frac{1}{2}.\nonumber
\end{eqnarray}
Applying once more (\ref{p24}) to $u$ and $u_R$ separately, we obtain (\ref{p25}).

\medskip

We now pass from the resolvent to the semi-group. We fix a function $f$ with $\int\eta_r |f|^2<\infty$,
and for $z$ within the sector (\ref{p26}) consider the solution $u_z$ of $-zu_z-\nabla\cdot a\nabla u_z=f$
in $\mathbb{R}^d$. Let $\Gamma$ be the set ${\calR}(z+2)=\frac{\lambda}{4}|z+2|$ oriented in such
a way that it positively circles around the positive axis, which is well within the sector
(\ref{p26}). Hence in view of (\ref{p24}), for $t>0$ the integral $\int_\Gamma e^{-tz} u_zdz$ converges
absolutely with values in the weighted $L^2$-space $(\int\eta_r(|\cdot|^2+|\nabla\cdot|^2))^\frac{1}{2}$. It is then
an easy consequence of $\frac{1}{2\pi i}\int_\Gamma e^{-tz}dz=0$ that
\begin{equation}\nonumber
u(t)=\frac{1}{2\pi i}\int_\Gamma e^{-tz}u_zdz
\end{equation}
solves the equation $\partial_tu-\nabla\cdot a\nabla u=0$, and of Cauchy's integral theorem
in form of $\lim_{t\downarrow0}\frac{1}{2\pi i}\int_\Gamma \frac{e^{-tz}}{\mu-z}dz=1$ for $\mu>-2$ 
that $u(t=0)=f$. In particular, we have
\begin{equation}\nonumber
{\textstyle\int}\tilde h\cdot(\nabla u(t),a\nabla u(t))=\frac{1}{2\pi i}\int_\Gamma e^{-tz}
{\textstyle\int}\tilde h\cdot(\nabla u_z,a\nabla u_z)\;dz,
\end{equation}
converging absolutely for any $\tilde h$ with $\int\frac{1}{\eta_r}|\tilde h|^2<\infty$. Hence by
Lebesgue's theory, the H-continuity of $a\mapsto\int\tilde h\cdot (\nabla u_z,a\nabla u_z)$ transfers
to $a\mapsto\int\tilde h\cdot (\nabla u(t),a\nabla u(t))$.

\medskip

Finally, the H-continuity of $(\nabla\int_0^1ud\tau,q(1))_\delta$ is a consequence of the
previous general statement. Indeed, we consider $U(a;t,x)=\int_0^tud\tau+x\cdot e$ and note that
it solves the homogeneous equation $\partial_tU-\nabla\cdot a\nabla U=0$ with 
initial data $U(t=0)=f$ where $f(x)=x\cdot e$. Since clearly $\int\eta_r f^2<\infty$ and
$\int\frac{1}{\eta_r}|G_\delta|^2<\infty$, where $G_\delta$ denotes the Gaussian of variance 
$\delta>0$, we have that $(\int_0^1\nabla ud\tau+e, q(1))_\delta=\big(\int G_\delta \nabla U(1),
\int G_\delta a \nabla U(1)\big)$ is an H-continuous function.

\medskip

\step{3} We now argue that a stationary ensemble $\langle\cdot\rangle$ that is of range $r$ for all
$r>0$ is deterministic in the sense that
\begin{equation}\label{p12}
a=\langle a\rangle\quad\mbox{ae in}\;\mathbb{R}^d\quad\mbox{almost-surely}.
\end{equation}
Here comes the argument: For a scale $\rho$, we consider the convolution with a Dirac sequence
$\tilde\eta_\rho(x)=\frac{1}{\rho^d}\tilde\eta_1(\frac{x}{\rho})$, where $\tilde\eta_1$ is supported
in $B_1$, $\tilde\eta_1\ge 0$, and $\int\tilde\eta_1=1$:
\begin{equation}\label{p14}
a_{\rho}(y)=\int\tilde\eta_\rho(x-y)a(x)dx=\int\tilde\eta_\rho(x)a(x+y)dx.
\end{equation}
We first note that for any $y\in\mathbb{R}^d$ and $\rho>0$, $a\mapsto a_{\rho}(y)$ is continuous
wrt H-convergence: Indeed, $e_i\cdot a_{\rho}(y)e_j$ is of the form $\int\tilde h\cdot a\nabla u$, provided we 
make the following choices.
\begin{itemize}
\item For the ball: $B=B_{2\rho}(y)$, 
\item for the square integrable vector field $\tilde h$ testing the weak convergence: $\tilde h(x)=\tilde\eta_{\rho}(x-y)e_i$,
\item and for the square integrable vector field $h$ defining the equation for $u$: $h=-a\nabla(\zeta x_j)$ where
$\zeta$ is a cut-off for $B_{\rho}(y)$ in $B_{2\rho}(y)$ so that $u=\zeta x_j$ and thus
$\nabla u=e_j$ in $B_{\rho}(y)$, the support of $x\mapsto \tilde\eta_\rho(x-y)$.
\end{itemize}
By definition of stationarity of $\langle\cdot\rangle$, 
we learn from the second representation in (\ref{p14})
\begin{equation}\nonumber
\langle a_\rho\rangle:=\langle a_\rho(y)\rangle\quad\mbox{does not depend on}\;y\in\mathbb{R}^d.
\end{equation}
By definition of range $r$, we learn from the fact that in view of (\ref{p14}), $a_\rho(y)$
depends on $a$ only through $a_{|B_\rho(y)}$:
\begin{equation}\nonumber
\langle a_\rho(y)\otimes a_\rho(y')\rangle=\langle a_\rho(y)\rangle\otimes \langle a_\rho(y')\rangle
\quad\mbox{for}\;|y-y'|\ge r+2\rho.
\end{equation}
Both last statements combine to
\begin{equation}\label{p15}
\langle (a_\rho(y)-\langle a_\rho\rangle):(a_\rho(y')-\langle a_\rho\rangle)\rangle=0
\quad\mbox{for}\;|y-y'|\ge r+2\rho.
\end{equation}
Denoting by the subscript $R$ another convolution, we compute
\begin{equation}\nonumber
\langle|(a_\rho)_R-\langle a_\rho\rangle|^2\rangle
=\int\int\tilde\eta_R(y)\tilde\eta_R(y')
\langle (a_\rho(y)-\langle a_\rho\rangle):(a_\rho(y')-\langle a_\rho\rangle)\rangle.
\end{equation}
By the upper bound in (\ref{1.3}), which implies $|a(x)|\le d$ for the Frobenius norm,
this yields in conjunction with (\ref{p15})
\begin{equation}\nonumber
\langle|(a_\rho)_R-\langle a_\rho\rangle|^2\rangle
\le d\sup|\tilde\eta_1||B_1|(\frac{r+2\rho}{R})^d,
\end{equation}
which we integrate over an arbitrary ball $B$ to obtain 
\begin{equation}\label{p16}
\langle\int_B|(a_\rho)_R-\langle a_\rho\rangle|^2\rangle\lesssim |B| (\frac{r+2\rho}{R})^d.
\end{equation}

\medskip

We first let $\rho$ tend to zero; since $a_\rho$ converges to the Lebesgue-measurable and bounded $a$
strongly in $L^2_{loc}$, we learn that also $\langle a_\rho\rangle$ converges to a limit we call $\langle a\rangle$,
so that (\ref{p16}) turns into $\langle\int_B|a_R-\langle a\rangle|^2\rangle\lesssim |B|(\frac{r}{R})^d$.
Since $r>0$ was arbitrary, this yields $\langle\int_B|a_R-\langle a\rangle|^2\rangle=0$. 
We finally let $R$ go to zero to obtain 
$\langle\int_B|a-\langle a\rangle|^2\rangle=0$. Since this holds for any ball $B$, we first
learn that $\langle a\rangle$ does not depend on $B$ and then deduce (\ref{p12}).

\medskip

\step{4} Indirect argument for \eqref{Laux17.1} by compactness. 

We assume that \eqref{Laux17.1} fails,
which means that there exists a sequence $\langle\cdot\rangle_n$ of admissible ensembles
and a sequence of times $T_n\uparrow\infty$ such that 
\begin{equation}\label{p13bis}
\liminf_{n\uparrow\infty}
\big\langle\big|(\int_0^{T_n}\nabla ud\tau, q(T_n)-\langle q(T_n)\rangle_n)_{\delta\sqrt{T_n}}\big|^2
\big\rangle_n>0.
\end{equation}
We write for brevity $\langle\cdot\rangle_{\hat n}$ for the rescaled ensemble $\langle\cdot\rangle_{n,\sqrt{T_n}}$
introduced in Step 1. We note that the space $\Omega$ of $\lambda$-uniformly elliptic coefficient fields
$a=a(x)$ on $\mathbb{R}^d$ is compact under H-convergence: For $\mathbb{R}^d$ replaced by a fixed bounded
open set this is explicit in the work of Murat \& Tartar, see \cite[Theorem 6.5]{Tartar}. In our case
like in theirs, compactness relies on the fact that the topology is characterized by the continuity of the maps
$F$ made explicit in Step~3, which are indexed by vector fields $h,\tilde h\in L^2(\mathbb{R}^d)$
and, in our case, by balls $B\subset\mathbb{R}^d$;
and that these maps can be approximated in the uniform topology by a countable set of maps,
which are generated by a countable dense set in $L^2(\mathbb{R}^d)$, and in our case,
by the balls with rational center and radius. The H-limits coming from different balls are compatible thanks to
the locality of that notion, see \cite[Lemma 10.5]{Tartar}. 
Since $\Omega$ is compact, the space of probability measures
on the topological space $\Omega$ (by which we of course mean those that respect the topology) is also compact.
Hence after extraction of a subsequence, which we don't indicate in our notation, we may assume that
our sequence of ensembles $\langle\cdot\rangle_{\hat n}$ weak-* converges to a probability measure $\langle\cdot\rangle$,
which means that for all continuous functions $F$ on $\Omega$, 
we have $\lim_{n\uparrow\infty}\langle F\rangle_{\hat n}=\langle F\rangle$.
By definition of the latter, this weak-* convergence preserves the stationarity and the property of 
having a given range. Since $\langle\cdot\rangle_{\hat n}$ has range $\frac{1}{\sqrt{T_n}}\downarrow 0$,
we have that the limiting measure is both stationary and of range $r$ for any $r>0$. Hence by Step~4,
it is deterministic in the sense of (\ref{p12}). The latter implies $u=0$ and thus
\begin{equation}\nonumber
(\nabla u(1),q(1))=(0,\langle a\rangle e)\quad\mbox{ae in}\;\mathbb{R}^d\quad\mbox{almost-surely},
\end{equation}
and thus in particular
\begin{equation}\nonumber
\langle|(\int_0^1\nabla u d\tau,q(1)-\langle q(1)\rangle)_{\delta}|^2\rangle=0.
\end{equation}
Since by Step~2, $(\int_0^1\nabla u d\tau,q(1))_{\delta}$ is a continuous function on $\Omega$,
this implies
\begin{equation}\nonumber
\lim_{n\uparrow\infty}\langle|(\int_0^1\nabla u d\tau,q(1)-\langle q(1)\rangle_{\hat n})_{\delta}|^2\rangle_{\hat n}=0,
\end{equation}
which by Step~1 is in contradiction with (\ref{p13bis}).

\section*{Acknowledgements}
We thank Peter Bella, Mitia Duerinckx, and Julian Fischer for suggestions on the manuscript.
AG acknowledges financial support from the European Research Council under
the European Community's Seventh Framework Programme (FP7/2014-2019 Grant Agreement
QUANTHOM 335410).

\bibliographystyle{plain}

\bibliography{lit.bib}

\end{document}